\documentclass[11pt]{article}
\usepackage[letterpaper,margin=1in]{geometry}
\usepackage{amsmath,amssymb,amsfonts,mathtools}
\usepackage{amsthm}
\usepackage{booktabs,array,longtable}
\usepackage{graphicx}
\usepackage{float}
\usepackage{enumitem}
\usepackage{microtype}
\usepackage{colortbl}
\usepackage{caption}
\usepackage{bm}
\usepackage{titletoc}
\usepackage{multirow}
\usepackage{tabularx}
\usepackage[authoryear,round]{natbib}
\usepackage[colorlinks=true,linkcolor=blue,citecolor=blue,urlcolor=blue]{hyperref}
\usepackage{xcolor}
\usepackage{etoolbox}
\definecolor{navy}{RGB}{0,51,180}

\setlist[itemize]{leftmargin=2em,itemsep=0.25em,topsep=0.25em}
\setlist[enumerate]{leftmargin=2.2em,itemsep=0.3em,topsep=0.3em}
\allowdisplaybreaks
\numberwithin{equation}{section}

\newtheorem{theorem}{Theorem}[section]
\newtheorem{proposition}[theorem]{Proposition}
\newtheorem{corollary}[theorem]{Corollary}
\newtheorem{lemma}[theorem]{Lemma}
\newtheorem{definition}[theorem]{Definition}
\theoremstyle{remark}
\newtheorem{remark}[theorem]{Remark}

\theoremstyle{plain}
\newtheorem{example}[theorem]{Example}

\newcommand{\Id}{I}

\newcommand{\norm}[1]{\left\lVert #1\right\rVert}

\hypersetup{
  pdftitle={The Statistical Geometry of Latent Recovery in Hidden Markov Models},
  pdfauthor={Trang Pham, Dung Le, Alessandro Rinaldo, Nhat Ho},
  pdfsubject={Separation-dependent inverse inequalities and latent-law recovery},
  pdfkeywords={hidden Markov models, Gaussian mixtures, separation, overfitting, Wasserstein distance}
}
\title{\bfseries The Statistical Geometry of Latent Recovery in Hidden Markov Models}

\begin{document}

\begingroup
\renewcommand{\thefootnote}{\fnsymbol{footnote}}
\footnotetext[1]{Equal contribution.}
\endgroup

\begin{center}

{\bfseries\Large
The Statistical Geometry of Latent Recovery in \\ Hidden Markov Models
}

\vspace*{.2in}

{\large
\begin{tabular}{cccc}
Trang Pham$^{*}$ &
Dung Le$^{*}$ &
Alessandro Rinaldo &
Nhat Ho
\end{tabular}
}

\vspace*{.2in}

\begin{tabular}{c}
Department of Statistics and Data Science \\
The University of Texas at Austin
\end{tabular}

\vspace*{.2in}
\today

\begin{abstract}
Hidden Markov models (HMMs) offer an important capability to interpret fitted latent states and transitions. When emission distributions nearly coincide, a small observation-density error does not ensure accurate estimation of the latent quantities needed for interpretation. This motivates us to characterize convergence rates for latent recovery that account for emission separation. Our main contribution lies in deriving deterministic inverse inequalities that convert Hellinger error between fitted and true densities into bounds on latent recovery errors. We develop these bounds for multivariate Gaussian HMMs with a common, known emission covariance matrix and combine them with a Hellinger estimation guarantee for the global maximum likelihood estimator to obtain latent recovery rates. Although the joint observation distribution is a Gaussian mixture with exponentially many path components, they share the same emission centers across time. Leveraging this shared structure, we show that the separation exponents are governed by the maximum number of simultaneously approaching centers rather than the number of path components. Under exact specification, the global bounds on state-sequence probability discrepancies and the 1-Wasserstein loss for latent location distributions exhibit sharp separation exponents that do not increase with the number of selected times. Moreover, sufficiently small observation-density error guarantees unique matching between fitted and true emission centers, which yields error bounds for individual emission centers and transition rows. Under over-specification, individual split-state transition rows are not necessarily identifiable, so we quantify aggregate latent recovery errors and show that state-sequence probability bounds incur larger separation exponents compared to exact specification. We conduct synthetic experiments to validate our theoretical results.

\end{abstract}
\end{center}
\section{Introduction}
Hidden Markov models (HMMs) serve as a standard approach to modeling sequential data with underlying state changes and temporal dependence \citep{baum1966statistical, rabiner1989tutorial, ephraim2002hidden, cappe2005inference, JSSv043b04}. In this framework, the hidden states evolve according to a finite-state Markov chain, and each observation is generated from the emission distribution associated with its current state. Owing to these modeling capabilities, HMMs and closely related Markov-switching models have been employed extensively in speech recognition \citep{rabiner1989tutorial,mark2024application}, genomics \citep{fridlyand2004hidden,guha2008bayesian,mahmud2011fast}, single-molecule biophysics \citep{mckinney2006analysis,bronson2009learning,van2014empirical}, ecological modeling \citep{langrock2012flexible,mcclintock2020uncovering}, neuroscience \citep{baker2014fast,vidaurre2017brain}, and a broad range of time-series applications \citep{fruhwirth2006finite,JSSv043b04}. In many of these applications, the central goal is to interpret the underlying process through fitted states and transitions. However, a close fit to the observations alone is insufficient for reliable interpretation unless the fitted emission distributions, state probabilities, and transition flows also accurately reflect their true counterparts. For Gaussian emissions with a common, known covariance, this difficulty becomes especially pronounced when emission centers approach one another. Specifically, models with substantially different state probabilities or transition matrices can produce nearly identical observation densities.

This gap between observation fit and latent recovery calls for deriving convergence rates that explicitly account for emission separation. To obtain these rates, we establish inverse inequalities that bound latent recovery error in terms of observation-density error and emission separation. Formally, let $\theta$ denote a candidate HMM parameter, $\theta^\star$ the true parameter, and $p_{\theta,T}$ the density of a length-$T$ observation sequence under $\theta$. Our inverse inequalities take the following form
\begin{equation}
\mathcal D_r(\theta,\theta^\star)
\leq C_r\Delta_{\mathrm{sep}}^{-\alpha_r}
h(p_{\theta,T},p_{\theta^\star,T})^{\nu_r},
\label{eq:problem}
\end{equation}
where $h$ denotes the Hellinger distance defined in equation~\eqref{eq:hellinger-definition}, $\Delta_{\mathrm{sep}}$ is the truncated minimum separation between true emission centers, as defined in Section~\ref{prelim:hmm}, and $\mathcal D_r$ measures error in a latent target involving $r$ time coordinates. The constant $C_r$ is independent of $\Delta_{\mathrm{sep}}$, while the exponents $\alpha_r$ and $\nu_r$ quantify the dependence of the recovery bound on emission separation and observation-density error, respectively. Once the Hellinger error of an estimator is controlled, these inequalities translate its density estimation rate into latent recovery rates. The recovery targets considered in this paper are model parameters and unconditional distributions of latent emission-center sequences, rather than the realized hidden path or its posterior distribution given the observations.

The separation exponent $\alpha_r$ in~\eqref{eq:problem} also depends on how many true emission centers approach one another, which is not captured by their minimum separation alone. We illustrate this dependence in Figure~\ref{fig:cluster-concept}, where four emission centers form either two separated pairs or a single cluster, with the same minimum separation in both configurations. In the first configuration, the four centers form two pairs of nearby centers that remain well separated from each other, whereas in the second, all four centers lie close together in a single cluster. This difference changes the separation factor in the exact-specification emission-center bound from $\Delta_{\mathrm{sep}}^{-2}$ to $\Delta_{\mathrm{sep}}^{-6}$, as established in Theorem~\ref{thm:globalexact}. Driven by this example, we analyze how the number of nearby centers, in addition to their minimum separation, determines the sensitivity of latent recovery to observation error.

\subsection{Latent Recovery in Hidden Markov Models}
\label{sec:intro-hmm-recovery}
Latent recovery in HMMs has long been a subject of statistical research, with classical work on likelihood-based estimation and model identification~\citep{ephraim2002hidden,cappe2005inference}. Over time, this literature has expanded to include finite-sample recovery guarantees, with more recent work examining how emission separation affects estimation rates.

\paragraph{Estimation and identifiability.}
Likelihood-based inference provides a classical approach to estimating HMM parameters from observed sequences. Consistency and asymptotic normality describe the behavior of these estimators as the sample size increases under suitable assumptions \citep{baum1966statistical,leroux1992maximum,bickel1998asymptotic,Douc_2004,Douc_2011}, while \citet{lehericy2021nonasymptotic} studies nonasymptotic guarantees for maximum likelihood under model misspecification. In practical applications, likelihood-based estimation also requires algorithms for computing the estimators. The expectation-maximization (EM) algorithm provides an iterative approach to this computation with latent variables \citep{dempster1977maximum}, and \citet{10.1214/16-AOS1435} develop guarantees that connect its optimization and statistical errors under suitable initialization and local stability conditions. Apart from computation, the specification of the fitted model requires choosing the number of hidden states. When this number is unknown, state-count estimation is addressed through penalized likelihood and information-theoretic procedures \citep{gassiat2003optimal} and nonparametric order estimators \citep{lehericy2019consistent}. From a Bayesian perspective, posterior analyses also account for an unknown state count or flexible emission distributions \citep{10.3150/13-BEJ550,10.1214/15-EJS1017}.

Underlying these estimation approaches is the question of whether the latent parameters can be uniquely determined from the observation distribution. Since hidden-state labels are arbitrary, uniqueness is understood up to a common state permutation. General identifiability results for latent-variable models, including HMMs, are developed by \citet{Allman_2009} through algebraic arguments based on conditional independence. More specifically, \citet{gassiat2016inference} establish identifiability for finite-state nonparametric HMMs with a full-rank transition matrix and linearly independent emission distributions. The linear-independence requirement can be relaxed to distinct emissions when the transition matrix is full-rank and ergodic, as shown by \citet{alexandrovich2016nonparametric}. These results clarify when the observation distribution distinguishes different latent models, but do not by themselves quantify how a small error in that distribution affects latent recovery. This quantitative relation is needed for recovery guarantees from finitely many observations, where the observation distribution must itself be estimated.

\begin{figure}[t]
\centering
\includegraphics[width=\textwidth]{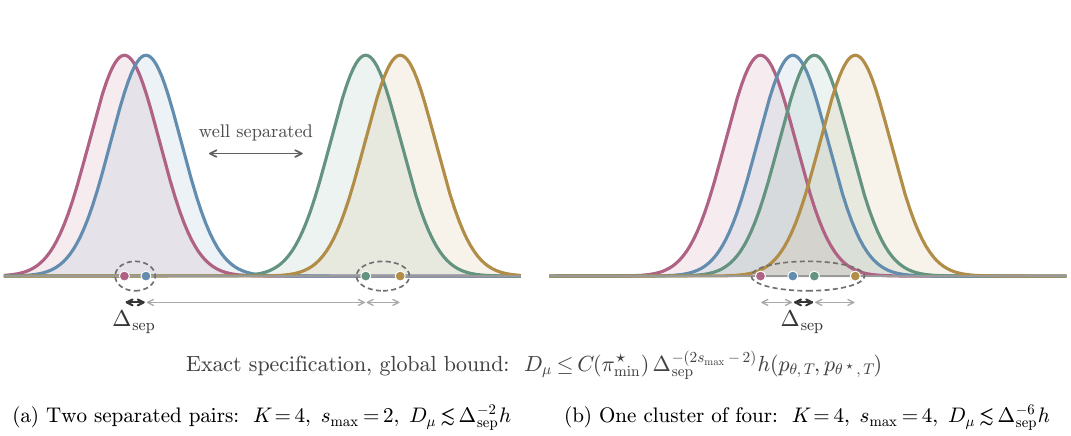}
\captionsetup{font=small}
\caption{\textbf{Same state count and minimum separation, different separation exponents.} Both configurations have $K=4$ Gaussian emission distributions and the same minimum separation $\Delta_{\mathrm{sep}}$. Colored dots mark emission centers, and dashed outlines group centers in the same cluster. Configuration (a) consists of two well-separated pairs with $s_{\max}=2$, whereas configuration (b) consists of one cluster of four with $s_{\max}=4$. Theorem~\ref{thm:globalexact} gives separation factors $\Delta_{\mathrm{sep}}^{-2}$ and $\Delta_{\mathrm{sep}}^{-6}$ in the global emission-center recovery bounds under exact specification for configurations (a) and (b), respectively.}
\label{fig:cluster-concept}
\end{figure}

\paragraph{Finite-sample recovery.}
Spectral and method-of-moments approaches provide finite-sample guarantees by using observable moments to construct matrix or tensor decompositions of the latent model \citep{HSU20121460,pmlr-v9-siddiqi10a,pmlr-v28-huang13,JMLR:v15:anandkumar14b}. Building on these approaches, subsequent work improves the statistical efficiency of moment-based estimation \citep{pmlr-v51-tran16} and extends spectral inference to hidden semi-Markov models \citep{JMLR:v18:15-468}. The recovery guarantees for such procedures depend on how accurately the latent quantities can be reconstructed from estimated moments. In the spectral analysis of \citet{HSU20121460}, this reconstruction is controlled through lower bounds for singular values of the model matrices. These conditions ensure that errors in observable quantities do not become arbitrarily amplified during recovery. The propagation of estimation error can also be studied directly through inverse inequalities that bound latent recovery error in terms of observation-density error. In the nonparametric HMM setting, \citet{JMLR:v17:15-381} use this approach to derive emission-estimation rates. Specifically, they first establish an oracle inequality for estimating the joint density of three consecutive observations, then show under an additional nondegeneracy condition that the squared emission-recovery error is bounded by the squared joint-density error. This relation converts their joint-density guarantee into minimax-adaptive rates for emission estimation, up to logarithmic factors. However, these rates depend on the lowest regularity among the emission densities, so they do not give faster rates for smoother densities within the same model. \citet{JMLR:v19:17-345} addresses this limitation through state-specific rates that adapt to the regularity of each emission density. In addition to emission-density estimation, parameter-recovery guarantees also support inference on hidden states conditional on the observations. For conditional inference, \citet{7907196} control filtering and marginal smoothing errors in terms of parameter-estimation errors and provide bounds for nonparametric spectral estimators. Nevertheless, the finite-sample results discussed above do not explicitly quantify how latent recovery rates depend on emission separation.

\paragraph{Recovery as emissions approach one another.}
\citet{9917566} explicitly quantify how parameter-estimation rates depend on emission separation in two-state HMMs with multinomial emissions. More precisely, they establish matching nonasymptotic minimax upper and lower bounds over classes in which the separation between the two emission probability vectors may decrease with the sample size. Consequently, their bounds quantify the additional data needed for recovery as the emissions become harder to distinguish. With nonparametric emissions, recovery additionally requires estimating the emission densities themselves, so their smoothness also affects estimation accuracy. Moreover, temporal dependence can provide information for distinguishing the hidden states, so recovery becomes more difficult as the observations approach an i.i.d. sequence. In this setting, \citet{abraham2025frontiers} obtain minimax bounds for two-state HMMs that account jointly for emission smoothness and how close the observations are to independence. Despite these differences in emission assumptions, both studies consider two-state HMMs, where separation involves only a single pair of emission distributions. With more states, however, the same minimum separation can arise from one cluster of approaching Gaussian emission centers or several smaller, separated clusters, as illustrated in Figure~\ref{fig:cluster-concept}. Thus, the two-state results leave unresolved how the number of approaching centers within each cluster affects latent recovery rates.

\subsection{Connections to Finite-Mixture Theory}
\label{sec:intro-mixture-theory}

For Gaussian HMMs with an arbitrary number of hidden states, finite-mixture theory provides a natural starting point for studying latent recovery. Indeed, over a fixed trajectory length, the joint observation density is a finite Gaussian mixture whose mixing measure describes the distribution of latent emission-center sequences. Therefore, recovering this measure amounts to recovering the latent location distributions of interest.

\paragraph{Mixing-distribution estimation.}
The mixture representation above connects latent recovery in Gaussian HMMs to the identification and estimation of mixing distributions. Both problems have been studied in finite and infinite mixture models \citep{10.1214/aoms/1177703862,10.1214/aoms/1177698520,le2026pde,bariletto2026convergence,le2026heterogeneous}, with broader discussions of mixture inference in \citet{titterington_85,7b17fdcf-c8b1-398f-aeab-286cfaff53a3,mclachlan2000finite}. At the population level, identifiability establishes when the observation density uniquely determines the mixing distribution. With finitely many observations, estimation additionally requires quantifying how accurately that distribution can be recovered from data. Likelihood-based methods address this estimation problem, and their computation is reviewed by \citet{8ca25de6-5fc1-3d5f-b37b-f82f74ac6984}. However, their statistical behavior can depart from regular parametric theory near singular or nonidentifiable configurations, as studied by \citet{10.1214/aos/1176324464,10.1214/aos/1017938921,10.1214/aos/1056562463}. Near such configurations, differences in component locations and probabilities can produce much smaller differences in the observation density. As a result, accurate density estimation need not imply equally accurate recovery of the mixing distribution. On the density-estimation side, \citet{10.1214/aos/1013203452} derive Hellinger convergence rates for maximum-likelihood estimators of Gaussian mixtures through bracketing-entropy bounds and empirical process tools \citep{geer2000empirical}. These guarantees control the observation-density error, but do not by themselves bound errors in component locations and probabilities. Accordingly, latent recovery requires an inverse bound that controls both types of error in terms of observation-density error.

\paragraph{Wasserstein recovery and separation geometry.}
Wasserstein distance measures errors in both component locations and probability masses, making it suitable for assessing mixing-distribution recovery \citep{villani_topics_2003}. As a link between density estimation and mixing-distribution recovery, \citet{10.1214/12-AOS1065} establishes quantitative relations between Wasserstein distances on mixing measures and divergences between the corresponding mixture densities for independent observations. Within this framework, the recovery rates depend on the identifiability properties of the mixture family and on whether the fitted model contains excess components. In particular, \citet{10.1214/16-AOS1444,10.1214/16-EJS1105} show how strong and weak identifiability lead to different convergence behavior in exactly specified and overfitted finite mixtures. At the local minimax level, \citet{10.1214/17-AOS1641} further establish how the number of excess components determines the optimal recovery rate near a fixed mixing distribution.

Through moment comparisons, Wasserstein recovery bounds can also be obtained for Gaussian location mixtures. Specifically, \citet{10.1214/19-AOS1873} derive polynomial interpolation inequalities that convert moment discrepancies into Wasserstein bounds and use them to obtain optimal estimation guarantees. In higher dimensions, \citet{10.1214/22-AOS2207} exploit moment tensors to establish minimax-optimal Wasserstein rates without separation assumptions. These minimax guarantees describe the worst-case difficulty over classes allowing arbitrarily close components, rather than the dependence of recovery accuracy on a particular configuration of component locations. From a complementary perspective, \citet{nguyen2026geometryseparationfinitegaussian} develop polynomial constructions that make the dependence on component separation and minimum mixture weight explicit in recovery bounds. 

However, directly applying generic mixture bounds to latent recovery in Gaussian HMMs poses two distinct difficulties. First, component-count-dependent bounds can make recovery guarantees depend on up to $K^T$ hidden paths, even though these paths reuse the same $K$ emission centers across time. Furthermore, even within this shared set, the total number of centers fails to reflect whether they all approach one another or instead concentrate into smaller, well-separated clusters. Second, unlike the freely varying weights inherent to generic mixtures, HMM path weights are constrained by the Markov factorization into initial-state and transition probabilities. In turn, generic bounds for unrestricted weights leave these structural dependencies unexploited. 

\subsection{Our Contributions}
\label{sec:intro-contributions}

This paper addresses these difficulties through inverse inequalities that exploit the shared emission centers and Markov structure of Gaussian HMMs. These inequalities connect observation-density error to latent recovery error under both exact specification, where the fitted and true models have the same number of states, and over-specification, where the fitted model has more states than the true model. Combined with density-estimation guarantees, they yield statistical recovery rates that quantify how estimation accuracy depends on sample size, emission separation, and cluster size. Our main contributions are summarized as follows.

\begin{enumerate}
\item \textbf{HMM-aware inverse bounds for latent recovery.} We establish deterministic inverse inequalities that bound latent recovery errors in terms of the Hellinger distance between fitted and true observation densities. Our analysis exploits the shared emission centers across time and the Markov factorization of path probabilities, avoiding a treatment of the observation density as an arbitrary $K^T$-component mixture. This structure allows us to express the dependence on emission separation through explicit powers, with constants independent of that separation. When the true centers form clusters whose diameters are bounded by a fixed multiple of the minimum separation and whose intercluster distances remain bounded away from zero, these powers depend on the largest cluster size $s_{\max}$. Without cluster assumptions, the corresponding bounds hold with $s_{\max}=K$.

\item \textbf{Sharp separation exponents under exact specification and over-specification.} We derive separation exponents under both exact specification and over-specification and prove their sharpness through explicit examples, with the powers of Hellinger error held fixed. Under exact specification, our global bounds have separation exponents $2s_{\max}-1$ for state-sequence probability recovery and $2s_{\max}-2$ for $1$-Wasserstein recovery of joint emission-center distributions, neither of which increases with the number of selected times. The constants in these global bounds may depend on the minimum true initial-state probability. When the observation-density error is sufficiently small, we prove that fitted and true emission centers admit a one-to-one matching and obtain inverse bounds with the same separation exponents but constants independent of that probability. Under over-specification, however, several fitted states may represent the same true state, so their individual transition rows need not be identifiable. Therefore, we control probabilities aggregated over these states, with separation exponents that now depend on the number of selected times. We further extend these aggregate probability bounds to grouping based on sufficiently accurate reference centers without changing their separation exponents.

\item {\textbf{Statistical recovery rates with explicit separation dependence.} We combine a nearly parametric Hellinger bound for the global MLE with the inverse inequalities above to obtain latent recovery rates that depend explicitly on emission separation and cluster size. For instance, emission-center displacement, weighted by fitted state probabilities, has rate $\mathcal{O}_{\mathbb P}\!\left(\Delta_{\mathrm{sep}}^{-(2s_{\max}-2)}(d/n)^{1/2}\right)$ under exact specification and $\mathcal{O}_{\mathbb P}\!\left(\Delta_{\mathrm{sep}}^{-(s_{\max}-1)}(d/n)^{1/4}\right)$ under over-specification, up to logarithmic factors for fixed trajectory length. In both cases, the separation factors increase as the true emission centers approach one another, so a larger sample size is needed to obtain the same upper bound. At the same minimum emission separation, these recovery rates have smaller inverse powers of separation when fewer emission centers approach one another within each cluster. This results in different statistical recovery guarantees for HMMs with the same state count and minimum emission separation, according to how many emission centers approach one another within the largest cluster.}
\end{enumerate}

\noindent In this paper, we consider independent trajectories observed over a fixed time horizon. Such settings arise in longitudinal studies and repeated experiments, where each subject or experimental run is observed at a fixed number of time points. In both cases, the amount of data increases through the number of trajectories rather than their length. As more independent trajectories become available, our recovery rates describe how accurately individual-state and aggregate latent quantities can be estimated, providing a statistical basis for interpreting fitted HMMs in applications such as disease progression, human activity recognition, and economic modeling.

\paragraph{Paper Organization.} The remainder of the paper is organized as follows. Section~\ref{sec:preliminaries} introduces the model, latent recovery losses, emission separation and clustering conditions, and establishes the density-estimation bound used to obtain statistical recovery rates. With these definitions, we develop the main recovery bounds under exact specification and over-specification in Sections~\ref{sec:exact_specification} and~\ref{sec:over_specification}, respectively. We empirically assess our recovery bounds and their corresponding statistical rates through synthetic experiments under exact specification and over-specification in Section~\ref{sec:simulation}. Section~\ref{sec:practical} discusses what these guarantees imply for the interpretation and use of fitted HMMs, and Section~\ref{sec:conclusion} concludes the paper. Proofs, supporting technical lemmas, and additional results are provided in the appendices.

\section{Notation and Preliminaries}
\label{sec:preliminaries}
In this section, we first introduce the Gaussian HMM and define loss functions for comparing fitted and true emission centers, state probabilities, transition flows, and emission-center distributions. Subsequently, we specify conditions on the separation of true emission centers within and between clusters that determine how observation-density error translates into latent recovery bounds. Finally, we establish a finite-sample Hellinger bound for the global maximum likelihood estimator, which we use in later sections to derive latent recovery rates from the inverse inequalities.

\subsection{Gaussian Hidden Markov Models}
\label{prelim:hmm}

\paragraph{Model parameters.} Let $K\ge2$ be the number of true hidden states and let $k\ge K$ be the number of states in the fitted model. A Gaussian HMM is parameterized by $\theta=(\pi,A,\mu_1,\ldots,\mu_k)$, where $\pi$ is the row vector of initial state probabilities, $A$ is the transition matrix, and $\mu_1,\ldots,\mu_k\in\mathbb R^d$ are the emission centers. We write $\mathbb P_\theta$ for probabilities under the model with parameter $\theta$. The latent states form a Markov chain with initial probability $\pi$ and transition matrix $A$
\[
Z_1\sim\pi,
\qquad
\mathbb P_\theta(Z_{t+1}=j\mid Z_t=i)=A_{ij}.
\]
Conditioned on the latent states, the observations are
independent and satisfy $X_t\mid Z_{1:T}\sim N_d(\mu_{Z_t},\Sigma)$, where the common covariance matrix $\Sigma$ is known and positive definite.
\paragraph{Observed data.} We observe $n$ independent trajectories $X_{1:T}^{(1)},\ldots,X_{1:T}^{(n)}$ generated by the same Gaussian HMM with true parameter $\theta^\star=(\pi^\star,A^\star,\mu_1^\star,\ldots,\mu_K^\star)$. Each trajectory consists of $T$ observations in $\mathbb{R}^d$ with its hidden-state sequence unobserved. Within each trajectory, however, the observations are dependent through the hidden Markov chain. Throughout our statistical analysis, $T\geq2$ remains fixed as $n$ tends to infinity, so additional data come from more independent trajectories rather than longer sequences.

\paragraph{Normalization of the emission covariance.} For simplicity, we reduce the emission covariance to the identity matrix by letting $S=\Sigma^{-1/2}$ and transforming each observation $X_t$ into $SX_t$. Accordingly, the transformed emission density is $\varphi_d(x-S\mu_{Z_t})$, where $\varphi_d$ denotes the standard Gaussian density on $\mathbb R^d$. The initial distribution $\pi$ and transition matrix $A$ remain unchanged. This transformation also preserves
the Hellinger distance in equation~\eqref{eq:problem}. For any $u\in\mathbb R^d$,
\[
\sqrt{\lambda_{\min}(\Sigma)}\,\|Su\|_2
\le \|u\|_2
\le \sqrt{\lambda_{\max}(\Sigma)}\,\|Su\|_2,
\]
where $\lambda_{\min}(\Sigma)$ and $\lambda_{\max}(\Sigma)$ denote the smallest and largest eigenvalues of $\Sigma$, respectively.
For fixed $\Sigma$, these inequalities imply that distances between emission centers before and after the transformation differ by at most constant factors depending on $\Sigma$. Consequently, the recovery bounds established in subsequent sections for the transformed observations and emission centers can be expressed in the original coordinates without changing their orders. We provide the detailed justification of the reduction to identity emission covariance and its effect on the recovery bounds in Appendix~\ref{app:whitening}.

Throughout the sequel, we retain the notation $X_t$ and $\mu_j$
for the transformed observations and centers.
All radius constraints, separation and clustering conditions,
and nearest-center assignments below are understood in
these transformed coordinates.
\paragraph{Observation density.}
We write $p_{\theta,T}$ for the density of a length-$T$
observation sequence under parameter $\theta$. For every positive integer $m$, write $[m]=\{1,\ldots,m\}$.
For any latent state sequence $z\in[k]^T$, its probability is
\[
\mathbb{P}_\theta(Z_{1:T}=z)
=\pi_{z_1}\prod_{t=1}^{T-1}A_{z_tz_{t+1}},
\]
by the Markov property. Conditional on $Z_{1:T}=z$, the
observations are independent, with joint density
\[
p_\theta(x_{1:T}\mid Z_{1:T}=z)
=\prod_{t=1}^{T}\varphi_d(x_t-\mu_{z_t}).
\]
The unconditional density of the length-$T$ observation
sequence $X_{1:T}$ is obtained by summing the conditional
densities over all latent paths $z\in[k]^T$, weighted by
their respective probabilities:
\begin{equation}
p_{\theta,T}(x_{1:T})
=
\sum_{z\in[k]^T}
\pi_{z_1}
\prod_{t=1}^{T-1}A_{z_tz_{t+1}}
\prod_{t=1}^{T}\varphi_d(x_t-\mu_{z_t}).
\end{equation}
\begin{definition}[Parameter space and model assumptions]
For $R\geq 1$, let $\Theta_k(R)$ be the class of parameters $\theta$ satisfying
\[
\pi_i\geq 0,\qquad \sum_{i=1}^k\pi_i=1,\qquad
A_{ij}\geq 0,\qquad \sum_{j=1}^k A_{ij}=1,\qquad
\mu_j\in B_d(0,R),
\]
for all $i,j\in[k]$, where $B_d(0,R)$ is the closed Euclidean ball of radius $R$ centered at the origin. Zero initial and transition probabilities and repeated fitted emission centers are allowed.

The true parameter $\theta^\star\in\Theta_K(R)$ has strictly positive initial state probabilities and pairwise distinct emission centers. Writing $a_j:=\mu_j^\star$, define
\[
\pi_{\min}^\star:=\min_{j\in[K]}\pi_j^\star>0,\qquad
\delta_{\mathrm{sep}}
:= \min_{\substack{j,\ell\in[K]\\ j\neq\ell}}
\lVert a_j-a_\ell\rVert_2 > 0,\qquad
\Delta_{\mathrm{sep}}:=1\wedge\delta_{\mathrm{sep}}.
\]
\end{definition}

The probability and radius conditions make $\Theta_k(R)$ a closed and bounded subset of a finite-dimensional Euclidean space, and hence compact. Since the log-likelihood is continuous on this space for any fixed observed sample, a global maximum exists. The radius bound also ensures that any two emission centers are at most $2R$ apart, which provides a uniform bound on transport costs in later recovery analysis. These properties hold without requiring the transition matrix to be invertible or the hidden chain to be stationary or irreducible. In particular, positive initial state probabilities do not ensure that every state remains occupied at later times. The fitted parameter space also permits states that are never visited, since initial and transition probabilities may be zero. Consequently, the true observation density can be represented in every class $\Theta_k(R)$ with $k\ge K$ by retaining the $K$ true states and adding $k-K$ unvisited states.\footnote{Assign zero initial probability to each added state and zero transition probability from every original state to each added state.}

\subsection{Latent Recovery Loss}
\label{sec:latent_recovery_loss}
State labels are arbitrary, and under over-specification several fitted states can represent a single true state. To compare state probabilities between the fitted and the true models, we assign each fitted state to its nearest true emission center and aggregate probabilities according to this assignment. For a fitted parameter
$\theta\in\Theta_k(R)$, define\footnote{Ties are resolved by choosing the smallest index.}
\begin{equation}
c(i):=\arg\min_{j\in[K]}\lVert\mu_i-a_j\rVert_2,
\qquad
V_j:=\{i\in[k]:c(i)=j\},
\qquad
r_i:=\lVert\mu_i-a_{c(i)}\rVert_2,
\end{equation}
The set $V_j$ contains the fitted states assigned to the true center $a_j$, and $r_i$ measures the displacement of the fitted center $\mu_i$ from its assigned true center. This construction allows multiple fitted states to belong to the same group without requiring a one-to-one matching.

\paragraph{State-sequence probability errors.} We write $q_t=\pi A^{t-1}$ and $q_t^\star=\pi^\star(A^\star)^{t-1}$ for the row vectors of state probabilities at time $t$ under the fitted and true models, respectively. Thus, $q_t(i)=\mathbb{P}_\theta(Z_t=i)$. To characterize temporal dependence beyond a single time point, consider indices $1\leq t_1<\cdots<t_r\leq T$ and write $I=(t_1,\ldots,t_r)$. We identify each nonempty subset of $[T]$ with its increasing enumeration; the full set $I=[T]$ is allowed. To compare state-sequence probabilities on the common index set $[K]^r$, define
\begin{equation}
\bar P_\theta^I(j_1,\ldots,j_r)
:=
\mathbb{P}_\theta
\bigl(c(Z_{t_1})=j_1,\ldots,c(Z_{t_r})=j_r\bigr)
\end{equation}
and
\begin{equation}
P_\star^I(j_1,\ldots,j_r)
:=
\mathbb{P}_{\theta^\star}
\bigl(Z_{t_1}=j_1,\ldots,Z_{t_r}=j_r\bigr).
\end{equation}
The distribution $\bar P_\theta^I$ combines the probabilities of fitted state sequences assigned to the same sequence of true-state labels. It can therefore be compared with $P_\star^I$ even when $k>K$. We measure their discrepancy by
\begin{equation}
D_I
:=
\lVert\bar P_\theta^I-P_\star^I\rVert_1
=
\sum_{(j_1,\ldots,j_r)\in[K]^r} \left|\bar P_\theta^I(j_1,\ldots,j_r)-P_\star^I(j_1,\ldots,j_r)\right|
=
2\operatorname{TV}(\bar P_\theta^I,P_\star^I).
\end{equation}
For a single time point, $D_{\{t\}}$ compares the aggregate fitted state probabilities with the true state probabilities. In particular, $D_\pi:=D_{\{1\}}$ measures the discrepancy between the aggregate initial probabilities $\sum_{i\in V_j}\pi_i$ and $\pi_j^\star$.

For consecutive times $t$ and $t+1$, the joint probability of $Z_t = i$ and $Z_{t+1} = b$ is $q_t(i)A_{ib}$. Summing over fitted states assigned to groups $V_j$ and $V_\ell$ therefore gives
\[
\bar P_\theta^{\{t,t+1\}}(j,\ell)
=
\sum_{i\in V_j}\sum_{b\in V_\ell}q_t(i)A_{ib},
\qquad
P_\star^{\{t,t+1\}}(j,\ell)
=
q_t^\star(j)A_{j\ell}^\star.
\]
We measure the largest discrepancy over consecutive time
pairs by
\begin{equation}
\begin{aligned}
D_A
&:= \max_{1\le t<T} D_{\{t,t+1\}} \\
&= \max_{1\le t<T}
\sum_{j,\ell=1}^{K}
\left|
\underbrace{
\sum_{i\in V_j}\sum_{b\in V_\ell}q_t(i)A_{ib}
}_{\text{fitted group-to-group flow}}
-
\underbrace{
q_t^\star(j)A_{j\ell}^\star
}_{\text{true state-to-state flow}}
\right|.
\end{aligned}
\end{equation}
Under state splitting, different allocations of flow among fitted states can leave the observation density unchanged. The fitted flows are therefore summed within each pair of groups before taking the absolute difference from the corresponding true flow.
Thus, $D_A$ measures discrepancies in \emph{occupancy-weighted joint probabilities}. Although its notation refers to the transition matrix, it is not a norm of $A-A^\star$: a transition row contributes according to the probability of occupying its source state. Recovering an individual row therefore requires the additional normalization described in Corollary~\ref{cor:rows}. As with the other losses in this section, the probabilities entering $D_A$ are unconditional under the fitted and true models.

\paragraph{Emission-center displacement.} The aggregate probability losses do not measure how far fitted centers lie from their assigned true centers. To quantify these location errors, we further define
\begin{equation}
D_\mu
:=
\max_{1\leq t\leq T}\sum_{i=1}^k q_t(i)r_i,
\qquad
V
:=
\max_{1\leq t\leq T}\sum_{i=1}^k q_t(i)r_i^2.
\end{equation}
At each time, the sums are the expected distance and expected squared distance from a fitted emission center to its assigned true center. The weights $q_t(i)$ account for the probability of occupying each fitted state, and the maximum gives the largest discrepancy over the observation horizon.

To measure the net displacement within each group, we also define
\begin{equation}
B
:=
\max_{1\leq t\leq T}
\sum_{j=1}^K
\left\|
\sum_{i\in V_j}q_t(i)(\mu_i-a_j)
\right\|_2.
\end{equation}
The displacement vectors within a group are added before taking the norm, so deviations in opposite directions can cancel. For example, suppose two fitted states in $V_j$ have emission centers $a_j-v$ and $a_j+v$, where $v\in\mathbb R^d$ is nonzero, and equal positive state probabilities at a given time. Although both centers differ from $a_j$, their net displacement is zero. However, their contributions to $D_\mu$ and $V$ remain positive. The quadratic quantity $V$ measures this spread even when the net displacement vanishes. In Section~\ref{sec:over_specification}, we bound $V$ under over-specification and use this bound to control both $D_\mu$ and squared Wasserstein distances.
The three displacement losses are related by elementary inequalities. For every fitted parameter, the triangle inequality, Cauchy--Schwarz, and the bound $r_i\le2R$ give
\begin{equation}
 B\le D_\mu\le\sqrt V,
 \qquad V\le2R D_\mu.
 \label{eq:basic-displacement-comparisons}
\end{equation}
These comparisons hold without matching or a density-error assumption, since they use only the definitions and bounded support. They show how a bound for $V$ controls the first-order displacement $D_\mu$ and the net displacement $B$. The converse implication from $B$ to $V$ is unavailable because opposite displacements within a group can cancel, as in the preceding example.

\paragraph{Wasserstein error for emission-center distributions.} The grouped state-sequence probabilities can be correct even when the fitted emission centers differ from the true centers. To assess errors in both locations and probabilities, we define \begin{equation} G_\theta^I := \mathcal{L}_\theta \bigl(\mu_{Z_{t_1}},\ldots,\mu_{Z_{t_r}}\bigr), \end{equation} where $\mathcal{L}_\theta$ denotes distribution under $\theta$, and write $G_\star^I:=G_{\theta^\star}^I$ for its true counterpart. This distribution assigns probabilities to sequences of actual emission-center locations, so it distinguishes fitted centers even when they belong to the same group. We compare $G_\theta^I$ and $G_\star^I$ through Wasserstein distances. For $I=(t_1,\ldots,t_r)$, both are supported on $(\mathbb{R}^d)^r$, whose elements represent sequences of emission centers at the selected times. For $p\in\{1,2\}$ and $u=(u_1,\ldots,u_r),v=(v_1,\ldots,v_r)\in(\mathbb{R}^d)^r$, define
\begin{equation}
d_{p,r}(u,v)
:=
\left(
\frac{1}{r}\sum_{\ell=1}^r
\|u_\ell-v_\ell\|_2^p
\right)^{1/p}.
\label{eq:distance}
\end{equation}
The associated Wasserstein distance is
\[
W_{p,I}(G_\theta^I,G_\star^I)
:=
\left(
\inf_{\gamma\in\Pi(G_\theta^I,G_\star^I)}
\int d_{p,r}(u,v)^p\,d\gamma(u,v)
\right)^{1/p},
\]
where $\Pi(G_\theta^I,G_\star^I)$ denotes the set of
couplings with marginals $G_\theta^I$ and $G_\star^I$.
This definition allows different numbers of support points and does not require one-to-one matching between fitted and true states. The factor $1/r$ in equation~\eqref{eq:distance} normalizes the transport cost by the number of selected time points. Table~\ref{tab:loss-guide} summarizes the quantities measured by these losses and the scope of the corresponding recovery guarantees.

\begin{table}[t]
\centering
\small
\captionsetup{font=small}
\setlength{\tabcolsep}{6pt}
\renewcommand{\arraystretch}{1.08}
\arrayrulecolor{black}
\begin{tabularx}{\textwidth}{@{}c>{\raggedright\arraybackslash}Xl@{}}
\toprule
\textbf{Notation} & \textbf{Recovery error} & \textbf{Time points} \\
\midrule[0.7pt]
\addlinespace[4pt]
& \multicolumn{2}{l}{\textbf{\textit{State probabilities and transition flows}}} \\
\addlinespace[4pt]
$D_\pi$ & $\ell_1$ error in grouped initial probabilities & $t=1$ \\
\arrayrulecolor{black!12}\midrule[0.25pt]\arrayrulecolor{black}
$D_A$ & $\ell_1$ error in group-to-group transition flows & Maximum over $1\le t<T$ \\
\arrayrulecolor{black!12}\midrule[0.25pt]\arrayrulecolor{black}
$D_I$ & $\ell_1$ error in grouped joint state probabilities & Selected times $I$ \\
\midrule[0.7pt]
\addlinespace[4pt]
& \multicolumn{2}{l}{\textbf{\textit{Emission-center displacement}}} \\
\addlinespace[4pt]
$D_\mu$ & Weighted sum of distances to assigned true centers & Maximum over $1\le t\le T$ \\
\arrayrulecolor{black!12}\midrule[0.25pt]\arrayrulecolor{black}
$B$ & Sum of norms of within-group weighted displacement sums & Maximum over $1\le t\le T$ \\
\arrayrulecolor{black!12}\midrule[0.25pt]\arrayrulecolor{black}
$V$ & Weighted sum of squared distances to assigned true centers & Maximum over $1\le t\le T$ \\
\midrule[0.7pt]
\addlinespace[4pt]
& \multicolumn{2}{l}{\textbf{\textit{Emission-center distributions}}} \\
\addlinespace[4pt]
$W_{p,I}$ & Wasserstein distance between fitted and true joint emission-center distributions & Selected times $I$ \\
\bottomrule
\end{tabularx}
\caption{\textbf{Summary of latent recovery losses.} Fitted states are grouped by their nearest true emission centers. Displacement weights are the fitted state probabilities $q_t(i)$. Bounds on $D_A$, $D_\mu$, $B$, and $V$ control the corresponding errors uniformly over their indicated time ranges within the observed block.}
\label{tab:loss-guide}
\end{table}

\subsection{Emission Separation and Clustering}
The inverse bounds in Sections~\ref{sec:exact_specification} and~\ref{sec:over_specification} depend on both the separation of the true emission centers and the number of centers in each cluster. Nevertheless, the separation quantity $\Delta_{\mathrm{sep}}$ alone does not describe these cluster sizes. We capture the cluster structure through a partition in which centers within each cluster are close, while centers in different clusters are separated by a positive distance independent of $\Delta_{\mathrm{sep}}$.

\paragraph{Structured regimes.} In the structured regimes, assume that the true-state index set
admits a partition $[K]=C_1\cup\cdots\cup C_M$ into nonempty sets
such that
\begin{equation}
\label{eqn:dung_cluster_assumption}
\begin{aligned}
\|a_j-a_\ell\|_2 &\le C_0\Delta_{\mathrm{sep}},
&& \forall m\in[M],\ \forall j,\ell\in C_m,\\
\|a_j-a_\ell\|_2 &\ge D_0,
&& \forall m,b\in[M],\ m\ne b,\
   \forall j\in C_m,\ \ell\in C_b,
\end{aligned}
\end{equation}
where $C_0,D_0>0$ are fixed constants independent of $\Delta_{\mathrm{sep}}$.
The diameter condition bounds every pairwise distance within a cluster by $C_0\Delta_{\mathrm{sep}}$, whereas the second condition separates centers in different clusters by at least $D_0$. Thus, as $\Delta_{\mathrm{sep}}$ decreases, centers within each cluster approach one another while distinct clusters remain separated. For $M\geq2$, we additionally assume $\Delta_{\mathrm{sep}}\leq \frac{D_0}{4C_0}$. This ensures that each cluster has diameter at most $D_0/4$, so its internal scale is smaller than the distance to any
other cluster.

To express the dependence of the recovery bounds on cluster size,
we define
\[
s_{\max} :=\max_{m\in[M]}|C_m|.
\]
This quantity denotes the largest number of true emission centers in a single cluster. Our proofs use polynomials that equal one at a chosen true center and zero at the other true centers. The expectation of such a polynomial under the true emission-center distribution is therefore the probability assigned to the chosen center. Comparing this expectation between the fitted and true models requires a bound on the polynomial's magnitude, which involves products of inverse distances between the centers. Within each cluster, these factors can contribute powers of $\Delta_{\mathrm{sep}}^{-1}$. Between clusters, they are bounded by $D_0^{-1}$ and affect only the constants in the bounds. This distinction allows the separation exponents to be expressed in terms of $s_{\max}$, although the constants can still depend on the total number of states $K$. When $M=1$, all $K$ true centers belong to a single cluster of diameter at most $C_0\Delta_{\mathrm{sep}}$, and $s_{\max}=K$. The condition on distances between distinct clusters is then vacuous.
\paragraph{Unstructured regime.} We also consider an unstructured regime, in which no cluster conditions are imposed beyond the boundedness and distinctness assumptions of Section~\ref{prelim:hmm}. In this regime, all the main inverse bounds use $s_{\max}=K$. Unlike the single-cluster regime, this convention does not require the full set of centers to have diameter $\mathcal O(\Delta_{\mathrm{sep}})$.

\subsection{Finite-Sample Density Estimation}
\label{sec:density}
To quantify the observation-density error, we use the Hellinger distance. For probability densities $p,q$ on a common observation space, it is defined by
\begin{equation}
h^2(p,q)
:=
\frac{1}{2}\int
\bigl(\sqrt{p(x)}-\sqrt{q(x)}\bigr)^2\,dx
=
1-\int\sqrt{p(x)q(x)}\,dx.
\label{eq:hellinger-definition}
\end{equation}
Under this convention, $0\le h(p,q)\le1$.

The inverse inequalities in the following sections relate latent recovery errors to discrepancies between observation densities, independently of how the fitted model is obtained. They therefore yield statistical recovery guarantees for any estimator with a suitable observation-density error bound. Let $\widetilde\theta_n$ denote any data-dependent fitted parameter and suppose that, for some radius $e_{n,T}(x)$,
\begin{equation}
\mathbb P_{\theta^\star}\!\left\{
 h(p_{\widetilde\theta_n,T},p_{\theta^\star,T})
 \le e_{n,T}(x)
\right\}\ge 1-e^{-x}.
\label{eq:generic-density-event}
\end{equation}
Then every applicable deterministic result in Sections~\ref{sec:exact_specification} and~\ref{sec:over_specification} immediately becomes a statistical guarantee after replacing $h(p_{\theta,T},p_{\theta^\star,T})$ by $e_{n,T}(x)$.

Throughout the paper, we obtain explicit statistical recovery rates by applying this principle to the global MLE. We first bound its observation-density error over the full fitted model class and then use the inverse inequalities to recover latent quantities. This density-estimation step does not require separated emission centers; the dependence on emission geometry enters through the subsequent inverse bounds.

\paragraph{Maximum likelihood estimation (MLE).} Given the $n$ independent trajectories introduced in Section~\ref{prelim:hmm}, we define the MLE $\widehat{\theta}_n$ over $\Theta_k(R)$ by
\[
\widehat{\theta}_n\in\arg\max_{\theta\in\Theta_k(R)}\ell_n(\theta),
\qquad
\ell_n(\theta):=\sum_{\ell=1}^n
\log p_{\theta,T}\bigl(X_{1:T}^{(\ell)}\bigr).
\]
We now establish a finite-sample Hellinger bound for this estimator
in the following proposition.

\begin{proposition}[Finite-sample Hellinger bound]
\label{prop:density}
There exists a constant $C>0$, depending only on $k$ and $R$, such that, for every $n\geq2$, the radius $\varepsilon_{n,T}:=C\sqrt{d\log(nTd)/n}$ satisfies
\[
\mathbb{P}_{\theta^\star}\left\{h\bigl(p_{\widehat{\theta}_n,T},p_{\theta^\star,T}\bigr)>\varepsilon_{n,T}\right\}\leq\frac{1}{n}.
\]
The bound holds for any measurable choice of a global maximizer, including when the maximizer is not unique.
\end{proposition}

The proof is given in Appendix~\ref{app:likelihood}. For fixed $k,d,T,R$, Proposition~\ref{prop:density} gives an $n^{-1/2}$ rate up to a logarithmic factor, with probability at least $1-1/n$. The density-estimation bound is independent of emission separation and remains valid at parameter configurations with singular Fisher information. Repeated fitted centers and zero fitted initial or transition probabilities are also allowed. The dependence on emission separation enters through the inverse inequalities used to recover latent quantities. We apply these inequalities on the common event $\mathcal E_n:=\{h(p_{\widehat{\theta}_n,T},p_{\theta^\star,T})\leq\varepsilon_{n,T}\}$, which has probability at least $1-1/n$. On this event, substituting $\varepsilon_{n,T}$ for the observation-density error gives the corresponding MLE recovery bounds simultaneously for all applicable latent losses.

Throughout the following inverse inequalities and MLE rate propositions, factors depending only on $K,k,R$ are absorbed into the constants. These constants may also depend on the cluster parameters in the structured regimes. Unless dependence on $\pi_{\min}^\star$ is displayed or explicitly permitted, they do not depend on $n,d,T,\Delta_{\mathrm{sep}},\pi_{\min}^\star$. A subscript $r$ additionally permits dependence on the number of selected time coordinates. For a sequence of random variables $Y_n$ and a positive deterministic sequence $a_n$, we write $Y_n=\mathcal{O}_{\mathbb{P}}(a_n)$ if $Y_n/a_n$ is bounded in probability.

The recovery bounds below also apply when emission separation decreases as the sample size grows. Smaller separation increases the recovery bounds, so the number of independent trajectories must grow fast enough for these upper bounds to tend to zero. This guarantees consistent latent recovery provided that the constants remain uniformly bounded.

\section{Exact Specification}
\label{sec:exact_specification}
In this section, we study latent recovery under exact specification, where the fitted and true models have the same number of states. However, equal state counts alone do not ensure that each true emission center has a fitted center nearby, since several fitted centers can be close to the same true center, while another true center has none nearby. To cover these configurations, we first establish inverse inequalities that apply to all fitted parameters in $\Theta_K(R)$, which we consider as \emph{global inverse bounds}. Using these results, we show that sufficiently small observation-density error guarantees a one-to-one matching between fitted and true centers. This matching allows us to obtain recovery bounds with constants independent of the minimum true initial probability in Section~\ref{sec:sharper_bound_exact_fit}. Combined with Proposition~\ref{prop:density}, these inverse bounds provide convergence rates for the maximum likelihood estimator in exact-specified setting.
\subsection{Separation-Dependent Bounds and Center Matching}
\label{sec:exact_specification_global}
To quantify latent recovery error, we consider both emission-center displacement and discrepancies in state probabilities and transition flows. Indeed, proximity of fitted centers to the true centers does not ensure that the probabilities assigned to these locations are accurate. As a result, location bounds must be accompanied by control of probability errors for individual states and state sequences. In addition to these separate location and probability losses, we consider Wasserstein distance to jointly assess errors in emission-center locations and their associated probabilities. The following theorem bounds these quantities in terms of the Hellinger distance between the fitted and true observation densities for a length-$T$ sequence.

\begin{theorem}[Global inverse bounds under exact specification]\label{thm:globalexact}
Let $\theta\in\Theta_K(R)$. In the structured regimes,
\begin{align}
D_\mu
&\le C(\pi_{\min}^\star)\Delta_{\mathrm{sep}}^{-(2s_{\max}-2)}h(p_{\theta,T},p_{\theta^\star,T}),
\label{eq:globalexact-location}\\
D_\pi\vee\max_{t\le T}D_{\{t\}}
&\le C(\pi_{\min}^\star)\Delta_{\mathrm{sep}}^{-(2s_{\max}-1)}h(p_{\theta,T},p_{\theta^\star,T}),
\label{eq:globalexact-mass}\\
D_A
&\le C(\pi_{\min}^\star)\Delta_{\mathrm{sep}}^{-(2s_{\max}-1)}h(p_{\theta,T},p_{\theta^\star,T}),
\label{eq:globalexact-flow}\\
\max_{t\le T}W_1(G_\theta^{\{t\}},G_\star^{\{t\}})
&\le C(\pi_{\min}^\star)\Delta_{\mathrm{sep}}^{-(2s_{\max}-2)}h(p_{\theta,T},p_{\theta^\star,T}).
\label{eq:globalexact-state-wasserstein}
\end{align}
More generally, for every nonempty $I\subseteq[T]$ with $|I|=r$,
\begin{equation}
D_I\le C_r(\pi_{\min}^\star)\Delta_{\mathrm{sep}}^{-(2s_{\max}-1)}h(p_{\theta,T},p_{\theta^\star,T}).
\label{eq:globalexact-joint}
\end{equation}
The displacement and one-time probability bounds hold without cluster assumptions with $s_{\max}=K$, and
\begin{align}
\max_{t\le T}W_1(G_\theta^{\{t\}},G_\star^{\{t\}})
&\le C(\pi_{\min}^\star)\Delta_{\mathrm{sep}}^{-(2K-2)}h(p_{\theta,T},p_{\theta^\star,T}),
\label{eq:globalexact-state-wasserstein-unstructured}\\
D_I&\le C_r(\pi_{\min}^\star)\Delta_{\mathrm{sep}}^{-(2K-1)}h(p_{\theta,T},p_{\theta^\star,T}),\qquad r\ge2.
\label{eq:globalexact-joint-unstructured}
\end{align}
The latter inequality also bounds $D_A$. All inequalities hold for every fitted parameter in $\Theta_K(R)$, with constants allowed to depend on the fixed positive quantity $\pi_{\min}^\star$.
\end{theorem}

The proof is given in Appendix~\ref{app:proof_3_1}. Equation~\eqref{eq:globalexact-location} controls the probability-weighted displacement of fitted centers, while equation~\eqref{eq:globalexact-mass} controls the probabilities assigned to the corresponding nearest-center groups. As $B\le D_\mu$, the displacement bound also controls the net displacement $B$ with the same separation exponent. To account for dependence between consecutive states, equation~\eqref{eq:globalexact-flow} further bounds the aggregate transition flows. These comparisons use the aggregated quantities defined in Section~\ref{sec:latent_recovery_loss}, so they remain valid before any one-to-one matching of fitted and true centers. In particular, several fitted states may be assigned to the same true state, and fitted states with zero probability are allowed.

\paragraph{Separation exponents and dependence on the smallest initial-state probability.}
In the structured regimes, the state-probability errors $D_\pi$, $D_{\{t\}}$, and $D_I$ and the transition-flow error $D_A$ have separation exponent $2s_{\max}-1$, regardless of the number of selected times. The exponent for weighted emission-center displacement and one-time Wasserstein distance is $2s_{\max}-2$. Without cluster assumptions, the corresponding exponents are obtained by setting $s_{\max}=K$. Thus, the bounds distinguish the difficulty of recovering probability masses from that of recovering their locations, while accounting for the number of centers that approach one another.

The constants in Theorem~\ref{thm:globalexact} may depend on $\pi_{\min}^\star$. For emission-center displacement and one-time Wasserstein distance, we can make this dependence explicit by increasing the separation exponent. More concretely, the resulting bound in the structured regimes is
\begin{equation}
D_\mu\vee\max_{t\le T}W_1(G_\theta^{\{t\}},G_\star^{\{t\}})
\le\frac{C}{\pi_{\min}^\star}\Delta_{\mathrm{sep}}^{-(2s_{\max}-1)}h(p_{\theta,T},p_{\theta^\star,T}).
\label{eq:globalexact-explicit}
\end{equation}
Compared with Theorem~\ref{thm:globalexact}, equation~\eqref{eq:globalexact-explicit} shows explicitly how small initial-state probabilities affect the recovery bound, at the cost of one additional inverse power of separation. Without cluster assumptions, the displacement bound remains valid with $s_{\max}=K$, while the one-time Wasserstein bound becomes
\[
\max_{t\le T}W_1(G_\theta^{\{t\}},G_\star^{\{t\}})
\le C\left\{(\pi_{\min}^\star)^{-1}\Delta_{\mathrm{sep}}^{-(2K-1)}
+\Delta_{\mathrm{sep}}^{-2K}\right\}h(p_{\theta,T},p_{\theta^\star,T}).
\]
Both alternatives are established in the proof of Theorem~\ref{thm:globalexact}. In the unstructured Wasserstein bound, the additional term $\Delta_{\mathrm{sep}}^{-2K}$ produces a loss of two separation powers relative to the sharp bound in that theorem. For displacement and structured transport, the loss is one power. These explicit alternatives are useful when the initial probabilities vary, since they quantify the effect of small $\pi_{\min}^\star$. For fixed positive $\pi_{\min}^\star$, the bounds of Theorem~\ref{thm:globalexact} have the better separation exponents. 

Applying the inverse bounds of Theorem~\ref{thm:globalexact} to the Hellinger error bound in Proposition~\ref{prop:density}, we obtain the following finite-sample rates for latent recovery by the global maximum likelihood estimator.

\begin{proposition}[Global MLE recovery rates]
\label{prop:global-mle-rates}
Let $\widehat{\theta}_n$ be a measurable global MLE over $\Theta_K(R)$. All losses below are evaluated at $(\widehat{\theta}_n,\theta^\star)$. In the structured regimes,
\begin{align*}
D_\mu &= \mathcal{O}_{\mathbb{P}}\left(\Delta_{\mathrm{sep}}^{-(2s_{\max}-2)}\sqrt{\frac{d\log(nTd)}{n}}\right),\\
D_\pi\vee\max_{t\leq T}D_{\{t\}} &= \mathcal{O}_{\mathbb{P}}\left(\Delta_{\mathrm{sep}}^{-(2s_{\max}-1)}\sqrt{\frac{d\log(nTd)}{n}}\right),\\
D_A &= \mathcal{O}_{\mathbb{P}}\left(\Delta_{\mathrm{sep}}^{-(2s_{\max}-1)}\sqrt{\frac{d\log(nTd)}{n}}\right),\\
\max_{t\leq T}W_1(G_{\widehat{\theta}_n}^{\{t\}},G_\star^{\{t\}}) &= \mathcal{O}_{\mathbb{P}}\left(\Delta_{\mathrm{sep}}^{-(2s_{\max}-2)}\sqrt{\frac{d\log(nTd)}{n}}\right).
\end{align*}
The same probability rate holds for $D_I$, with a constant depending on $r=|I|$. Without cluster assumptions, all these rates hold with $s_{\max}=K$. All constants may depend on $\pi_{\min}^\star$. These global rates require no sample-size condition for matching.
\end{proposition}

The Wasserstein rates above describe recovery of the emission-center distribution at each individual time. However, these marginal guarantees alone do not quantify recovery of dependence across time. To address this limitation, the following corollary extends the Wasserstein bounds to the joint emission-center distributions $G_\theta^I$ and $G_\star^I$, from which we derive the corresponding MLE rates.

\begin{corollary}[Global Wasserstein bounds for emission-center distributions]
\label{cor:globalpath}
Let $\theta\in\Theta_K(R)$. For nonempty $I\subseteq[T]$ with $|I|=r$, the structured regimes satisfy
\begin{equation}\label{eq:globalexact-path-wasserstein}
W_{1,I}(G_\theta^I,G_\star^I)
\le C_r(\pi_{\min}^\star)\Delta_{\mathrm{sep}}^{-(2s_{\max}-2)}h(p_{\theta,T},p_{\theta^\star,T}).
\end{equation}
In the unstructured regime, for every nonempty $I$,
\begin{equation}\label{eq:globalexact-path-wasserstein-unstructured}
W_{1,I}(G_\theta^I,G_\star^I)
\le C_r(\pi_{\min}^\star)\Delta_{\mathrm{sep}}^{-(2K-2)}h(p_{\theta,T},p_{\theta^\star,T}).
\end{equation}
These inequalities do not require one-to-one matching between fitted and true emission centers.
\end{corollary}
The proof is given in Appendix~\ref{app:proof_cor_global_path}. Compared with the one-time Wasserstein bounds in Theorem~\ref{thm:globalexact}, Corollary~\ref{cor:globalpath} controls the joint emission-center distributions with the same separation exponents. Consequently, the bounds account for dependence across several selected times without introducing additional inverse powers of separation. The constants may still depend on the number of selected times. The transport interpretation below explains why controlling the joint distribution of latent emission centers does not increase the separation exponent under exact specification.

\paragraph{Transport within and between cluster sequences.}
The distinction between the structured and unstructured bounds depends on both the transport distances and the probability mass moved over those distances. To describe this distinction, a cluster sequence records the cluster containing the center at each selected time. Two location sequences therefore share a cluster sequence when their corresponding centers belong to the same cluster at every selected time. For example, if $a_1,a_2$ belong to one cluster and $a_3,a_4$ to another, then $(a_1,a_3)$ and $(a_2,a_4)$ share a cluster sequence.

\begin{enumerate}
\renewcommand{\labelenumi}{(\roman{enumi})}

\item \emph{Within a cluster sequence.} Transport between location sequences that share the same cluster sequence involves shorter distances than transport between different cluster sequences, so we bound these two contributions separately. After fitted locations are moved to their assigned true centers, redistributing probability within a cluster sequence costs at most $C_0\Delta_{\mathrm{sep}}$ per unit mass. This factor offsets one inverse power of separation in the probability bound~\eqref{eq:globalexact-joint}, giving Wasserstein exponent $2s_{\max}-2$.

\item \emph{Between cluster sequences.} Transport between different cluster sequences can involve distances that do not shrink with $\Delta_{\mathrm{sep}}$. The transport bound therefore relies on controlling how much probability must move between them. This comparison considers the total probability of each cluster sequence, without distinguishing individual centers within the same cluster. Since different clusters remain separated by a fixed gap, the proof can control these totals more accurately than individual state-sequence probabilities. More precisely, Lemma~\ref{lem:cluster-mass} bounds the discrepancy between cluster-sequence probabilities by $C_r\Delta_{\mathrm{sep}}^{-(2s_{\max}-2)}h(p_{\theta,T},p_{\theta^\star,T})$. Since each unit of this discrepant mass can be transported over a distance of at most $2R$, the resulting transport cost retains the same separation exponent.

\item \emph{Unstructured regime.} Without cluster assumptions, distances between true centers need not shrink with $\Delta_{\mathrm{sep}}$, so the preceding within-cluster argument no longer applies. The largest possible transport distance is instead the full support diameter $D_\star:=\max_{i,j}\|a_i-a_j\|_2$. A direct bound using this diameter and the probability inequality~\eqref{eq:globalexact-joint-unstructured} would retain separation exponent $2K-1$, rather than the desired $2K-2$. The proof avoids this loss through a stronger probability bound that accounts for the actual distances between centers. More precisely, the proof of Corollary~\ref{cor:globalpath} gives
\[
D_I\le C_r(\pi_{\min}^\star)D_\star^{-2}\Delta_{\mathrm{sep}}^{-(2K-3)}h(p_{\theta,T},p_{\theta^\star,T}).
\]
After fitted locations are moved to their assigned true centers, the remaining probability mass can be redistributed at a cost of at most $D_\star$ per unit mass. The resulting cost is therefore bounded by $D_\star D_I/2$, whose separation-dependent factor is $D_\star^{-1}\Delta_{\mathrm{sep}}^{-(2K-3)}$. Because $D_\star\ge\Delta_{\mathrm{sep}}$, this factor is at most $\Delta_{\mathrm{sep}}^{-(2K-2)}$. The redistribution cost thus has the same separation exponent as the emission-center displacement, giving the joint Wasserstein bound without a cluster-diameter assumption.

\end{enumerate}

Combining Corollary~\ref{cor:globalpath} with the density-estimation bound in Proposition~\ref{prop:density} gives the following recovery rates for joint emission-center distributions.

\begin{proposition}[Global MLE Wasserstein rates]
\label{prop:global-mle-wasserstein}
Let $\widehat{\theta}_n$ be a measurable global MLE over $\Theta_K(R)$. For any nonempty $I\subseteq[T]$, let $r=|I|$. In the structured regimes,
\[
W_{1,I}(G_{\widehat{\theta}_n}^I,G_\star^I)=\mathcal{O}_{\mathbb{P}}\left(\Delta_{\mathrm{sep}}^{-(2s_{\max}-2)}\sqrt{\frac{d\log(nTd)}{n}}\right).
\]
In the unstructured regime, the separation exponent is $2K-2$ for every nonempty $I$. Constants may depend on $\pi_{\min}^\star$ and $r$.
\end{proposition}
Because the same full-block Hellinger error controls every nonempty $I\subseteq[T]$, the recovery inequalities hold simultaneously for all these choices of $I$ on the event in Proposition~\ref{prop:density}. These global guarantees allow several fitted centers to be associated with the same true center. They therefore establish recovery without first requiring a one-to-one correspondence between fitted and true states. To obtain such a correspondence, we next impose a sufficiently small observation-density error. The following corollary shows that this condition places exactly one fitted center within distance $\Delta_{\mathrm{sep}}/4$ of each true emission center. Whenever this correspondence holds, we write $b_j$ for the fitted center relabeled to match $a_j$.

\begin{corollary}[Density-error threshold for emission-center matching]
\label{cor:localization}
For $\theta\in\Theta_K(R)$, there exists $c(\pi_{\min}^\star)>0$ such that, if
\begin{equation}\label{eq:exact-localization-threshold}
h(p_{\theta,T},p_{\theta^\star,T})\le c(\pi_{\min}^\star)\Delta_{\mathrm{sep}}^{2s_{\max}-1},
\end{equation}
then a unique permutation pairs every fitted center $b_j$ with a true center $a_j$ such that $\norm{b_j-a_j}\le\Delta_{\mathrm{sep}}/4$. After relabeling, $\pi_j\ge\pi_j^\star/2$ for every $j$. In the unstructured regime, take $s_{\max}=K$.
\end{corollary}

The proof of Corollary~\ref{cor:localization} is given in Appendix~\ref{app:proof_corollary3.3}. The key step is to show that the density-error condition ensures a fitted center within distance $\Delta_{\mathrm{sep}}/4$ of every true center. Since these neighborhoods are disjoint and there are exactly $K$ fitted centers for $K$ true centers, each neighborhood contains exactly one fitted center, giving a unique matching. We use this matching to relabel the fitted states, their initial probabilities, and both indices of the transition matrix, so that fitted state $j$ corresponds to true state $j$. After relabeling, the corollary also gives $\pi_j\ge\pi_j^\star/2$. Thus, the density-error condition in equation~\eqref{eq:exact-localization-threshold} provides both the center matching and the probability lower bound needed for the subsequent matched-center recovery bounds.

\paragraph{Sample size for guaranteed matching.}
For the global maximum likelihood estimator, Corollary~\ref{cor:localization} turns the Hellinger error bound of Proposition~\ref{prop:density} into a sufficient sample-size condition for matching. Specifically, Proposition~\ref{prop:density} gives an event $\mathcal E_n:=\{h(p_{\widehat{\theta}_n,T},p_{\theta^\star,T})\le\varepsilon_{n,T}\}$ with probability at least $1-1/n$. On this event, the density-error threshold in Corollary~\ref{cor:localization} is satisfied whenever
\begin{equation}
\label{eq:exact-localization-sample-size}
n\ge C(\pi_{\min}^\star)\Delta_{\mathrm{sep}}^{-(4s_{\max}-2)}d\log(nTd).
\end{equation}
For $n\ge2$, this condition guarantees one-to-one emission-center matching with probability at least $1-1/n$. In particular, closer emission centers require more independent trajectories to guarantee matching, with the separation exponent determined by the largest cluster size. Without cluster assumptions, the same conclusion holds with $s_{\max}=K$.

The constant in~\eqref{eq:exact-localization-sample-size} may depend on $\pi_{\min}^\star$, so this condition does not specify how the sample size changes when some true initial-state probabilities become small. To quantify this dependence, we also give an alternative sufficient condition for the conclusions of Corollary~\ref{cor:localization}:
\begin{equation}
h(p_{\theta,T},p_{\theta^\star,T})\le c\pi_{\min}^\star\Delta_{\mathrm{sep}}^{2s_{\max}},
\label{eq:alternative-localization}
\end{equation}
On the same event $\mathcal E_n$, condition~\eqref{eq:alternative-localization} holds whenever
\[
n\ge \frac{C}{(\pi_{\min}^\star)^2}\Delta_{\mathrm{sep}}^{-4s_{\max}}d\log(nTd).
\]
This alternative condition shows that smaller minimum true initial-state probabilities require a larger sample size to guarantee matching, with an inverse-square dependence on $\pi_{\min}^\star$, but has two additional inverse powers of emission separation compared with~\eqref{eq:exact-localization-sample-size}. For fixed positive $\pi_{\min}^\star$, condition~\eqref{eq:exact-localization-sample-size} gives the better separation exponent. Without cluster assumptions, both alternative conditions hold with $s_{\max}=K$.

\subsection{Inverse Bounds under Emission-Center Matching}
\label{sec:sharper_bound_exact_fit}
Once fitted and true emission centers admit a one-to-one matching, we can compare the probabilities assigned to individual states and to sequences of states under a common labeling. Corollary~\ref{cor:localization} provides one sufficient condition for obtaining this correspondence. The following theorem gives recovery bounds under this matching, with $\|b_j-a_j\|_2\le\Delta_{\mathrm{sep}}/4$ for every $j$ after relabeling. Compared with the global bounds, the following result retains the same separation exponents for state and state-sequence probabilities, transition flows, probability-weighted emission-center displacement, and Wasserstein distances, while removing dependence on $\pi_{\min}^\star$ from their constants. Thus, once the centers are matched, the same observation-density error controls these latent discrepancies without additional amplification from small initial-state probabilities. As in the global bounds, recovery across more selected times does not increase the separation exponents.

\begin{theorem}[Paired-center inverse bounds]\label{thm:paired}
Suppose $\theta\in\Theta_K(R)$ is labeled so that $\norm{b_j-a_j}\le\Delta_{\mathrm{sep}}/4$ for every $j$. Then
\begin{align}
D_\mu&\le C\Delta_{\mathrm{sep}}^{-(2s_{\max}-2)}
h(p_{\theta,T},p_{\theta^\star,T}),\\
\max_{t\le T}D_{\{t\}}&\le C\Delta_{\mathrm{sep}}^{-(2s_{\max}-1)}
h(p_{\theta,T},p_{\theta^\star,T}). \label{eq:one-time-recovery}
\end{align}
In the structured regimes, for every $I\subseteq[T]$ with $|I|=r\ge1$,
\begin{align}
D_I&\le C_r\Delta_{\mathrm{sep}}^{-(2s_{\max}-1)}
h(p_{\theta,T},p_{\theta^\star,T}),
\label{eq:paired-path-improved}\\
W_{1,I}(G_\theta^I,G_\star^I)
&\le C_r\Delta_{\mathrm{sep}}^{-(2s_{\max}-2)}
h(p_{\theta,T},p_{\theta^\star,T}).
\label{eq:paired-transport-improved}
\end{align}
In particular, $D_A\le C\Delta_{\mathrm{sep}}^{-(2s_{\max}-1)}h(p_{\theta,T},p_{\theta^\star,T})$.
In the unstructured regime, the first two bounds hold with $s_{\max}=K$, and
\begin{align}
D_I&\le C_r\Delta_{\mathrm{sep}}^{-(2K-1)}
h(p_{\theta,T},p_{\theta^\star,T}),
\label{eq:paired-unstructured-probability}\\
\max_{t\le T}W_1(G_\theta^{\{t\}},G_\star^{\{t\}})
&\le C\Delta_{\mathrm{sep}}^{-(2K-2)}h(p_{\theta,T},p_{\theta^\star,T}),
\label{eq:paired-unstructured-marginal}\\
W_{1,I}(G_\theta^I,G_\star^I)
&\le C_r\Delta_{\mathrm{sep}}^{-(2K-2)}h(p_{\theta,T},p_{\theta^\star,T}),
\qquad r\ge2.
\label{eq:paired-unstructured-path}
\end{align}

The constants in these bounds are independent of $\pi_{\min}^\star$, and no lower bound on fitted state probabilities is required.
\end{theorem}
We defer the proof to Appendix~\ref{app:proof_paired_theorem}. Under emission-center matching, we bound $D_\mu$ directly with a constant independent of $\pi_{\min}^\star$. Since this loss already weights each center displacement by its fitted state probability, the bound does not require division by individual state probabilities. For joint state probabilities, the proof uses the Markov structure to reduce the comparison across all selected times to comparisons at successive selected times. These pairs need not be consecutive in the original trajectory. Under the common labeling supplied by matching, Lemma~\ref{lem:markov-path-comparison} bounds the joint probability discrepancy by a sum of the corresponding two-time discrepancies. As each two-time discrepancy has separation exponent $2s_{\max}-1$, the joint bound retains this exponent, while the number of selected times affects only the multiplicative constant. 

\paragraph{Sharpness of the structured separation exponents.}
Theorem~\ref{thm:paired} attains the sharp separation exponents for recovery under emission-center matching over structured classes containing collinear clusters of size $s_{\max}\ge2$. Specifically, the exponents are $2s_{\max}-1$ for state-probability and transition-flow recovery and $2s_{\max}-2$ for weighted emission-center displacement and path Wasserstein recovery. Proposition~\ref{prop:matched-sharpness} establishes their sharpness through examples satisfying the matching condition with true initial-state probabilities bounded away from zero. The construction uses two mixing distributions with identical low-order moments but different mixing weights and center locations. These moment equalities cancel the corresponding low-order terms in the Gaussian observation densities, making the observation error small relative to the latent discrepancies. Each mixing distribution defines an HMM whose initial probabilities and transition rows all equal its mixture weights, so the same comparison holds for HMMs with independent hidden states. These examples also establish sharpness of the global structured bounds in Theorem~\ref{thm:globalexact} for fixed positive $\pi_{\min}^\star$. The separation exponents are therefore optimal in the worst case over these model classes.

\begin{table}[t]
\centering
\setlength{\tabcolsep}{8pt}
\renewcommand{\arraystretch}{1.6}
\setlength{\extrarowheight}{3pt}
\arrayrulecolor{black}
\begin{tabular*}{\textwidth}{@{\extracolsep{\fill}}lcc@{}}
\toprule
True emission geometry
& \shortstack{State-probability and\\transition-flow errors}
& \shortstack{Weighted center-displacement\\and Wasserstein errors} \\
\midrule[0.7pt]
Single cluster
& $\Delta_{\mathrm{sep}}^{-(2K-1)}\sqrt{\frac{d\log(nTd)}{n}}$
& $\Delta_{\mathrm{sep}}^{-(2K-2)}\sqrt{\frac{d\log(nTd)}{n}}$ \\
\arrayrulecolor{black!12}\midrule[0.25pt]\arrayrulecolor{black}
Separated clusters
& $\Delta_{\mathrm{sep}}^{-(2s_{\max}-1)}\sqrt{\frac{d\log(nTd)}{n}}$
& $\Delta_{\mathrm{sep}}^{-(2s_{\max}-2)}\sqrt{\frac{d\log(nTd)}{n}}$ \\
\arrayrulecolor{black!12}\midrule[0.25pt]\arrayrulecolor{black}
Unstructured
& $\Delta_{\mathrm{sep}}^{-(2K-1)}\sqrt{\frac{d\log(nTd)}{n}}$
& $\Delta_{\mathrm{sep}}^{-(2K-2)}\sqrt{\frac{d\log(nTd)}{n}}$ \\
\bottomrule
\end{tabular*}
\caption{\textbf{Exact-fit MLE rates.} Each entry gives an $\mathcal{O}_{\mathbb P}$ upper rate. The second column covers $D_\pi$, $D_{\{t\}}$, $D_I$, and $D_A$; the third covers $D_\mu$ and $W_{1,I}$. Here, $K$ is the number of true states and $s_{\max}$ is the largest cluster size.}
\label{tab:exactfit}
\end{table}

The following proposition states the paired-center MLE rates. We use the sample-size condition in equation~\eqref{eq:exact-localization-sample-size} to ensure that the fitted centers satisfy the matching hypothesis of Theorem~\ref{thm:paired} on the density event of Proposition~\ref{prop:density}.

\begin{proposition}[MLE recovery rates under emission-center matching]
\label{prop:mle-matched-rates}
Let $\widehat{\theta}_n$ be a measurable global MLE over $\Theta_K(R)$, and suppose that the sample-size condition in \eqref{eq:exact-localization-sample-size} holds. After relabeling the fitted states according to the matching in Corollary~\ref{cor:localization}, the following rates hold as $n\to\infty$. All losses are evaluated at $(\widehat{\theta}_n,\theta^\star)$:
\begin{align*}
D_\mu &= \mathcal{O}_{\mathbb{P}}\left(\Delta_{\mathrm{sep}}^{-(2s_{\max}-2)}\sqrt{\frac{d\log(nTd)}{n}}\right),\\
D_\pi\vee\max_{t\leq T}D_{\{t\}} &= \mathcal{O}_{\mathbb{P}}\left(\Delta_{\mathrm{sep}}^{-(2s_{\max}-1)}\sqrt{\frac{d\log(nTd)}{n}}\right).
\end{align*}

For structured regimes,
\begin{align*}
D_A &= \mathcal{O}_{\mathbb{P}}\left(\Delta_{\mathrm{sep}}^{-(2s_{\max}-1)}\sqrt{\frac{d\log(nTd)}{n}}\right),\\
D_I &= \mathcal{O}_{\mathbb{P}}\left(C_r\Delta_{\mathrm{sep}}^{-(2s_{\max}-1)}\sqrt{\frac{d\log(nTd)}{n}}\right),\qquad |I|=r,\\
W_{1,I}(G_{\widehat{\theta}_n}^I,G_\star^I)
&=\mathcal{O}_{\mathbb{P}}\left(C_r\Delta_{\mathrm{sep}}^{-(2s_{\max}-2)}\sqrt{\frac{d\log(nTd)}{n}}\right).
\end{align*}
In the unstructured regime, the sample-size condition and the first two rates use $s_{\max}=K$, while
\begin{align*}
D_A &= \mathcal{O}_{\mathbb{P}}\left(\Delta_{\mathrm{sep}}^{-(2K-1)}\sqrt{\frac{d\log(nTd)}{n}}\right),\\
D_I &= \mathcal{O}_{\mathbb{P}}\left(C_r\Delta_{\mathrm{sep}}^{-(2K-1)}\sqrt{\frac{d\log(nTd)}{n}}\right),\\
\max_{t\le T}W_1(G_{\widehat{\theta}_n}^{\{t\}},G_\star^{\{t\}})
&=\mathcal{O}_{\mathbb{P}}\left(\Delta_{\mathrm{sep}}^{-(2K-2)}\sqrt{\frac{d\log(nTd)}{n}}\right),\\
W_{1,I}(G_{\widehat{\theta}_n}^I,G_\star^I)
&=\mathcal{O}_{\mathbb{P}}\left(C_r\Delta_{\mathrm{sep}}^{-(2K-2)}\sqrt{\frac{d\log(nTd)}{n}}\right),\qquad r\ge2.
\end{align*}
These rates hold simultaneously for all nonempty $I\subseteq[T]$.
\end{proposition}

We defer the proof to Appendix~\ref{app:corollary_3_5_3_6_proof}. Compared with the global rates in Proposition~\ref{prop:global-mle-rates}, these rates retain the same dependence on sample size and emission separation. Once the centers are matched, the constants in the underlying probability-weighted inverse bounds are independent of $\pi_{\min}^\star$. Table~\ref{tab:exactfit} summarizes how these probability and location recovery rates depend on the emission geometry.

\subsection{Emission-Center and Transition-Row Recovery}
The bounds in Section~\ref{sec:sharper_bound_exact_fit} measure emission-center displacement and transition-flow error through probability-weighted quantities. A state with small probability can therefore contribute little to these losses even when its emission center or transition probabilities are poorly estimated. To assess recovery of every state, we instead consider the largest emission-center error and the largest transition-row error. Using the labeling established in Corollary~\ref{cor:localization}, the following corollary bounds both errors in terms of the observation-density error.

\begin{corollary}[Emission-center and transition-row bounds]\label{cor:rows}
Under the localization condition of Corollary~\ref{cor:localization}, use its aligned labeling. In the structured regimes,
\begin{align}
\max_j\norm{b_j-a_j} &\le \frac{C}{\pi_{\min}^\star}\Delta_{\mathrm{sep}}^{-(2s_{\max}-2)}h(p_{\theta,T},p_{\theta^\star,T}),\\
\max_j\norm{A_{j\cdot}-A^\star_{j\cdot}}_1
&\le \frac{C}{\pi_{\min}^\star}\Delta_{\mathrm{sep}}^{-(2s_{\max}-1)}h(p_{\theta,T},p_{\theta^\star,T})
\end{align}
Without cluster assumptions, both bounds hold with $s_{\max}=K$.
\end{corollary}

The proof of Corollary~\ref{cor:rows} is deferred to Appendix~\ref{app:corollary_3_5_3_6_proof}. The corollary shows that controlling every center and transition row incurs a dependence on $(\pi_{\min}^\star)^{-1}$ without increasing the separation exponents of Theorem~\ref{thm:paired}. This additional dependence arises because passing from probability-weighted errors to errors for each state requires division by the corresponding initial probability. For centers, this division removes the displacement weight. For transition rows, division by the initial probability requires control of both flow and initial-probability errors. In either case, smaller initial probabilities amplify the resulting error bounds, so the guarantee over all states depends on $\pi_{\min}^\star$. Since recovering transition rows requires only normalization by initial probabilities, neither invertibility of the transition matrix nor a positive lower bound on its individual entries is required.

\noindent For the global maximum likelihood estimator, these bounds result in the following rates for recovering individual emission centers and transition rows.

\begin{proposition}[MLE rates for emission centers and transition rows]
\label{prop:mle-rows}
Let $\widehat{\theta}_n$ be a measurable global MLE over
$\Theta_K(R)$, and suppose that the sample-size condition in \eqref{eq:exact-localization-sample-size} holds. After relabeling the fitted states according to the matching in
Corollary~\ref{cor:localization},
\begin{align*}
\max_{j\in[K]}\|\widehat{\mu}_j-\mu_j^\star\|_2
&=
\mathcal{O}_{\mathbb{P}}\left(
(
\pi_{\min}^\star)^{-1}\Delta_{\mathrm{sep}}^{-(2s_{\max}-2)}
\sqrt{\frac{d\log(nTd)}{n}}
\right),\\
\max_{j\in[K]}\|\widehat{A}_{j\cdot}-A_{j\cdot}^\star\|_1
&=
\mathcal{O}_{\mathbb{P}}\left(
(
\pi_{\min}^\star)^{-1}\Delta_{\mathrm{sep}}^{-(2s_{\max}-1)}
\sqrt{\frac{d\log(nTd)}{n}}
\right).
\end{align*}
Without cluster assumptions, the sample-size condition and both rates hold with $s_{\max}=K$.
\end{proposition}

The result follows by applying Corollary~\ref{cor:rows} on the density event of Proposition~\ref{prop:density}, with emission-center matching ensured by~\eqref{eq:exact-localization-sample-size}.
\begin{remark}[Two observations already identify the exact-fit model]
\label{rem:two-observation-identifiability}
At the population level, two consecutive observations suffice to identify an exactly specified $K$-state HMM under our Gaussian location assumptions. Specifically, if two exact-fit $K$-state HMMs satisfy $p_{\theta,2}=p_{\theta^\star,2},$
then one common permutation aligns all emission centers, initial probabilities, and transition matrices. No invertibility assumption on the transition matrix is required. The formal statement is Corollary~\ref{cor:identifiability-two}, whose proof in Appendix~\ref{app:corollary_3_5_3_6_proof} applies the localization and matched one- and two-time inverse bounds at zero Hellinger distance. Thus, two observation coordinates already contain enough information for exact identification in the present model. The inverse inequalities quantify the stability of that identification as emission centers approach one another.
\end{remark}

\section{Over-Specification}
\label{sec:over_specification}
In this section, we study latent recovery when the fitted HMM has more states than the true model. Several fitted states may then represent the same true state, so a comparison of individual fitted and true parameters is no longer sufficient to describe what the observation distribution determines. We therefore consider recovery after grouping fitted states associated with the same true emission center. Our analysis focuses on the following objectives:

\begin{enumerate}
\item[\textit{(i)}] \emph{Recover aggregate state probabilities and transition flows.} Splitting a true state into several fitted states with the same emission center can preserve the observation distribution while changing their individual initial and transition probabilities. We therefore compare the total fitted probability assigned to each true center, or each sequence of true centers, with its true counterpart. Unlike the exactly specified case, the grouped process $c(Z_t)$ need not be Markov, as illustrated in Example~\ref{ex:aggregate-memory}. Consequently, we bound aggregate state-sequence probabilities directly from the observation density rather than applying the Markov comparison through two-time probabilities used in Section~\ref{sec:sharper_bound_exact_fit}.

\item[\textit{(ii)}] \emph{Control emission-center displacement without losing opposing shifts.} Fitted centers assigned to the same true center may move in opposite directions, so their contributions to the net displacement $B$ can cancel. We also control the quadratic displacement $V$, which records squared distances without cancellation, and recover a bound for $D_\mu$ through $D_\mu\le\sqrt{V}$. Section~\ref{sec:overfit_global} establishes bounds for $B$, $V$, and the aggregate probabilities considered in~(i).

\item[\textit{(iii)}] \emph{Bound Wasserstein errors and derive MLE recovery rates.} Section~\ref{sec:overfit_latent_path_transport} combines the displacement and aggregate-probability bounds to control Wasserstein distances between fitted and true emission-center distributions. Applying the Hellinger error bound of Proposition~\ref{prop:density} then gives recovery rates for the global maximum likelihood estimator.
\end{enumerate}

The following example illustrates the need for aggregation in (i) by showing that individual fitted transition rows can vary while the observation distribution remains unchanged.

\begin{example}[Arbitrary internal transitions]
\label{ex:state-splitting}
Split each true state $j$ into substates $(j,1),\ldots,(j,m_j)$ sharing the same emission center $a_j$.  For arbitrary probability vectors $r_{j\ell,u}(\cdot)$ on the substates of state $\ell$, set
\[
\widetilde A_{(j,u),(\ell,v)}=A^\star_{j\ell}r_{j\ell,u}(v).
\]
If the initial substate masses sum to $\pi_j^\star$ within each true state, then the fitted and true HMMs have identical joint observation densities for every sequence length. If $m_\ell\ge 2$ and $A^\star_{j\ell}>0$ for some $j,\ell$, the individual substate transition rows can vary continuously without changing these observation densities. Thus, individual split-state transition rows are not meaningful estimands without an additional identifying convention, whereas aggregate probabilities and transition flows remain meaningful recovery targets.
\end{example}
We present the detailed verification of this example in Appendix~\ref{app:state-splitting}. The construction shows that individual split-state transition rows cannot generally be recovered even if the observation distribution is known exactly. The relevant quantities are instead the total probabilities and transition flows associated with each true center, together with the distributions of latent emission-center sequences. Throughout this section, aggregation refers to summing fitted-state probabilities according to the assignment map.

\subsection{Recovery after Aggregating Fitted States}
\label{sec:overfit_global}
We begin by extending the global inverse bounds of
Section~\ref{sec:exact_specification_global} to over-specification. Two changes arise when several fitted states can correspond to one true state. First, the aggregate probability bounds have larger separation exponents than their exact-fit counterparts. Second, emission-center displacement can no longer be controlled linearly by the observation-density error in general. For instance, splitting a true center into two equally weighted fitted centers at distance $\varepsilon$ on opposite sides gives a displacement proportional to $\varepsilon$. Relative to the true emission density, the two opposite shifts produce first-order changes of opposite signs. These changes cancel when the shifted densities are averaged, leaving a Hellinger error proportional to $\varepsilon^2$.
This first-order cancellation prevents a bound for $D_\mu$ that is linear in $h(p_{\theta,T},p_{\theta^\star,T})$ from holding uniformly over all over-specified configurations. We instead
bound $V$ and $B$ directly and control $D_\mu$
through $D_\mu\le\sqrt{V}$. The following theorem bounds these displacement quantities and aggregate state-sequence probability errors by the full-block Hellinger error, uniformly over fitted parameters
in $\Theta_k(R)$.
\begin{theorem}[Overfitted location and aggregate-flow bounds]\label{thm:overglobal}
For any $k\ge K$ and $\theta\in\Theta_k(R)$, in the structured regimes, put $\beta_r=\max\{2s_{\max},r(s_{\max}-1)+s_{\max}\}$. Then
\begin{align}
V &\le C\Delta_{\mathrm{sep}}^{-(2s_{\max}-2)}h(p_{\theta,T},p_{\theta^\star,T}),
\label{eq:overglobal-V}\\
B &\le C\Delta_{\mathrm{sep}}^{-(2s_{\max}-1)}h(p_{\theta,T},p_{\theta^\star,T}),\\
D_I &\le C_r\Delta_{\mathrm{sep}}^{-\beta_r}h(p_{\theta,T},p_{\theta^\star,T}),
\qquad I\subseteq[T],\ |I|=r\ge1.
\label{eq:overglobal-DI}
\end{align}
In particular,
\[
D_\pi\le C\Delta_{\mathrm{sep}}^{-2s_{\max}}h(p_{\theta,T},p_{\theta^\star,T}),
\qquad
D_A\le C\Delta_{\mathrm{sep}}^{-\beta_2}h(p_{\theta,T},p_{\theta^\star,T}).
\]
Without a cluster assumption, the $V$ and $B$ bounds hold with $s_{\max}=K$, while $D_I$ has separation exponent
\[
\gamma_{K,r}:=\max\{2K,\ r(K-1)+K\}.
\]
In particular, $\gamma_{K,1}=2K$ and $\gamma_{K,2}=3K-2$. No lower bound on fitted state probabilities and no factor $(\pi_{\min}^\star)^{-1}$ is required.
\end{theorem}
The proof of Theorem~\ref{thm:overglobal} is given in Appendix~\ref{proof:global_overift}. For aggregate state-sequence probabilities, the theorem gives bounds linear in the observation-density error even when the probabilities of individual split states cannot be identified. These bounds compare the total probability assigned to each true state sequence, regardless of how that probability is divided among fitted states. Unlike the exact-fit bounds, however, the separation exponents depend on the number of selected times $r$. For two selected times and $s_{\max}\ge2$, the structured exponent is $\beta_2=3s_{\max}-2$, compared with $2s_{\max}$ for a single time. When $s_{\max}>2$, this larger exponent means that the guarantee for joint probabilities deteriorates faster as emission centers approach one another.
The exponent $\beta_r$ accounts for both the polynomial-moment comparison and the approximation error caused by displaced fitted centers. We first construct a polynomial that equals one at the target true-center sequence and zero at every other true-center sequence. Its expectation under the true model is therefore the probability of that sequence. Lemma~\ref{lem:moment-extraction} then relates the difference between its fitted and true expectations to the observation Hellinger error, with a multiplier determined by the polynomial's magnitude. The bound on this magnitude contributes the separation exponent $(2s_{\max}-1)+(r-1)(s_{\max}-1)=r(s_{\max}-1)+s_{\max}$ for $s_{\max}\ge2$. However, fitted centers may differ from their assigned true centers, so the fitted expectation only approximates the corresponding grouped probability. The polynomial is chosen to have zero first derivatives at the true centers, which makes this approximation error quadratic in the center displacements. The proof controls this second contribution with separation exponent $\beta_r=\max\{2s_{\max},r(s_{\max}-1)+s_{\max}\}$.

For emission-center recovery, however, the dependence on observation error changes. Since $D_\mu\le\sqrt{V}$, the quadratic displacement bound gives $D_\mu\le C\Delta_{\mathrm{sep}}^{-(s_{\max}-1)}\sqrt{h(p_{\theta,T},p_{\theta^\star,T})}$, with $C$ independent of $\pi_{\min}^\star$. In comparison, the aggregate probability bound is $D_I\le C_r\Delta_{\mathrm{sep}}^{-\beta_r}h(p_{\theta,T},p_{\theta^\star,T})$. For fixed emission geometry, the distinction is therefore between the powers $h^{1/2}$ for probability-weighted emission-center displacement and $h$ for aggregate state-sequence probabilities. If separation also varies, the corresponding factors $\Delta_{\mathrm{sep}}^{-(s_{\max}-1)}$ and $\Delta_{\mathrm{sep}}^{-\beta_r}$ must also be taken into account. This difference reflects the cancellation illustrated by symmetric splitting, where opposite center shifts can produce an observation error of quadratic order in their displacement. 

We next relate these deterministic comparisons into statistical recovery rates. Combining Theorem~\ref{thm:overglobal} with the Hellinger bound in Proposition~\ref{prop:density}, we obtain the following rates for the global maximum likelihood estimator.

\begin{proposition}[MLE recovery rates under over-specification]
\label{prop:mle-overglobal}
Let $\widehat{\theta}_n$ be a measurable global MLE over $\Theta_k(R)$, where $k\geq K$. All losses below are evaluated at $(\widehat{\theta}_n,\theta^\star)$. In the structured regimes,
\begin{align*}
V &= \mathcal{O}_{\mathbb{P}}\left(\Delta_{\mathrm{sep}}^{-(2s_{\max}-2)}\sqrt{\frac{d\log(nTd)}{n}}\right),\\
B &= \mathcal{O}_{\mathbb{P}}\left(\Delta_{\mathrm{sep}}^{-(2s_{\max}-1)}\sqrt{\frac{d\log(nTd)}{n}}\right),\\
D_\pi &= \mathcal{O}_{\mathbb{P}}\left(\Delta_{\mathrm{sep}}^{-2s_{\max}}\sqrt{\frac{d\log(nTd)}{n}}\right),\\
D_A &= \mathcal{O}_{\mathbb{P}}\left(\Delta_{\mathrm{sep}}^{-\beta_2}\sqrt{\frac{d\log(nTd)}{n}}\right),\\
D_I &= \mathcal{O}_{\mathbb{P}}\left(C_r\Delta_{\mathrm{sep}}^{-\beta_r}\sqrt{\frac{d\log(nTd)}{n}}\right),\qquad |I|=r,\\
D_\mu &= \mathcal{O}_{\mathbb{P}}\left(\Delta_{\mathrm{sep}}^{-(s_{\max}-1)}\left[\frac{d\log(nTd)}{n}\right]^{1/4}\right).
\end{align*}
The rates hold simultaneously for all nonempty $I\subseteq[T]$. Without cluster assumptions, take $s_{\max}=K$ in the $V,B,D_\pi,D_\mu$ rates, and use exponents $\gamma_{K,2}$ and $\gamma_{K,r}$ for $D_A$ and $D_I$, respectively. No lower bound on fitted state probabilities is required, and none of the rates contains a factor $(\pi_{\min}^\star)^{-1}$.
\end{proposition}

The first five rates follow by applying Theorem~\ref{thm:overglobal} on the density event of Proposition~\ref{prop:density}. For fixed model dimensions and emission geometry, these quantities, including the quadratic displacement $V$, have an $n^{-1/2}$ upper rate, up to logarithmic factors. Since $D_\mu\leq\sqrt{V}$, taking the square root of the bound for $V$ gives the $n^{-1/4}$ upper rate for probability-weighted emission-center displacement. This slower displacement rate reflects the symmetric splitting described above, where opposite center shifts can produce an observation-density error of quadratic order in their displacement.

\paragraph{Exact recovery of emission-center distributions.} At zero density error, the emission-center distributions are recovered exactly within the observed block. This recovery requires both the emission-center locations and the probabilities of the corresponding location sequences to equal their true counterparts. Since $V$ is a sum of nonnegative weighted squared distances at each time, $V=0$ forces every fitted center occupied with positive probability to equal its assigned true center, although several fitted states may share that center. The equality $D_I=0$ then ensures that, after these states are grouped, the joint probabilities of the corresponding location sequences equal those of the true model. It follows that the fitted and true latent location measures coincide within the observed block, even though their representations by individual states need not be unique.
\subsection{Wasserstein Recovery of Emission-Center Sequence}
We now use the bounds of Section~\ref{sec:overfit_global} to assess recovery
of emission-center distributions. This comparison depends on
both the displacement of fitted centers and the distances
over which probability discrepancies must be transported.
Because $V$ controls squared displacement, the natural
transport loss here is squared Wasserstein distance.
The following theorem gives global bounds for squared Wasserstein error at any selection of times within the observed block.
\label{sec:overfit_latent_path_transport}
\begin{theorem}[Wasserstein Recovery of Emission-Center Sequence under over-specification]\label{thm:overtransport}
Let $\theta\in\Theta_k(R)$ with $k\ge K$. Let $I\subseteq[T]$ be nonempty, with $|I|=r$, and let $\beta_r$ be as in Theorem~\ref{thm:overglobal}. In the structured regimes,
\begin{equation}
W_{2,I}^2(G_\theta^I,G_\star^I)
\le C_r\Delta_{\mathrm{sep}}^{-(\beta_r-2)}h(p_{\theta,T},p_{\theta^\star,T}).
\end{equation}
Without cluster assumptions, the one-time bound is
\begin{equation}
\max_{t\le T}W_2^2(G_\theta^{\{t\}},G_\star^{\{t\}})
\le C\Delta_{\mathrm{sep}}^{-(2K-2)}h(p_{\theta,T},p_{\theta^\star,T}).
\label{eq:overtransport-unstructured-marginal}
\end{equation}
For $r\ge2$, the unstructured bound is
\begin{equation}
W_{2,I}^2(G_\theta^I,G_\star^I)
\le C_r\Delta_{\mathrm{sep}}^{-(\gamma_{K,r}-2)}h(p_{\theta,T},p_{\theta^\star,T}).
\end{equation}
\end{theorem}

The proof of Theorem~\ref{thm:overtransport} is given in
Appendix~\ref{app:proof_43}.
The theorem combines location and probability errors
to control the joint distribution of latent locations,
without requiring a one-to-one correspondence between
fitted and true states. Because this control is linear
in $h(p_{\theta,T},p_{\theta^\star,T})$ for squared Wasserstein distance, taking square
roots gives the corresponding bound for $W_{2,I}$,
which is proportional to $\sqrt{h(p_{\theta,T},p_{\theta^\star,T})}$.

\begin{proposition}[Sharpness of the transport separation exponents]
\label{prop:transport-sharpness}
Fix $s\ge2$, $r\ge1$, and $T\ge\max\{2,r\}$. There exist over-specified Gaussian HMMs with $K=s$ true states in a single collinear cluster and true initial-state probabilities bounded away from zero such that, as $\Delta_{\mathrm{sep}}\to0$,
\[
\frac{W_{2,I}^2(G_\theta^I,G_\star^I)}
{h(p_{\theta,T},p_{\theta^\star,T})}
\asymp
\Delta_{\mathrm{sep}}^{-(\beta_r-2)},
\qquad I=\{1,\ldots,r\},
\]
where $\beta_r=\max\{2s,r(s-1)+s\}$. The fitted-state count can be chosen independently of $\Delta_{\mathrm{sep}}$. Hence, the structured exponent $\beta_r-2$ in Theorem~\ref{thm:overtransport} is sharp over model classes containing these collinear configurations. Taking $s=K$ also establishes sharpness of the unstructured exponent $\gamma_{K,r}-2$.
\end{proposition}

The result follows from Proposition~\ref{prop:overfit-sharpness}, where we construct pairs of true and over-specified Gaussian HMMs with collinear emission centers that attain the separation exponents in Theorem~\ref{thm:overtransport}. These joint constructions establish sharpness over model classes containing the corresponding collinear configurations. For each choice of $s$ and $r$, both the observation horizon and the number of fitted states remain fixed as emission separation tends to zero.

\paragraph{Cluster geometry and transport cost.}
The Wasserstein bounds account for both the displacement of fitted centers and the redistribution of probability between true centers. The first contribution is controlled by the bound for $V$ in Theorem~\ref{thm:overglobal}. To control the second, we move fitted locations to their assigned true centers and account for the distances over which probability mass is redistributed. The argument differs according to whether cluster assumptions are available.

\begin{enumerate}
\item[\textit{(i)}] \emph{Under cluster assumptions.} We separate the mass redistributed within each cluster sequence from the remaining mass transported between different cluster sequences. Within a cluster sequence, the two locations at each selected time belong to the same cluster and are at most $C_0\Delta_{\mathrm{sep}}$ apart. Each unit of redistributed mass therefore incurs a squared cost of at most $C_0^2\Delta_{\mathrm{sep}}^2$. Multiplying the aggregate probability-error bound in Theorem~\ref{thm:overglobal} by this cost reduces its separation exponent from $\beta_r$ to $\beta_r-2$.

For transport between different cluster sequences, however, the distances need not shrink with $\Delta_{\mathrm{sep}}$. We instead control the amount of mass moved through the bound on cluster-sequence probability discrepancies in Lemma~\ref{lem:cluster-mass}, of order $C_r\Delta_{\mathrm{sep}}^{-(2s_{\max}-2)}h(p_{\theta,T},p_{\theta^\star,T}).$ Since the squared transport distances are bounded, the resulting contribution has separation exponent $2s_{\max}-2$. This exponent is no larger than $\beta_r-2$, so transport between cluster sequences does not worsen the separation dependence obtained within cluster sequences. The displacement contribution also has exponent $2s_{\max}-2$ by the bound for $V$. Combining these contributions gives the squared Wasserstein exponent $\beta_r-2$.

\item[\textit{(ii)}] \emph{Without cluster assumptions.} The true centers need not form clusters whose diameters shrink with $\Delta_{\mathrm{sep}}$, so we cannot use the within-cluster cost reduction in~(i). We instead retain the diameter of the full true support, $D_\star:=\max_{i,j}\|a_i-a_j\|_2$, in the probability comparison. The interpolation bounds in Lemma~\ref{lem:diameter-joint-selectors} give $D_I\le C_rD_\star^{-2}\Delta_{\mathrm{sep}}^{-(\gamma_{K,r}-2)}h(p_{\theta,T},p_{\theta^\star,T}).$ Redistributing probability between true centers incurs squared cost at most $D_\star^2$ per unit mass. Multiplying by this cost cancels the factor $D_\star^{-2}$, leaving separation exponent $\gamma_{K,r}-2$. The displacement contribution has exponent $2K-2$, which is no larger because $\gamma_{K,r}\ge2K$. The combined bound therefore has squared Wasserstein exponent $\gamma_{K,r}-2$. In particular, $\gamma_{K,1}=2K$ gives the one-time exponent $2K-2$.
\end{enumerate}

Having accounted for both displacement and probability redistribution, we can now obtain statistical rates for the latent emission-center distributions. Applying Theorem~\ref{thm:overtransport} on the density event of Proposition~\ref{prop:density} gives the following result for the global maximum likelihood estimator.

\begin{proposition}[MLE Wasserstein rates under over-specification]
\label{prop:mle-overtransport}
Let $\widehat{\theta}_n$ be a measurable global MLE over $\Theta_k(R)$, where $k\geq K$. For any nonempty $I\subseteq[T]$, let $r=|I|$. In the structured regimes,
\begin{align*}
W_{2,I}(G_{\widehat{\theta}_n}^I,G_\star^I)
&=
\mathcal{O}_{\mathbb{P}}\left(
\Delta_{\mathrm{sep}}^{-(\beta_r-2)/2}
\left[\frac{d\log(nTd)}{n}\right]^{1/4}
\right).
\end{align*}
Without cluster assumptions,
\begin{align*}
\max_{t\le T}W_2(G_{\widehat{\theta}_n}^{\{t\}},G_\star^{\{t\}})
&=\mathcal{O}_{\mathbb{P}}\left(
\Delta_{\mathrm{sep}}^{-(K-1)}
\left[\frac{d\log(nTd)}{n}\right]^{1/4}\right),\\
W_{2,I}(G_{\widehat{\theta}_n}^I,G_\star^I)
&=\mathcal{O}_{\mathbb{P}}\left(
\Delta_{\mathrm{sep}}^{-(\gamma_{K,r}-2)/2}
\left[\frac{d\log(nTd)}{n}\right]^{1/4}\right),\qquad r\ge2.
\end{align*}
The rates hold simultaneously for all nonempty $I\subseteq[T]$. These bounds require neither emission-center matching nor a positive lower bound on fitted state probabilities.
\end{proposition}

The result follows by applying Theorem~\ref{thm:overtransport} on the density event of Proposition~\ref{prop:density} and taking square roots. For fixed model dimensions and emission geometry, both single-time and joint emission-center distributions have an $n^{-1/4}$ upper rate in $W_2$, up to logarithmic factors. The corresponding squared Wasserstein errors have an $n^{-1/2}$ upper rate. 

Table~\ref{tab:mle-overfit-rates} summarizes the MLE convergence rates and their separation dependence under over-specification. The common density-estimation scale is $\sqrt{\frac{d\log(nTd)}{n}}$. The upper bounds for aggregate probabilities, net displacement $B$, and quadratic displacement $V$ are linear in this scale, whereas the bounds for probability-weighted emission-center displacement $D_\mu$ and Wasserstein distance involve its square root. The table therefore distinguishes recovery of squared displacement from recovery in $W_2$ itself.

\begin{table}[t]
\centering
\small
\setlength{\tabcolsep}{5pt}
\renewcommand{\arraystretch}{1.6}
\setlength{\extrarowheight}{3pt}
\arrayrulecolor{black}
\begin{tabular*}{\textwidth}{@{\extracolsep{\fill}}lcc@{}}
\toprule
Recovery loss & Structured geometry & Unstructured geometry \\
\midrule[0.7pt]
$V$
& $\Delta_{\mathrm{sep}}^{-(2s_{\max}-2)}\sqrt{\frac{d\log(nTd)}{n}}$
& $\Delta_{\mathrm{sep}}^{-(2K-2)}\sqrt{\frac{d\log(nTd)}{n}}$ \\
\arrayrulecolor{black!12}\midrule[0.25pt]\arrayrulecolor{black}
$B$
& $\Delta_{\mathrm{sep}}^{-(2s_{\max}-1)}\sqrt{\frac{d\log(nTd)}{n}}$
& $\Delta_{\mathrm{sep}}^{-(2K-1)}\sqrt{\frac{d\log(nTd)}{n}}$ \\
\arrayrulecolor{black!12}\midrule[0.25pt]\arrayrulecolor{black}
$D_\pi\vee\max_t D_{\{t\}}$
& $\Delta_{\mathrm{sep}}^{-2s_{\max}}\sqrt{\frac{d\log(nTd)}{n}}$
& $\Delta_{\mathrm{sep}}^{-2K}\sqrt{\frac{d\log(nTd)}{n}}$ \\
\arrayrulecolor{black!12}\midrule[0.25pt]\arrayrulecolor{black}
$D_A$
& $\Delta_{\mathrm{sep}}^{-\beta_2}\sqrt{\frac{d\log(nTd)}{n}}$
& $\Delta_{\mathrm{sep}}^{-\gamma_{K,2}}\sqrt{\frac{d\log(nTd)}{n}}$ \\
\arrayrulecolor{black!12}\midrule[0.25pt]\arrayrulecolor{black}
$D_I$
& $\Delta_{\mathrm{sep}}^{-\beta_r}\sqrt{\frac{d\log(nTd)}{n}}$
& $\Delta_{\mathrm{sep}}^{-\gamma_{K,r}}\sqrt{\frac{d\log(nTd)}{n}}$ \\
\midrule[0.7pt]
$D_\mu\vee\max_t W_2$
& $\Delta_{\mathrm{sep}}^{-(s_{\max}-1)}\left(\frac{d\log(nTd)}{n}\right)^{1/4}$
& $\Delta_{\mathrm{sep}}^{-(K-1)}\left(\frac{d\log(nTd)}{n}\right)^{1/4}$ \\
\arrayrulecolor{black!12}\midrule[0.25pt]\arrayrulecolor{black}
$W_{2,I}$, $r\ge2$
& $\Delta_{\mathrm{sep}}^{-(\beta_r-2)/2}\left(\frac{d\log(nTd)}{n}\right)^{1/4}$
& $\Delta_{\mathrm{sep}}^{-(\gamma_{K,r}-2)/2}\left(\frac{d\log(nTd)}{n}\right)^{1/4}$ \\
\bottomrule
\end{tabular*}
\caption{\textbf{Overfitted MLE rates for $k\ge K$.} Each entry is an $\mathcal O_{\mathbb P}$ upper rate. Here $\beta_r=\max\{2s_{\max},r(s_{\max}-1)+s_{\max}\}$ and $\gamma_{K,r}=\max\{2K,r(K-1)+K\}$. No center matching or fitted-weight lower bound is required, and the constants are independent of $\pi_{\min}^\star$. Joint constants may depend on $r=|I|$. The notation $\max_t W_2$ abbreviates $\max_t W_2(G_{\widehat\theta_n}^{\{t\}},G_\star^{\{t\}})$.}
\label{tab:mle-overfit-rates}
\end{table}

The following example compares the exact-fit and over-specified bounds for the same emission geometry, showing how both change when four centers form two separated pairs rather than a single cluster.

\begin{example}[Four centers in two dimensions]
For $0<\Delta_{\mathrm{sep}}\le1/4$, let
\[
a_1=(-2,0),\quad a_2=(-2,\Delta_{\mathrm{sep}}),\quad
 a_3=(2,0),\quad a_4=(2,\Delta_{\mathrm{sep}}).
\]
Under this configuration, the centers form two clusters of size two separated by a macroscopic gap, with $K=4$ but $s_{\max}=2$. The matched one-time $W_1$ bound has separation factor $\Delta_{\mathrm{sep}}^{-2}$, while the matched transition-flow bound has factor $\Delta_{\mathrm{sep}}^{-3}$. The overfitted one-time $W_2$ bound has separation factor $\Delta_{\mathrm{sep}}^{-1}$. The first two factors multiply the Hellinger error, whereas the $W_2$ factor multiplies its square root. If all four centers instead lie in one $\mathcal O(\Delta_{\mathrm{sep}})$ cluster, these multipliers become $\Delta_{\mathrm{sep}}^{-6}$, $\Delta_{\mathrm{sep}}^{-7}$, and $\Delta_{\mathrm{sep}}^{-3}$.
\end{example}
This example shows why the largest cluster size, rather than the total number of states, determines the structured separation exponents. Although both configurations have four true states and the same minimum separation, two separated pairs give smaller recovery-error multipliers than a single cluster of four centers. The distances between clusters enter only through the constants in the bounds, while the within-cluster distances determine the inverse powers of $\Delta_{\mathrm{sep}}$. For recovery at two selected times in the two-pair configuration, the overfitted joint probability bound has separation exponent $\beta_2=4$, while the corresponding $W_{2,I}$ bound has exponent $(\beta_2-2)/2=1$. The same comparison holds in higher dimensions after adding zero coordinates to the centers, although the density-estimation bound still depends on $d$.

\subsection{Aggregate Recovery with Approximate Reference Centers}
\label{sec:estimated-reference-grouping}

The recovery losses introduced above assign each fitted state to its nearest \emph{true} emission center. This assignment provides a precise target for the inverse inequalities, but the true centers are unavailable in an application. A natural alternative is to group fitted states using estimated reference centers. To relate this construction to the theoretical losses, we must control how errors in the reference locations affect the assignment of fitted states and the resulting aggregate probabilities. The following result provides this comparison without imposing a Markov property on the grouped process or specifying the algorithm used to obtain the references.

Let $\xi_1,\ldots,\xi_K\in\mathbb R^d$ be reference centers with labels aligned to the true centers, and define
\[
\widehat c(i):=\arg\min_{j\in[K]}\|\mu_i-\xi_j\|_2,
\qquad
\widehat P_\theta^I(\boldsymbol j)
:=\mathbb P_\theta\{\widehat c(Z_{t_\ell})=j_\ell,
\ 1\le\ell\le r\},
\]
with smallest-index tie breaking. Thus, $\widehat P_\theta^I$ is obtained by summing the probabilities of fitted state sequences assigned to the same reference-label sequence. Write $\widehat D_I=\|\widehat P_\theta^I-P_\star^I\|_1$, and define $\widehat D_\pi=\widehat D_{\{1\}}$ and $\widehat D_A=\max_{t<T}\widehat D_{\{t,t+1\}}$ analogously. 

\begin{proposition}[Stability of grouped distributions under reference error]
\label{prop:reference-grouping-stability}
Suppose that, after a common relabeling of the reference centers,
\[
\max_{j\in[K]}\|\xi_j-a_j\|_2\le\tau
\le\frac14\delta_{\mathrm{sep}}.
\]
For each fitted state $i$, let
\[
m_i(\boldsymbol\xi):=\min_{\ell\ne\widehat c(i)}
\left\{\|\mu_i-\xi_\ell\|_2-\|\mu_i-\xi_{\widehat c(i)}\|_2\right\},
\qquad
\mathcal B_\tau(\boldsymbol\xi):=\{i:m_i(\boldsymbol\xi)\le2\tau\}.
\]
The margin $m_i(\boldsymbol\xi)$ is the difference between the distances to the second-nearest and nearest reference centers and is computed without the true centers. Under the reference accuracy condition, each distance changes by at most $\tau$ when a reference center is replaced by its corresponding true center, so the difference between two distances changes by at most $2\tau$. A margin greater than $2\tau$ therefore guarantees an unchanged assignment, and only states in $\mathcal B_\tau(\boldsymbol\xi)$ can contribute to the additional grouping error. For any fitted parameter $\theta\in\Theta_k(R)$ with $k\ge K$ and every nonempty $I\subseteq[T]$, $r=|I|$, one has
\begin{equation}
\widehat D_I
\le
\underbrace{D_I}_{\text{true-center grouping error}}
+
\underbrace{
\min\!\left\{
2\sum_{t\in I}\sum_{i\in\mathcal B_\tau(\boldsymbol\xi)}q_t(i),
\frac{8r}{\delta_{\mathrm{sep}}}D_\mu,
\frac{32r}{\delta_{\mathrm{sep}}^2}V
\right\}
}_{\text{additional error from reference grouping}}.
\label{eq:reference-grouping-stability}
\end{equation}
In particular, if no state with positive occupancy at a selected time lies in $\mathcal B_\tau(\boldsymbol\xi)$, then $\widehat P_\theta^I=\bar P_\theta^I$ exactly. Under the reference accuracy condition, the margin term can be computed from the fitted model and reference centers to bound the additional error caused by reference-based grouping. Under exact specification, we control the additional grouping error through $D_\mu/\delta_{\mathrm{sep}}$, and under over-specification through $V/\delta_{\mathrm{sep}}^2$. Substituting the preceding displacement bounds shows that these additional errors have no larger separation exponents than the corresponding bound for $D_I$.
Consequently, the global exact-fit probability bounds of Theorem~\ref{thm:globalexact} and the overfitted probability bounds of Theorem~\ref{thm:overglobal} remain valid with $D_I$ replaced by $\widehat D_I$, with the \emph{same separation exponents}. The same conclusion holds for $\widehat D_\pi$ and $\widehat D_A$, and in the unstructured regimes with $s_{\max}=K$. In the matched exact-fit regime the corresponding constants remain independent of $\pi_{\min}^\star$. In the overfitted regime no factor $(\pi_{\min}^\star)^{-1}$ is introduced.
\end{proposition}

\begin{remark}[Why a reference accuracy condition is necessary]
The following example shows why controlling reference-center error is necessary for aggregate recovery. Consider two exactly fitted states with true centers $a_1=0$, $a_2=\delta>0$ and initial masses $(1/2,1/2)$. If the two supplied reference centers both equal $\delta/2$, their maximum true-center error is $\delta/2$. Smallest-index tie breaking then assigns both fitted states to the first reference group. Thus $\widehat D_\pi=1$ even though $D_\pi=D_\mu=V=h(p_{\theta,T},p_{\theta^\star,T})=0$. Therefore no reference-based density-to-group-recovery implication can hold without some condition excluding such ambiguity.
\end{remark}

The proof of Proposition~\ref{prop:reference-grouping-stability} is given in Appendix~\ref{app:proof-reference-grouping-stability}. The bound separates two sources of error: the discrepancy under the true-center assignment and the additional discrepancy caused by changing that assignment. The margin term measures the fitted occupancy of states whose assignments may change, and it can be computed from the fitted model and supplied references once an accuracy bound $\tau$ has been justified.

To obtain a statistical guarantee, suppose the reference accuracy event holds with probability at least $1-\gamma_n$ and the density event underlying a displayed recovery rate holds with probability at least $1-\alpha_n$. On their intersection, the same recovery rate applies to the reference-grouped probabilities. This intersection has probability at least $1-\alpha_n-\gamma_n$, without requiring independence between the reference estimates and the fitted model.

\paragraph{Obtaining reference centers from a pilot estimator.}
The exact-fit theory provides a way to construct sufficiently accurate reference centers for grouping. When $K$ is known, we first fit a global $K$-state MLE to $m$ independent trajectories and use its emission centers as references. If $m$ satisfies the sample-size condition~\eqref{eq:exact-localization-sample-size}, Corollary~\ref{cor:localization} guarantees that these reference centers can be matched to the true centers within distance $\tau=\Delta_{\mathrm{sep}}/4\le\delta_{\mathrm{sep}}/4$ with probability at least $1-1/m$. This preliminary fit, which we call the pilot estimator, therefore supplies reference centers satisfying the accuracy condition required for aggregate recovery.

We then use these reference centers to group the states of a $k$-state estimator whose full-block Hellinger error is at most $\varepsilon_n$ with probability at least $1-\alpha_n$. On the intersection of the reference-accuracy and density-accuracy events, the reference-grouped probabilities satisfy the corresponding recovery bounds with $h$ replaced by $\varepsilon_n$. This gives a probability guarantee of at least $1-1/m-\alpha_n$, without requiring independence between the two estimators. In particular, Proposition~\ref{prop:density} gives $\alpha_n=1/n$ when the second estimator is the global MLE. Since the pilot fit determines reference centers only up to a permutation, comparing the grouped probabilities with their true counterparts also requires aligning the labels. The next corollary shows that one common relabeling suffices for all marginal and joint state-sequence distributions.

\begin{corollary}[Reference recovery up to one common permutation]
\label{cor:reference-common-permutation}
Suppose the $K$ reference centers are unlabeled and some permutation
$\sigma$ of $[K]$ satisfies
\begin{equation}\label{eq:reference-permuted-centers}
\max_{\ell\in[K]}\|\xi_\ell-a_{\sigma(\ell)}\|_2
 \le\tau\le\delta_{\mathrm{sep}}/4.
\end{equation}
For a reference-label distribution $Q$ on $[K]^r$, put
$(\sigma_\#Q)(j_1,\ldots,j_r)
=Q(\sigma^{-1}(j_1),\ldots,\sigma^{-1}(j_r))$.
There is a unique permutation $\sigma$ satisfying
\eqref{eq:reference-permuted-centers}, and this \emph{same}
permutation satisfies, simultaneously for all nonempty $I\subseteq[T]$,
\[
\|\sigma_\#\widehat P_\theta^I-P_\star^I\|_1
\le D_I+\min\!\left\{
2\sum_{t\in I}\sum_{i\in\mathcal B_\tau(\boldsymbol\xi)}q_t(i),
\frac{8|I|D_\mu}{\delta_{\mathrm{sep}}},
\frac{32|I|V}{\delta_{\mathrm{sep}}^2}\right\}.
\]
Consequently, the reference-grouped probability bounds retain the
separation exponents obtained using the true emission centers for grouping under one common relabeling of every
time block.

If $\boldsymbol\xi$ consists of the centers from a measurable
exactly specified global $K$-state MLE on $m$ independent
trajectories, with $K$ known and $m$ satisfying
\eqref{eq:exact-localization-sample-size}, and a second (possibly
dependent) estimator $\widetilde\theta_n$ satisfies
$\mathbb P_\star\{h(p_{\widetilde\theta_n,T},p_{\theta^\star,T})
\le\varepsilon_n\}\ge1-\alpha_n$,
then with probability at least $1-m^{-1}-\alpha_n$ one common
permutation yields the corresponding $D_I$ rates obtained using the true emission centers for grouping for
all nonempty $I\subseteq[T]$.
\end{corollary}

The corollary uses the same permutation across all marginal and joint state-sequence distributions to preserve a consistent correspondence between reference groups and true states. Relabeling these distributions separately could conceal differences in the underlying dynamics. This common permutation is needed only to compare the grouped probabilities with their true counterparts, since constructing the reference groups requires only the reference centers.

Proposition~\ref{prop:reference-grouping-stability} and Corollary~\ref{cor:reference-common-permutation} therefore extend the recovery bounds to sufficiently accurate reference centers while retaining the separation exponents obtained using the true emission centers for grouping. To examine how reference errors affect aggregate recovery, we perturb reference centers within the permitted radius and measure the resulting changes in assignments and grouped probability errors in Appendix~\ref{app:synthetic-reference-sensitivity}. Beyond these controlled perturbations, we also examine grouping by complete-linkage clustering of the fitted centers in Appendix~\ref{app:synthetic-grouping}, where we compare the resulting aggregate-flow errors with those obtained using the true emission centers for grouping.

\section{Empirical Analysis}
\label{sec:simulation}
In the preceding sections, we have established latent recovery bounds with explicit dependence on emission separation and derived the corresponding convergence rates under exact specification and over-specification. Guided by these results, we empirically investigate how emission geometry, sample size, and model specification affect observation fit and latent recovery through synthetic experiments with known generating parameters. More concretely, our experiments address three main questions: (i) how emission separation and local cluster size determine the extent to which observation error can translate into errors in state probabilities and latent location distributions; (ii) how splitting a true state into multiple fitted states affects center displacement and the recovery of aggregate transition flows versus individual transition rows; and (iii) how increasing the number of independent trajectories reduces recovery errors and promotes one-to-one matching between fitted and true emission centers in exact-sepecified setting. To this end, we implement two complementary experimental designs:

\begin{itemize}
    \item The first setup is intended to assess the separation powers in our inverse inequalities independently of estimation and optimization error. In particular, we construct pairs of HMMs, with one serving as the true model and the other as a candidate model. For each pair, we specify both sets of parameters directly to satisfy the moment-matching or joint cancellation conditions used in the sharpness constructions of Propositions~\ref{prop:matched-sharpness} and~\ref{prop:overfit-sharpness}. For each fixed cluster size, we evaluate the ratios of latent recovery errors to Hellinger error over a range of emission separations and compare the resulting separation powers with the theoretical predictions. This comparison numerically illustrates the sharpness of the bounds, as ratios growing at the predicted powers show that the separation exponents cannot be reduced uniformly over the corresponding model classes. We repeat this comparison across cluster sizes to examine their effect on these powers. The model constructions and separation grids are detailed in Sections~\ref{sec:synthetic-exact} and~\ref{sec:synthetic-overconditioning}.
    \item The second setup examines how latent recovery errors and emission-center matching behave when HMMs are estimated from data. Specifically, we generate independent observation trajectories from a known true HMM and fit models with either the true number of states or a larger number. Next, we compare the fitted models with the generating model to evaluate observation and latent recovery errors. At each fixed sample size, we repeat this procedure on independently generated data sets and summarize the errors and the frequency of successful one-to-one center matching under exact specification. The sample-size and matching experiments are presented in Section~\ref{sec:synthetic-statistical}, which specifies the generating models, sample-size grids, and numbers of repetitions. The fitting and evaluation procedures are detailed in Appendix~\ref{app:synthetic-protocol}.
\end{itemize}

\noindent To evaluate aggregate recovery, we assign each fitted state to the group associated with its nearest true emission center, following Section~\ref{sec:latent_recovery_loss}. We sum state probabilities within each group and transition flows between each pair of groups, and compare these quantities with their true counterparts. This grouping is used only to calculate the recovery losses, leaving the fitted emission centers and transition parameters unchanged. Since the true emission centers are unavailable in practice, Appendix~\ref{app:synthetic-grouping} compares grouping based on fitted centers with the nearest-true-center assignment. Appendix~\ref{app:synthetic-reference-sensitivity} then examines the effect of reference-center error under the stability result in Section~\ref{sec:estimated-reference-grouping}.
\paragraph{Data generation and model fitting.}  Throughout the estimation experiments in the second setup described above, we generate $n$ independent trajectories of fixed length $T=6$ and hold the common emission covariance fixed at its known value during fitting. In these experiments, we generate and fit the models with identity emission covariance and constrain the emission centers to the Euclidean ball $B_d(0,6)$. With $T$ fixed, increasing $n$ adds independent trajectories rather than extending individual trajectories, which follows the sampling scheme presented in Section~\ref{prelim:hmm}. For model fitting, we use the Baum--Welch algorithm with four initializations, except for the globally optimized benchmark described in Appendix~\ref{app:synthetic-single}. The stopping criteria for each run are an absolute change in log-likelihood divided by $nT$ of at most $10^{-7}$ or a maximum of $400$ iterations, whichever is reached first. Among the four runs, we select the fit with the highest final training likelihood and include it in the reported results regardless of which stopping condition is met. Table~\ref{tab:synthetic-map} summarizes the experiments and the quantities evaluated.

\begin{table}[t]
\centering
\small
\captionsetup{font=small,margin=0pt,width=\linewidth}
\setlength{\tabcolsep}{6pt}
\renewcommand{\arraystretch}{1.12}
\arrayrulecolor{black}
\begin{tabularx}{\linewidth}{
  @{}>{\raggedright\arraybackslash}p{.28\linewidth}
  >{\raggedright\arraybackslash}X@{}
}
\toprule
\textbf{Experiment} & \textbf{Quantities evaluated} \\
\midrule[0.7pt]
\addlinespace[4pt]
\multicolumn{2}{@{}l}{
  \textbf{\textit{Separation dependence}}
} \\
\addlinespace[5pt]
Exact specification
& $D_I/h$, $D_\mu/h$, and $W_{1,I}/h$ across emission separations and cluster sizes \\
\arrayrulecolor{black!15}\midrule[0.25pt]\arrayrulecolor{black}
Over-specification at one time
& $D_\pi/h$, $B/h$, $V/h$, and $W_2^2/h$ across emission separations \\
\arrayrulecolor{black!15}\midrule[0.25pt]\arrayrulecolor{black}
Over-specification across multiple times
& Joint probability and squared Wasserstein errors relative to observation error across emission separations, with $|I|=2,3,4$ selected times \\
\midrule[0.7pt]
\addlinespace[4pt]
\multicolumn{2}{@{}l}{
  \textbf{\textit{Statistical recovery}}
} \\
\addlinespace[5pt]
Increasing the sample size
& Observation and latent recovery errors as the number of independent trajectories increases \\
\arrayrulecolor{black!15}\midrule[0.25pt]\arrayrulecolor{black}
Globally optimized benchmark
& Exact-fit $W_1$ and split $W_2$ errors compared with the $n^{-1/2}$ and $n^{-1/4}$ reference scales \\
\arrayrulecolor{black!15}\midrule[0.25pt]\arrayrulecolor{black}
Matching fitted and true centers
& Frequency of one-to-one matching within $\Delta_{\mathrm{sep}}/4$ across sample sizes and emission separations \\
\midrule[0.7pt]
\addlinespace[4pt]
\multicolumn{2}{@{}l}{
  \textbf{\textit{State splitting}}
} \\
\addlinespace[5pt]
Aggregate flows and individual rows
& Aggregate transition-flow error $D_A$ and individual-row discrepancy $R_{\mathrm{row}}$ across sample sizes \\
\arrayrulecolor{black!15}\midrule[0.25pt]\arrayrulecolor{black}
Exact constructions
& Individual-row ambiguity despite identical observation distributions, and higher-order grouped-path discrepancies despite agreement of lower-order marginals \\
\bottomrule
\end{tabularx}
\caption{\textbf{Overview of the synthetic experiments in the main text.} Here, $h$ denotes the Hellinger distance between the compared observation distributions, and $I$ indexes the selected observation times. The recovery losses and Wasserstein distances are defined in Section~\ref{sec:latent_recovery_loss}.}
\label{tab:synthetic-map}
\end{table}

\paragraph{Evaluation protocol.} We evaluate errors in both the observation distribution and the latent structure. The former is measured by Hellinger distance, and the latter by the probability and displacement losses $D_\pi$, $\max_t D_{\{t\}}$, $D_A$, $D_\mu$, $V$, and $B$. For $T=6$, we evaluate $D_A$ as the maximum over all five consecutive time pairs. In addition, we consider the selected-time sets $\{1\}$, $\{1,2\}$, $\{1,6\}$, and $\{1,3,6\}$ to examine marginal, consecutive, nonconsecutive, and three-time recovery. For each selected-time set $I$, we compute $D_I$ and the corresponding Wasserstein distances between the fitted and true joint latent distributions. These calculations account for dependence across the selected times through the joint probabilities specified by each HMM. We use the normalized block costs $d_{1,r}$ and $d_{2,r}$, with $r=|I|$, to average the transport cost over the selected times.

For the directly specified pairs of HMMs with one-dimensional observations, we compute Hellinger distance using Gauss--Hermite quadrature \citep{golub1969calculation} with 80-digit arithmetic. We repeat the calculation with more quadrature points to check that the computed values remain stable. Regarding the estimation experiments, we estimate Hellinger distance using $4096$ new trajectories from the true model and another $4096$ from the fitted model, independently of the training data. The Monte Carlo and quadrature calculations are detailed in Appendices~\ref{app:synthetic-losses} and~\ref{app:synthetic-finite}, respectively. 

\subsection{Separation Exponents under Exact Specification}
\label{sec:synthetic-exact}
We first examine how cluster size affects the separation powers in the exact-specification bounds and whether these powers change with the number of selected times. For this purpose, we compare pairs of HMMs with the same number of states, constructed from the moment-matched measures in Proposition~\ref{prop:matched-sharpness}.

\paragraph{Setup.} For each $s\in\{2,3,4\}$, we construct a true HMM and a candidate HMM, each with $s$ states and one-dimensional Gaussian emissions of variance one. Let $a_1,\ldots,a_s$ and $b_1,\ldots,b_s$ denote their respective unscaled emission centers. The true model has centers $a_j=j-(s+1)/2$ and state probabilities $w_j=1/s$. The candidate model has positive state probabilities $v_1,\ldots,v_s$ satisfying $\sum_{j=1}^s v_j=1$. Each probability vector serves as both the initial distribution and every row of the corresponding transition matrix. Thus, the hidden states are independent across time within each model, and state $j$ has probability $w_j$ in the true model and $v_j$ in the candidate model at every time.

We choose different centers and probabilities for the two models but require their center moments to satisfy
\[
\underbrace{\sum_{j=1}^{s} w_j a_j^\ell}_{\text{true center moment}}
=
\underbrace{\sum_{j=1}^{s} v_j b_j^\ell}_{\text{candidate center moment}},
\qquad \ell=1,\ldots,2s-2,
\]
with unequal sums at $\ell=2s-1$. These moment equalities make the observation densities difficult to distinguish when the centers approach one another, despite the differences in their latent parameters.

Subsequently, we vary emission separation by setting the true and candidate centers to $\delta a_j$ and $\delta b_j$, respectively, keeping the probabilities and emission variance fixed. Smaller values of $\delta$ bring the centers closer together in both models. The moment equalities are preserved because scaling multiplies both moments of order $\ell$ by $\delta^\ell$. As a result, the difference between the one-time observation densities first appears at order $\delta^{2s-1}$. Evaluating the latent-to-Hellinger error ratios across values of $\delta$ allows us to compare their separation dependence with the powers predicted by the theory. We use $\delta\in\{0.80,0.55,0.38,0.26,0.18,0.12,0.08,0.05\}$ and fix the trajectory length at $T=6$. Because the unscaled true centers have minimum separation one and $\delta<1$ throughout this grid, the scaled centers satisfy $\delta_{\mathrm{sep}}=\Delta_{\mathrm{sep}}=\delta$.

We also compare four-state models whose emission centers form two separated pairs. Starting from the two-state construction above, we define the four true center--probability pairs as
\[
\begin{array}{c|cc}
& \text{Center} & \text{Probability} \\
\hline
\text{First pair}
& -2.5+\delta a_1 & w_1/2 \\
& -2.5+\delta a_2 & w_2/2 \\[3pt]
\text{Second pair}
& 2.5+\delta a_1 & w_1/2 \\
& 2.5+\delta a_2 & w_2/2
\end{array}
\]
The candidate model replaces $a_1,a_2$ by $b_1,b_2$ and $w_1,w_2$ by $v_1,v_2$ in this construction. In each model, these four probabilities form the initial distribution and every transition row. As $\delta$ decreases, the centers within each pair approach one another, while the two pairs remain separated. 

\paragraph{Evaluation.} At each value of $\delta$, we compute the Hellinger distance $h$ between the length-$T$ observation distributions of the true and candidate models. We then evaluate the ratios $D_I/h$, $D_\mu/h$, and $W_{1,I}/h$, which measure the discrepancies in state probabilities and emission locations relative to the observation error. We consider sets $I$ containing $r=1,2,3$ distinct times to examine whether the separation powers of $D_I/h$ and $W_{1,I}/h$ change when recovery is evaluated jointly across more times. Since the hidden states are independent and identically distributed in this construction, these losses depend only on the number of selected times, not on their positions. This independence also allows us to obtain the full-block Hellinger distance from the one-time distance $h_1$ through $h^2=1-(1-h_1^2)^T$. We then evaluate $h_1$ numerically from the two one-time observation densities. Similarly, the joint state probabilities are products of the one-time probabilities, allowing us to compute $D_I$ directly from the vectors $w$ and $v$. Under the normalized block cost, $W_{1,I}$ equals the one-time Wasserstein distance, which we compute by monotone matching between the two discrete distributions of emission centers.

From these ratios, we estimate the separation powers by fitting a straight line to the logarithm of each ratio against $\log\delta$ over the four smallest separation values, $\delta\in\{0.18,0.12,0.08,0.05\}$. If a ratio behaves as $C\delta^{-\alpha}$, its logarithm behaves as $\log C-\alpha\log\delta$, so the negative fitted slope estimates the exponent $\alpha$. We compare the estimated exponents with the values $2s-1$ for $D_I/h$ and $2s-2$ for both $D_\mu/h$ and $W_{1,I}/h$, predicted by Proposition~\ref{prop:matched-sharpness}.

\begin{figure}[t]
\centering
\includegraphics[width=0.87\textwidth]{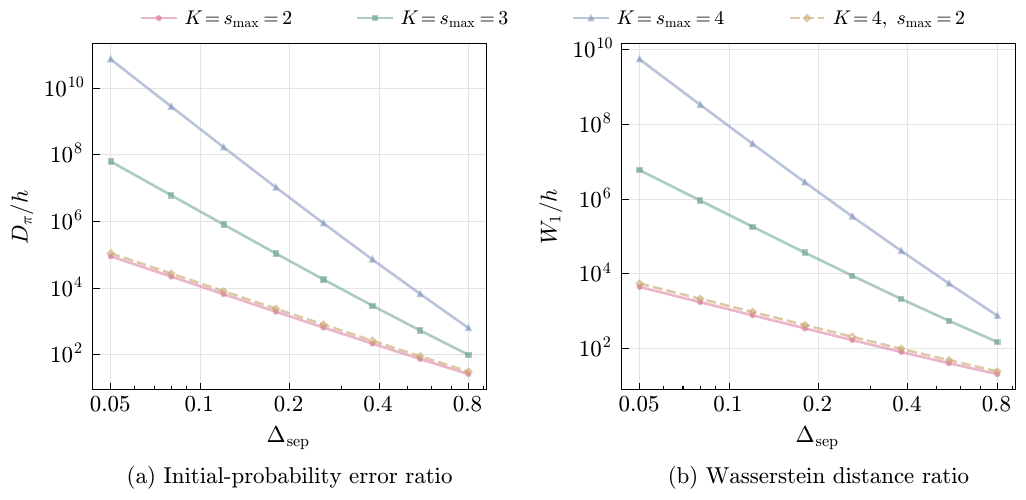}
\caption{\textbf{Exact specification with different emission-cluster sizes.} Left: $D_\pi/h$ for a single cluster of size $s=2,3,4$ and for four states arranged in two separated pairs. Right: the corresponding $W_1/h$ ratios. Each point is computed from a finite moment-matched pair of Gaussian HMMs. The separation slopes are estimated from the four smallest grid values after evaluating the densities.}
\label{fig:synthetic-exact-conditioning}
\end{figure}

\paragraph{Results.} Figure~\ref{fig:synthetic-exact-conditioning} shows that the ratios $D_\pi/h$ and $W_1/h$ increase as the emission separation decreases, with a steeper increase for larger emission clusters. This increase illustrates that, for HMM pairs satisfying the assumptions of Theorem~\ref{thm:globalexact}, a small discrepancy between observation distributions can correspond to much larger discrepancies in latent probabilities and emission locations as the emission centers approach one another. For example, for a single cluster with $s_{\max}=3$, reducing the separation from $0.8$ to $0.05$ increases $D_\pi/h$ from approximately $10^2$ to $6\times10^7$, while $W_1/h$ rises from approximately $1.5\times10^2$ to $6\times10^6$. Appendix~\ref{app:synthetic-geometry} complements these comparisons with a separate estimation experiment using four-state true models whose centers form either one cluster or two separated pairs, with the same minimum separation. At the smallest separation considered, the two-pair configuration yields smaller aggregate transition-flow errors despite larger observation-density errors, under both exact specification and over-specification.

For the directly specified pairs in Figure~\ref{fig:synthetic-exact-conditioning}, Table~\ref{tab:exact_separation_exponents} quantifies how the growth of the error ratios depends on cluster size by comparing the estimated separation exponents with their theoretical values. For single clusters of size $s_{\max}=2,3,4$, the estimated exponents for $D_\pi/h$ are $2.992$, $4.965$, and $6.916$, respectively, close to $2s_{\max}-1$. The corresponding estimates for both $D_\mu/h$ and $W_1/h$ are $1.992$, $3.965$, and $5.916$, close to $2s_{\max}-2$. These estimates show that the probability ratio grows by one additional inverse power of separation compared with the displacement and Wasserstein ratios. In each case, the observed growth is close to the power in the corresponding upper bound, numerically illustrating the sharpness established by Proposition~\ref{prop:matched-sharpness}. The dependence of these powers on cluster size can also be examined while keeping the total state count fixed. In Figure~\ref{fig:synthetic-exact-conditioning}, the four-state model with two separated pairs follows the two-state single-cluster curves closely. Table~\ref{tab:exact_separation_exponents} quantifies this similarity, with probability and Wasserstein exponents of $2.993$ and $1.993$ for the two-pair configuration, compared with $2.992$ and $1.992$ for the single two-state cluster. This agreement supports the role of the largest cluster size $s_{\max}$ in Theorem~\ref{thm:globalexact}, whose structured bounds distinguish these configurations from a single four-state cluster despite their common total state count.

\begin{figure}[!tp]
    \centering
    \includegraphics[width=0.94\textwidth]{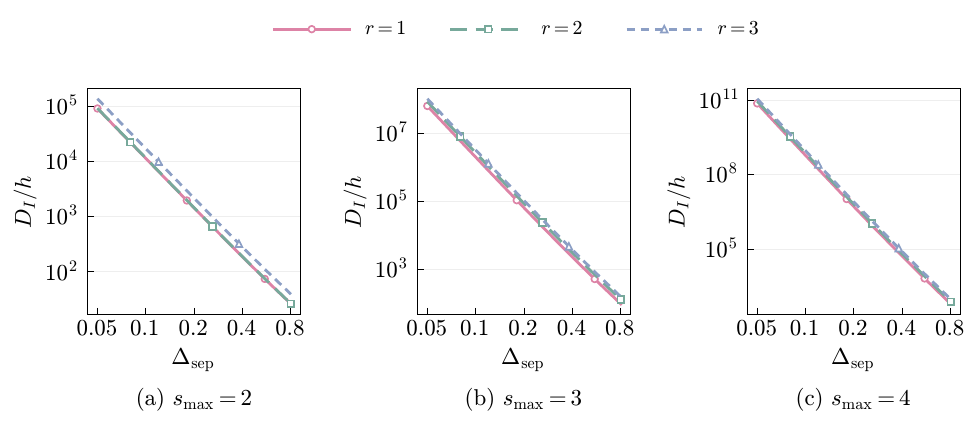}
    \caption{\textbf{Exact specification across selected times.} The ratios $D_I/h$ are plotted against $\Delta_{\mathrm{sep}}$ for $r=|I|\in\{1,2,3\}$ and a single emission cluster of size $s_{\max}=2,3,4$. Here, $h$ denotes the Hellinger distance between the length-$T$ observation distributions, with $T=6$. For each cluster size, the separation slopes remain unchanged across $r$, while the multiplicative constants differ.
    Some lines nearly overlap because these constants are close.}
    \label{fig:exact_selected_times}
\end{figure}

The dependence on cluster size raises a further question of whether joint recovery across more times incurs additional separation powers. Figure~\ref{fig:exact_selected_times} shows that, for each fixed cluster size, the curves for $r=1,2,3$ are nearly parallel on logarithmic axes, with some curves nearly overlapping. The similar slopes indicate that increasing the number of selected times does not introduce the stronger separation dependence observed when the cluster size increases. Instead, the differences appear in the multiplicative constants, as reflected in the joint probability bound of Theorem~\ref{thm:globalexact}, where $r$ enters $C_r$ while the separation exponent remains $2s_{\max}-1$.

\subsection{Separation Exponents under Over-Specification}
\label{sec:synthetic-overconditioning}
We next examine how the separation dependence of latent recovery changes under over-specification, where several fitted states may represent the same true state. We first consider recovery at a single time to examine the effect of opposing center shifts within these groups, and then study whether recovering aggregate probabilities across several times introduces additional separation powers.

\paragraph{Setup.}
To distinguish the effect of state splitting from that of recovering probabilities across several times, we consider the following two experiments.
\begin{enumerate}[label=(\alph*),font=\itshape]
\item \emph{Recovery at a single time.}
For each $s\in\{2,3,4\}$, we construct a true HMM with $K=s$ states in one emission cluster and a candidate HMM with $k=s+1$ states. As in Section~\ref{sec:synthetic-exact}, both models have one-dimensional Gaussian emissions of variance one, and each model uses its state-probability vector as both the initial distribution and every row of its transition matrix. With this choice, the hidden states are independent across time, so we can examine the effect of splitting without introducing temporal dependence.

To represent one true state by two candidate states, we place their centers at $c_1-\epsilon$ and $c_1+\epsilon$, where $c_1$ is their group center and $\epsilon=0.01$, and assign half of the group probability to each. Since the two shifts have equal weights and opposite signs, their average location remains $c_1$. Their first-order contributions to the emission density also cancel, so changes in the individual centers need not produce observation discrepancies of the same order.

To extend this cancellation across the cluster, we adjust the candidate group probabilities and centers to match the moments of the true emission-center distribution through degree $2s-1$. In particular, for true centers and probabilities $(a_j,w_j)_{j=1}^{s}$ and candidate counterparts $(b_i,v_i)_{i=1}^{s+1}$ after splitting, we require
\[
\underbrace{\sum_{j=1}^{s}w_j a_j^\ell}_{\text{true center moment}}
=
\underbrace{\sum_{i=1}^{s+1}v_i b_i^\ell}_{\text{candidate center moment}},
\qquad \ell=1,\ldots,2s-1.
\]
Compared with the exact-specification construction, this choice matches one additional moment. To examine its effect as the centers approach one another, we multiply all centers by $\delta$ while keeping their probabilities fixed. Under this scaling, the moment equalities remain valid and the distance between the split centers becomes $2\delta\epsilon$. We use $\delta\in\{0.80,0.55,0.38,0.26,0.18,0.12,0.08,0.05\}$ and fix $T=6$.

\begin{table}[!tp]
\centering
\small
\captionsetup{font=small,margin=0pt,width=\linewidth}
\setlength{\tabcolsep}{7pt}
\renewcommand{\arraystretch}{1.12}
\arrayrulecolor{black}
\begin{tabular}{@{}cccccccc@{}}
\toprule
& & \multicolumn{6}{c}{\textbf{Separation exponent}} \\
\cmidrule(lr){3-8}
& & \multicolumn{2}{c}{$D_\pi/h$}
  & \multicolumn{2}{c}{$D_\mu/h$}
  & \multicolumn{2}{c}{$W_1/h$} \\
\cmidrule(lr){3-4}
\cmidrule(lr){5-6}
\cmidrule(lr){7-8}
$K$ & $s_{\max}$
& Theoretical & Estimated
& Theoretical & Estimated
& Theoretical & Estimated \\
\midrule[0.7pt]
2 & 2 & 3 & 2.992 & 2 & 1.992 & 2 & 1.992 \\
\arrayrulecolor{black!15}\midrule[0.25pt]\arrayrulecolor{black}
3 & 3 & 5 & 4.965 & 4 & 3.965 & 4 & 3.965 \\
\arrayrulecolor{black!15}\midrule[0.25pt]\arrayrulecolor{black}
4 & 4 & 7 & 6.916 & 6 & 5.916 & 6 & 5.916 \\
\midrule[0.5pt]
\rowcolor{black!5}
4 & \textbf{2}
& 3 & \textbf{2.993}
& 2 & \textbf{1.993}
& 2 & \textbf{1.993} \\
\bottomrule
\end{tabular}
\caption{\textbf{Separation exponents under exact specification.} Estimated exponents are the negative log--log slopes over the four smallest separation values. The theoretical exponents are $2s_{\max}-1$ for $D_\pi/h$ and $2s_{\max}-2$ for $D_\mu/h$ and $W_1/h$. The shaded row uses two separated pairs; the other rows use a single cluster. Bold estimates highlight the smaller exponents for two separated pairs than for a single four-center cluster at the same $K=4$.}
\label{tab:exact_separation_exponents}
\end{table}

\item \emph{Recovery across several times.}
To examine joint recovery, we use a true three-state HMM with emission centers $(-\delta,0,\delta)$ and independent hidden states, each uniformly distributed over the three labels. We compare its latent emission-center distribution at $I=\{1,\ldots,r\}$, for $r=2,3,4$, with that of an over-specified candidate model. The candidate first selects a sequence of $r$ labels, then follows that sequence to generate the first $r$ observations. Each possible sequence has its own chain of $r$ hidden states, one for each position. This realizes the joint construction in Proposition~\ref{prop:overfit-sharpness} as a time-homogeneous HMM.

Each label sequence initially has probability $3^{-r}$. We first change the probabilities of the label sequences, creating differences in the joint latent probabilities. To keep these differences difficult to detect from observations, we also shift the first emission center of each sequence so that the probability changes and center shifts partly compensate in the observation density. These changes follow the formulas in Proposition~\ref{prop:overfit-sharpness}, which cancel the low-degree density terms at first order in the perturbation size. We then examine how much the joint latent probabilities and emission-center distributions can differ relative to the resulting observation error. After the first $r$ observations, the candidate generates independent observations from the true model, so all observation discrepancy comes from the first $r$ times. For each choice of $r$, we vary the emission separation over the same grid as in comparison (a). Throughout each scan, the candidate state count remains fixed at $k=r3^r+3$, although this count differs across the three choices of $r$.
\end{enumerate}

\paragraph{Evaluation.}
For experiment (a), we first assign each candidate state to its nearest true center and sum the probabilities within each group. After this aggregation, we compute $D_\pi$ to assess whether the total probability assigned to each true state is recovered. To relate this discrepancy to the observation error, we divide it by the Hellinger distance $h$ between the length-$T$ observation distributions. In addition to these probability discrepancies, we examine the locations of the candidate centers within each group. Since opposing shifts can cancel when averaged, we retain both the net displacement $B$ and the quadratic displacement $V$. For example, if a true center at zero is represented by two equally weighted centers at $-u$ and $u$, then
\[
\begin{aligned}
B
&=\underbrace{\left|\tfrac12(-u)+\tfrac12u\right|}_{\text{magnitude of the average shift}}
=0,\\
V
&=\underbrace{\tfrac12(-u)^2+\tfrac12u^2}_{\text{average squared displacement}}
=u^2.
\end{aligned}
\]
Although the average location equals the true center, $V$ remains positive because both fitted centers are displaced from it. To account also for differences in the mass assigned to emission locations, we compute the squared Wasserstein distance $W_2^2$.

For each separation, we therefore evaluate $D_\pi/h$, $B/h$, $V/h$, and $W_2^2/h$ from the recomputed observation densities and losses. We then estimate their separation exponents over the four smallest grid values, following the logarithmic regression in Section~\ref{sec:synthetic-exact}. According to Theorems~\ref{thm:overglobal} and~\ref{thm:overtransport}, the corresponding powers are $2s$ for $D_\pi/h$, $2s-1$ for $B/h$, and $2s-2$ for both $V/h$ and $W_2^2/h$.

For experiment (b), let $\eta$ control the size of the probability changes and center shifts. At $\eta=0$, the candidate and true models have the same observation and latent emission-center distributions. We measure how rapidly the observation error, joint probability error, and squared Wasserstein error grow as $\eta$ increases from zero. Comparing these initial rates of change isolates the effect of small perturbations before we examine their dependence on separation. For each fixed $\delta$, we compute
\[
\begin{aligned}
c_h(\delta)
&=\lim_{\eta\to0}\frac{h}{\eta},
&&\text{for observation error},\\
c_D(\delta)
&=\lim_{\eta\to0}\frac{D_I}{\eta},
&&\text{for joint probability error},\\
c_W(\delta)
&=\lim_{\eta\to0}\frac{W_{2,I}^2}{\eta},
&&\text{for squared Wasserstein error}.
\end{aligned}
\]
For each fixed separation $\delta$ with $c_h(\delta)>0$, the coefficient ratios $c_D(\delta)/c_h(\delta)$ and $c_W(\delta)/c_h(\delta)$ are the limits of $D_I/h$ and $W_{2,I}^2/h$, respectively, as $\eta$ tends to zero. We therefore compare the same latent-to-observation error ratios as before, but in the limit of small perturbations.

To obtain these ratios, we evaluate $c_h$ by numerical integration and $c_D$ directly from the changes in joint probabilities. For $c_W$, we solve a finite linear program that redistributes probability mass among emission-location sequences under the normalized block cost, as detailed in Appendix~\ref{app:synthetic-joint}. We then vary $\delta$ and compare the resulting ratios with the theoretical powers $\beta_r$ and $\beta_r-2$, where $\beta_r=\max\{2s,r(s-1)+s\}$. In the displayed cases with $s=3$ and $r=2,3,4$, these powers are $7,9,11$ for joint probability error and $5,7,9$ for squared Wasserstein error. Throughout this comparison, the limit in $\eta$ is taken at each fixed separation before different separations are compared.

\begin{table}[!tp]
\centering
\captionsetup{font=small,margin=0pt,width=\linewidth}
\fontsize{9}{11}\selectfont
\setlength{\tabcolsep}{4.5pt}
\renewcommand{\arraystretch}{1.15}
\arrayrulecolor{black}
\begin{tabular}{@{}cccccccccc@{}}
\toprule
& & \multicolumn{8}{c}{\textbf{Separation exponent}} \\
\cmidrule(lr){3-10}
& & \multicolumn{2}{c}{\boldmath$D_\pi/h$}
  & \multicolumn{2}{c}{\boldmath$B/h$}
  & \multicolumn{2}{c}{\boldmath$V/h$}
  & \multicolumn{2}{c}{\boldmath$W_2^2/h$} \\
\cmidrule(lr){3-4}
\cmidrule(lr){5-6}
\cmidrule(lr){7-8}
\cmidrule(lr){9-10}
$K=s_{\max}$ & $k$
& Theoretical & Estimated
& Theoretical & Estimated
& Theoretical & Estimated
& Theoretical & Estimated \\
\midrule[0.7pt]
2 & 3 & 4 & 3.989 & 3 & 2.989 & 2 & 1.989 & 2 & 1.989 \\
\arrayrulecolor{black!15}\midrule[0.25pt]\arrayrulecolor{black}
3 & 4 & 6 & 5.960 & 5 & 4.960 & 4 & 3.960 & 4 & 3.960 \\
\arrayrulecolor{black!15}\midrule[0.25pt]\arrayrulecolor{black}
4 & 5 & 8 & 7.908 & 7 & 6.908 & 6 & 5.908 & 6 & 5.908 \\
\bottomrule
\end{tabular}
\caption{\textbf{Separation exponents under over-specification.} Each true model has $K=s_{\max}$ states in one emission cluster, and the fitted model has $k=K+1$ states. Estimated exponents are the negative log--log slopes over the four smallest separation values. The theoretical exponents for $D_\pi/h$, $B/h$, $V/h$, and $W_2^2/h$ are $2s_{\max}$, $2s_{\max}-1$, $2s_{\max}-2$, and $2s_{\max}-2$, respectively. The loss $D_\pi$ compares initial probabilities after aggregating fitted states within each group.}
\label{tab:overfit_separation_exponents}
\end{table}

\paragraph{Results.}
As shown in Table~\ref{tab:overfit_separation_exponents}, permitting a state to split increases the separation exponent for probability recovery. In particular, for $s_{\max}=2,3,4$, the estimated exponents for $D_\pi/h$ are $3.989$, $5.960$, and $7.908$, respectively, close to the predicted values $4$, $6$, and $8$. Compared with Table~\ref{tab:exact_separation_exponents}, each estimate is approximately one larger than its exact-specification counterpart. For $s_{\max}=3$, for example, the exponent increases from $4.965$ to $5.960$, so the probability discrepancy grows more rapidly relative to Hellinger error as separation decreases.

Within the same configuration, however, the location ratios exhibit smaller separation powers. For $s_{\max}=3$, the estimated exponent decreases from $5.960$ for $D_\pi/h$ to $4.960$ for $B/h$ and $3.960$ for both $V/h$ and $W_2^2/h$. Across the other cluster sizes, the estimates show the same ordering, close to $2s_{\max}$, $2s_{\max}-1$, and $2s_{\max}-2$. This agreement with Theorems~\ref{thm:overglobal} and~\ref{thm:overtransport} numerically illustrates the sharpness of their one-time separation exponents, with the constructed pairs exhibiting the respective powers for probability, net displacement, and quadratic displacement. Although $W_2^2/h$ has approximately the same exponent as $W_1/h$ under exact specification, it controls a squared distance and therefore describes a different dependence on observation error.

\begin{figure}[t]
\centering
\includegraphics[width=0.89\textwidth]{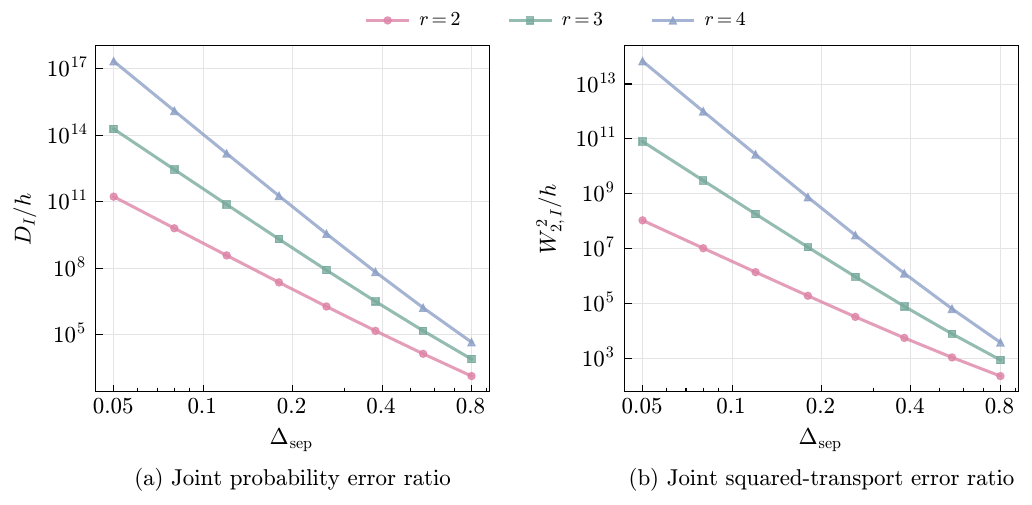}
\caption{\textbf{Joint recovery under over-specification.} (a) Joint probability-error ratios $c_D/c_h$. (b) Squared Wasserstein-error ratios $c_W/c_h$. The true model has three states with emission centers in one cluster. Curves correspond to $r=2,3,4$ selected times. These ratios are the limits of $D_I/h$ and $W_{2,I}^2/h$, respectively, as the perturbation size tends to zero at each fixed separation.}
\label{fig:synthetic-joint-conditioning}
\end{figure}

When recovery is extended across several times, Figure~\ref{fig:synthetic-joint-conditioning} shows a further increase in sensitivity to separation. For both joint probability and squared Wasserstein error, the curves become steeper as $r$ increases from $2$ to $4$. In particular, with $s=3$, the probability ratios exhibit powers close to $7$, $9$, and $11$, while the squared Wasserstein ratios exhibit powers close to $5$, $7$, and $9$. Each additional selected time therefore contributes approximately two inverse powers of separation in these configurations. Compared with the nearly parallel curves in Figure~\ref{fig:exact_selected_times}, this increase shows that joint recovery under over-specification can change the separation exponent itself. The observed powers are close to $\beta_r$ and $\beta_r-2$ in Theorems~\ref{thm:overglobal} and~\ref{thm:overtransport}, respectively, which numerically illustrates the sharpness of their joint separation exponents. The additional separation cost with increasing $r$ is thus visible in these constructions, whereas $r$ enters only the multiplicative constant in the exact-specification bounds.

\subsection{Recovery Errors and Center Matching across Sample Sizes}
\label{sec:synthetic-statistical}

The preceding comparisons show how latent discrepancies can grow relative to observation error as emission separation decreases. We now keep the true model fixed and increase the number of independent trajectories to examine whether improved observation fit is accompanied by improved latent recovery. We first evaluate models fitted by multistart EM, then examine location recovery in a simpler setting where global likelihood maximization can be verified. Finally, we investigate whether increasing the sample size makes it possible to pair every fitted emission center with a nearby true center.

\paragraph{Recovery as sample size increases.}
We generate data from a Gaussian HMM with three states, emission centers $(-2,0,2)$, variance one, and
\[
\pi^\star=(0.2,0.5,0.3),
\qquad
A^\star=
\begin{pmatrix}
0.85&0.10&0.05\\
0.10&0.80&0.10\\
0.05&0.15&0.80
\end{pmatrix}.
\]
For each $n\in\{128,512,2048\}$, we generate 24 independent data sets, each containing $n$ trajectories of length $T=6$. To compare exact specification and over-specification under the same observations, we fit models with $k=3,4,5$ to every data set. For each model, we use four initializations based on the observations and retain the parameter with the largest final training likelihood.

\begin{figure}[t]
\centering
\includegraphics[width=0.95\textwidth]{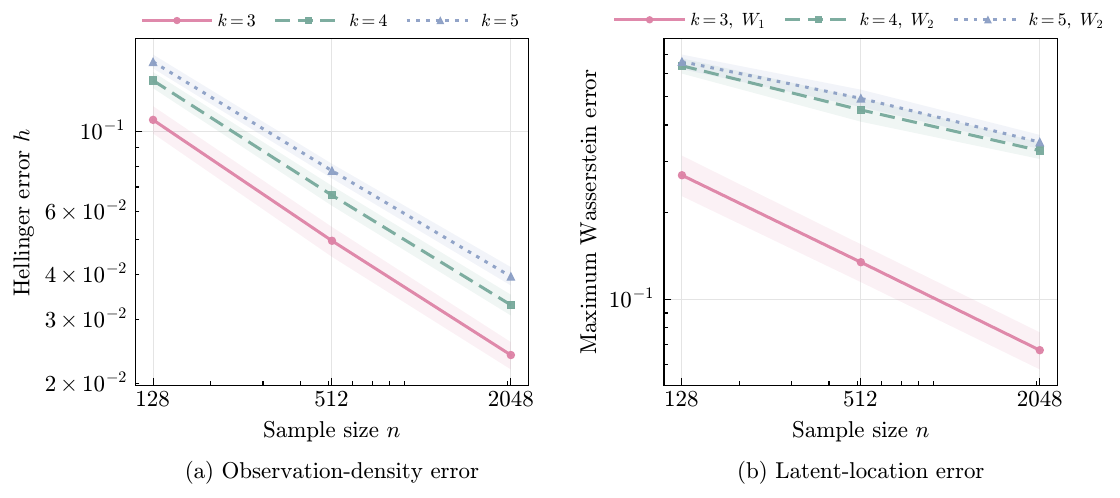}
\caption{\textbf{Observation and location errors as sample size increases.} Data are generated from a three-state HMM with $n$ independent trajectories of length $T=6$. (a) Mean Hellinger error for $k=3,4,5$. (b) Mean maximum Wasserstein error over the six times, using $W_1$ for $k=3$ and $W_2$ for $k=4,5$. The maximum is computed within each fit before averaging over $24$ independent data sets. Shaded bands show pointwise $95\%$ bootstrap confidence intervals.}
\label{fig:synthetic-statistical}
\end{figure}

For these fitted models, we evaluate the Hellinger error over the full trajectory and the aggregate transition-flow error $D_A$ over all five consecutive time pairs. To assess emission-center recovery, we compute the Wasserstein distances between the fitted and true emission-center distributions at each of the six times and take the maximum within each fit. We also report $R_{\mathrm{row}}$, the largest $\ell_1$ discrepancy between a fitted state's transition probabilities into the true-center groups and the transition row of its assigned true state. Different from $D_A$, this diagnostic retains each fitted source state separately. Figure~\ref{fig:synthetic-statistical} shows how the location errors decrease with sample size, using $W_1$ under exact specification and $W_2$ under over-specification to follow the corresponding theoretical bounds. Table~\ref{tab:statistical-recovery} additionally reports both distances for every fitted state count, allowing comparisons across state counts using the same loss. All summaries are means over the 24 data sets, and the figure includes pointwise $95\%$ bootstrap intervals for these means.

As shown in Figure~\ref{fig:synthetic-statistical}, both observation error and Wasserstein error decrease as more trajectories become available. Under exact specification, Table~\ref{tab:statistical-recovery} shows that the mean Hellinger error decreases from $0.107$ at $n=128$ to $0.024$ at $n=2048$, while the mean maximum $W_1$ error decreases from $0.268$ to $0.067$. Over the same sample sizes, the mean aggregate-flow error falls from $0.196$ to $0.047$. The improved observation fit is therefore accompanied by more accurate recovery of emission locations and transitions between the true states. Under over-specification, the errors also decrease, although the additional states incur a recovery cost in this experiment. For $k=5$, the mean Hellinger error falls from $0.155$ to $0.040$, while $D_A$ falls from $0.290$ to $0.074$ and the mean maximum $W_2$ error from $0.661$ to $0.349$. At $n=2048$, the mean maximum $W_1$ error is $0.150$, compared with $0.067$ under exact specification, and the aggregate-flow error also remains larger. These comparisons show that allowing more fitted states can affect latent recovery even when all models achieve small observation errors. The remaining individual-row discrepancies in Table~\ref{tab:statistical-recovery} motivate the constructions in Section~\ref{sec:synthetic-identification}, which show that these discrepancies can persist even when aggregate transition flows are recovered exactly.
\begin{table}[t]
\centering
\small
\captionsetup{font=small,margin=0pt,width=\linewidth}
\setlength{\tabcolsep}{9pt}
\renewcommand{\arraystretch}{1.12}
\arrayrulecolor{black}
\begin{tabular}{@{}ccccccc@{}}
\toprule
& & \textbf{Observation fit}
& \multicolumn{2}{c}{\textbf{Transition recovery}}
& \multicolumn{2}{c}{\textbf{Location recovery}} \\
\cmidrule(lr){3-3}
\cmidrule(lr){4-5}
\cmidrule(lr){6-7}
$n$ & $k$ & $h$
& $D_A$ & $R_{\mathrm{row}}$
& $\max_t W_1$ & $\max_t W_2$ \\
\midrule[0.7pt]
\multirow{3}{*}{128}
& 3 & 0.107 & 0.196 & 0.139 & 0.268 & 0.613 \\
& 4 & 0.138 & 0.255 & 0.441 & 0.363 & 0.640 \\
& 5 & 0.155 & 0.290 & 0.663 & 0.421 & 0.661 \\
\arrayrulecolor{black!15}\midrule[0.3pt]\arrayrulecolor{black}
\multirow{3}{*}{512}
& 3 & 0.050 & 0.094 & 0.061 & 0.135 & 0.442 \\
& 4 & 0.066 & 0.114 & 0.439 & 0.196 & 0.451 \\
& 5 & 0.078 & 0.136 & 0.576 & 0.228 & 0.494 \\
\arrayrulecolor{black!15}\midrule[0.3pt]\arrayrulecolor{black}
& 3 & 0.024 & 0.047 & 0.031 & 0.067 & 0.305 \\
\rowcolor{black!5}
& 4 & 0.033
& \textbf{0.061} & \textbf{0.319}
& 0.122 & 0.326 \\
\rowcolor{black!5}
\multirow{-3}{*}{2048}
& 5 & 0.040
& \textbf{0.074} & \textbf{0.421}
& 0.150 & 0.349 \\
\bottomrule
\end{tabular}
\caption{\textbf{Statistical recovery as the number of trajectories increases.} The true HMM has three states, and each data set contains $n$ independent trajectories of length $T=6$. Entries are means over $24$ independent data sets, shared across $k=3,4,5$. The loss $D_A$ measures aggregate transition-flow error, while $R_{\mathrm{row}}$ measures the largest individual-row discrepancy after destination states are grouped. Each Wasserstein error is maximized over the six times within each fit before averaging over data sets. Bold entries in the shaded rows show that, at $n=2048$, individual transition-row errors remain much larger than aggregate transition-flow errors when fitting $k=4,5$ states to the true three-state model.}
\label{tab:statistical-recovery}
\end{table}

\paragraph{Sample-size rates for emission-center recovery.}
To validate the $n^{-1/2}$ and $n^{-1/4}$ rates for emission-center recovery in Remark~\ref{remark:single_true_state}, we next compare the global maximum likelihood estimators in two simple Gaussian models. More concretely, we generate independent observations from $\mathcal N(0,1)$ and compare fitting a single Gaussian $\mathcal N(m,1)$, with $m\in[-6,6]$, against fitting the symmetric mixture $\frac12\mathcal N(-a,1)+\frac12\mathcal N(a,1)$, with $a\in[0,6]$. The first model is exactly specified, whereas the second permits the true center at zero to split into two fitted centers. Both models can be represented as HMMs with independent hidden states. For each $n\in\{128,512,2048,8192\}$, we generate 256 independent data sets, each containing $n$ trajectories of length $T=6$, or equivalently $6n$ independent scalar observations. The global MLE in the single-Gaussian model is the sample mean projected onto $[-6,6]$. In the symmetric-mixture model, we obtain the global maximizer by optimizing over the scalar $a$, as detailed in Appendix~\ref{app:synthetic-single}.Since the true emission-center distribution is concentrated at zero, the corresponding recovery errors are $W_1=|\widehat m|$ under exact specification and $W_2=\widehat a$ under over-specification. We compare their empirical means with reference curves proportional to $n^{-1/2}$ and $n^{-1/4}$, respectively. The curves are scaled to pass through the first empirical mean in each case, allowing us to assess the predicted decay rates without specifying their multiplicative constants.

\begin{figure}[t]
\centering
\includegraphics[width=.94\textwidth]{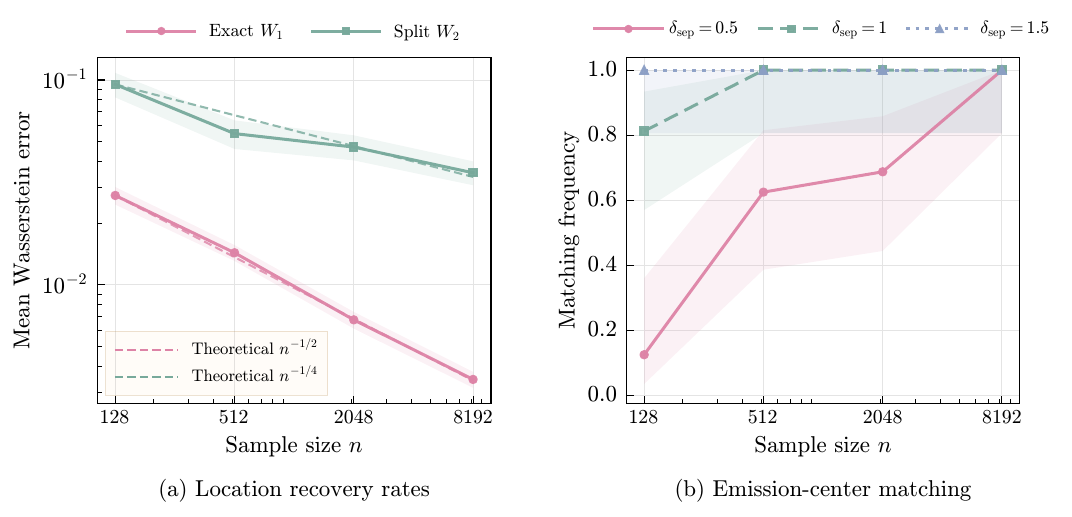}
\caption{\textbf{Location recovery and emission-center matching.} (a) Mean $W_1$ error for a single-Gaussian fit and mean $W_2$ error for a symmetric two-Gaussian fit to data from $\mathcal N(0,1)$, using global MLEs over $256$ independent data sets at each $n$. Dashed curves show the $n^{-1/2}$ and $n^{-1/4}$ reference rates, respectively, scaled to the first empirical means. (b) Matching success rates over $16$ independent data sets at each $n$ and separation, with tolerance $\Delta_{\mathrm{sep}}/4$. Shaded bands give $95\%$ bootstrap intervals in (a) and $95\%$ Wilson intervals in (b).}
\label{fig:synthetic-rates-matching}
\end{figure}

 Figure~\ref{fig:synthetic-rates-matching}(a) shows that the single-Gaussian error closely follows the $n^{-1/2}$ reference, while the symmetric-mixture error decreases more slowly, broadly following the $n^{-1/4}$ reference. These results support the theoretical distinction between exact specification and over-specification: allowing a true center to split can slow emission-center recovery even when the fitted parameters globally maximize the likelihood.

\paragraph{Emission-center matching.}
While the preceding comparison measures the size of the location error, we also examine whether every true center has a nearby fitted counterpart. We generate two-state Gaussian HMMs with variance one, initial probabilities $(0.4,0.6)$, transition matrix
\[
A^\star=
\begin{pmatrix}
0.9&0.1\\
0.15&0.85
\end{pmatrix},
\]
and emission centers $(-\delta_{\mathrm{sep}}/2,\delta_{\mathrm{sep}}/2)$, where $\delta_{\mathrm{sep}}\in\{0.5,1,1.5\}$ is their minimum separation. For each gap and $n\in\{128,512,2048,8192\}$, we fit two-state models to 16 independent data sets of trajectories with length $T=6$. For each fit, we check whether a permutation pairs both fitted centers with the true centers within distance $\Delta_{\mathrm{sep}}/4$, where $\Delta_{\mathrm{sep}}=\min\{1,\delta_{\mathrm{sep}}\}$. 

\noindent We report the proportion satisfying this location condition on the right of Figure~\ref{fig:synthetic-rates-matching}, with Wilson intervals over the 16 repetitions. \footnote{The reported frequency measures how often the fitted emission centers satisfy the matching condition.} At $n=128$, only $2$ of $16$ fits satisfy the location condition for gap $0.5$, compared with $13$ of $16$ for gap $1$ and all $16$ for gap $1.5$. For the smallest gap, the frequency increases to $10/16$ at $n=512$, $11/16$ at $n=2048$, and $16/16$ at $n=8192$. The closely spaced centers therefore require more observations before matching becomes reliable in this experiment.

\subsection{Aggregate-Flow and Individual-Row Recovery}
\label{sec:synthetic-identification}

The experiment presented in Section~\ref{sec:synthetic-statistical} shows that aggregate transition-flow error decreases with sample size. Under over-specification, however, several fitted states may represent the same true state, and accurate transitions between their groups need not determine the transition probabilities of each fitted state. We examine this distinction using the same estimated models and then construct an exact example to explain why individual-row ambiguity can persist even when observation error vanishes.

\noindent\emph{Empirical comparison.} 
We reuse the $k=4,5$ fits from Section~\ref{sec:synthetic-statistical} and assign each fitted state to its nearest true emission center. To compare individual transition rows, we combine destination states within each group while retaining every fitted source state separately. With the nearest-center assignment $c(i)$ defined in Section~\ref{sec:latent_recovery_loss}, we compute
\[
R_{\mathrm{row}}
=
\max_{i\in[k]}
\sum_{j=1}^{K}
\Bigl|
\underbrace{
\sum_{\ell:\,c(\ell)=j} A_{i\ell}
}_{\substack{
\text{transition probability}\\
\text{from fitted state }i\text{ to group }j
}}
-
\underbrace{
A^\star_{c(i),j}
}_{\substack{
\text{corresponding true}\\
\text{transition probability}
}}
\Bigr|.
\]

This diagnostic considers whether every fitted source state has the correct probabilities of moving to the true-center groups. To assess transitions at the group level, $D_A$ instead combines both source and destination states, weighting each source state by its occupancy probability. It compares the resulting matrices of consecutive group probabilities, with the maximum taken over all five consecutive time pairs. Since states in the same source group can have opposing row discrepancies, these discrepancies may average out in the aggregate flow while remaining visible in $R_{\mathrm{row}}$. Both quantities are evaluated under the original fitted HMM.

\begin{figure}[t]
\centering
\includegraphics[width=.80\textwidth]{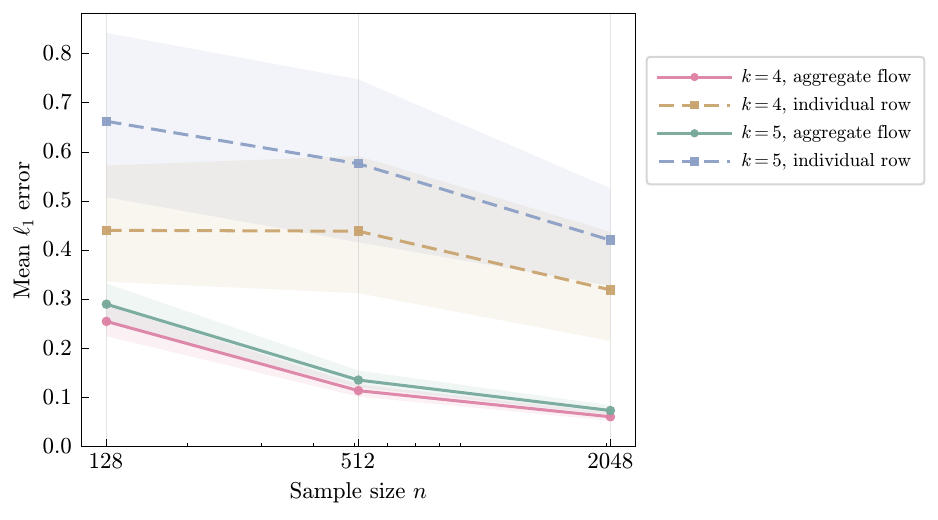}
\caption{\textbf{Aggregate-flow and individual-row recovery under over-specification.} Results use the $k=4,5$ fits from the three-state experiment in Section~\ref{sec:synthetic-statistical}. Each curve reports a mean over the same 24 independent data sets, with pointwise $95\%$ bootstrap confidence intervals. The loss $D_A$ combines source and destination states into groups and weights sources by their occupancies. In contrast, $R_{\mathrm{row}}$ combines destinations but retains each fitted source state separately, taking the largest row discrepancy. All fitted source states enter this diagnostic.}
\label{fig:synthetic-overfit}
\end{figure}

Figure~\ref{fig:synthetic-overfit} shows that the aggregate-flow error decreases substantially with sample size, whereas the individual-row discrepancy remains appreciable. For $k=4$, the mean $D_A$ decreases from $0.255$ to $0.061$ as $n$ increases from $128$ to $2048$, while the mean $R_{\mathrm{row}}$ decreases from $0.441$ to $0.319$. For $k=5$, the corresponding changes are from $0.290$ to $0.074$ and from $0.663$ to $0.421$. By comparison, the exact-fit row discrepancy in Table~\ref{tab:statistical-recovery} falls from $0.139$ to $0.031$. More observations therefore improve aggregate recovery without producing a comparable improvement in the worst individual row of the over-specified models.

\noindent\emph{Individual-row ambiguity.} In Example~\ref{ex:state-splitting}, summing transition probabilities over destination substates gives the correct group-transition probabilities for every fitted source state. The following construction shows that these probabilities can also differ from their true counterparts even when the fitted and true observation distributions are identical. Consider a true two-state HMM with distinct centers $a,b$, initial probabilities $(1/2,1/2)$, and both transition rows equal to $(1/2,1/2)$. We split the state at $a$ into two candidate states and set
\[
\widetilde{\pi}=(1/4,1/4,1/2),
\qquad
\widetilde{A}=
\begin{pmatrix}
1/2&1/2&0\\
0&0&1\\
1/4&1/4&1/2
\end{pmatrix},
\qquad
\widetilde{\mu}=(a,a,b).
\]
The candidate model replaces true state 1, whose emission center is $a$, by candidate states 1 and 2 with the same center. Candidate state 3 has center $b$ and corresponds to true state 2. We therefore place candidate states 1 and 2 in group 1 and candidate state 3 in group 2. Although candidate states 1 and 2 have the same emission distribution, their transition probabilities differ. State 1 moves to group 1 with probability one, whereas state 2 moves to group 2 with probability one. Their group-transition rows are $(1,0)$ and $(0,1)$, both different from the true row $(1/2,1/2)$.

However, the transition probabilities between groups are still $(1/2,1/2)$, exactly as in the true model. The two candidate states in group 1 have equal initial probabilities, and every row of $\widetilde A$ assigns equal probabilities to entering them. Conditional on any group history up to the current time that ends in group 1, the two states therefore have probabilities $1/2$ each. The next-group probabilities are $\frac12(1,0)+\frac12(0,1)=(1/2,1/2).$ Candidate state 3 also moves to the two groups with probabilities $(1/2,1/2)$. Since the initial group probabilities are $(1/2,1/2)$, the grouped process has exactly the same distribution as the true hidden-state process, despite the differences in the individual candidate transition rows.

Because the grouped process has the same distribution as the true hidden-state process and each group has the corresponding true emission distribution, the two models induce the same joint observation distribution for every sequence length. Hence $h=0$, $D_A=0$, and all aggregate path discrepancies vanish, but $R_{\mathrm{row}}=1$. The transition probabilities of individual split states therefore remain unidentifiable, even though both models have identical observation distributions and aggregate state-sequence probabilities. Moreover, every candidate state has positive occupancy at every time. The row discrepancy therefore reflects ambiguity in how a true state is represented by several fitted states, even when neither observation error nor vanishing occupancy is present.

The previous example shows that aggregate transition flows can be recovered even when individual fitted rows remain ambiguous. A separate question is whether these flows determine aggregate state-sequence probabilities at more than two times. In the eight-state construction of Example~\ref{ex:aggregate-memory}, direct path enumeration gives zero discrepancy on every subset of at most three of the first four times, while $D_{\{1,2,3,4\}}=1$. The discrepancy is therefore visible only when all four times are considered jointly. This explains why the over-specification analysis controls $D_I$ directly for the selected times, rather than reconstructing all joint probabilities from aggregate transition flows.

\section{Practical Implications}
\label{sec:practical}
Our recovery bounds distinguish the accuracy of the observation density from that of the latent quantities used to interpret a fitted Gaussian HMM. The relation between them depends on emission geometry and on whether the recovery targets are individual-state parameters or aggregate probabilities and transition flows. These results provide a theoretical basis for interpreting individual fitted states and transition probabilities in practical applications, and for assessing whether groups of states with similar emissions retain stable aggregate probabilities and transition flows across fits.

\begin{enumerate}
\item \textbf{Interpreting individual states and transition probabilities.} In practice, emission geometry, matching stability across fits, and state probabilities guide the interpretation of individual fitted states and the estimated transition probabilities between them.

\emph{Emission geometry.} Since a close fit to the observations can coexist with substantial latent recovery error, practitioners should examine the geometry of the fitted emission centers before assigning separate interpretations to their states. As illustrated by Figure~\ref{fig:cluster-concept} and quantified in Theorems~\ref{thm:globalexact} and~\ref{thm:paired}, configurations with the same minimum separation can have different separation exponents when different numbers of centers approach one another. Accordingly, this examination should consider both the distances between centers and the sizes of groups of close centers. For an empirical summary of these groups, complete-linkage clustering provides control of the within-group diameter. By contrast, a distance-threshold graph can join a chain of close pairs whose overall diameter is large, so its connected components need not satisfy the theoretical diameter condition.

\emph{Matching across fits.} Under exact specification, Corollary~\ref{cor:localization} gives a sufficient density-error condition for matching every fitted center to a true center. Since the true centers are unavailable, practitioners can examine whether the correspondence between fitted centers remains stable across bootstrap resamples, random initializations, or small changes in the training data, after accounting for label permutations. If this correspondence repeatedly changes, or centers repeatedly merge and split despite comparable observation fit, separate state interpretations require additional evidence.

\emph{Transition probabilities.} In addition to matching emission centers across fits, practitioners should examine estimated state occupancies and transition flows before interpreting individual transition rows. This is especially important for rarely occupied states, since Corollary~\ref{cor:rows} shows that recovering individual rows from flows introduces an inverse initial-state probability factor. Consequently, accurate flow recovery can provide a weaker guarantee for the outgoing transition probabilities of such states. To assess whether these probabilities support a stable transition pattern, practitioners can compare them across the same fits used to examine center matching and report the corresponding occupancies and transition flows alongside the rows.

\item \textbf{Interpreting groups of similar states.} Under over-specification, several fitted states may represent the same true emission center, so their individual transition rows can differ even when the observation distributions are identical. As Section~\ref{sec:synthetic-identification} shows, this ambiguity persists even when every fitted state has positive occupancy, so interpreting each fitted state separately may obscure the behavior shared by the group. Instead, practitioners can examine aggregate probabilities, transition flows, and emission-center distributions, whose recovery is controlled by Theorems~\ref{thm:overglobal} and~\ref{thm:overtransport}. For interpretations involving several times, this examination should also include the joint probabilities of the corresponding group sequences: grouping can destroy the Markov property, so transition flows between groups alone need not determine these probabilities, as illustrated by Example~\ref{ex:aggregate-memory}.

To compare these aggregate quantities, our theoretical losses group fitted states by assigning each to its nearest true emission center. As these centers are unavailable in practice, Proposition~\ref{prop:reference-grouping-stability} bounds the additional grouping error caused by replacing them with sufficiently accurate estimated reference centers. Under its reference accuracy condition, the distance margins identify fitted states whose group assignments may change under reference perturbations. Practitioners can therefore perturb the proposed reference centers and compare the resulting aggregate probabilities and transition flows to assess how sensitive their interpretation is to the chosen grouping.
\end{enumerate}

These implications are not specific to a single fitting procedure. Because the inverse inequalities apply to any fitted estimator with suitable observation-density error control, the same reasoning can inform the interpretation of models obtained by maximum likelihood, Baum--Welch, spectral or moment-based methods, and Bayesian procedures.

\section{Conclusion}
\label{sec:conclusion}
In this paper, we establish latent recovery guarantees for multivariate Gaussian HMMs with a common, known emission covariance matrix. The guarantees relate observation-density error to latent recovery error through inverse powers of emission separation, making explicit how recovery becomes more demanding as emission centers approach one another. For mutually separated clusters, the shared emission centers across time allow us to obtain bounds whose separation exponents depend on the largest cluster of approaching centers rather than the total number of states or the exponentially many hidden paths. These bounds therefore distinguish configurations with the same minimum separation and quantify the advantage of smaller, mutually separated clusters. Under exact specification, direct transport interpolation and the Markov factorization yield global bounds with sharp separation exponents $2s_{\max}-1$ for state-sequence probability errors and $2s_{\max}-2$ for $W_1$ recovery of joint emission-center distributions, neither of which increases with the number of selected times. Without cluster assumptions, the same bounds hold with $s_{\max}=K$. When the observation-density error is sufficiently small, these bounds also guarantee emission-center matching, after which the probability and Wasserstein recovery bounds have constants independent of the minimum true initial probability. Under over-specification, several fitted states may represent the same true state, so we instead control aggregate probabilities, transition flows, and emission-center distributions that accommodate state splitting. In this case, the separation exponents for joint probability and $W_2$ recovery can increase with the number of selected times. Since these inverse bounds apply to any estimator with suitable observation-density error control, they also provide a basis for deriving statistical recovery rates. For the global MLE, the Hellinger estimation guarantee yields latent recovery rates of the form $\mathcal{O}_{\mathbb P}\!\left(\Delta_{\mathrm{sep}}^{-\alpha_r}[d\log(nTd)/n]^{\nu_r/2}\right)$, with explicit separation exponents determined by cluster size, the recovery target, and model specification. In these rates, the sample-size dependence is $n^{-1/2}$ for errors in the probabilities of individual states and state sequences and in transition flows under exact specification, up to logarithmic factors. Under over-specification, the corresponding grouped probability and transition-flow bounds retain this dependence. By comparison, the bounds for probability-weighted emission-center displacement and Wasserstein error involve $n^{-1/2}$ under exact specification using $W_1$, and $n^{-1/4}$ under over-specification using $W_2$. Moreover, the underlying density guarantee remains valid with coincident fitted centers and zero fitted probabilities, so the statistical recovery analysis also covers these degenerate fitted models. To examine the separation dependence of the recovery guarantees, we first compare controlled model pairs with directly specified parameters. These comparisons illustrate how the predicted separation powers depend on cluster size and the number of selected times. We then examine recovery from data as the sample size increases, including emission-center matching under exact specification. Under over-specification, however, the fitted models show that decreasing aggregate-flow error can coexist with substantial individual-row discrepancies. Indeed, our constructions show that these discrepancies can persist even when the fitted and true observation distributions are identical. Beyond these statistical guarantees and empirical findings, our analysis offers a quantitative foundation for interpreting latent regimes in applications such as genomics, biophysics, and neuroscience. In such settings, the scientific value of a fitted HMM lies in what its states reveal about underlying observation patterns and their temporal organization. This perspective supports more informed interpretation of fitted models and lays the groundwork for future research on prediction and conditional inference with explicit latent recovery guarantees.

\clearpage
\appendix
\startcontents[appendices]

\begin{center}
  {\Large\bfseries
  Supplement for\\[0.3em]
  ``The Statistical Geometry of Latent Recovery in \\ Hidden Markov Models''\par}
\end{center}

\vspace{1em}

\begingroup

\hypersetup{linkcolor=navy}

\fontsize{11}{15}\selectfont
\contentsmargin{2.5em} 

\titlecontents{section}
  [1.6em]
  {\addvspace{0.9em}\color{navy}\bfseries}
  {\contentslabel{1.6em}}
  {}
  {\hfill\contentspage}

\titlecontents{subsection}
  [4em]
  {\addvspace{0.15em}\color{navy}}
  {\contentslabel{2.4em}}
  {}
  {\titlerule*[0.75em]{.}\contentspage}

\titlecontents{subsubsection}
  [7.2em]
  {\addvspace{0.15em}\color{navy}}
  {\contentslabel{3.2em}}
  {}
  {\titlerule*[0.75em]{.}\contentspage}

\noindent{\large\bfseries Table of Contents\par}
\vspace{0.3em}
\hrule height 0.4pt
\vspace{0.6em}

\printcontents[appendices]{}{1}{\setcounter{tocdepth}{2}}

\vspace{0.6em}
\hrule height 0.4pt

\endgroup

\clearpage

\section{Gaussian moment extraction and elementary transport}
\label{app:gaussian-moments}
This appendix section establishes the analytic tools used in the inverse bounds. We first justify the reduction to identity emission covariance and compare transport distances before and after the transformation. We then show that the full-block Hellinger error controls polynomial moments of the latent locations, with constants independent of the ambient observation dimension and the number of possible latent paths. Finally, we relate errors in locations and aggregate probabilities to Wasserstein distances through a coupling on the true support.

\subsection{Covariance normalization and comparison of transport metrics}
\label{app:whitening}
We first verify that the covariance normalization in Section~\ref{prelim:hmm} preserves Hellinger distance and likelihood maximizers. The following lemma also compares Wasserstein distances before and after the transformation.
\begin{lemma}[Whitening reduction for Gaussian HMMs]
\label{lem:whitening_reduction}
Let $\Sigma\in\mathbb R^{d\times d}$ be a fixed symmetric
positive-definite matrix. Consider Gaussian HMMs with parameter
$\theta=(\pi,A,\mu_1,\ldots,\mu_k)$ and common emission
covariance $\Sigma$. Define $S=\Sigma^{-1/2}$, $\widetilde\theta=(\pi,A,S\mu_1,\ldots,S\mu_k)$,
and, for every $r\geq1$, let $S_r(x_1,\ldots,x_r)=(Sx_1,\ldots,Sx_r)$. Write $p_{\theta,r}^{\Sigma}$ and
$p_{\widetilde\theta,r}^{I_d}$ for the corresponding
length-$r$ observation densities. Then $Y_t=SX_t$ is a
Gaussian HMM with parameter $\widetilde\theta$ with common
emission covariance $I_d$, and
\[
p_{\widetilde\theta,r}^{I_d}(y)
=
|\det S|^{-r}
p_{\theta,r}^{\Sigma}(S_r^{-1}y).
\]
Moreover, the following conclusions hold.
\begin{enumerate}
    \item \textnormal{(Hellinger invariance)}
    For any two parameters $\theta,\theta'$ and every $r\geq1$,
    \[
    h\!\left(p_{\theta,r}^{\Sigma},
             p_{\theta',r}^{\Sigma}\right)
    =
    h\!\left(p_{\widetilde\theta,r}^{I_d},
             p_{\widetilde\theta',r}^{I_d}\right).
    \]

    \item \textnormal{(Correspondence of likelihood maximizers)}
    For $n$ independent length-$T$ trajectories, let
    $L_X(\theta)$ and $L_Y(\widetilde\theta)$ denote the
    likelihoods based on the original and transformed
    observations, respectively. Then
    \[
    L_Y(\widetilde\theta)
    =
    |\det S|^{-nT}L_X(\theta).
    \]
    Consequently, for any parameter class $\Theta$ and its
    transformed image
    $\widetilde\Theta=\{\widetilde\theta:\theta\in\Theta\}$,
    the map $\theta\mapsto\widetilde\theta$ induces a
    bijection between the sets of global likelihood
    maximizers over $\Theta$ and $\widetilde\Theta$.

    \item \textnormal{(Transport identity and comparison)}
    Fix $p\in\{1,2\}$ and $r\geq1$. Throughout this lemma, $W_{p,r}$ is the Wasserstein distance associated with
    $d_{p,r}(x,y)=\{r^{-1}\sum_{\ell=1}^r\|x_\ell-y_\ell\|_2^p\}^{1/p}$, exactly as in equation~\eqref{eq:distance}.
    Let $G,H$ be probability measures on $(\mathbb R^d)^r$
    with finite $p$th moments, and define
    \[
    \widetilde G=(S_r)_\#G,
    \qquad
    \widetilde H=(S_r)_\#H.
    \]
    Then, by defining 
    $\lambda_-=\lambda_{\min}(\Sigma)$ and
    $\lambda_+=\lambda_{\max}(\Sigma)$, we have
    \begin{equation}
    \label{eqn:W_2_equivalent}
    \frac{1}{\sqrt{\lambda_+}}W_{p,r}(G,H)
    \leq
    W_{p,r}(\widetilde G,\widetilde H)
    \leq
    \frac{1}{\sqrt{\lambda_-}}W_{p,r}(G,H).
    \end{equation}
\end{enumerate}
\end{lemma}

\begin{proof}
The transformation leaves the hidden chain unchanged and takes each conditional
emission $N_d(\mu_i,\Sigma)$ to $N_d(S\mu_i,I_d)$. The block
change of variables has Jacobian $|\det S|^r$, giving the asserted
density identity. Applying this identity to both densities in the Hellinger
affinity proves invariance; applying it to $n$ independent length-$T$
trajectories multiplies every likelihood by the same factor
$|\det S|^{-nT}$. Hence maximizers correspond bijectively over a
parameter class and its \emph{transformed image}.

For the transport comparison, the eigenvalue bounds give, for any $x,y$,
\[
\lambda_+^{-p/2}d_{p,r}(x,y)^p
\le d_{p,r}(S_rx,S_ry)^p
\le \lambda_-^{-p/2}d_{p,r}(x,y)^p.
\]
The map $\gamma\mapsto(S_r,S_r)_\#\gamma$ is a bijection from
$\Pi(G,H)$ onto $\Pi((S_r)_\#G,(S_r)_\#H)$, because $S_r$
is invertible. Integrate the pointwise inequality, take infima over
couplings and then take $p$th roots. This proves
\eqref{eqn:W_2_equivalent}, with the same normalized factor $1/r$
for $p=1$ and $p=2$.
\end{proof}

For fixed $\Sigma$, the lemma transfers transport bounds
between the original and whitened coordinates with only
a change in constants. Throughout the subsequent proofs,
separation, cluster geometry, radius constraints, and
nearest-center assignments refer to the whitened centers.
In the original coordinates, these correspond to the
Mahalanobis geometry induced by $\Sigma^{-1}$.
In particular, nearest-center assignments after whitening
need not agree with those defined by the original
Euclidean metric.

\subsection{Hellinger distance decreases under marginalization}
The inverse bounds use observation distributions at selected times within the observed block. To this end, the following lemma shows that the Hellinger error at any selection of times is bounded above by the Hellinger error of the full observation block.

\begin{lemma}[Marginalization]\label{lem:marginalization}
Let $p$ and $q$ be joint densities, and let $p_I$ and $q_I$ be their marginals on a subset $I$ of coordinates. Then
\[
h(p_I,q_I)\le h(p,q).
\]
\end{lemma}

\begin{proof}
Write the retained and discarded coordinates as $(x,y)$. Cauchy--Schwarz
implies, for almost every $x$,
\[
\int\sqrt{p(x,y)q(x,y)}\,dy
\le\sqrt{\int p(x,y)\,dy\int q(x,y)\,dy}
=\sqrt{p_I(x)q_I(x)}.
\]
Integrating in $x$ and using the affinity identity
$h^2(p,q)=1-\int\sqrt{pq}$ gives $h^2(p_I,q_I)\le h^2(p,q)$.
\end{proof}

\subsection{An inverse heat identity with a uniform moment bound}
The preceding lemma controls the Hellinger error of any retained observation marginal by the Hellinger error of full block. To recover information about the latent locations, we next relate their polynomial moments to these observation densities. For a given polynomial $U$ of the latent locations, we construct a polynomial $g$ of the observations whose conditional expectation equals $U$. Applying Cauchy--Schwarz to the difference between the expectations of $g$ produces a Hellinger factor and a factor involving its second moments under the two observation densities. We therefore bound these second moments uniformly over the mixing distributions to obtain the desired inequality.

\begin{lemma}[Polynomial moments of Gaussian location mixtures]\label{lem:moment-extraction}
Let $S$ be a Euclidean space of dimension $q\le K+k$. Let $H,H'$ be probability measures supported on $B_S(0,R)^r$, and let $p_H,p_{H'}$ be the corresponding Gaussian location-mixture densities on $S^r$ with identity covariance. If $U$ is a polynomial of total degree at most $D$, then
\begin{equation}
\left|\int U\,d(H-H')\right|
\le C(q,r,R,D)\,\|U\|_{\infty,B_S(0,R)^r}\,h(p_H,p_{H'}).
\label{proof:A1}
\end{equation}
The constant does not depend on the numbers of support points of $H$ and $H'$. For a HMM block with $r$ retained times, one may take
\[
S=\operatorname{span}\{\mu_1,\ldots,\mu_k,a_1,\ldots,a_K\},
\]
so the constant can be bounded by a quantity depending only on $K,k,r,R,D$.
\end{lemma}

\begin{proof}
Choose orthonormal coordinates on $S^r$, put $m=qr$, and write $\mathcal K=B_S(0,R)^r$. If $q=0$ or $R=0$, the mixing measures coincide; if $D=0$, their difference integrates the constant polynomial to zero. We may therefore assume $q,D\ge1$ and $R>0$.

\paragraph{The inverse heat identity.}
Let $Z\sim N_m(0,I_m)$ and $\Delta_L=\sum_{b=1}^m\partial_b^2$. For a polynomial $F$ of total degree at most $D$, its finite Taylor expansion and the Gaussian moments give
\begin{equation}
\mathbb EF(u+Z)
=\sum_{2|\alpha|\le D}\frac{\partial^{2\alpha}F(u)}{2^{|\alpha|}\alpha!}
=\sum_{j=0}^{J}\frac{\Delta_L^jF(u)}{2^j j!},
\qquad J=\lfloor D/2\rfloor.
\label{proof:A2}
\end{equation}
Here $\alpha$ is a multi-index, and the second equality is the multinomial formula. On this finite-dimensional polynomial space, $N=\Delta_L/2$ satisfies $N^{J+1}=0$. Define
\[
g=\sum_{j=0}^{J}\frac{(-1)^j}{j!}N^jU.
\]
For $1\le a\le J$, the coefficient of $N^a$ in the product of the two truncated operator sums is
\[
\sum_{j=0}^{a}\frac{(-1)^j}{j!(a-j)!}
=\frac{(1-1)^a}{a!}=0.
\]
The constant coefficient is one, and all higher powers annihilate $U$. Thus equation~\eqref{proof:A2} gives $\mathbb Eg(u+Z)=U(u)$. No identity between unrestricted, truncated formal power series is being asserted.

\paragraph{A coefficient bound and uniform integrability.}
For $F(x)=\sum_\alpha f_\alpha x^\alpha$, write $\|F\|_{\mathrm{coef},1}=\sum_\alpha|f_\alpha|$. Set $a=R/\sqrt q$ and
$\xi_i=-a+2ai/D$ for $0\le i\le D$. The product grid $\{\xi_0,\ldots,\xi_D\}^m$ lies in $\mathcal K$, since each of its $r$ blocks has Euclidean norm at most $R$. Let
\[
L_i(t)=\prod_{\substack{0\le b\le D\\b\ne i}}
       \frac{t-\xi_b}{\xi_i-\xi_b}.
\]
The numerator has coefficient norm at most $(1+a)^D$. The denominator has absolute value $(2a/D)^D i!(D-i)!$. Consequently,
\[
A_{a,D}:=\sum_{i=0}^{D}\|L_i\|_{\mathrm{coef},1}
\le\frac{(1+a)^D D^D}{a^D D!}.
\]
Since the degree of $U$ in each coordinate is at most $D$, tensor-product interpolation is exact and implies
\[
\|U\|_{\mathrm{coef},1}
\le A_{a,D}^{m}\|U\|_{\infty,\mathcal K}.
\]
For every polynomial of total degree at most $D$,
$\|\Delta_L F\|_{\mathrm{coef},1}\le mD^2\|F\|_{\mathrm{coef},1}$.
Iteration and the finite definition of $g$ therefore give
\[
\|g\|_{\mathrm{coef},1}
\le e^{mD^2/2}A_{a,D}^{m}\|U\|_{\infty,\mathcal K}.
\]
It follows that $|g(x)|\le C\|U\|_{\infty,\mathcal K}(1+\|x\|_2^D)$, with $C$ depending only on $q,r,R,D$. Since $\|u\|_2\le R\sqrt r$ on $\mathcal K$ and
$\mathbb E\|Z\|_2^{2D}=\prod_{j=0}^{D-1}(m+2j)<\infty$, we obtain
\[
\sup_{u\in\mathcal K}\mathbb E|g(u+Z)|^2
\le B(q,r,R,D)\|U\|_{\infty,\mathcal K}^{2}.
\]
This bound proves, before any interchange of integrals, absolute integrability of $g$ under either observation mixture. Fubini's theorem and the inverse heat identity now yield
\[
\int g(x)p_H(x)\,dx=\int U\,dH,
\qquad
\int g(x)p_{H'}(x)\,dx=\int U\,dH'.
\]
Factoring the difference of the densities and applying Cauchy--Schwarz gives
\begin{align}
\left|\int U\,d(H-H')\right|
&\le\sqrt2\,h(p_H,p_{H'})
 \left\{2\int g^2p_H+2\int g^2p_{H'}\right\}^{1/2}\nonumber\\
&\le 2\sqrt{2B}\,\|U\|_{\infty,\mathcal K}h(p_H,p_{H'}).
\label{eqn:dung_hellinger_bound_for_u_difference}
\end{align}
This proves equation~\eqref{proof:A1} for arbitrary mixing measures on $\mathcal K$; their support cardinalities have not entered the argument.

\paragraph{Reduction of an HMM block to the center span.}
For the stated HMM application, choose $S$ to contain every fitted and true center. After an orthogonal change of coordinates, the $S^\perp$ coordinates of the observations have the same standard Gaussian law under both models. They are independent of the $S$ coordinates and of the latent path. Thus, for the selected time set $I$,
\[
p_\theta^I(x_S,x_\perp)=p_H(x_S)\varphi_{(d-q)r}(x_\perp),
\qquad
p_\star^I(x_S,x_\perp)=p_{H'}(x_S)\varphi_{(d-q)r}(x_\perp).
\]
The common factor integrates to one in the Hellinger affinity, so
\[
h(p_H,p_{H'})=h(p_\theta^I,p_\star^I)
\le h(p_{\theta,T},p_{\theta^\star,T}).
\]
If $q=d$, the orthogonal factor is absent; if $q=0$, the desired moment difference is zero. Finally, $q\le K+k$, so maximizing the constants over the finitely many admissible dimensions gives the asserted dependence on $K,k,r,R,D$, without dependence on $d$ or on the number of latent paths.
\end{proof}
\subsection{Transport through a common finite support}

To compare the fitted and true latent location distributions, we first move each fitted center to its assigned true center. After this, both distributions are supported on the same finite set, and the remaining transport concerns differences in their probabilities. Under the cluster assumptions, some of this probability mismatch can be transported within clusters at a cost determined by their diameter. We therefore distinguish this part from the mass that must move between clusters.

Let $I=\{t_1,\ldots,t_r\}$ with $t_1<\cdots<t_r$, and define the cluster map $\kappa:[K]\to[M]$ by $\kappa(j)=m$ for $j\in C_m$. For $\boldsymbol m\in[M]^r$, set
\[
\begin{aligned}
\bar P_\theta^{I,\mathrm{cl}}(\boldsymbol m)
&=\mathbb P_\theta\{\kappa(c(Z_{t_\ell}))=m_\ell\text{ for every }\ell\},\\
P_\star^{I,\mathrm{cl}}(\boldsymbol m)
&=\mathbb P_{\theta^\star}\{\kappa(Z_{t_\ell})=m_\ell\text{ for every }\ell\},\\
J_I&=\sum_{\boldsymbol m\in[M]^r}
\left|\bar P_\theta^{I,\mathrm{cl}}(\boldsymbol m)
             -P_\star^{I,\mathrm{cl}}(\boldsymbol m)\right|.
\end{aligned}
\]
Each cluster-sequence difference is a sum of label-sequence differences. The triangle inequality therefore gives $0\le J_I\le D_I$. For a single cluster both cluster-sequence laws are the same point mass, so $J_I=0$.

\begin{lemma}[Hierarchical transport]\label{lem:hierarchical-transport}
For $p\in\{1,2\}$, let
\[
V_{p,I}=\frac1r\sum_{t\in I}\sum_{i=1}^k q_t(i)r_i^p,
\]
where $q_t(i)=\mathbb P_\theta(Z_t=i)$ and $r_i=\|\mu_i-a_{c(i)}\|$. In the structured regimes,
\begin{equation}
W_{p,I}^p(G_\theta^I,G_\star^I)
\le 2^{p-1}\left\{
V_{p,I}
+\frac{(C_0\Delta_{\mathrm{sep}})^p}{2}D_I
+\frac{(2R)^p}{2}J_I
\right\}.
\label{proof:A5}
\end{equation}
Otherwise, for unstructured regime, we have 
\begin{equation}
W_{p,I}^p(G_\theta^I,G_\star^I)
\le 2^{p-1}\left\{
V_{p,I}+\frac{(2R)^p}{2}D_I
\right\}.
\label{proof:A6}
\end{equation}
\end{lemma}

\begin{proof}
\emph{Projection to the true support.}
Draw the fitted location tuple $X=(\mu_{Z_{t_1}},\ldots,\mu_{Z_{t_r}})$ and put
$Y=(a_{c(Z_{t_1})},\ldots,a_{c(Z_{t_r})})$. Its law on the true grid
$\mathcal J=[K]^r$ is $u(j)=\bar P_\theta^I(j)$; write
$v(j)=P_\star^I(j)$ for the true grid law. This coupling has exact cost
\[
\mathbb E d_{p,r}(X,Y)^p
=\frac1r\sum_{t\in I}\sum_iq_t(i)r_i^p=V_{p,I}.
\]
No independence among the selected hidden states is used.

\emph{Redistribution on the true grid.}
Write $w(j)=\min\{u(j),v(j)\}$ and let
$u'(j)=u(j)-w(j)$, $v'(j)=v(j)-w(j)$.
Their common total is $t=D_I/2$. Couple $w(j)$ diagonally, at zero cost.
Under cluster structure, partition $\mathcal J$ into the cells
$\mathcal J_m$ of equal cluster sequences $m\in[M]^r$. Define
\[
 a_m=\sum_{j\in\mathcal J_m}u'(j),\quad
 b_m=\sum_{j\in\mathcal J_m}v'(j),\quad
 c_m=\min\{a_m,b_m\}.
\]
For $a_mb_m>0$, put
\[
 \Gamma_m(j,l)=c_m\frac{u'(j)}{a_m}\frac{v'(l)}{b_m}
 \quad(j,l\in\mathcal J_m),
\]
and put $\Gamma_m=0$ if $a_mb_m=0$. Its row sums are
$c_mu'(j)/a_m\le u'(j)$ and column sums are
$c_mv'(l)/b_m\le v'(l)$. Thus it couples exactly $c_m$ units
within the cell, at cost per unit at most $(C_0\Delta_{\mathrm{sep}})^p$.
Let $U_m=\sum_{j\in\mathcal J_m}u(j)$ and
$V_m=\sum_{j\in\mathcal J_m}v(j)$. Because diagonal matching removes
the same mass from both cell totals,
\[
 a_m-c_m=(U_m-V_m)_+,
 \qquad b_m-c_m=(V_m-U_m)_+.
\]
The remaining supply and demand therefore each have total
$t_0=\tfrac12\sum_m|U_m-V_m|=J_I/2$. Denote their pointwise
masses by $u''(j)$ and $v''(l)$. Couple them by
$\Gamma_0(j,l)=u''(j)v''(l)/t_0$ when $t_0>0$, and set
$\Gamma_0=0$ otherwise. Its marginals are the two residual arrays.
Consequently,
\[
 \Gamma(j,l)=w(j)\mathbf 1\{j=l\}
  +\sum_m\Gamma_m(j,l)+\Gamma_0(j,l)
\]
has marginals $u,v$, with total cost at most
\[
 (C_0\Delta_{\mathrm{sep}})^p(t-t_0)+(2R)^pt_0
 \le \frac{(C_0\Delta_{\mathrm{sep}})^p}{2}D_I
     +\frac{(2R)^p}{2}J_I.
\]
For unstructured centers, omit the cell couplings and couple the
entire residual mass $t=D_I/2$ by $u'(j)v'(l)/t$ if $t>0$
(and by zero otherwise). Its cost is at most $(2R)^pD_I/2$.

\emph{Gluing the couplings.}
Conditional on $Y=\mathbf a_j$ with $u(j)>0$, draw
$Y'=\mathbf a_l$ with probability $\Gamma(j,l)/u(j)$.
For $u(j)=0$ choose an arbitrary conditional probability: this event
is null and the corresponding row of $\Gamma$ vanishes. Then
$X\sim G_\theta^I$, $Y'\sim G_\star^I$ and
$(Y,Y')\sim\Gamma$. Applying
$d_{p,r}(X,Y')^p\le2^{p-1}\{d_{p,r}(X,Y)^p+d_{p,r}(Y,Y')^p\}$
and taking expectations gives exactly \eqref{proof:A5} and
\eqref{proof:A6}, respectively.
\end{proof}

\paragraph{Lower transport bounds and matching.}
The preceding coupling gives upper transport bounds. We also record the reverse comparisons used to extract probability errors and matching from transport. They require no Markov property of the grouped fitted process.

\begin{lemma}[Transport, assigned labels, and positive-mass matching]
\label{lem:transport-assignment-matching}
Let $0<\delta\le\delta_{\mathrm{sep}}$ be any lower bound on the true-center separation. For $p\in\{1,2\}$ and every nonempty $I\subseteq[T]$ with $|I|=r$,
\begin{equation}
 V_{p,I}\le W_{p,I}^p(G_\theta^I,G_\star^I),
 \qquad
 D_I\le \frac{2^{p+1}r}{\delta^p}
              W_{p,I}^p(G_\theta^I,G_\star^I).
 \label{eq:assignment-transport-lower}
\end{equation}
The same statements hold for arbitrary fitted and true location-block laws with the indicated supports. Zero atom masses and repeated fitted centers are allowed.

For one-time measures $H=\sum_{i=1}^K v_i\delta_{b_i}$ and $H^\star=\sum_{j=1}^K w_j\delta_{a_j}$, assume additionally that $w_j\ge p_0>0$. If
\[
 W_1(H,H^\star)<p_0\delta/4,
\]
then each open ball $B(a_j,\delta/4)$ contains exactly one fitted label and that label has positive mass. After this unique relabeling,
\begin{equation}
 \|v-w\|_1\le\frac{8}{3\delta}W_1(H,H^\star).
 \label{eq:matching-transport-masses}
\end{equation}
In particular, $W_1(H,H^\star)\le p_0\delta/32$ implies $v_j\ge w_j/2$ for every $j$. Every assertion is valid both for $\delta=\delta_{\mathrm{sep}}$ and for $\delta=\Delta_{\mathrm{sep}}$.
\end{lemma}

\begin{proof}
We first compare the cost of an arbitrary coupling, then specialize to one-time exact specification.

\paragraph{Nearest-support cost.}
Let $(X,Y)$ have any coupling of the fitted and true location-block laws. Let $\chi(x)$ denote the smallest-index nearest true-center label, so that $\chi(\mu_i)=c(i)$, and put $C(X)=(\chi(X_1),\ldots,\chi(X_r))$. For every coordinate, $Y_\ell$ is a true center, so
\[
 \|X_\ell-Y_\ell\|_2^p
 \ge \min_j\|X_\ell-a_j\|_2^p.
\]
Averaging over coordinates and taking expectations gives a lower bound of $V_{p,I}$ for the cost of every coupling. This proves the first inequality in equation~\eqref{eq:assignment-transport-lower}.

\paragraph{Mass assigned to different labels.}
Because the true centers are distinct, the labels of $Y$ are uniquely determined by its locations. The two label marginals are $\bar P_\theta^I$ and $P_\star^I$. The coupling can put at most
$\sum_{\boldsymbol j}\min\{\bar P_\theta^I(\boldsymbol j),P_\star^I(\boldsymbol j)\}=1-D_I/2$
on identical label sequences. On the complementary event, at least one coordinate has $\chi(X_\ell)=i$ and $Y_\ell=a_j$ with $i\ne j$. The nearest-center property and the triangle inequality give
\[
 \delta\le\|a_i-a_j\|_2
 \le\|X_\ell-a_i\|_2+\|X_\ell-a_j\|_2
 \le2\|X_\ell-a_j\|_2.
\]
The normalized cost on this event is therefore at least $\delta^p/(2^p r)$. Every coupling has expected cost at least $\delta^pD_I/(2^{p+1}r)$, proving the second inequality. This argument acts on locations, so repeated fitted labels at the same location have the same assigned label; no division by an atom mass is involved.

\paragraph{Matching and fitted masses.}
If an open ball $B(a_j,\delta/4)$ contained no positive-mass fitted atom, all mass $w_j$ arriving at $a_j$ would travel at least $\delta/4$, contradicting $W_1(H,H^\star)<p_0\delta/4$. The $K$ balls are disjoint and there are only $K$ fitted labels. Thus they contain exactly one label each, all with positive mass; no zero-weight or repeated extra label can remain. Within these balls, each matched true center is the unique nearest center, so the relabeling agrees with $c$. The fitted centers are now distinct, so an optimal location coupling is a matrix $\Gamma$ on the aligned labels. Put $O=\sum_{i\ne j}\Gamma_{ij}$. For $i\ne j$,
$\|b_i-a_j\|_2\ge\delta-\|b_i-a_i\|_2>3\delta/4$.
Also,
\[
 |v_i-w_i|
 =\left|\sum_{j\ne i}\Gamma_{ij}-\sum_{j\ne i}\Gamma_{ji}\right|
 \le\sum_{j\ne i}(\Gamma_{ij}+\Gamma_{ji}).
\]
Summation yields $\|v-w\|_1\le2O\le8W_1(H,H^\star)/(3\delta)$, which is equation~\eqref{eq:matching-transport-masses}. At $W_1\le p_0\delta/32$, this norm is at most $p_0/12$, and hence $v_j\ge w_j-p_0/12\ge w_j/2$. All these comparisons are deterministic.
\end{proof}

\subsection{Stability of aggregation with estimated reference centers}
\label{app:proof-reference-grouping-stability}
\begin{proof}[Proof of Proposition~\ref{prop:reference-grouping-stability}]
Let $\mathcal E=\{i\in[k]:\widehat c(i)\ne c(i)\}$ be the set of labels whose assignment changes. Put $\delta=\delta_{\mathrm{sep}}$ and $\eta=\max_j\|\xi_j-a_j\|_2\le\tau$. For any $i\in\mathcal E$, write $j=c(i)$ and $\ell=\widehat c(i)\ne j$. Since $\ell$ is closest to $\mu_i$ among the reference centers, including at ties,
\[
\begin{aligned}
\delta
&\le\|a_j-a_\ell\|_2
\le\|a_j-\mu_i\|_2+\|\mu_i-a_\ell\|_2\\
&\le r_i+\|\mu_i-\xi_\ell\|_2+\eta
\le r_i+\|\mu_i-\xi_j\|_2+\eta
\le 2r_i+2\eta.
\end{aligned}
\]
Thus $r_i\ge\delta/2-\eta\ge\delta/4$. To prove the margin claim, fix any label $i$ whose reference assignment is $\ell=\widehat c(i)$, and suppose $m_i(\boldsymbol\xi)>2\tau$. For every $j\ne\ell$,
\[
\|\mu_i-a_j\|_2-\|\mu_i-a_\ell\|_2
\ge\|\mu_i-\xi_j\|_2-\|\mu_i-\xi_\ell\|_2-2\eta
\ge m_i(\boldsymbol\xi)-2\eta>0.
\]
Hence $c(i)=\widehat c(i)=\ell$, and by contraposition $\mathcal E\subseteq\mathcal B_\tau(\boldsymbol\xi)$. For every $t\le T$ the probability of occupying a relabeled state satisfies both
\begin{equation}
\mathbb P_\theta\{Z_t\in\mathcal E\}
=\sum_{i\in\mathcal E}q_t(i)
\le\frac4\delta\sum_iq_t(i)r_i
\le\frac{4D_\mu}{\delta},
\label{eq:estimated-reference-bad-state-linear}
\end{equation}
\begin{equation}
\mathbb P_\theta\{Z_t\in\mathcal E\}
\le\frac{16}{\delta^2}\sum_iq_t(i)r_i^2
\le\frac{16V}{\delta^2}.
\label{eq:estimated-reference-bad-state-quadratic}
\end{equation}

On the same fitted hidden-state path, compare the two sequences of labels $(c(Z_{t_\ell}))_{\ell=1}^r$ and $(\widehat c(Z_{t_\ell}))_{\ell=1}^r$. They agree unless at least one of the selected states belongs to $\mathcal E$. The elementary coupling inequality for total variation and a union bound therefore give
\[
\begin{aligned}
\|\widehat P_\theta^I-\bar P_\theta^I\|_1
&\le2\,\mathbb P_\theta\bigl\{\exists\ell\le r:
Z_{t_\ell}\in\mathcal E\bigr\}\\
&\le2\sum_{\ell=1}^r\mathbb P_\theta\{Z_{t_\ell}\in\mathcal E\}
\le\min\!\left\{
2\sum_{t\in I}\sum_{i\in\mathcal B_\tau(\boldsymbol\xi)}q_t(i),
\frac{8rD_\mu}{\delta},\frac{32rV}{\delta^2}\right\}.
\end{aligned}
\]
The triangle inequality with the true label-sequence law $P_\star^I$ proves equation~\eqref{eq:reference-grouping-stability}. This coupling works even if fitted states have zero masses or coincident emission centers, and it never uses a Markov assumption for the grouped process.

For the asserted exponents, write $s=s_{\max}$ in a structured regime, and $s=K$ otherwise. Under exact specification, Theorem~\ref{thm:globalexact} supplies $D_I\le C_r(\pi_{\min}^\star)\Delta_{\mathrm{sep}}^{-(2s-1)}h$ and $D_\mu\le C(\pi_{\min}^\star)\Delta_{\mathrm{sep}}^{-(2s-2)}h$. Since $\delta\ge\Delta_{\mathrm{sep}}$, the linear alternative in equation~\eqref{eq:reference-grouping-stability} preserves exponent $2s-1$. If the paired-center hypothesis of Theorem~\ref{thm:paired} is available, the same argument uses its constants independent of $\pi_{\min}^\star$. Under over-specification, Theorem~\ref{thm:overglobal} gives $D_I\le C_r\Delta_{\mathrm{sep}}^{-b_r}h$, where $b_r=\beta_r$ or $\gamma_{K,r}$, and $V\le C\Delta_{\mathrm{sep}}^{-(2s-2)}h$. The quadratic alternative contributes at most $C_r\Delta_{\mathrm{sep}}^{-2s}h$, which is bounded by $C_r\Delta_{\mathrm{sep}}^{-b_r}h$ because $b_r\ge2s$ and $\Delta_{\mathrm{sep}}\le1$. For one time and consecutive pairs, the definitions of $\widehat D_\pi$ and $\widehat D_A$ give the same conclusions. All constants have the stated dependence, completing the proof.
\end{proof}

\begin{proof}[Proof of Corollary~\ref{cor:reference-common-permutation}]
The balls $B(a_j,\tau)$ are disjoint because
$2\tau\le\delta_{\mathrm{sep}}/2<\delta_{\mathrm{sep}}$.
Thus the $K$ reference centers determine a unique bijection
$\sigma$ with the true centers. Apply $\sigma$ to the labels
\emph{after} grouping with the originally indexed references.
This transports the original smallest-index tie rule; it does not
recompute a new tie order after permuting the centers. Every
assigned reference remains a nearest reference under the aligned
indices. Hence the coupling proof of
Proposition~\ref{prop:reference-grouping-stability}, which only
uses this minimizing property, applies to $\sigma_\#\widehat P_\theta^I$.
The changed-label event is the same for all $I$, and the same
permutation yields all the stated bounds. Since $T$ is fixed, all
$2^T-1$ inequalities hold on the same center-accuracy event
without any additional union bound.

For the statistical assertion, on the global-MLE matching event
Corollary~\ref{cor:localization} supplies
\eqref{eq:reference-permuted-centers} with
$\tau=\Delta_{\mathrm{sep}}/4$. Intersect it with the stated
density event for $\widetilde\theta_n$ and apply the deterministic
bounds. The union bound gives failure probability at most
$m^{-1}+\alpha_n$, regardless of dependence between the estimators.
For measurability, choose the first admissible permutation in a
fixed ordering on the matching event, and the identity otherwise.
The finitely many matching conditions are Borel in the reference
centers, so this convention defines a measurable relabeling.
\end{proof}

\section{Finite-sample likelihood analysis}
\label{app:likelihood}
We prove Proposition~\ref{prop:density}, which bounds the Hellinger error of the global maximum likelihood estimator. The likelihood maximization is over the entire parameter space, so the proof must control the likelihood comparison simultaneously for all fitted parameters. Lemma~\ref{lem:upper-cover} constructs finitely many upper functions, each of which dominates a group of observation densities and has integral close to one. Lemma~\ref{lem:likelihood-tail} uses the domination and integral bounds in~\eqref{eq:upper-cover-domination}, combined with the collection-size bound~\eqref{eq:upper-cover-cardinality}, to control the probability of a large Hellinger error. Appendix~\ref{app:density-proof} finally establishes the existence of a measurable global maximizer and chooses the radius in Proposition~\ref{prop:density}.

We construct the finite collection by grouping parameters around finitely many representatives. For each representative, we need a function that dominates every observation density in its group. Its integral must remain close to one, since an upper function with large excess mass would give a poor likelihood comparison. We must also bound the number of functions, since the probability bound will sum over the collection. The following lemma provides these three properties.

\begin{lemma}[Upper likelihood cover]
\label{lem:upper-cover}
For every $0<\eta\le1/16$, there is a finite collection of nonnegative integrable functions $\{u_b:b\in\mathcal B_\eta\}$ such that every $b\in\mathcal B_\eta$ satisfies the integral bounds below, and every $\theta\in\Theta_k(R)$ has at least one $b\in\mathcal B_\eta$ satisfying the pointwise domination:
\begin{equation}
\begin{aligned}
p_{\theta,T}(x_{1:T})&\le u_b(x_{1:T})
\qquad\text{for all }x_{1:T}\in(\mathbb R^d)^T,\\
1&\le\int_{(\mathbb R^d)^T}u_b(x_{1:T})\,dx_{1:T}\le1+\eta.
\end{aligned}
\label{eq:upper-cover-domination}
\end{equation}
Moreover, with $m_k=k+k^2+kd$,
\begin{equation}
\log|\mathcal B_\eta|
\le m_k\log\!\left\{\frac{CT(k+1)(R+1)(d+1)}{\eta}\right\}.
\label{eq:upper-cover-cardinality}
\end{equation}
\end{lemma}

\begin{proof}
Fix $0<\rho\le 1/2$. Partition each probability simplex into
coordinate cells of side at most $\rho$, taking a simplex-valued
representative from every intersecting cell, and take a $\rho$-net of
$B_d(0,R)$ for each emission center. This constructs a finite map
$\theta\mapsto\theta^0$ satisfying
\[
\|\pi-\pi^0\|_\infty\vee\|A-A^0\|_\infty
\vee\max_i\|\mu_i-\mu_i^0\|_2\le\rho,
\qquad
N_\rho\le\left\{\frac{2(R+1)}{\rho}\right\}^{m_k},
\]
where $m_k=k+k^2+kd$. The covering count follows from the
coordinate-simplex grids and the volume bound for Euclidean ball nets.

For a representative center, put
$f_{\mu^0,\rho}(x)=\varphi_d(x-\mu^0)e^{\rho\|x-\mu^0\|_2}$.
The Gaussian likelihood ratio gives
$\varphi_d(x-\mu)\le f_{\mu^0,\rho}(x)$ whenever
$\|\mu-\mu^0\|_2\le\rho$. Its integral does not depend on $\mu^0$:
\[
M_d(\rho)=\mathbb E e^{\rho\|Z\|_2}
\le e^{\rho\sqrt d/2}(1-\rho/\sqrt d)^{-d/2}
\le e^{\rho\sqrt d+\rho^2/2},\qquad Z\sim N_d(0,I_d),
\]
provided $\rho/\sqrt d\le1/2$. The first inequality follows from
$\|Z\|_2\le\sqrt d/2+\|Z\|_2^2/(2\sqrt d)$ and the
chi-square moment-generating function; the second uses
$-\log(1-t)\le t+t^2$ for $0\le t\le1/2$.

For each representative $\theta^0$, inflate the discrete path weights
coordinatewise, setting $\pi_i^+=\pi_i^0+\rho$ and
$A_{ij}^+=A_{ij}^0+\rho$, and define
\[
u_{\theta^0,\rho}(x_{1:T})=
\sum_{z\in[k]^T}\pi^+_{z_1}
\prod_{t=1}^{T-1}A^+_{z_tz_{t+1}}
\prod_{t=1}^T f_{\mu^0_{z_t},\rho}(x_t).
\]
It dominates $p_{\theta,T}$ for every parameter mapped to $\theta^0$.
Each inflated initial or transition row sums to $1+k\rho$, so
integration and successive summation over the hidden path give
\[
1\le\int u_{\theta^0,\rho}
=M_d(\rho)^T(1+k\rho)^T
\le\exp\{T\rho(k+\sqrt d+1)\}.
\]
Choose $\rho=\eta/[4T(k+\sqrt d+1)]$; it satisfies the restrictions
above, and $e^{\eta/4}\le1+\eta$ for $0<\eta\le1/16$.
Thus every upper function has integral in $[1,1+\eta]$ and the
family dominates every fitted density. Finally,
\[
\log|\mathcal B_\eta|\le\log N_\rho
\le m_k\log\left\{
\frac{8T(R+1)(k+\sqrt d+1)}{\eta}\right\}
\le m_k\log\left\{\frac{CT(k+1)(R+1)(d+1)}{\eta}\right\}.
\]
No separation or positive fitted-weight assumption is used.
\end{proof}

\begin{lemma}[Hellinger tail bound]
\label{lem:likelihood-tail}
Let $\widehat p=p_{\widehat\theta,T}$ be any measurable global maximum likelihood estimator over $\Theta_k(R)$, with $k\ge K$, and let $p^\star=p_{\theta^\star,T}$. For $n$ independent length-$T$ training trajectories and every $0<r\le1$,
\begin{equation}
\mathbb P_\star\{h(\widehat p,p^\star)>r\}
\le
\exp\!\left[
m_k\log\!\left\{
\frac{CT(k+1)(R+1)(d+1)}{r^4}
\right\}
-\frac34nr^2
\right],
\label{eq:likelihood-tail}
\end{equation}
where $m_k=k+k^2+kd$ and $C>0$ is a numerical constant.
\end{lemma}

\begin{proof}
Fix $0<r\le1$ and apply Lemma~\ref{lem:upper-cover} with
$\eta=r^4/16$. Retain precisely those upper functions $u_b$ whose
domination cells contain some density $p$ with $h(p,p^\star)>r$.
For each retained $u_b$, since $u_b\ge p$ and
$\int(u_b-p)\le\eta$, Cauchy--Schwarz yields
\begin{equation}\label{eq:likelihood-upper-affinity}
\int\sqrt{u_b p^\star}
\le\int\sqrt{pp^\star}
+\left\{\int(\sqrt{u_b}-\sqrt p)^2\right\}^{1/2}
\le1-r^2+\sqrt\eta=1-3r^2/4.
\end{equation}
The retained collection is deterministic (depends on $p^\star$ and
$r$, not on the sample). By the unused-state embedding in the model
definition, $p^\star$ belongs to the fitted density class for
$k\ge K$. Hence a global MLE satisfies
$\prod_{\ell=1}^n\widehat p(X^{(\ell)})
\ge\prod_{\ell=1}^np^\star(X^{(\ell)})$.
If $h(\widehat p,p^\star)>r$, a retained $u_b\ge\widehat p$
therefore satisfies
$\prod_\ell u_b(X^{(\ell)})/p^\star(X^{(\ell)})\ge1$.
The denominators are strictly positive, because every nonempty
Gaussian mixture density is positive everywhere. For each fixed $u_b$,
independence, Markov's inequality applied to the square-root likelihood
ratio, and \eqref{eq:likelihood-upper-affinity} give
\[
\mathbb P_\star\left\{
\prod_{\ell=1}^n\frac{u_b(X^{(\ell)})}{p^\star(X^{(\ell)})}\ge1
\right\}
\le\left(\int\sqrt{u_b p^\star}\right)^n
\le e^{-3nr^2/4}.
\]
A union bound over at most $|\mathcal B_\eta|$ functions and
\eqref{eq:upper-cover-cardinality}, with $\eta=r^4/16$,
prove \eqref{eq:likelihood-tail}; absorb the factor $16$ into $C$.
If the retained family is empty, the bad event is empty.
\end{proof}

\subsection{Proof of Proposition~\ref{prop:density}}
\label{app:density-proof}
Lemma~\ref{lem:likelihood-tail} bounds the probability of a large Hellinger error for a measurable global maximum likelihood estimator. We complete the proof of Proposition~\ref{prop:density} by establishing the existence of such an estimator and choosing the radius that gives the stated probability bound.
\begin{proof}
We first prove that (i) the likelihood attains a global maximum for every sample and (ii) a global maximizer can be selected as a Borel measurable function of the sample. We then derive the Hellinger probability bound in Proposition~\ref{prop:density} from Lemma~\ref{lem:likelihood-tail}.

\medskip
\noindent\emph{(i) Existence of a global maximizer.} The parameter space $\Theta_k(R)$ is a finite product of closed probability simplexes and closed Euclidean balls, and is therefore a nonempty compact subset of $\mathbb R^{m_k}$. The finite hidden-path expansion shows that $(\theta,y_{1:T})\mapsto p_{\theta,T}(y_{1:T})$ is jointly continuous. As verified in the proof of Lemma~\ref{lem:likelihood-tail}, this density is strictly positive for every admissible parameter and every observation vector. Thus, for a sample $y=(y^{(1)},\ldots,y^{(n)})$, the likelihood $L(\theta;y)=\prod_{\ell=1}^n p_{\theta,T}(y^{(\ell)})$ is positive and jointly continuous in $(\theta,y)$. Compactness ensures that it attains a global maximum over $\Theta_k(R)$. Its maximizers also maximize the log likelihood, since the logarithm is strictly increasing.

\noindent\emph{(ii) A measurable global maximizer.} We now construct a global maximizer $\widehat\theta_n(y)$ that is a Borel measurable function of the sample $y$. Enclose $\Theta_k(R)$ in a fixed closed coordinate cube and repeatedly subdivide each cube into its smaller closed cubes obtained by bisecting each side. Fix an ordering of the children of every cube. For each cube $Q$ meeting $\Theta_k(R)$, define
\[
M_Q(y)=\max_{\theta\in Q\cap\Theta_k(R)}L(\theta;y),
\qquad
M(y)=\max_{\theta\in\Theta_k(R)}L(\theta;y).
\]
Both functions are continuous in $y$. Indeed, joint continuity makes $L$ uniformly continuous on $\Theta_k(R)$ times any compact neighborhood of a fixed sample. Hence, if $y_j\to y$, then
\[
|M_Q(y_j)-M_Q(y)|
\le \sup_{\theta\in\Theta_k(R)}|L(\theta;y_j)-L(\theta;y)|
\longrightarrow0,
\]
and the same bound holds for $M$.

Starting from the outer cube, select the first child $Q$ meeting $\Theta_k(R)$ for which $M_Q(y)=M(y)$, and repeat this selection inside the chosen cube. Such a child exists at every depth because the parent contains a global maximizer and its children cover it. Since $M_Q$ and $M$ are continuous, the set $\{y:M_Q(y)=M(y)\}$ is Borel. The fixed ordering therefore makes the selected cube at every finite depth a Borel measurable function of $y$.

For every cube meeting $\Theta_k(R)$, fix a representative in its intersection with $\Theta_k(R)$, independently of the sample. The representative of the selected cube at depth $j$ is a finite-valued measurable function of $y$. The selected cubes are nested and their diameters tend to zero, so these representatives converge to a measurable limit $\widehat\theta_n(y)\in\Theta_k(R)$. For each fixed $y$, a global maximizer in the selected cube lies within its diameter of the representative. These maximizing points therefore converge to $\widehat\theta_n(y)$ as well, and continuity gives $L(\widehat\theta_n(y);y)=M(y)$. This constructs a measurable global maximizer even when the maximum is attained at more than one parameter.

\noindent\emph{(iii) Hellinger probability bound.} We consider \(n\geq2\) and set \(L_n=\log(nTd)\) and \(r=C\sqrt{dL_n/n}\), where \(C\geq(\log2)^{-1/2}\) will be chosen below depending only on \(k\) and \(R\). If \(r\geq1\), the claimed bound follows immediately from \(h\leq1\). We therefore consider \(r<1\).

Let \(C_{\mathrm{tail}}\) denote the numerical constant in Lemma~\ref{lem:likelihood-tail}, enlarged if necessary so that \(C_{\mathrm{tail}}\geq1\), and write \(a_T=C_{\mathrm{tail}}T(k+1)(R+1)(d+1)\). Since \(d\geq1\), we have \(a_T\leq B_{k,R}Td\), where \(B_{k,R}=2C_{\mathrm{tail}}(k+1)(R+1)\). Moreover, \(L_n\geq\log2\) and \(C\geq(\log2)^{-1/2}\) imply \(r\geq n^{-1/2}\). Consequently,
\[
\log(a_T/r^4)\leq\log B_{k,R}+\log(Td)+2\log n\leq\left(2+\frac{\log B_{k,R}}{\log2}\right)L_n.
\]
Combining this inequality with \(m_k=kd+k^2+k\leq(k^2+2k)d\) gives \(m_k\log(a_T/r^4)\leq A_{k,R}dL_n\) for a constant \(A_{k,R}>0\) depending only on \(k\) and \(R\). Lemma~\ref{lem:likelihood-tail} now yields
\[
\mathbb P_\star\left\{h(p_{\widehat\theta_n,T},p_{\theta^\star,T})>r\right\}\leq\exp\left\{\left(A_{k,R}-\frac34C^2\right)dL_n\right\}.
\]
Choose \(C\) sufficiently large that \(\frac34C^2\geq A_{k,R}+1\). The right-hand side is then at most \(e^{-dL_n}\leq1/n\). Since \(r=\varepsilon_{n,T}\), this proves Proposition~\ref{prop:density}. The argument applies to every measurable global maximizer, including when the maximizing parameter is not unique.

\end{proof}

\paragraph{Measurability of the statistical losses.}
\label{app:loss-measurability}
The measurable selection just constructed must also be compatible with the losses and the matching rule used later. We verify those properties for each fixed true parameter.

First, the map $\theta\mapsto p_{\theta,T}(x)$ is continuous for every $x$, and its values are probability densities. If $\theta_m\to\theta$, dominated convergence applied to $\min\{p_{\theta_m,T},p_{\theta,T}\}\le p_{\theta,T}$ gives
\[
\int|p_{\theta_m,T}-p_{\theta,T}|
=2-2\int\min\{p_{\theta_m,T},p_{\theta,T}\}\longrightarrow0.
\]
Because $h^2(p,q)\le\tfrac12\|p-q\|_1$, the observation density is continuous in Hellinger distance. The triangle inequality then makes
$\theta\mapsto h(p_{\theta,T},p_{\theta^\star,T})$ continuous. In particular, the density events used in the probability bounds are measurable.

For each fitted label $i$, the smallest-index nearest-center rule has the Borel cells
\[
\{c(i)=j\}
=\bigcap_{\ell<j}\{\|\mu_i-a_j\|_2<\|\mu_i-a_\ell\|_2\}
 \cap\bigcap_{\ell>j}\{\|\mu_i-a_j\|_2\le\|\mu_i-a_\ell\|_2\}.
\]
State and path probabilities are polynomial functions of $\pi,A$. Hence each grouped probability is a finite sum of continuous path probabilities multiplied by these Borel indicators. This proves that $D_I,D_\pi,D_A,D_\mu,V,B$, and $J_I$ are Borel, including at nearest-center ties.

The Wasserstein losses are in fact continuous despite repeated centers and zero path probabilities. To see this, put
$\lambda_z(\theta)=\pi_{z_1}\prod_{t<T}A_{z_tz_{t+1}}$ for $z\in[k]^T$.
Couple common mass $\min\{\lambda_z(\theta_m),\lambda_z(\theta)\}$ at identical hidden paths. The remaining mass is
$\frac12\sum_z|\lambda_z(\theta_m)-\lambda_z(\theta)|$ and can be coupled arbitrarily; its normalized cost is at most $(2R)^p$. Thus, for $p\in\{1,2\}$ and any nonempty $I$,
\[
\begin{aligned}
W_{p,I}^p(G_{\theta_m}^I,G_\theta^I)
&\le \max_{i\in[k]}\|\mu_i(\theta_m)-\mu_i(\theta)\|_2^p\\
&\quad+\frac{(2R)^p}{2}
       \sum_{z\in[k]^T}|\lambda_z(\theta_m)-\lambda_z(\theta)|
\longrightarrow0.
\end{aligned}
\]
The triangle inequality proves continuity of the distance to the fixed true location law. The construction does not require distinct support points or strictly positive path weights.

Finally, when $k=K$, enumerate the finitely many permutations in a fixed order and select the first satisfying all matching inequalities $\|b_j-a_j\|_2\le\Delta_{\mathrm{sep}}/4$, if one exists; select the identity permutation otherwise. Each feasibility set is closed, so the selected permutation is Borel. On the matching event it is the unique admissible permutation, because the true-center neighborhoods are disjoint. Off this event the convention still defines all relabeled center and transition-row losses. Consequently, the probability and $\mathcal O_{\mathbb P}$ statements in the paper concern well-defined random variables on the entire sample space, rather than quantities defined only on a high-probability event.

\subsection{Likelihood errors and approximate maximizers}
\label{app:likelihood-tolerance}

The likelihood comparison in the preceding proof is the only step that uses global maximization. The following lemma allows a quantified likelihood shortfall and states the weaker comparison with the true likelihood that is actually sufficient. Write
\[
\ell_n(\theta)=\sum_{i=1}^n\log p_{\theta,T}(X_{1:T}^{(i)}).
\]

\begin{lemma}[Density control with a likelihood tolerance]
\label{lem:likelihood-tolerance}
Suppose that the assumptions of Proposition~\ref{prop:density} hold, and let $\widetilde\theta_n$ be a measurable estimator taking values in $\Theta_k(R)$. Let $\tau_n\ge0$ and $\alpha_n\in[0,1]$ be deterministic. If
\begin{equation}
\mathbb P_\star(\mathcal G_n)\ge1-\alpha_n,
\qquad
\mathcal G_n:=\bigl\{\ell_n(\widetilde\theta_n)
                 \ge\ell_n(\theta^\star)-n\tau_n\bigr\},
\label{eq:likelihood-tolerance-event}
\end{equation}
then a constant $C=C(k,R)$ satisfies, for every $n\ge2$,
\begin{equation}
\mathbb P_\star\!\left\{
 h(p_{\widetilde\theta_n,T},p_{\theta^\star,T})
 >C\sqrt{\frac{d\log(nTd)}{n}+\tau_n}
 \right\}
\le \frac1n+\alpha_n.
\label{eq:likelihood-tolerance-radius}
\end{equation}
A sufficient condition for equation~\eqref{eq:likelihood-tolerance-event} is that
$\ell_n(\widetilde\theta_n)\ge\sup_{\theta\in\Theta_k(R)}\ell_n(\theta)-n\tau_n$
with probability at least $1-\alpha_n$. The probability statements may also include independent auxiliary randomness used to construct the estimator.
\end{lemma}

\begin{proof}
Put $p^\star=p_{\theta^\star,T}$ and $m_k=kd+k^2+k$. We first establish the explicit tail comparison
\begin{equation}
\begin{split}
&\mathbb P_\star\bigl\{
 h(p_{\widetilde\theta_n,T},p^\star)>r,\ \mathcal G_n\bigr\}\\
&\qquad\le
 \exp\!\left\{
 m_k\log\!\left(\frac{C_{\mathrm{cov}}T(k+1)(R+1)(d+1)}{r^4}\right)
 -\frac34nr^2+\frac12n\tau_n
 \right\},\qquad 0<r\le1,
\end{split}
\label{eq:likelihood-tolerance-tail}
\end{equation}
where $C_{\mathrm{cov}}\ge1$ is numerical. To see this, use the upper cover in Lemma~\ref{lem:upper-cover} with $\eta=r^4/16$, and retain exactly the cells used in Lemma~\ref{lem:likelihood-tail}. These cells are fixed before the sample is observed. For every retained upper function $u_b$, equation~\eqref{eq:likelihood-upper-affinity} gives
\[
\int\sqrt{u_bp^\star}\le1-\frac34r^2.
\]
Define the nonnegative random variable
\[
Z_b:=\prod_{i=1}^n
\sqrt{\frac{u_b(X_{1:T}^{(i)})}{p^\star(X_{1:T}^{(i)})}}.
\]
Its expectation is at most $(1-3r^2/4)^n\le e^{-3nr^2/4}$ by independence of the trajectories. On the event on the left of equation~\eqref{eq:likelihood-tolerance-tail}, some retained $u_b$ dominates $p_{\widetilde\theta_n,T}$, and the likelihood comparison gives $Z_b\ge e^{-n\tau_n/2}$. Markov's inequality and a union bound therefore yield the right-hand side of equation~\eqref{eq:likelihood-tolerance-tail}. If the estimator uses auxiliary randomness, this event inclusion still holds pointwise, while $Z_b$ depends only on the observations; the same bound follows without any conditioning on a selected fit.

Set $L_n=\log(nTd)$ and $r=C\sqrt{dL_n/n+\tau_n}$, with $C\ge(\log2)^{-1/2}$. When $r\ge1$, equation~\eqref{eq:likelihood-tolerance-radius} follows from $h\le1$. Otherwise $r\ge n^{-1/2}$, and the covering logarithm is bounded exactly as in the proof of Proposition~\ref{prop:density}:
\[
 m_k\log\!\left(\frac{C_{\mathrm{cov}}T(k+1)(R+1)(d+1)}{r^4}\right)
 \le A_{k,R}dL_n.
\]
Inserting the chosen $r$ into equation~\eqref{eq:likelihood-tolerance-tail} bounds its exponent by
\[
\left(A_{k,R}-\frac34C^2\right)dL_n
+\left(\frac12-\frac34C^2\right)n\tau_n.
\]
Choose $C$ so that $3C^2/4\ge A_{k,R}+1$ and $3C^2/4\ge1/2$. The exponent is then at most $-dL_n$, whose exponential is at most $1/n$. Adding $\mathbb P_\star(\mathcal G_n^c)\le\alpha_n$ proves the claim. Finally, the true observation density belongs to the fitted class by the unused-state embedding already established, so a global likelihood-gap bound implies equation~\eqref{eq:likelihood-tolerance-event}.
\end{proof}

A tolerance $\tau_n=O(d\log(nTd)/n)$ preserves the density order used in the recovery corollaries. More generally, outside the exceptional event bounded in equation~\eqref{eq:likelihood-tolerance-radius}, every deterministic inverse inequality applies with its density error replaced by the displayed radius. This observation requires a likelihood-gap or true-likelihood comparison with controlled probability. A small increase between successive EM iterates, or selection of the best among finitely many restarts, does not by itself supply that control.

\section{Multivariate witness polynomials}
\label{app:witness-polynomials}

This appendix section constructs polynomials that convert observation-density error into bounds on latent recovery. We begin with a nonnegative polynomial that vanishes at every true center and controls the squared distance to the true support. Throughout, $x\in B_d(0,R)$, and $c(x)$ denotes the index of a nearest true center, with ties resolved by choosing the smallest index. When applying Lemma~\ref{lem:moment-extraction}, we restrict the polynomials to the span of the true and fitted centers and take their supremum norms over the intersection of this span with $B_d(0,R)$.

The first subsections provide positive witnesses, flat selectors, and interpolation tools. Joint selectors then control multitime probabilities with explicit separation and diameter factors. The final subsection constructs polynomials that interpolate both true and fitted supports under an assumed pairing, for use in the matched bounds under exact specification.

\subsection{A positive witness for all fitted centers}

We first construct a polynomial that detects displacement from the true support without cancellation between fitted states. Define
\begin{equation}
P(x)=\prod_{j=1}^K\|x-a_j\|^2.
\label{proof:B1}
\end{equation}
Since $P$ is nonnegative and vanishes at every true center, its expectation under the true latent distribution is zero, whereas its expectation under the fitted distribution is a sum of nonnegative terms. The following lemma shows that this sum controls the probability-weighted squared distances of fitted centers from their nearest true centers.

\begin{lemma}[Global positive coercivity]
\label{lem:positive-coercivity}
In a structured configuration, there is a constant $c_{\mathrm{pos}}>0$, independent of $\Delta_{\mathrm{sep}}$, such that
\begin{equation}
P(x)\ge c_{\mathrm{pos}}\Delta_{\mathrm{sep}}^{2s_{\max}-2}\|x-a_{c(x)}\|^2,
\qquad
\|P\|_\infty\le(2R)^{2K}.
\label{proof:B2}
\end{equation}
Without a cluster assumption, the same statement holds with $s_{\max}=K$. Consequently, for every $t\le T$,
\begin{equation}
\begin{aligned}
0\le \sum_i q_t(i)P(\mu_i)
&\le C h(p_{\theta,T},p_{\theta^\star,T}),\\
\sum_i q_t(i)r_i^2
&\le C\Delta_{\mathrm{sep}}^{-(2s_{\max}-2)}h(p_{\theta,T},p_{\theta^\star,T}).
\end{aligned}
\label{proof:B3}
\end{equation}
\end{lemma}

\begin{proof}
We first establish the pointwise lower bound for $P$ and its supremum norm. We then apply Lemma~\ref{lem:moment-extraction} to bound the fitted expectation of $P$ and obtain the quadratic displacement bound in Equation~\eqref{proof:B3}.

\paragraph{Separating the nearest-center distance.}
Fix $x\in B_d(0,R)$, and write $j_0=c(x)$ and $\rho=\|x-a_{j_0}\|$. We isolate the factor $\rho^2$ in $P(x)$ and bound the remaining factors using distances between true centers. For each $\ell\ne j_0$, the nearest-center assignment gives $\rho\le\|x-a_\ell\|$, so
\[
\begin{aligned}
d_{j_0\ell}:=\|a_{j_0}-a_\ell\|
&\le \|a_{j_0}-x\|+\|x-a_\ell\|
&&\text{(by triangle inequality)}\\
&\le 2\|x-a_\ell\|
&&\text{(nearest-center assignment)}.
\end{aligned}
\]
It follows that
\[
\begin{aligned}
P(x)
&=\rho^2\prod_{\ell\ne j_0}\|x-a_\ell\|^2 \ge 4^{-(K-1)}\rho^2
\underbrace{\prod_{\ell\ne j_0}d_{j_0\ell}^2}_{=:D_{j_0}}.
\end{aligned}
\]
This inequality holds throughout the ball, including points tied between nearest centers; it does not require $x$ to be close to the true support.

\paragraph{Accounting for the cluster sizes.}
Suppose $a_{j_0}$ belongs to a cluster $\mathcal C$ containing $s_{j_0}\le s_{\max}$ true centers. The product $D_{j_0}$ separates into distances to centers inside and outside this cluster:
\[
D_{j_0}
=
\underbrace{\prod_{\substack{\ell\in\mathcal C\\\ell\ne j_0}}d_{j_0\ell}^2}_{s_{j_0}-1\text{ within-cluster factors  }}
\underbrace{\prod_{\ell\notin\mathcal C}d_{j_0\ell}^2}_{K-s_{j_0}\text{   between-cluster factors}}
\ge
\Delta_{\mathrm{sep}}^{2(s_{j_0}-1)}D_0^{2(K-s_{j_0})}.
\]
Only the within-cluster factors require the lower bound $\Delta_{\mathrm{sep}}$; the remaining distances are bounded below by the fixed constant $D_0$. Let $\lambda=\min\{1,D_0\}$. Since $0<\Delta_{\mathrm{sep}}\le1$ and $s_{j_0}\le s_{\max}$, we have
\[
\Delta_{\mathrm{sep}}^{2(s_{j_0}-1)}
\ge \Delta_{\mathrm{sep}}^{2(s_{\max}-1)},
\qquad
D_0^{2(K-s_{j_0})}
\ge \lambda^{2(K-s_{j_0})}
\ge \lambda^{2K}.
\]
Consequently,
\[
P(x)\ge
\underbrace{4^{-(K-1)}\lambda^{2K}}_{=:c_{\mathrm{pos}}}
\Delta_{\mathrm{sep}}^{2s_{\max}-2}\|x-a_{c(x)}\|^2.
\]
Without a cluster assumption, every factor in $D_{j_0}$ is at least $\Delta_{\mathrm{sep}}^2$, which gives
\[
D_{j_0}\ge\Delta_{\mathrm{sep}}^{2(K-1)}.
\]
Thus the same argument applies with $s_{\max}=K$ and $c_{\mathrm{pos}}=4^{-(K-1)}$.

\paragraph{Bounding the polynomial norm.}
The polynomial $P$ has degree $2K$. Since $x,a_j\in B_d(0,R)$, the triangle inequality gives $\|x-a_j\|\le2R$ for every $j$, and hence
\[
0\le P(x)=\prod_{j=1}^K\|x-a_j\|^2\le(2R)^{2K}.
\]
This proves the supremum bound in~\eqref{proof:B2}, which also holds after restricting $P$ to the span of the true and fitted centers.

\paragraph{Controlling the fitted polynomial expectation.}
Fix $t\le T$. The true one-time latent location distribution is supported on $\{a_1,\ldots,a_K\}$, where $P$ vanishes. Therefore,
\[
\begin{aligned}
\int P\,d\bigl(G_\theta^{\{t\}}-G_\star^{\{t\}}\bigr)
&=\sum_i q_t(i)P(\mu_i)
-\underbrace{\int P\,dG_\star^{\{t\}}}_{=\,0\;\text{ since }P(a_j)=0\text{ for all }j}\\
&=\sum_i q_t(i)P(\mu_i)\\
&\ge0 \qquad\text{since }q_t(i)\ge0\text{ and }P(\mu_i)\ge0\text{ for all }i.
\end{aligned}
\]
Let $p_\theta^{\{t\}}$ and $p_{\theta^\star}^{\{t\}}$ denote the corresponding one-time observation densities. Applying Lemma~\ref{lem:moment-extraction} on the span of the true and fitted centers gives
\[
\begin{aligned}
\sum_i q_t(i)P(\mu_i)
&=\left|\int P\,d\bigl(G_\theta^{\{t\}}-G_\star^{\{t\}}\bigr)\right|\\
&\le C\|P\|_\infty h\bigl(p_\theta^{\{t\}},p_{\theta^\star}^{\{t\}}\bigr)
&&\text{(Lemma~\ref{lem:moment-extraction})}\\
&\le C\|P\|_\infty h(p_{\theta,T},p_{\theta^\star,T})
&&\text{(Hellinger contraction under marginalization)}\\
&\le C h(p_{\theta,T},p_{\theta^\star,T})
&&\text{(the supremum bound in~\eqref{proof:B2})}.
\end{aligned}
\]
The span has dimension at most $K+k$, and the polynomial has degree $2K$. Thus the constant supplied by the moment-extraction lemma is independent of emission separation and the ambient dimension. The fixed factor $(2R)^{2K}$ is absorbed into $C$.

\paragraph{Recovering the quadratic displacement bound.}
For each fitted center, the pointwise bound in~\eqref{proof:B2} gives
\[
P(\mu_i)\ge c_{\mathrm{pos}}\Delta_{\mathrm{sep}}^{2s_{\max}-2}r_i^2,
\qquad
r_i=\|\mu_i-a_{c(\mu_i)}\|.
\]
Multiplying by $q_t(i)\ge0$ and summing over fitted states yields
\[
c_{\mathrm{pos}}\Delta_{\mathrm{sep}}^{2s_{\max}-2}\sum_i q_t(i)r_i^2
\le \sum_i q_t(i)P(\mu_i)
\le C h(p_{\theta,T},p_{\theta^\star,T}).
\]
Dividing by $c_{\mathrm{pos}}\Delta_{\mathrm{sep}}^{2s_{\max}-2}$ proves the second inequality in~\eqref{proof:B3}. No fitted probability is divided out, so the argument also applies to zero fitted probabilities and repeated fitted centers.
\end{proof}

\subsection{Flat mass selectors and directional mean selectors}

Lemma~\ref{lem:positive-coercivity} bounds the probability-weighted squared distances of fitted centers from the true support. We next construct the directional polynomial $H_{j,u}$ to approximate $\mathbf1\{c(x)=j\}\langle u,x-a_j\rangle$. Its expectation will recover the displacement sum within each assigned group. We also construct the flat polynomial $E_j$, whose values and local second-order behavior will be used in Lemma~\ref{lem:diameter-joint-selectors}.

The mass polynomial should equal one at $a_j$ and zero at the other true centers. We additionally require its gradient to vanish at every true center, so that small location displacements contribute only a second-order error. The directional polynomial should vanish at every true center, with gradient $u$ at $a_j$ and zero at the others. Its first-order term then agrees with the directional displacement we want to recover.

For $j\in[K]$, consider
\begin{equation}
Q_j(x)=\prod_{l\ne j}\frac{\|x-a_l\|^2}{\|a_j-a_l\|^2}.
\label{B.4}
\end{equation}
Each factor equals one at $a_j$, so $Q_j(a_j)=1$. At any other true center $a_l$, the squared factor $\|x-a_l\|^2$ makes both the value and gradient vanish. The remaining gradient at $a_j$ is
\begin{equation}
v_j:=\nabla Q_j(a_j)
=2\sum_{l\ne j}\frac{a_j-a_l}{\|a_j-a_l\|^2}.
\label{B.5}
\end{equation}
The polynomial $Q_j$ already has the required values at the true centers, but its gradient at $a_j$ is $v_j$, which need not be zero. Multiplication by $1-\langle v_j,x-a_j\rangle$ cancels this gradient without changing those values. Multiplication by $\langle u,x-a_j\rangle$ instead gives a polynomial with value zero and gradient $u$ at $a_j$. We therefore define
\[
\left\{
\begin{aligned}
E_j(x)&=Q_j(x)\{1-\langle v_j,x-a_j\rangle\},\\
H_{j,u}(x)&=Q_j(x)\langle u,x-a_j\rangle,
\qquad \|u\|\le1.
\end{aligned}
\right.
\]
The following lemma verifies the interpolation values and gradients, bounds the supremum norms of $Q_j$ and $H_{j,u}$, and controls the directional approximation error by $P(x)$. Its proof also establishes the local Hessian and Taylor bounds for $E_j$ needed later.

\begin{lemma}[Values, norms, and pointwise leakage]
\label{lem:flat-selectors}
We have
\begin{equation}
E_j(a_l)=\mathbf1\{j=l\},\qquad \nabla E_j(a_l)=0,
\label{B.7}
\end{equation}
\begin{equation}
H_{j,u}(a_l)=0,\qquad \nabla H_{j,u}(a_l)=\mathbf1\{j=l\}u.
\label{B.8}
\end{equation}

Their degrees are at most $2K-1$, and
\begin{equation}
\|Q_j\|_\infty+\|H_{j,u}\|_\infty
\le C\Delta_{\mathrm{sep}}^{-(2s_{\max}-2)}.
\label{B.9}
\end{equation}
For every $x\in B_d(0,R)$,
\begin{equation}
|H_{j,u}(x)-\mathbf1\{c(x)=j\}\langle u,x-a_j\rangle|
\le C\Delta_{\mathrm{sep}}^{-(2s_{\max}-1)}P(x).
\label{B.11}
\end{equation}
\end{lemma}

\begin{proof}
We first verify the interpolation values and gradients and bound the supremum norms of $Q_j$ and $H_{j,u}$. We then prove the directional approximation bound by treating points away from $a_j$ through direct product comparisons and points near $a_j$ through Taylor expansion. The local derivative calculation also provides the second-order bound for $E_j$ used in Lemma~\ref{lem:diameter-joint-selectors}.

\paragraph{Values and gradients at the true centers.}
The definition gives $Q_j(a_j)=1$. Differentiating the product at $a_j$ yields
\[
\begin{aligned}
\nabla Q_j(a_j)
&=\sum_{l\ne j}\frac{2(a_j-a_l)}{\|a_j-a_l\|^2}
\prod_{\substack{m\ne j\\m\ne l}}
\frac{\|a_j-a_m\|^2}{\|a_j-a_m\|^2}\\
&=2\sum_{l\ne j}\frac{a_j-a_l}{\|a_j-a_l\|^2}
=v_j.
\end{aligned}
\]
For $l\ne j$, we can instead isolate the factor that vanishes at $a_l$:
\[
Q_j(x)=\frac{\|x-a_l\|^2}{\|a_j-a_l\|^2}
\prod_{\substack{m\ne j\\m\ne l}}
\frac{\|x-a_m\|^2}{\|a_j-a_m\|^2}.
\]
Both $\|x-a_l\|^2$ and its gradient vanish at $a_l$. The product rule therefore gives $Q_j(a_l)=0$ and $\nabla Q_j(a_l)=0$.

Differentiating the selectors, we obtain
\[
\begin{aligned}
\nabla E_j(x)
&=\{1-\langle v_j,x-a_j\rangle\}\nabla Q_j(x)-Q_j(x)v_j,\\
\nabla H_{j,u}(x)
&=\langle u,x-a_j\rangle\nabla Q_j(x)+Q_j(x)u.
\end{aligned}
\]
At $a_j$, these identities give
\[
\begin{aligned}
E_j(a_j)&=1,
\qquad &\nabla E_j(a_j)&=v_j-v_j=0,\\
H_{j,u}(a_j)&=0,
\qquad &\nabla H_{j,u}(a_j)&=u.
\end{aligned}
\]
At every other true center, both values and gradients vanish because $Q_j(a_l)=0$ and $\nabla Q_j(a_l)=0$. This proves~\eqref{B.7} and~\eqref{B.8}. The degree of $Q_j$ is $2K-2$, so multiplication by either linear factor gives degree at most $2K-1$.

\paragraph{Supremum bounds for moment extraction.}
The product-distance estimate in the proof of Lemma~\ref{lem:positive-coercivity} gives
\[
\prod_{l\ne j}\|a_j-a_l\|^2
\ge c\Delta_{\mathrm{sep}}^{2s_{\max}-2}.
\]
Without cluster assumptions, the same estimate holds with $s_{\max}=K$.

For $x\in B_d(0,R)$, every numerator distance is at most $2R$. Therefore,
\[
\begin{aligned}
|Q_j(x)|
&=\frac{\prod_{l\ne j}\|x-a_l\|^2}
{\prod_{l\ne j}\|a_j-a_l\|^2}\\
&\le\frac{(2R)^{2K-2}}
{\prod_{l\ne j}\|a_j-a_l\|^2}\\
&\le C\Delta_{\mathrm{sep}}^{-(2s_{\max}-2)}.
\end{aligned}
\]
The directional selector is obtained by multiplying $Q_j$ by $\langle u,x-a_j\rangle$. Since $\|u\|\le1$ and both $x$ and $a_j$ lie in $B_d(0,R)$,
\[
|\langle u,x-a_j\rangle|
\le\|u\|\|x-a_j\|
\le2R.
\]
Consequently,
\[
\begin{aligned}
|H_{j,u}(x)|
&=|Q_j(x)|\,|\langle u,x-a_j\rangle|\\
&\le2R\,|Q_j(x)|\\
&\le C\Delta_{\mathrm{sep}}^{-(2s_{\max}-2)}.
\end{aligned}
\]
Taking suprema and absorbing the fixed factors into $C$ proves~\eqref{B.9}.

We also record a bound for $v_j$, which will enter the local Hessian calculation for $E_j$. Its definition and the separation of the true centers give
\[
\begin{aligned}
\|v_j\|
&=\left\|2\sum_{l\ne j}
\frac{a_j-a_l}{\|a_j-a_l\|^2}\right\|\\
&\le2\sum_{l\ne j}\frac{1}{\|a_j-a_l\|}\\
&\le2(K-1)\Delta_{\mathrm{sep}}^{-1}.
\end{aligned}
\]

\paragraph{Comparison with the positive witness away from the selected center.}
We next prove~\eqref{B.11} away from $a_j$. If $P(x)=0$, then $x$ is a true center, and the interpolation values make the directional approximation error zero. We may therefore assume that $P(x)>0$.

Write $r=\|x-a_j\|$. The product $Q_j$ contains every numerator factor of $P$ except $\|x-a_j\|^2$, so
\[
\begin{aligned}
\frac{Q_j(x)}{P(x)}
&=\frac{1}{r^2\prod_{l\ne j}\|a_j-a_l\|^2}\\
&\le C\Delta_{\mathrm{sep}}^{-(2s_{\max}-2)}r^{-2}.
\end{aligned}
\]
Whenever $r>\Delta_{\mathrm{sep}}/4$, this gives
\begin{equation}
\begin{aligned}
\frac{|H_{j,u}(x)|}{P(x)}
&=\frac{Q_j(x)|\langle u,x-a_j\rangle|}{P(x)}\\
&\le r\frac{Q_j(x)}{P(x)}
\qquad\text{(because $\|u\|\le1$)}\\
&\le C\Delta_{\mathrm{sep}}^{-(2s_{\max}-2)}r^{-1}\\
&\le C\Delta_{\mathrm{sep}}^{-(2s_{\max}-1)}
\qquad\text{(because $r>\Delta_{\mathrm{sep}}/4$)}.
\end{aligned}
\label{eq:direction_positive}
\end{equation}

Suppose first that $c(x)\ne j$. The assigned center is distinct from $a_j$, and nearest-center assignment gives
\[
\begin{aligned}
\Delta_{\mathrm{sep}}
&\le\|a_j-a_{c(x)}\|\\
&\le\|a_j-x\|+\|x-a_{c(x)}\|
\qquad\text{(by the triangle inequality)}\\
&\le2r
\qquad\text{(because $\|x-a_{c(x)}\|\le\|x-a_j\|=r$)}.
\end{aligned}
\]
Thus $r\ge\Delta_{\mathrm{sep}}/2$, so~\eqref{eq:direction_positive} applies. The target $\mathbf1\{c(x)=j\}\langle u,x-a_j\rangle$ is zero in this case. Hence
\[
|H_{j,u}(x)-\mathbf1\{c(x)=j\}\langle u,x-a_j\rangle|
=|H_{j,u}(x)|
\le C\Delta_{\mathrm{sep}}^{-(2s_{\max}-1)}P(x).
\]

Suppose instead that $c(x)=j$ and $r>\Delta_{\mathrm{sep}}/4$. We must now bound the target term as well as the polynomial. Since $a_j$ is the assigned nearest center, Lemma~\ref{lem:positive-coercivity} gives
\[
P(x)\ge c\Delta_{\mathrm{sep}}^{2s_{\max}-2}r^2.
\]
Therefore,
\[
\begin{aligned}
\frac{|\langle u,x-a_j\rangle|}{P(x)}
&\le\frac{r}{P(x)}
\qquad\text{(because $\|u\|\le1$)}\\
&\le C\Delta_{\mathrm{sep}}^{-(2s_{\max}-2)}r^{-1}\\
&\le C\Delta_{\mathrm{sep}}^{-(2s_{\max}-1)}
\qquad\text{(because $r>\Delta_{\mathrm{sep}}/4$)}.
\end{aligned}
\]
We can now bound the approximation error by adding the polynomial and target contributions:
\[
\begin{aligned}
\frac{|H_{j,u}(x)-\langle u,x-a_j\rangle|}{P(x)}
&\le\frac{|H_{j,u}(x)|}{P(x)}
+\frac{|\langle u,x-a_j\rangle|}{P(x)}\\
&\le C\Delta_{\mathrm{sep}}^{-(2s_{\max}-1)}
+C\Delta_{\mathrm{sep}}^{-(2s_{\max}-1)}\\
&\le C\Delta_{\mathrm{sep}}^{-(2s_{\max}-1)}.
\end{aligned}
\]
The first contribution is bounded by~\eqref{eq:direction_positive}, and the second by the preceding calculation. This proves~\eqref{B.11} in the second case. It remains to treat $c(x)=j$ and $r\le\Delta_{\mathrm{sep}}/4$.

\paragraph{Derivative bounds near the selected center.}
Suppose $c(x)=j$ and $r=\|x-a_j\|\le\Delta_{\mathrm{sep}}/4$. Directly bounding the selector and target terms separately would lose their cancellation at $a_j$. We therefore use their matching values and gradients and bound the second-order remainders. This requires Hessian bounds along the segment joining $a_j$ to $x$.

Let $y=a_j+t(x-a_j)$, where $0\le t\le1$. For every $l\ne j$,
\[
\|y-a_j\|\le\frac{\Delta_{\mathrm{sep}}}{4}
\le\frac{\|a_j-a_l\|}{4}.
\]
The triangle and reverse triangle inequalities imply
\[
\frac34\|a_j-a_l\|\le\|y-a_l\|\le\frac54\|a_j-a_l\|.
\]
Substituting into the product defining $Q_j$ gives
\[
(3/4)^{2(K-1)}\le Q_j(y)\le(5/4)^{2(K-1)}.
\]
All factors are positive on this segment. Differentiating $\log Q_j$ expresses the product derivatives through sums:
\[
g_j(y):=\nabla\log Q_j(y)
=2\sum_{l\ne j}\frac{y-a_l}{\|y-a_l\|^2},
\]
\[
B_j(y):=\nabla^2\log Q_j(y)
=\sum_{l\ne j}\left\{
\frac{2I_d}{\|y-a_l\|^2}
-\frac{4(y-a_l)(y-a_l)^\top}{\|y-a_l\|^4}
\right\}.
\]
Since $\|y-a_l\|\ge3\Delta_{\mathrm{sep}}/4$, we obtain
\[
\begin{aligned}
\|g_j(y)\|
&\le2\sum_{l\ne j}\frac{1}{\|y-a_l\|}\le\frac{8(K-1)}{3\Delta_{\mathrm{sep}}}.
\end{aligned}
\]
For a matrix $M$, we write $\|M\|_{\mathrm{op}}=\sup_{\|v\|=1}\|Mv\|$ for its operator norm, with Euclidean norms on the vectors. Using $\|I_d\|_{\mathrm{op}}=1$ and $\|zz^\top\|_{\mathrm{op}}=\|z\|^2$ also gives
\[
\begin{aligned}
\|B_j(y)\|_{\mathrm{op}}
&\le6\sum_{l\ne j}\frac{1}{\|y-a_l\|^2}\le\frac{96(K-1)}{9\Delta_{\mathrm{sep}}^2}.
\end{aligned}
\]
The identities
\[
\nabla Q_j=Q_jg_j,
\qquad
\nabla^2Q_j=Q_j(B_j+g_jg_j^\top)
\]
therefore imply
\[
\|\nabla Q_j(y)\|\le C\Delta_{\mathrm{sep}}^{-1},
\qquad
\|\nabla^2Q_j(y)\|_{\mathrm{op}}\le C\Delta_{\mathrm{sep}}^{-2}.
\]
The operator-norm estimates introduce no factor depending on the ambient dimension.

We now apply these estimates to the selectors. On the same segment,
\[
|1-\langle v_j,y-a_j\rangle|
\le1+\|v_j\|\|y-a_j\|
\le1+\frac{K-1}{2}.
\]
Differentiating twice gives
\[
\begin{aligned}
\nabla^2E_j(y)
&=\{1-\langle v_j,y-a_j\rangle\}\nabla^2Q_j(y)-\nabla Q_j(y)v_j^\top-v_j\nabla Q_j(y)^\top,
\end{aligned}
\]
and
\[
\begin{aligned}
\nabla^2H_{j,u}(y)
&=\langle u,y-a_j\rangle\nabla^2Q_j(y)+\nabla Q_j(y)u^\top+u\nabla Q_j(y)^\top.
\end{aligned}
\]
Thus, for every $y\in[a_j,x]$,
\begin{equation}
\begin{aligned}
\|\nabla^2E_j(y)\|_{\mathrm{op}}
&\le
|1-\langle v_j,y-a_j\rangle|
\|\nabla^2Q_j(y)\|_{\mathrm{op}}
+2\|\nabla Q_j(y)\|\|v_j\|\\
&\le C\Delta_{\mathrm{sep}}^{-2}.
\end{aligned}
\label{eq:flat-mass-hessian}
\end{equation}
For the directional selector, $\|u\|\le1$ and $\|y-a_j\|\le\Delta_{\mathrm{sep}}/4$ imply $|\langle u,y-a_j\rangle|\le\Delta_{\mathrm{sep}}/4$. Consequently,
\begin{equation}
\begin{aligned}
\|\nabla^2H_{j,u}(y)\|_{\mathrm{op}}
&\le
|\langle u,y-a_j\rangle|
\|\nabla^2Q_j(y)\|_{\mathrm{op}}
+2\|\nabla Q_j(y)\|\|u\|\\
&\le
\frac{\Delta_{\mathrm{sep}}}{4}
C\Delta_{\mathrm{sep}}^{-2}
+2C\Delta_{\mathrm{sep}}^{-1}\\
&\le C\Delta_{\mathrm{sep}}^{-1}.
\end{aligned}
\label{eq:flat-direction-hessian}
\end{equation}
These Hessian bounds hold uniformly along $[a_j,x]$. Taylor's formula gives local quadratic error bounds for both selectors. We then use the lower bound for $P(x)$ to convert the directional error into~\eqref{B.11}.

\paragraph{From Taylor remainders to the positive witness.}
The Hessian bounds control the approximation errors because the constant and linear terms already agree at $a_j$. Set $h=x-a_j$, so that $\|h\|=r$. For a twice continuously differentiable function $F$, put $f(t)=F(a_j+th)$. Applying the fundamental theorem of calculus twice gives
\[
\begin{aligned}
f(1)-f(0)-f'(0)
&=\int_0^1\{f'(t)-f'(0)\}\,dt=\int_0^1\int_0^t f''(z)\,dz\,dt.
\end{aligned}
\]
The integration region is $0\le z\le t\le1$. For each fixed $z\in[0,1]$, the variable $t$ ranges from $z$ to $1$. Reversing the order of integration therefore gives
\[
\begin{aligned}
\int_0^1\int_0^t f''(z)\,dz\,dt
&=\int_0^1\int_z^1 f''(z)\,dt\,dz\\[4pt]
&=\int_0^1 f''(z)\underbrace{\left(\int_z^1dt\right)}_{=\,1-z}\,dz
\qquad\text{(because }f''(z)\text{ is independent of }t\text{)}\\[4pt]
&=\int_0^1(1-z)f''(z)\,dz.
\end{aligned}
\]
We have therefore proved
\[
f(1)=f(0)+f'(0)+\int_0^1(1-t)f''(t)\,dt,
\]
where the integration variable has been renamed from $z$ to $t$. By the definition $f(t)=F(a_j+th)$ and the chain rule,
\[
\begin{aligned}
f(1)&=F(a_j+h),\\
f(0)&=F(a_j),\\
f'(0)&=\langle\nabla F(a_j),h\rangle,\\
f''(t)&=h^\top\nabla^2F(a_j+th)h.
\end{aligned}
\]
Substituting these expressions into the preceding identity gives
\begin{equation}
F(a_j+h)=F(a_j)+\langle\nabla F(a_j),h\rangle
+\int_0^1(1-t)h^\top\nabla^2F(a_j+th)h\,dt.
\label{B.12}
\end{equation}
For $E_j$, Equation~\eqref{B.7} gives
$E_j(a_j)=1$ and $\nabla E_j(a_j)=0$.
Substituting these values into Equation~\eqref{B.12} yields
\[
E_j(x)-1
=
\int_0^1(1-t)
h^\top\nabla^2E_j(a_j+th)h\,dt.
\]
To bound this remainder, recall that, for any matrix $M$,
\[
\begin{aligned}
|h^\top Mh|
&=|\langle h,Mh\rangle|\\
&\le \|h\|\|Mh\|
\qquad\text{(by Cauchy--Schwarz)}\\
&\le \|h\|^2\|M\|_{\mathrm{op}}
\qquad\text{(by the definition of the operator norm)}\\
&=r^2\|M\|_{\mathrm{op}}.
\end{aligned}
\]
Taking absolute values in the remainder identity and applying
this inequality gives
\begin{equation}
\begin{aligned}
|E_j(x)-1|
&\le
\int_0^1(1-t)
|h^\top\nabla^2E_j(a_j+th)h|\,dt\\
&\le
r^2\int_0^1(1-t)
\|\nabla^2E_j(a_j+th)\|_{\mathrm{op}}\,dt\\
&\le
C\Delta_{\mathrm{sep}}^{-2}r^2
\int_0^1(1-t)\,dt
\qquad\text{(by \eqref{eq:flat-mass-hessian})}\\
&=\frac C2\Delta_{\mathrm{sep}}^{-2}r^2\\
&\le C\Delta_{\mathrm{sep}}^{-2}r^2.
\end{aligned}
\label{eq:flat-mass-taylor-error}
\end{equation}
Here Equation~\eqref{eq:flat-mass-hessian} applies at every
point $a_j+th$, since $0\le t\le1$.

For $H_{j,u}$, Equation~\eqref{B.8} instead gives
$H_{j,u}(a_j)=0$ and $\nabla H_{j,u}(a_j)=u$.
Equation~\eqref{B.12} therefore becomes
\[
H_{j,u}(x)-\langle u,h\rangle
=
\int_0^1(1-t)
h^\top\nabla^2H_{j,u}(a_j+th)h\,dt.
\]
Using the same operator-norm inequality and the Hessian bound
in Equation~\eqref{eq:flat-direction-hessian}, we obtain
\begin{equation}
\begin{aligned}
|H_{j,u}(x)-\langle u,x-a_j\rangle|
&\le
\int_0^1(1-t)
|h^\top\nabla^2H_{j,u}(a_j+th)h|\,dt\\
&\le
r^2\int_0^1(1-t)
\|\nabla^2H_{j,u}(a_j+th)\|_{\mathrm{op}}\,dt\\
&\le
C\Delta_{\mathrm{sep}}^{-1}r^2
\int_0^1(1-t)\,dt\\
&=\frac C2\Delta_{\mathrm{sep}}^{-1}r^2\\
&\le C\Delta_{\mathrm{sep}}^{-1}r^2.
\end{aligned}
\label{eq:flat-direction-taylor-error}
\end{equation}

It remains to express the directional error in terms of $P(x)$.
Since $c(x)=j$, the positive-coercivity bound
in Equation~\eqref{proof:B2} gives
\[
P(x)
\ge
c_{\mathrm{pos}}\Delta_{\mathrm{sep}}^{2s_{\max}-2}
\|x-a_j\|^2
=
c_{\mathrm{pos}}\Delta_{\mathrm{sep}}^{2s_{\max}-2}r^2.
\]
Dividing by the positive coefficient on the right yields
\begin{equation}
r^2
\le
c_{\mathrm{pos}}^{-1}
\Delta_{\mathrm{sep}}^{-(2s_{\max}-2)}P(x).
\label{eq:flat-local-distance-witness}
\end{equation}

Equation~\eqref{eq:flat-direction-taylor-error} gives
\[
\begin{aligned}
|H_{j,u}(x)-\langle u,x-a_j\rangle|
&\le
C\Delta_{\mathrm{sep}}^{-1}r^2\\
&\le
C\Delta_{\mathrm{sep}}^{-1}
\Delta_{\mathrm{sep}}^{-(2s_{\max}-2)}P(x)\\
&=
C\Delta_{\mathrm{sep}}^{-(2s_{\max}-1)}P(x),
\end{aligned}
\]
where $c_{\mathrm{pos}}^{-1}$ is absorbed into $C$.

Because $c(x)=j$, the last inequality is precisely~\eqref{B.11} in the remaining case $r\le\Delta_{\mathrm{sep}}/4$. The cases $c(x)\ne j$ and $c(x)=j$ with $\|x-a_j\|>\Delta_{\mathrm{sep}}/4$ were established above. Thus~\eqref{B.11} holds for every $x\in B_d(0,R)$.
\end{proof}
\subsection{Polynomial approximation of cluster indicators}

The polynomials in the preceding subsection distinguish individual true centers, and their bounds involve inverse powers of $\Delta_{\mathrm{sep}}$. For cluster probabilities, we only need to distinguish centers belonging to different clusters. We therefore construct a polynomial that takes the same value at every center in a given cluster, with bounds that remain uniform as centers within that cluster approach one another. These polynomials will control the cluster-sequence probability error $J_I$ appearing in the transport bounds of Appendix~\ref{app:gaussian-moments}.

Our construction reduces interpolation to one dimension. We first choose a projection that preserves pairwise distances up to a factor depending only on $K$. We then establish an interpolation bound for smooth functions that does not involve inverse node spacings. This allows us to prescribe cluster values and zero gradients without introducing additional powers of $\Delta_{\mathrm{sep}}^{-1}$.

\begin{lemma}[A projection preserving pairwise separation]\label{lem:separation-projection}
Let $a_1,\ldots,a_K$ be distinct points in a Euclidean space, with $K\ge2$. There exists a unit vector $v$ in the span of their differences such that
\begin{equation}
|\langle v,a_j-a_l\rangle|
\ge c_K\|a_j-a_l\|,
\qquad j\ne l,
\label{eq:projection-separation}
\end{equation}
where $c_K>0$ depends only on $K$.
\end{lemma}

\begin{proof}
Let $S_0=\operatorname{span}\{a_j-a_1:2\le j\le K\}$ and
$q=\dim S_0\in[1,K-1]$. If $q=1$, either unit direction in $S_0$
works with $c_K=1$. Otherwise, take a uniform random unit vector
$V\in S_0$. For every unit $u\in S_0$, the density of
$\langle V,u\rangle$ is bounded on $[-1/2,1/2]$ by a constant
$M_q<\infty$, and thus
$\mathbb P(|\langle V,u\rangle|<\eta)\le2M_q\eta$
for $\eta\le1/2$. There are $N_K=K(K-1)/2$ normalized difference
directions $u_{jl}=(a_j-a_l)/\|a_j-a_l\|_2$. Since
$q\le K-1$, choose
$\eta_K=\min\{1/2,[4N_K\max_{2\le q'\le K-1}M_{q'}]^{-1}\}$.
A union bound shows that with probability at least $1/2$,
$|\langle V,u_{jl}\rangle|\ge\eta_K$ for all pairs.
Fix such a direction $v$ and set $c_K=\eta_K$.
\end{proof}

The projection preserves the gap between different clusters, but projected centers within a cluster may still be arbitrarily close. We therefore need interpolation bounds that remain valid as nodes approach one another. Repeated nodes will impose derivative conditions, allowing us to prescribe both the value and the first derivative at each projected center.

\begin{lemma}[Interpolation bounds independent of node spacings]\label{lem:confluent-interpolation}
Let $z_1,\ldots,z_N\in[-R,R]$ be interpolation nodes, with repetitions allowed, and let $f\in C^N([-R,R])$. A node occurring $m$ times prescribes the value of $f$ and its derivatives through order $m-1$ at that node. There is a unique polynomial $H_f$ of degree at most $N-1$ satisfying these conditions. It has the Newton representation
\begin{equation}
H_f(t)
=
\sum_{b=0}^{N-1}
f[z_1,\ldots,z_{b+1}]
\prod_{a=1}^{b}(t-z_a),
\label{eq:hermite-newton}
\end{equation}
where an empty product equals one, and
\begin{equation}
|f[z_1,\ldots,z_{b+1}]|
\le \frac{\|f^{(b)}\|_{\infty,[-R,R]}}{b!},
\qquad 0\le b\le N-1.
\label{eq:divided-difference-bound}
\end{equation}
Consequently, bounds on $f$ and its derivatives through order $N-1$ give bounds on $\|H_f\|_{\infty,[-R,R]}$ and $\|H_f''\|_{\infty,[-R,R]}$ that do not depend on the node spacings.
\end{lemma}

\begin{proof}
For $0\le b\le N-1$, the Hermite conditions at the first $b+1$
nodes (counted with multiplicity) determine a unique degree-$b$
polynomial $H_b$: a polynomial in the kernel would be divisible by
a degree-$(b+1)$ product of node factors, hence vanish identically.
As $H_b-H_{b-1}$ satisfies the first $b$ homogeneous conditions,
$H_b-H_{b-1}=c_b\omega_b$, where
$\omega_b(t)=\prod_{a=1}^b(t-z_a)$ and $\omega_0=1$.
Summing this recursion proves the Newton formula
\eqref{eq:hermite-newton} with
$c_b=f[z_1,\ldots,z_{b+1}]$.

The function $f-H_b$ has $b+1$ zeros counted with multiplicity.
Repeated application of Rolle's theorem, including confluent nodes,
produces a point $\xi$ between the extreme nodes where
$(f-H_b)^{(b)}(\xi)=0$. Since
$H_b^{(b)}=b!c_b$, this gives
$|c_b|\le\|f^{(b)}\|_\infty/b!$, proving
\eqref{eq:divided-difference-bound} without division by node gaps.
Furthermore, for $0\le\ell\le b$,
\[
\|\omega_b^{(\ell)}\|_\infty
\le\frac{b!}{(b-\ell)!}(2R)^{b-\ell}.
\]
Therefore the Newton formula yields
\[
\|H_f^{(\ell)}\|_\infty\le
\sum_{b=\ell}^{N-1}
\frac{\|f^{(b)}\|_\infty}{(b-\ell)!}(2R)^{b-\ell},
\qquad 0\le\ell<N.
\]
Take $\ell=0,2$; if $N\le2$, then $H_f''=0$.
All constants are independent of the spacing and repetition of nodes.
\end{proof}

We now construct the smooth function to be interpolated. Its value will be constant near every projected center, with value one for the selected cluster and zero for the others. The gap between different clusters allows these values to be joined smoothly with uniformly bounded derivatives. Interpolating both values and first derivatives then produces the required polynomial.

\begin{lemma}[Polynomial approximation of cluster indicators]\label{lem:cluster-polynomials}
Under the structured geometry, for each cluster $C_m$ there exists a polynomial $F_m$ of degree at most $2K-1$ such that
\begin{equation}
F_m(a_j)=\mathbf1\{j\in C_m\},
\qquad
\nabla F_m(a_j)=0,
\qquad j\in[K],
\label{eq:cluster-polynomial-values}
\end{equation}
\begin{equation}
\|F_m\|_{\infty,B_d(0,R)}
+
\sup_{x\in B_d(0,R)}\|\nabla^2F_m(x)\|_{\mathrm{op}}
\le C,
\label{eq:cluster-polynomial-norms}
\end{equation}
and
\begin{equation}
|F_m(x)-\mathbf1\{c(x)\in C_m\}|
\le C\|x-a_{c(x)}\|^2,
\qquad x\in B_d(0,R).
\label{eq:cluster-polynomial-error}
\end{equation}
The constant depends only on $K,R,D_0$, independently of $\Delta_{\mathrm{sep}}$ and the ambient dimension $d$. For a single cluster, one may take $F_1\equiv1$.
\end{lemma}

\begin{proof}
If there is only one cluster, $F_1\equiv1$ satisfies all three claims. Assume that there are at least two clusters, so $K\ge2$.

\paragraph{Projected centers and the gap between clusters.}
Choose the unit vector $v$ from Lemma~\ref{lem:separation-projection} and set
\[
t_j=\langle v,a_j\rangle,
\qquad j\in[K].
\]
The projected centers are distinct and satisfy $|t_j|\le R$. For centers in different clusters,
\[
\begin{aligned}
|t_j-t_l|
&=|\langle v,a_j-a_l\rangle|\\
&\ge c_K\|a_j-a_l\|
\qquad \text{(by \eqref{eq:projection-separation})}\\
&\ge c_KD_0
\qquad \text{(by the between-cluster separation)}.
\end{aligned}
\]
Write $\gamma=c_KD_0$. Thus different clusters remain separated by at least $\gamma$ on the line, even though projected centers within a cluster may be arbitrarily close.

\paragraph{A smooth function with the required cluster values.}
Fix a cluster $C_m$. We construct a smooth function that equals one near its projected centers and zero near all other projected centers. Begin with a fixed smooth cutoff $\psi$ satisfying
\[
0\le\psi\le1,\qquad
\psi(t)=1\text{ for }|t|\le\frac12,\qquad
\psi(t)=0\text{ for }|t|\ge1.
\]
For an explicit choice, define
\[
\zeta(t)=
\begin{cases}
e^{-1/t},&t>0,\\
0,&t\le0,
\end{cases}
\qquad
\chi(t)=\frac{\zeta(t)}{\zeta(t)+\zeta(1-t)},
\qquad
\psi(t)=\chi\!\left(\frac{4(1-t^2)}{3}\right).
\]
The function $\chi$ equals zero on $(-\infty,0]$ and one on $[1,\infty)$, which gives the stated properties of $\psi$. Since $\psi$ is smooth and compactly supported, each derivative has a finite bound
\[
M_b:=\sup_{t\in\mathbb R}|\psi^{(b)}(t)|.
\]

For every $j\in C_m$, place a copy of this cutoff around $t_j$ by setting
\[
b_j(t)=\psi\!\left(\frac{3(t-t_j)}{\gamma}\right).
\]
Then
\[
b_j(t)=1\text{ if }|t-t_j|\le\frac{\gamma}{6},
\qquad
b_j(t)=0\text{ if }|t-t_j|\ge\frac{\gamma}{3}.
\]
To combine these cutoffs, define
\[
f_m(t)=1-\prod_{j\in C_m}\{1-b_j(t)\}.
\]
If any one of the cutoffs equals one, the product vanishes and $f_m=1$. If all cutoffs vanish, the product equals one and $f_m=0$. In particular, $f_m$ equals one in a neighborhood of every $t_j$ with $j\in C_m$.

For $l\notin C_m$ and $|t-t_l|<\gamma/3$, every $j\in C_m$ satisfies
\[
\begin{aligned}
|t-t_j|
&\ge |t_l-t_j|-|t-t_l|\\
&>\gamma-\frac{\gamma}{3}
=\frac{2\gamma}{3}.
\end{aligned}
\]
All $b_j(t)$ therefore vanish in this neighborhood of $t_l$, so $f_m=0$ there. We have obtained
\[
f_m(t_j)=\mathbf1\{j\in C_m\},
\qquad
f_m^{(b)}(t_j)=0
\quad\text{for every }b\ge1.
\]

The derivative bounds of $f_m$ depend on $\gamma$, not on distances between centers in the same cluster. To verify this, the chain rule gives
\[
\|b_j^{(b)}\|_{\infty,\mathbb R}
\le \left(\frac3\gamma\right)^b M_b.
\]
Put
\[
B_\gamma
=
\max\left\{
1,\,
\max_{1\le b\le2K}
\left(\frac3\gamma\right)^bM_b
\right\},
\qquad
s_m=|C_m|.
\]
All derivatives through order $2K$ of each factor $1-b_j$, including order zero, are bounded by $B_\gamma$. In each term of the product rule, the $b$ derivatives are distributed among the $s_m$ factors $1-b_j$. If these factors receive $\nu_1,\ldots,\nu_{s_m}$ derivatives, respectively, then $\nu_1+\cdots+\nu_{s_m}=b$, the coefficient is $b!/(\nu_1!\cdots\nu_{s_m}!)$, and the product of derivatives has absolute value at most $B_\gamma^{s_m}$. Thus, for $1\le b\le2K$,
\[
\begin{aligned}
\|f_m^{(b)}\|_{\infty,\mathbb R}
&\le
\sum_{\substack{\nu_1+\cdots+\nu_{s_m}=b\\\nu_i\ge0}}
\frac{b!}{\nu_1!\cdots\nu_{s_m}!}
B_\gamma^{s_m}\\
&=s_m^bB_\gamma^{s_m}
\qquad \text{(by the multinomial formula)}.
\end{aligned}
\]
Also, $0\le f_m\le1$. Since $s_m\le K$ and $\gamma=c_KD_0$, these bounds depend only on $K,D_0$.

\paragraph{Replacing the smooth function by a polynomial.}
Apply Lemma~\ref{lem:confluent-interpolation} to $f_m$ with each projected center repeated twice,
\[
(t_1,t_1,\ldots,t_K,t_K).
\]
The two occurrences of $t_j$ prescribe its function value and first derivative. The resulting polynomial $H_m$ has degree at most $2K-1$ and satisfies
\[
\underbrace{H_m(t_j)=f_m(t_j)=\mathbf1\{j\in C_m\}}
_{\text{prescribed cluster value}},
\qquad
\underbrace{H_m'(t_j)=f_m'(t_j)=0}
_{\text{prescribed zero derivative}}.
\]
Moreover, the derivative bounds just established give
\[
\begin{aligned}
\|H_m\|_{\infty,[-R,R]}
&\le
\sum_{b=0}^{2K-1}
\frac{\|f_m^{(b)}\|_{\infty,[-R,R]}}{b!}(2R)^b,\\
\|H_m''\|_{\infty,[-R,R]}
&\le
\sum_{b=2}^{2K-1}
\frac{\|f_m^{(b)}\|_{\infty,[-R,R]}}{(b-2)!}(2R)^{b-2}.
\end{aligned}
\]
Both sums are bounded by a constant depending only on $K,R,D_0$. In particular, no inverse within-cluster spacing enters these estimates.

Return to the original space by defining
\[
F_m(x)=H_m(\langle v,x\rangle).
\]
The chain rule gives
\[
\nabla F_m(x)=H_m'(\langle v,x\rangle)v,
\qquad
\nabla^2F_m(x)=H_m''(\langle v,x\rangle)vv^\top.
\]
Evaluating at $a_j$ proves \eqref{eq:cluster-polynomial-values}. Furthermore, $|\langle v,x\rangle|\le R$ for $x\in B_d(0,R)$ and $\|vv^\top\|_{\mathrm{op}}=1$, so
\[
\begin{aligned}
\|F_m\|_{\infty,B_d(0,R)}
&\le\|H_m\|_{\infty,[-R,R]},\\
\sup_{x\in B_d(0,R)}\|\nabla^2F_m(x)\|_{\mathrm{op}}
&\le\|H_m''\|_{\infty,[-R,R]}.
\end{aligned}
\]
This proves \eqref{eq:cluster-polynomial-norms}, with no factor depending on $d$.

\paragraph{Bounding the error away from the true centers.}
Fix $x\in B_d(0,R)$, let $a=a_{c(x)}$ be its assigned true center, and write $h=x-a$. The segment $a+th$, $0\le t\le1$, lies in $B_d(0,R)$ by convexity. Taylor's formula with integral remainder gives
\[
\begin{aligned}
F_m(x)
&=
\underbrace{F_m(a)}_{\mathbf1\{c(x)\in C_m\}}
+
\underbrace{\langle\nabla F_m(a),h\rangle}_{0}\\
&\qquad
+\int_0^1(1-t)\,
h^\top\nabla^2F_m(a+th)h\,dt.
\end{aligned}
\]
The zero gradient removes the linear term, leaving an error controlled by the squared displacement. Specifically,
\[
\begin{aligned}
|F_m(x)-\mathbf1\{c(x)\in C_m\}|
&\le
\int_0^1(1-t)\,
|h^\top\nabla^2F_m(a+th)h|\,dt\\
&\le
\|h\|^2
\int_0^1(1-t)
\|\nabla^2F_m(a+th)\|_{\mathrm{op}}\,dt\\
&\le
C\|h\|^2\int_0^1(1-t)\,dt
\qquad \text{(by \eqref{eq:cluster-polynomial-norms})}\\
&=\frac C2\|x-a_{c(x)}\|^2.
\end{aligned}
\]
After adjusting the constant, this proves \eqref{eq:cluster-polynomial-error}. The argument applies to every $x\in B_d(0,R)$ and requires no proximity assumption on the fitted centers.
\end{proof}

\subsection{Joint selectors with vanishing first derivatives}
\label{app:joint-first-order}

For a joint probability, it is enough to interpolate its values and first derivatives on the true support. Requiring a separate flat selector in every coordinate imposes additional conditions that are not needed for the quadratic approximation error. The following construction uses ordinary cardinal polynomials in the remaining coordinates.

\begin{lemma}[Joint probability selectors]
\label{lem:joint-first-order}
Suppose that the structured geometry holds, put $s=s_{\max}$ and $\Delta=\Delta_{\mathrm{sep}}$, and define
\[
d_r=r(s-1)+s,\qquad \beta_r=\max\{2s,d_r\}.
\]
For every $j=(j_1,\ldots,j_r)\in[K]^r$, there is a polynomial $U_j$ of degree bounded in terms of $K,r$ such that
\begin{align}
U_j(a_{i_1},\ldots,a_{i_r})&=\mathbf1\{i=j\},\qquad
\nabla U_j(a_{i_1},\ldots,a_{i_r})=0,\label{eq:joint-first-values}\\
\|U_j\|_{\infty,B_d(0,R)^r}&\le C_r\Delta^{-d_r},\label{eq:joint-first-norm}\\
\left|U_j(x_1,\ldots,x_r)-\mathbf1\{(c(x_1),\ldots,c(x_r))=j\}\right|
&\le C_r\Delta^{-\beta_r}\sum_{\ell=1}^r P(x_\ell).
\label{eq:joint-first-error}
\end{align}
Here $P(x)=\prod_{i=1}^K\|x-a_i\|_2^2$ is the positive polynomial in~\eqref{proof:B1}. Constants are independent of the ambient dimension and of $\Delta$.
\end{lemma}

\begin{proof}
We first construct localized one-coordinate polynomials, and then combine them so that their first derivatives cancel.

\paragraph{Cluster selectors with higher-order flatness.}
Choose $N=2s+2$. Repeat the construction of Lemma~\ref{lem:cluster-polynomials}, now repeating each projected true center $N+1$ times in Lemma~\ref{lem:confluent-interpolation}. The smooth cluster indicator used there is constant near every projected center. Its derivatives of every fixed order are bounded independently of within-cluster spacing. The resulting polynomials $F_m$ therefore satisfy
\[
F_m(a_i)=\mathbf1\{i\in C_m\},\qquad
D^vF_m(a_i)=0\quad(1\le v\le N),
\]
and have uniformly bounded derivatives of every fixed order on $B_d(0,R)$. Their degree is at most $(N+1)K-1$. These claims follow from the same Newton interpolation formula, with derivative orders up to $(N+1)K-1$ in place of $2K-1$. For a single cluster, take $F_1=1$.

For $j\in C_m$, define the ordinary within-cluster cardinal polynomial
\[
\ell_j(x)=\prod_{i\in C_m\setminus\{j\}}
\frac{\langle x-a_i,a_j-a_i\rangle}{\|a_j-a_i\|_2^2},
\]
and its flat counterpart
\[
e_j(x)=\{1-2\langle\nabla\ell_j(a_j),x-a_j\rangle\}\ell_j(x)^2.
\]
Their degrees are at most $s_m-1$ and $2s_m-1$, respectively, where $s_m=|C_m|$. Within $C_m$, both have the cardinal values, and $\nabla e_j$ vanishes at every center in that cluster. No such assertion is needed for the unlocalized $e_j$ at centers outside $C_m$. Set
\[
L_j=F_m\ell_j,\qquad E_j=F_me_j.
\]
Then $L_j(a_i)=E_j(a_i)=\mathbf1\{i=j\}$ for every true center, and $\nabla E_j(a_i)=0$. The additional flatness of $F_m$ also controls these polynomials near centers outside $C_m$.

\paragraph{Derivative bounds near the true support.}
Assume first $s\ge2$, and put $p=s-1$, $q=2s-1$. For any true center $a_i$, any $y\in B_d(0,R)$, and $v=0,1,2$, the preceding polynomials satisfy
\begin{equation}
\|D^vL_j(y)\|\le C\Delta^{-p}(\Delta+\|y-a_i\|)^{p-v},\qquad
\|D^vE_j(y)\|\le C\Delta^{-q}(\Delta+\|y-a_i\|)^{q-v}.
\label{eq:joint-local-derivatives}
\end{equation}
Here derivative norms are operator norms, with absolute value for $v=0$. To verify the bounds when $i\in C_m$, every linear factor in $\ell_j$ has numerator bounded by $C(\Delta+\|y-a_i\|)$ and denominator at least $\Delta$. Differentiation removes at most one such numerator per derivative. Also $\|\nabla\ell_j(a_j)\|\le C\Delta^{-1}$, giving the stated scaling for $e_j$. Multiplication by the uniformly smooth $F_m$ preserves these bounds. More explicitly, put $t=\Delta+\|y-a_i\|$ and $p_m=s_m-1$. A derivative of order $b\le p_m$ of the product of linear factors is bounded by $C\Delta^{-p_m}t^{p_m-b}$, and a derivative of order $b>p_m$ is identically zero. A term in $D^v(F_m\ell_j)$ with $v-b$ derivatives on $F_m$ is therefore bounded by the desired expression times $t^{v-b}$, which is uniformly bounded on the ball. This also covers singleton clusters, where $\ell_j=1$. The same product argument applies to $e_j$ with degree bound $2s_m-1$; a derivative beyond its actual degree is zero, not a term with a negative power. If $s_m<s$, increase the scaling degrees to $p,q$, using $\Delta/t\le1$.

When $i\notin C_m$, Taylor's formula and the vanishing derivatives of $F_m$ give
\[
\|D^vF_m(y)\|\le C\|y-a_i\|^{N+1-v},\qquad v=0,1,2.
\]
The derivatives through order two of $\ell_j,e_j$ on the ball are bounded by $C\Delta^{-p},C\Delta^{-q}$. The product rule proves~\eqref{eq:joint-local-derivatives}, since $N+1\ge q$ and all distances are bounded. This argument does not require the projected interval hulls of different clusters to be disjoint.

\paragraph{Cancelling the first derivatives.}
Define
\begin{equation}
U_j(x_1,\ldots,x_r)
=\sum_{\ell=1}^r E_{j_\ell}(x_\ell)
\prod_{u\ne\ell}L_{j_u}(x_u)
-(r-1)\prod_{u=1}^rL_{j_u}(x_u).
\label{eq:joint-first-construction}
\end{equation}
At a true grid point, each of the first $r$ terms has value $\mathbf1\{i=j\}$. For a derivative in coordinate $v$, the term with $E_{j_v}$ has zero derivative. Each of the other $r-1$ terms has exactly the same derivative as the last product, because $E_{j_u}$ and $L_{j_u}$ have the same values on the true support. These derivatives cancel, proving~\eqref{eq:joint-first-values}. The norm bounds for the factors give~\eqref{eq:joint-first-norm}, since $q+(r-1)p=d_r$ and $rp\le d_r$.

Let $a_{i_\ell}=a_{c(x_\ell)}$ and $M=\max_\ell\|x_\ell-a_{i_\ell}\|$. Along the segment from this true grid point to $(x_1,\ldots,x_r)$, the product rule and~\eqref{eq:joint-local-derivatives} bound the Hessian of the first sum by
$C_r\Delta^{-d_r}(\Delta+M)^{d_r-2}$.
The last product has the bound $C_r\Delta^{-rp}(\Delta+M)^{rp-2}$, which is no larger, since $\Delta/(\Delta+M)\le1$. When $p=1$ and two derivatives fall on one $L_j$, one cannot simply replace its individual distance by $M$ in a negative exponent. Instead, $\ell_j$ has degree at most one, so $D^2\ell_j=0$. The product rule gives
\[
D^2(F_m\ell_j)=\ell_jD^2F_m
 +(DF_m)\otimes D\ell_j+D\ell_j\otimes(DF_m),
\]
whose norm is at most $C\Delta^{-1}$ on the ball. Since $\Delta+M\le1+2R$, this is at most $C(1+2R)\Delta^{-1}(\Delta+M)^{-1}$. All mixed-coordinate Hessian blocks differentiate each affected factor once and have nonnegative remaining powers. For $r=1$, the final product in equation~\eqref{eq:joint-first-construction} has coefficient zero and $U_j=E_j$. Thus the displayed Hessian bound covers $s=2$, $r=1$, and singleton clusters without replacing an individual distance by a maximum in an unhandled negative power. Taylor's formula therefore yields
\begin{equation}
|U_j(x_1,\ldots,x_r)-\mathbf1\{i=j\}|
\le C_r\Delta^{-d_r}M^2(\Delta+M)^{d_r-2}.
\label{eq:joint-first-taylor}
\end{equation}
The factor $r$ from the squared length of the segment is absorbed into $C_r$.

\paragraph{Comparison with the positive polynomial.}
Writing $\rho=\|x-a_{c(x)}\|$, we have the stronger pointwise bound
\begin{equation}
P(x)\ge c\rho^2(\Delta+\rho)^{2s-2}.
\label{eq:positive-distance-refined}
\end{equation}
Indeed, for any other true center $a_i$, nearest-center assignment and separation give
$\|x-a_i\|\ge\rho$ and $\|x-a_i\|\ge\Delta/2$.
If $\rho<D_0/4$, the factors from other clusters are at least $3D_0/4$, while the $s_m-1$ remaining factors in the assigned cluster are at least $c(\Delta+\rho)$. Increasing the exponent from $2s_m-2$ to $2s-2$ only changes the constant on the bounded ball. If $\rho\ge D_0/4$, all factors are bounded below by a fixed positive number, which proves the same bound. In the single-cluster case, the first two inequalities prove it directly, without $D_0$.

Consequently, $\sum_\ell P(x_\ell)\ge cM^2(\Delta+M)^{2s-2}$. Dividing~\eqref{eq:joint-first-taylor} by this lower bound leaves
$C_r\Delta^{-d_r}(\Delta+M)^{d_r-2s}$.
If $d_r\ge2s$, boundedness of $M$ gives the multiplier $C_r\Delta^{-d_r}$; otherwise $\Delta+M\ge\Delta$ gives $C_r\Delta^{-2s}$. This proves~\eqref{eq:joint-first-error}.

Finally, if $s=1$, distinct true centers are separated by at least $D_0$, so $\Delta$ is bounded below. Take $U_j=\prod_\ell F_{\kappa(j_\ell)}(x_\ell)$ with the original flat cluster selectors, where $\kappa(j_\ell)$ is the cluster index, not the state label. Their norms and Hessians are uniformly bounded, and $P(x)\ge c\rho^2$. Taylor's formula proves all the stated bounds, with $d_r=1$ and $\beta_r=2$.
\end{proof}

\subsection{Joint selectors without cluster assumptions}
\label{app:unstructured-joint-selectors}

The product of distances from a true center to the other true centers contains at least one distance comparable to the diameter of the true support. Retaining this factor in the interpolation bounds sharpens both joint probability and transport recovery without cluster assumptions.

\begin{lemma}[Diameter-dependent joint selectors]
\label{lem:diameter-joint-selectors}
For $K\ge2$ arbitrary distinct true centers, put
\[
\Delta=\Delta_{\mathrm{sep}},\qquad
D_\star=\max_{i,j}\|a_i-a_j\|_2,\qquad
\mathcal D_j=\prod_{i\ne j}\|a_i-a_j\|_2.
\]
Then $\Delta\le D_\star\le2R$ and
$\mathcal D_j\ge(D_\star/2)\Delta^{K-2}$.
Let $P(x)=\prod_i\|x-a_i\|_2^2$, and define
\[
L_j(x)=\prod_{i\ne j}\frac{\langle x-a_i,a_j-a_i\rangle}{\|a_j-a_i\|_2^2},
\quad
Q_j(x)=\prod_{i\ne j}\frac{\|x-a_i\|_2^2}{\|a_j-a_i\|_2^2},
\quad
E_j(x)=\{1-\langle\nabla Q_j(a_j),x-a_j\rangle\}Q_j(x).
\]
For every $x\in B_d(0,R)$,
\begin{align}
|L_j(x)-\mathbf1\{c(x)=j\}|&\le\frac{C\sqrt{P(x)}}{\Delta\mathcal D_j},
\label{eq:diameter-linear-error}\\
|E_j(x)-\mathbf1\{c(x)=j\}|&\le\frac{CP(x)}{\Delta^2\mathcal D_j^2}.
\label{eq:diameter-flat-error}
\end{align}
For $r\ge1$, set
\[
d_r=r(K-1)+K,\qquad \gamma_{K,r}=\max\{2K,d_r\}.
\]
For each $j=(j_1,\ldots,j_r)\in[K]^r$, the construction in~\eqref{eq:joint-first-construction}, with the one-coordinate polynomials above, gives a polynomial $U_j$ of degree at most $d_r$, with the cardinal values and vanishing first derivatives on the true grid, such that
\begin{align}
\|U_j\|_\infty&\le C_rD_\star^{-2}\Delta^{-(\gamma_{K,r}-2)},
\label{eq:diameter-joint-norm}\\
|U_j(x)-\mathbf1\{c(x)=j\}|
&\le C_rD_\star^{-2}\Delta^{-(\gamma_{K,r}-2)}
\sum_{\ell=1}^rP(x_\ell).
\label{eq:diameter-joint-error}
\end{align}
Here $c(x)=(c(x_1),\ldots,c(x_r))$, and the norm is over $B_d(0,R)^r$. Constants depend only on $K,r,R$.
\end{lemma}

\begin{proof}
We first retain the distance products in the one-coordinate estimates, then compare the largest nearest-center distance with the diameter of the true support.

\paragraph{One-coordinate estimates.}
For every $j$, at least one true center is at distance at least $D_\star/2$ from $a_j$, by the triangle inequality applied to a pair attaining the diameter. Every other distance is at least $\Delta$, proving the bound on $\mathcal D_j$.

Write $\rho_j=\|x-a_j\|_2$. If $c(x)\ne j$, then $\rho_j\ge\Delta/2$, and
\[
|L_j(x)|\le\frac{\sqrt{P(x)}}{\rho_j\mathcal D_j},\qquad
Q_j(x)=\frac{P(x)}{\rho_j^2\mathcal D_j^2}.
\]
Also $\|\nabla Q_j(a_j)\|\le C/\Delta$. These inequalities give~\eqref{eq:diameter-linear-error}--\eqref{eq:diameter-flat-error} when $c(x)\ne j$.
If $c(x)=j$, nearest-center assignment gives
$P(x)\ge c\rho_j^2\mathcal D_j^2$.
For $\rho_j\ge\Delta/4$, the same formulas bound both the polynomial and indicator terms by the asserted right sides. For $\rho_j<\Delta/4$, the factors defining $L_j$ give
$|L_j(x)-1|\le C\rho_j/\Delta$.
The polynomial $E_j$ is the same flat selector considered in Lemma~\ref{lem:flat-selectors}. Its interpolation values give $E_j(a_j)=1$ and $\nabla E_j(a_j)=0$. Thus the constant and linear terms cancel when Taylor's formula is applied to $E_j(x)-1$. The local Taylor estimate~\eqref{eq:flat-mass-taylor-error}, proved using the Hessian bound~\eqref{eq:flat-mass-hessian}, gives
\[
|E_j(x)-1|\le C\rho_j^2/\Delta^2.
\]
To obtain the required comparison with $P(x)$, use the nearest-center bound $P(x)\ge c\rho_j^2\mathcal D_j^2$, which implies
\[
\rho_j^2\le \frac{CP(x)}{\mathcal D_j^2}.
\]
Substitution yields
\[
|E_j(x)-1|
\le \frac{CP(x)}{\Delta^2\mathcal D_j^2}.
\]
Since $c(x)=j$, this is exactly~\eqref{eq:diameter-flat-error} in the remaining case. These arguments also cover a tie in nearest-center assignment.

The cardinal values and cancellation of first derivatives of $U_j$ follow as in Lemma~\ref{lem:joint-first-order}; that cancellation uses only the equal cardinal values of $L_j,E_j$ and the vanishing gradients of $E_j$. For $r=1$, $U_j=E_j$, and~\eqref{eq:diameter-flat-error} already proves~\eqref{eq:diameter-joint-error}. We henceforth take $r\ge2$, so $d_r\ge2K$ and $\gamma_{K,r}=d_r$.

\paragraph{Points within the true-support diameter.}
Put $M=\max_\ell\|x_\ell-a_{c(x_\ell)}\|_2$ and suppose $M\le D_\star$. Every distance from $x_\ell$ to a true center is then at most $2D_\star$. Cauchy--Schwarz and the lower bound on $\mathcal D_j$ give
\[
|L_j(x_\ell)|\le C(D_\star/\Delta)^{K-2}=:B_\star,
\qquad B_\star\ge1,
\]
after increasing the constant. Write $v_\ell=\mathbf1\{c(x_\ell)=j_\ell\}$,
$\varepsilon_\ell=L_{j_\ell}(x_\ell)-v_\ell$, and
$\eta_\ell=E_{j_\ell}(x_\ell)-v_\ell$.
Expanding~\eqref{eq:joint-first-construction}, the constant term is $\prod_\ell v_\ell$ and all terms containing exactly one $\varepsilon_\ell$ cancel. Terms involving an $\eta_\ell$ are bounded by $|\eta_\ell|B_\star^{r-1}$; the remaining terms contain at least two $\varepsilon$ factors and are bounded by $C_rB_\star^{r-2}\sum_{\ell<u}|\varepsilon_\ell\varepsilon_u|$. By~\eqref{eq:diameter-linear-error}--\eqref{eq:diameter-flat-error} and $2\sqrt{P(x_\ell)P(x_u)}\le P(x_\ell)+P(x_u)$,
\begin{align*}
|U_j(x)-\mathbf1\{c(x)=j\}|
&\le C_rD_\star^{-2}\Delta^{-(2K-2)}
(D_\star/\Delta)^{(r-1)(K-2)}\sum_\ell P(x_\ell)\\
&\le C_rD_\star^{-2}\Delta^{-(d_r-2)}\sum_\ell P(x_\ell).
\end{align*}
For the last inequality, the ratio of the first multiplier to the second is
$D_\star^{(r-1)(K-2)}\Delta^{r-2}$, which is bounded by a constant depending only on $K,r,R$.

\paragraph{Points beyond the true-support diameter.}
Suppose $M>D_\star$. Every distance from a selected $x_\ell$ to a true center is at most $2M$. The product formulas and $\|\nabla Q_j(a_j)\|\le C/\Delta$ therefore give
\[
|L_j(x_\ell)|\le\frac{CM^{K-1}}{D_\star\Delta^{K-2}},\qquad
|E_j(x_\ell)|\le\frac{CM^{2K-1}}{D_\star^2\Delta^{2K-3}}.
\]
Consequently,
\[
|U_j(x)|\le
\frac{C_rM^{d_r}}{D_\star^{r+1}\Delta^{d_r-r-1}}.
\]
The last product in~\eqref{eq:joint-first-construction} is absorbed in this bound because $D_\star\Delta^{K-1}/M^K\le1$.
For a coordinate attaining $M$, all its distances to true centers are at least $M$, so $\sum_\ell P(x_\ell)\ge M^{2K}$. Dividing the preceding bound by this sum and using $M\le2R$, $d_r\ge2K$, and $\Delta\le D_\star$ yields
\[
\frac{|U_j(x)|}{\sum_\ell P(x_\ell)}
\le C_rD_\star^{-2}\Delta^{-(d_r-2)}.
\]
The indicator term satisfies the same bound, since
$D_\star^{-2}\Delta^{-(d_r-2)}M^{2K}\ge1$ when $M>D_\star\ge\Delta$ and $\Delta\le1$.
This proves~\eqref{eq:diameter-joint-error}. Finally, $P$ is uniformly bounded on the ball and $D_\star\le2R$, so~\eqref{eq:diameter-joint-error} also gives~\eqref{eq:diameter-joint-norm}.
\end{proof}

\subsection{Exact interpolation at paired centers}
\label{app:paired-interpolation}

The global selecting polynomials constructed above have the required values at true centers, but their values at fitted centers introduce approximation errors. For the matched bounds under exact specification, we assume a pairing satisfying the proximity condition below and impose the required values on both supports. We construct these polynomials by first making them vanish at every other true--fitted pair and then adjusting their values at the selected pair. This construction also covers pairs whose two centers coincide.

\subsubsection{A polynomial that vanishes at every other pair}

Assume that the fitted centers are labeled so that $\|b_j-a_j\|_2\le\Delta_{\mathrm{sep}}/4$ for every $j\in[K]$. Both $a_j$ and $b_j$ belong to $B_d(0,R)$. For $j\ne\ell$, write $d_{j\ell}=\|a_j-a_\ell\|_2$ and define
\begin{equation}
\label{eq:paired-product-definition}
u_{j\ell}=\frac{a_j-a_\ell}{d_{j\ell}},
\qquad
R_j(x)=
\prod_{\ell\ne j}
\frac{
\langle x-a_\ell,u_{j\ell}\rangle
\langle x-b_\ell,u_{j\ell}\rangle
}{
d_{j\ell}^2
}.
\end{equation}
Each factor indexed by $\ell$ introduces a zero at $a_\ell$ and a zero at $b_\ell$. We next verify that this product remains bounded away from zero on the selected segment $[a_j,b_j]$, so that its values there can be inverted.

\begin{lemma}[Conditioning of the paired product]
\label{lem:paired-product}
The polynomial $R_j$ has degree $2K-2$ and satisfies $R_j(a_\ell)=R_j(b_\ell)=0$ for every $\ell\ne j$. On the segment $[a_j,b_j]$,
\begin{equation}
\label{eq:paired-product-local-bounds}
0<c\le R_j(x)\le C,
\qquad
\|\nabla(1/R_j)(x)\|_2
\le C\Delta_{\mathrm{sep}}^{-1}.
\end{equation}
In the structured regimes,
\begin{equation}
\label{eq:paired-product-sup-bound}
\|R_j\|_{\infty,B_d(0,R)}
\le
C\Delta_{\mathrm{sep}}^{-(2s_{\max}-2)}.
\end{equation}
Without cluster assumptions, the same supremum bound holds with $s_{\max}=K$.
\end{lemma}

\begin{proof}
Fix $j\in[K]$. Distinctness of the true centers gives $d_{j\ell}>0$, and the definition of $u_{j\ell}$ gives $\|u_{j\ell}\|_2=1$. Every factor in Equation~\eqref{eq:paired-product-definition} is therefore a well-defined polynomial of degree two. Their product has degree $2K-2$. At the centers belonging to another pair,
\[
\langle a_\ell-a_\ell,u_{j\ell}\rangle=0,
\qquad
\langle b_\ell-b_\ell,u_{j\ell}\rangle=0,
\]
so $R_j(a_\ell)=R_j(b_\ell)=0$. If $a_\ell=b_\ell$, the two linear factors coincide, and the resulting squared factor still vanishes at that common center.

\paragraph{Bounds on the selected segment.}
Fix $x\in[a_j,b_j]$, and write $x=a_j+t(b_j-a_j)$ for some $t\in[0,1]$. Then
\begin{equation}
\begin{aligned}
\|x-a_j\|_2
&=t\|b_j-a_j\|_2\le \frac{\Delta_{\mathrm{sep}}}{4}
\qquad\text{(because $0\le t\le1$)}.
\end{aligned}
\label{eq:bound_x-aj_norm}
\end{equation}
For $\ell\ne j$, denote the two linear factors by
\[
L_{j\ell}(x)=\langle x-a_\ell,u_{j\ell}\rangle,
\qquad
M_{j\ell}(x)=\langle x-b_\ell,u_{j\ell}\rangle.
\]
To bound the first factor, separate the displacement from $a_j$:
\[
\begin{aligned}
L_{j\ell}(x)
&=\langle a_j-a_\ell,u_{j\ell}\rangle
  +\langle x-a_j,u_{j\ell}\rangle\\
&=d_{j\ell}+\langle x-a_j,u_{j\ell}\rangle
\qquad\text{(by the definition of $u_{j\ell}$)}.
\end{aligned}
\]
The second term satisfies
\[
\begin{aligned}
|\langle x-a_j,u_{j\ell}\rangle|
&\le \|x-a_j\|_2\|u_{j\ell}\|_2
\qquad\text{(by Cauchy--Schwarz)}\\
&\le \frac{\Delta_{\mathrm{sep}}}{4} \qquad \text{(because $\|x-a_j\|_2 \le \frac{\Delta_{\mathrm{sep}}}{4}$ by~\eqref{eq:bound_x-aj_norm} and $\|u_{jl}\|_2 = 1$ by its definition)}\\
&\le \frac{d_{j\ell}}{4}
\qquad\text{(because $d_{j\ell}\ge\Delta_{\mathrm{sep}}$)}.
\end{aligned}
\]
Consequently,
\begin{equation}
\label{eq:paired-first-factor-bounds}
\frac34d_{j\ell}
\le L_{j\ell}(x)
\le \frac54d_{j\ell}.
\end{equation}

For the second factor, the displacement of $b_\ell$ from its true center must also be included:
\[
M_{j\ell}(x)
=
d_{j\ell}
+\langle x-a_j,u_{j\ell}\rangle
-\langle b_\ell-a_\ell,u_{j\ell}\rangle.
\]
Its deviation from $d_{j\ell}$ is bounded by
\[
\begin{aligned}
|M_{j\ell}(x)-d_{j\ell}|
&\le
|\langle x-a_j,u_{j\ell}\rangle|
+
|\langle b_\ell-a_\ell,u_{j\ell}\rangle|
\qquad\text{(by the triangle inequality)}\\
&\le
\|x-a_j\|_2+\|b_\ell-a_\ell\|_2
\qquad\text{(by Cauchy--Schwarz)}\\
&\le
\frac{\Delta_{\mathrm{sep}}}{4}
+\frac{\Delta_{\mathrm{sep}}}{4}\\
&\le \frac{d_{j\ell}}{2}.
\end{aligned}
\]
It follows that
\begin{equation}
\label{eq:paired-second-factor-bounds}
\frac12d_{j\ell}
\le M_{j\ell}(x)
\le \frac32d_{j\ell}.
\end{equation}
Both linear factors are positive on the entire selected segment. We may therefore multiply their lower and upper bounds to obtain
\[
\frac38
=
\frac{(3d_{j\ell}/4)(d_{j\ell}/2)}{d_{j\ell}^2}
\le
\frac{L_{j\ell}(x)M_{j\ell}(x)}{d_{j\ell}^2}
\le
\frac{(5d_{j\ell}/4)(3d_{j\ell}/2)}{d_{j\ell}^2}
=
\frac{15}{8}.
\]
There are $K-1$ factors in $R_j$, so multiplying these inequalities over $\ell\ne j$ gives
\begin{equation}
\label{eq:paired-product-segment-bounds}
\left(\frac38\right)^{K-1}
\le R_j(x)
\le
\left(\frac{15}{8}\right)^{K-1},
\qquad x\in[a_j,b_j].
\end{equation}
In particular, $R_j(a_j)$ and $R_j(b_j)$ are strictly positive.

\paragraph{Derivative of the reciprocal.}
The positivity of all linear factors at a segment point persists in a neighborhood of that point. Thus $1/R_j$ is differentiable near the segment. Since $\nabla L_{j\ell}=\nabla M_{j\ell}=u_{j\ell}$, the product rule gives
\[
\nabla R_j(x)
=
R_j(x)
\sum_{\ell\ne j}
\left\{
\frac{u_{j\ell}}{L_{j\ell}(x)}
+
\frac{u_{j\ell}}{M_{j\ell}(x)}
\right\}.
\]
Applying the derivative formula for a reciprocal, we have
\[
\begin{aligned}
\nabla(1/R_j)(x)
&=-\frac{\nabla R_j(x)}{R_j(x)^2}=-\frac1{R_j(x)}
\sum_{\ell\ne j}
\left\{
\frac{u_{j\ell}}{L_{j\ell}(x)}
+
\frac{u_{j\ell}}{M_{j\ell}(x)}
\right\}.
\end{aligned}
\]
Taking norms, we have:
\[
\begin{aligned}
\|\nabla(1/R_j)(x)\|_2
&\le
\frac1{R_j(x)}
\sum_{\ell\ne j}
\left\{
\frac1{L_{j\ell}(x)}
+
\frac1{M_{j\ell}(x)}
\right\}
\qquad\text{(because $\|u_{j\ell}\|_2=1$)}\\
&\le
\frac1{R_j(x)}
\sum_{\ell\ne j}
\left\{
\frac4{3d_{j\ell}}+\frac2{d_{j\ell}}
\right\} \qquad \text{(by~\eqref{eq:paired-first-factor-bounds} and~\eqref{eq:paired-second-factor-bounds}}\\
&=
\frac{10}{3R_j(x)}
\sum_{\ell\ne j}\frac1{d_{j\ell}}\\
&\le
\frac{10(K-1)}{3R_j(x)\Delta_{\mathrm{sep}}}
\qquad\text{(because $d_{j\ell}\ge\Delta_{\mathrm{sep}}$)}\\
&\le
\left(\frac83\right)^{K-1}
\frac{10(K-1)}{3\Delta_{\mathrm{sep}}}
\qquad\text{(by Equation~\eqref{eq:paired-product-segment-bounds})}.
\end{aligned}
\]
This proves the reciprocal derivative bound. Its denominator involves distances between distinct true centers, and does not involve the distance $\|b_j-a_j\|_2$ within the selected pair.

\paragraph{Supremum bound on the full ball.}
For an arbitrary $x\in B_d(0,R)$, the signs of the linear factors need not be positive. Their absolute values nevertheless satisfy
\[
\begin{aligned}
|L_{j\ell}(x)|
&\le \|x-a_\ell\|_2
\qquad\text{(by Cauchy--Schwarz)}\\
&\le \|x\|_2+\|a_\ell\|_2
\le 2R,
\end{aligned}
\]
and the same calculation gives $|M_{j\ell}(x)|\le2R$. Hence, writing $D_j=\prod_{\ell\ne j}d_{j\ell}^2$, we obtain
\[
\begin{aligned}
|R_j(x)|
&=
\prod_{\ell\ne j}
\frac{|L_{j\ell}(x)||M_{j\ell}(x)|}{d_{j\ell}^2}\le
\frac{(2R)^{2K-2}}{D_j}.
\end{aligned}
\]

In the structured regimes, let $\mathcal C_j$ be the cluster containing $j$, and put $s_j=|\mathcal C_j|\le s_{\max}$. Distances within this cluster are at least $\Delta_{\mathrm{sep}}$, while distances to other clusters are at least $D_0$. With $\lambda_0=\min\{1,D_0\}$, the denominator satisfies
\[
\begin{aligned}
D_j
&=
\prod_{\ell\in\mathcal C_j\setminus\{j\}}d_{j\ell}^2
\prod_{\ell\notin\mathcal C_j}d_{j\ell}^2\\
&\ge
\Delta_{\mathrm{sep}}^{2(s_j-1)}
D_0^{2(K-s_j)}
\qquad\text{(by the distance bounds)}\\
&\ge
\Delta_{\mathrm{sep}}^{2(s_j-1)}
\lambda_0^{2(K-s_j)}
\qquad\text{(because $D_0\ge\lambda_0$)}\\
&\ge
\Delta_{\mathrm{sep}}^{2(s_j-1)}
\lambda_0^{2(K-1)}
\qquad\text{(because $0<\lambda_0\le1$)}\\
&\ge
\Delta_{\mathrm{sep}}^{2(s_{\max}-1)}
\lambda_0^{2(K-1)}
\qquad\text{(because $s_j\le s_{\max}$ and $\Delta_{\mathrm{sep}}\le1$)}.
\end{aligned}
\]
Substitution into the bound for $|R_j(x)|$ yields
\[
\|R_j\|_{\infty,B_d(0,R)}
\le
(2R)^{2K-2}\lambda_0^{-2(K-1)}
\Delta_{\mathrm{sep}}^{-(2s_{\max}-2)}
\le
C\Delta_{\mathrm{sep}}^{-(2s_{\max}-2)}.
\]
Without clustering, every factor in $D_j$ is at least $\Delta_{\mathrm{sep}}^2$, so
\[
D_j\ge\Delta_{\mathrm{sep}}^{2K-2}.
\]
The same numerator bound then proves the assertion with $s_{\max}=K$.
\end{proof}

\subsubsection{Exact mass and mean extraction}

The product $R_j$ already removes every pair other than the selected pair $(a_j,b_j)$. To extract its probability, we adjust the remaining factor so that the resulting polynomial equals one at both selected centers. For displacement, we instead require value zero at $a_j$ and value $\|b_j-a_j\|_2$ at $b_j$.

Write $h_j=\|b_j-a_j\|_2$ and $\alpha_j=R_j(a_j)^{-1}$. When $h_j>0$, define
\[
w_j=\frac{b_j-a_j}{h_j},
\qquad
\gamma_j=
\frac{R_j(b_j)^{-1}-R_j(a_j)^{-1}}{h_j}.
\]
When $h_j=0$, choose any unit vector $w_j$ and set
\[
\gamma_j=
\langle\nabla(1/R_j)(a_j),w_j\rangle.
\]
In both cases, define
\begin{equation}
\label{eq:paired-cardinal-definition}
\widetilde E_j(x)
=
R_j(x)
\left\{
\alpha_j+\gamma_j\langle w_j,x-a_j\rangle
\right\},
\end{equation}
and
\begin{equation}
\label{eq:paired-displacement-definition}
\widetilde H_j(x)
=
\begin{cases}
\displaystyle
\frac{R_j(x)}{R_j(b_j)}
\langle w_j,x-a_j\rangle,
& h_j>0,\\[4pt]
0,
& h_j=0.
\end{cases}
\end{equation}

\begin{lemma}[Paired cardinal and displacement polynomials]
\label{lem:paired-cardinal}
The polynomials $\widetilde E_j$ and $\widetilde H_j$ have degree at most $2K-1$ and satisfy
\begin{equation}
\label{eq:paired-cardinal-values}
\widetilde E_j(a_\ell)
=
\widetilde E_j(b_\ell)
=
\mathbf1\{j=\ell\},
\end{equation}
and
\begin{equation}
\label{eq:paired-displacement-values}
\widetilde H_j(a_\ell)=0,
\qquad
\widetilde H_j(b_\ell)=\mathbf1\{j=\ell\}h_j.
\end{equation}
In the structured regimes, their supremum norms on $B_d(0,R)$ satisfy
\begin{equation}
\label{eq:paired-selector-norms}
\|\widetilde E_j\|_\infty
\le C\Delta_{\mathrm{sep}}^{-(2s_{\max}-1)},
\qquad
\|\widetilde H_j\|_\infty
\le C\Delta_{\mathrm{sep}}^{-(2s_{\max}-2)}.
\end{equation}
Without cluster assumptions, these bounds hold with $s_{\max}=K$.
\end{lemma}

\begin{proof}
If $h_j>0$, the fundamental theorem of calculus on $[a_j,b_j]$ gives
\[
\gamma_j
=\int_0^1\big\langle
 \nabla(1/R_j)(a_j+t(b_j-a_j)),w_j\big\rangle\,dt.
\]
Lemma~\ref{lem:paired-product} therefore implies
$|\gamma_j|\le C\Delta_{\mathrm{sep}}^{-1}$, independently of
$h_j$. When $h_j=0$, the same estimate follows from the
directional-derivative definition. The lemma also gives
$|\alpha_j|+R_j(b_j)^{-1}\le C$ and
$\|R_j\|_\infty\le C\Delta_{\mathrm{sep}}^{-(2s_{\max}-2)}$.

For $\ell\ne j$, $R_j$ vanishes at both $a_\ell$ and $b_\ell$,
so both selectors vanish there. At $a_j$,
$\widetilde E_j(a_j)=R_j(a_j)\alpha_j=1$ and
$\widetilde H_j(a_j)=0$. If $h_j>0$, the definition of
$\gamma_j$ gives
\[
\widetilde E_j(b_j)
 =R_j(b_j)(\alpha_j+\gamma_jh_j)=1,
\qquad
\widetilde H_j(b_j)=h_j.
\]
If $h_j=0$, these assertions follow from $a_j=b_j$ and
$\widetilde H_j\equiv0$. This proves the interpolation
identities even for collapsed fitted pairs.

The degrees are at most $2K-1$. Since
$|\langle w_j,x-a_j\rangle|\le2R$ on $B_d(0,R)$,
\[
 \|\widetilde E_j\|_\infty
 \le\|R_j\|_\infty(|\alpha_j|+2R|\gamma_j|)
 \le C\Delta_{\mathrm{sep}}^{-(2s_{\max}-1)},
\]
and, for $h_j>0$,
\[
 \|\widetilde H_j\|_\infty
 \le 2R\|R_j\|_\infty/R_j(b_j)
 \le C\Delta_{\mathrm{sep}}^{-(2s_{\max}-2)}.
\]
The second norm is zero when $h_j=0$. Replacing $s_{\max}$
by $K$ proves the unstructured assertions.
\end{proof}

\section{Proofs under exact specification}
\label{app:exact-specification-proofs}

We first establish auxiliary displacement and probability bounds for any fitted state count $k\ge K$, which will be used in both the exact-specification and over-specification analyses. Under exact specification, direct transport interpolation and the positive initial weights give the global recovery and localization bounds. The counting argument also retains alternatives with explicit initial-probability dependence. The paired bounds are proved independently under an assumed matching. Their proof uses the polynomials in Appendix~\ref{app:paired-interpolation} for one-time displacement and probability errors, followed by lower-degree selectors, full products of center distances, and a Markov path comparison for joint recovery. The global argument invokes these paired bounds only after establishing matching. Throughout this proof, for simplicity, we denote $h_{T} = h(p_{\theta,T},p_{\theta^\star,T})$.

\paragraph{Dependency order.}
The exposition follows the main-text statements, so some ingredients are used before their proofs appear. The global theorem relies on the displacement, counting, moment-to-transport, and diameter-localization lemmas below, and then invokes the paired-center theorem. The paired-center proof uses the polynomial witnesses in Appendix~\ref{app:witness-polynomials}, Gaussian moment extraction, Lemmas~\ref{lem:cluster-mass} and~\ref{lem:direct-transport}, and the cardinal-polynomial and Markov-comparison lemmas in its own subsection. It does not use the global theorem or the center-matching corollary. Thus the forward reference is not circular. The overfitted proofs in Appendix~\ref{app:over-specification-proofs} use the global recovery \emph{lemma}, which allows arbitrary fitted state counts, not the exact-fit theorem.

\subsection{Proof of Theorem~\ref{thm:globalexact}}
\label{app:proof_3_1}
\label{app:global-bounds}
\subsubsection{Auxiliary recovery bounds}

The following lemma bounds the quadratic displacement $V$, the within-group displacement sums defining $B$, and the aggregate state-sequence probability errors $D_I$. The positive polynomial $P$ first controls $V$ and provides a fitted expectation bound that will also control the approximation errors of the other polynomials. We then use $H_{j,u}$ to recover displacement in a chosen direction. For state-sequence probabilities, we construct polynomials that equal one at the selected true-center tuple and zero at every other true-center tuple, with vanishing first derivatives at all such tuples. Their values extract the target probability, while the vanishing derivatives leave approximation errors that can be controlled by the quadratic displacement bound. In the unstructured regime, there is no fixed positive lower bound on distances between clusters; we therefore retain the full products of distances between true centers. This gives all three bounds without requiring a separate identification of fitted states assigned to the same true center.

\begin{lemma}[Global displacement and aggregate probability bounds]
\label{lem:global-recovery-bounds}
Under the standing model assumptions, for any $k\ge K$ and $\theta\in\Theta_k(R)$, in the structured regimes, with $\beta_r=\max\{2s_{\max},r(s_{\max}-1)+s_{\max}\}$,
\[
\begin{aligned}
V&\le C\Delta_{\mathrm{sep}}^{-(2s_{\max}-2)}h_{T},\\
B&\le C\Delta_{\mathrm{sep}}^{-(2s_{\max}-1)}h_{T},\\
D_I&\le C_r\Delta_{\mathrm{sep}}^{-\beta_r}h_{T},
\qquad I\subseteq[T],\quad |I|=r\ge1.
\end{aligned}
\]
Without a cluster assumption, the $V$ and $B$ bounds hold with $s_{\max}=K$, and $D_I$ has exponent $\gamma_{K,r}=\max\{2K,r(K-1)+K\}$. More precisely, writing $D_\star=\max_{i,j}\|a_i-a_j\|_2$, we have
\begin{equation}
D_I\le C_rD_\star^{-2}\Delta_{\mathrm{sep}}^{-(\gamma_{K,r}-2)}h_{T}.
\label{eq:overfit-diameter-probability}
\end{equation} 
No lower bound on fitted state probabilities and no factor $(\pi_{\min}^\star)^{-1}$ is required.
\end{lemma}

\begin{proof}
In this proof, all test functions are
polynomials on the span of the fitted and true centers, so
Lemma~\ref{lem:moment-extraction} and marginalization apply with
constants depending only on the standing fixed parameters.

\emph{Step 1: Quadratic displacement.} We show the bound of $V$ based on $\Delta_{\mathrm{sep}}$ and $h(p_{\theta,T},p_{\theta^*,T})$. 
Consider the test polynomial $P$ from Lemma~\ref{lem:positive-coercivity}
which vanishes on the true support and satisfies
$P(\mu_i)\ge c\Delta_{\mathrm{sep}}^{2s_{\max}-2}r_i^2$.
Therefore, by using Lemma~\ref{lem:positive-coercivity} and Lemma~\ref{lem:moment-extraction}, we have simultaneously for all $t\le T$,
\[
 c\Delta_{\mathrm{sep}}^{2s_{\max}-2}\sum_i q_t(i)r_i^2
 \le\mathbb E_\theta P(\mu_{Z_t})
 =\int P\,d(G_\theta^{\{t\}}-G_\star^{\{t\}})
 \le C h_{T}.
\]
Thus $V\le C\delta^{-(2s_{\max}-2)}h_{T}$; we retain the
stronger bound $\max_t\mathbb E_\theta P(\mu_{Z_t})\le Ch_{T}$
for the approximation errors below.

\emph{Step 2: Within-group first displacement.} We prove the bound of $B$ based on $\Delta_{\mathrm{sep}}$ and $h(p_{\theta,T},p_{\theta^*,T})$.
Fix $j,t$ and let
$b_{j,t}=\sum_{i\in V_j}q_t(i)(\mu_i-a_j)$.
Choose $u=b_{j,t}/\|b_{j,t}\|_2$ when nonzero and $u=0$
otherwise. According to Lemma~\ref{lem:flat-selectors} there exists a polynomial
$H_{j,u}$ such that $H_{j,u}(\mu_i) = 0$ and
\[
\|H_{j,u}\|_\infty\le C\delta^{-(2s_{\max}-2)},\qquad
|H_{j,u}(x)-\mathbf1\{c(x)=j\}\langle u,x-a_j\rangle|
\le C\delta^{-(2s_{\max}-1)}P(x).
\]
Consequently, by moment extraction and the preceding bound for $P$,
\[
\begin{aligned}
\|b_{j,t}\|_2
& = \left|\int\mathbf1\{c(x)=j\}\langle u,x-a_j\rangle| d(G_\theta^{\{t\}}-G_\star^{\{t\}})\right|\\
&\leq \left|\int H_{j,u}(x)\,d(G_\theta^{\{t\}}-G_\star^{\{t\}}) \right| + 
\left|\int\left(\mathbf1\{c(x)=j\}\langle u,x-a_j\rangle-H_{j,u}(x)\right)\,d(G_\theta^{\{t\}}-G_\star^{\{t\}})\right|\\
&\leq \left|\int H_{j,u}(x)\,d(G_\theta^{\{t\}}-G_\star^{\{t\}}) \right|
 +C\delta^{-(2s_{\max}-1)}\mathbb E_\theta P(\mu_{Z_t})\\
&\le C\left(\delta^{-(2s_{\max}-2)}
 +\delta^{-(2s_{\max}-1)}\right)h_{T}
 \le C\delta^{-(2s_{\max}-1)}h_{T}.
 \end{aligned}
\]
Summation over $j$ and maximization over $t$ prove the bound for $B$.

\emph{Step 3: Assigned path probabilities.} We show the bound for $D_I$ based on $\Delta_{\mathrm{sep}}$ and $h(p_{\theta,T},p_{\theta^*,T})$.
For $I=(t_1,\ldots,t_r)$ and $j\in[K]^r$, write
$\mathbf1_j(x)=\mathbf1\{(c(x_1),\ldots,c(x_r))=j\}$.
In the structured case, the polynomial $U_j$ of
Lemma~\ref{lem:joint-first-order} equals $\mathbf1_j$ on the true
grid and obeys, with
$d_r=r(s_{\max}-1)+s_{\max}$ and
$\beta_r=\max\{2s_{\max},d_r\}$,
\[
\|U_j\|_\infty\le C_r\delta^{-d_r},\qquad
|U_j(x)-\mathbf1_j(x)|
\le C_r\delta^{-\beta_r}\sum_{\ell=1}^r P(x_\ell).
\]
Since $\int U_j\,dG_\star^I=P_\star^I(j)$ and
$\int\mathbf1_j\,dG_\theta^I=\bar P_\theta^I(j)$,
we obtain
\[
\begin{aligned}
|\bar P_\theta^I(j)-P_\star^I(j)|
&\le\left|\int U_j\,d(G_\theta^I-G_\star^I)\right|
   +\int|U_j-\mathbf1_j|\,dG_\theta^I\\
&\le C_r\delta^{-d_r}h_{T}
  +C_r\delta^{-\beta_r}\sum_{\ell=1}^r
   \mathbb E_\theta P(\mu_{Z_{t_\ell}})
 \le C_r\delta^{-\beta_r}h_{T}.
\end{aligned}
\]
Summing over the $K^r$ label tuples proves the structured bound.
There is no independence assumption between selected hidden states.

Without cluster assumptions, Steps~1 and~2 apply with
$s_{\max}=K$, giving the stated bounds for $V$ and $B$.
For $D_I$, we retain the diameter dependence provided by
Lemma~\ref{lem:diameter-joint-selectors}. Write
\[
\gamma_{K,r}=\max\{2K,r(K-1)+K\},
\qquad
F_{K,r}=D_\star^{-2}
\delta^{-(\gamma_{K,r}-2)}.
\]
For each $j\in[K]^r$, that lemma supplies a polynomial
$U_j$ equal to $\mathbf1_j$ on the true grid and satisfying
\[
\|U_j\|_\infty\le C_rF_{K,r},
\qquad
|U_j(x)-\mathbf1_j(x)|
\le C_rF_{K,r}\sum_{\ell=1}^rP(x_\ell).
\]
The cardinal values give
$\int U_j\,dG_\star^I=P_\star^I(j)$. Hence moment extraction
and the expectation bound from Step~1 yield
\[
\begin{aligned}
|\bar P_\theta^I(j)-P_\star^I(j)|
&\le
\left|\int U_j\,d(G_\theta^I-G_\star^I)\right|
+\int|\mathbf1_j-U_j|\,dG_\theta^I\\
&\le
C_rF_{K,r}h_{T}
+C_rF_{K,r}\sum_{\ell=1}^r
\mathbb E_\theta P(\mu_{Z_{t_\ell}})\\
&\le
C_rF_{K,r}h_{T},
\end{aligned}
\]
where the fixed factor $r$ is absorbed into $C_r$.
Summing over $j\in[K]^r$ and absorbing $K^r$ into the
constant proves
\[
D_I\le
C_rD_\star^{-2}
\delta^{-(\gamma_{K,r}-2)}
h_{T},
\]
which is~\eqref{eq:overfit-diameter-probability}.
Since $D_\star\ge\delta$, this also implies
$D_I\le C_r\delta^{-\gamma_{K,r}}
h_{T}$.
All fitted probabilities enter only as nonnegative weights;
thus zero fitted probabilities and coincident fitted
centers are allowed.
\end{proof}

We next compare probabilities of cluster sequences. For these probabilities, Lemma~\ref{lem:cluster-polynomials} provides polynomials $F_m$ whose supremum norms are bounded independently of emission separation. Their approximation errors are controlled by squared distances to the assigned true centers. Products of these polynomials therefore give an observation-density contribution with no inverse separation factor and an approximation contribution controlled by $V$.

\begin{lemma}[Cluster-sequence mass error]
\label{lem:cluster-mass}
For every structured configuration and every $I\subseteq[T]$ with $|I|=r$,
\begin{equation}
J_I
\le C_r\bigl(h_{T}+V\bigr)
\le C_r\Delta_{\mathrm{sep}}^{-(2s_{\max}-2)}
h_{T}.
\label{C.2}
\end{equation}
This holds for any fitted number of states $k$.
\end{lemma}

\begin{proof}
If there is only one true cluster, both cluster-sequence laws are
point masses at the same label sequence, and $J_I=0$, equation~\eqref{C.2} automatically holds. Otherwise, for $m=(m_1,\ldots,m_r)$ define
\[
U_m(x)=\prod_{\ell=1}^rF_{m_\ell}(x_\ell),\qquad
B_m(x)=\prod_{\ell=1}^r
\mathbf1\{c(x_\ell)\in C_{m_\ell}\},
\]
where the cluster polynomials $F_a$ are those of
Lemma~\ref{lem:cluster-polynomials}. Their interpolation values
on true centers show
$\int U_m\,dG_\star^I=P_\star^{I,\mathrm{cl}}(m)$,
while $\int B_m\,dG_\theta^I=\bar P_\theta^{I,\mathrm{cl}}(m)$.
The uniform bounds $\|F_a\|_\infty\le C$ yield
$\|U_m\|_\infty\le C_r$, so moment extraction gives
\[
\left|\int U_m\,d(G_\theta^I-G_\star^I)\right|\le C_r h_{T}.
\]
The product telescoping identity and the cluster-polynomial error
$|F_a(x)-\mathbf1\{c(x)\in C_a\}|
\le C\|x-a_{c(x)}\|_2^2$ imply
\[
|U_m(x)-B_m(x)|\le C_r\sum_{\ell=1}^r
\|x_\ell-a_{c(x_\ell)}\|_2^2.
\]
Integrating this bound against $G_\theta^I$ gives
\[
\begin{aligned}
\int|U_m-B_m|\,dG_\theta^I
\le
C_r\sum_{\ell=1}^r
\mathbb E_\theta
\|\mu_{Z_{t_\ell}}-a_{c(\mu_{Z_{t_\ell}})}\|_2^2 =
C_r\sum_{\ell=1}^r V_{t_\ell}
\le C_r V,
\end{aligned}
\]
where $V=\max_{t\le T}V_t$ and the factor $r$ is absorbed
into $C_r$. Using the integral representations of the
cluster-sequence probabilities, we obtain
\[
\begin{aligned}
\left|
\bar P_\theta^{I,\mathrm{cl}}(m)
-P_\star^{I,\mathrm{cl}}(m)
\right|&=
\left|
\int B_m\,dG_\theta^I-\int U_m\,dG_\star^I
\right|\le
\left|\int U_m\,d(G_\theta^I-G_\star^I)\right|
+\int|B_m-U_m|\,dG_\theta^I\\
&\le C_r(h_{T}+V).
\end{aligned}
\]
There are at most $K^r$ cluster-label sequences. Summing
over them and absorbing this factor into $C_r$ yields
\[
J_I
=
\sum_m
\left|
\bar P_\theta^{I,\mathrm{cl}}(m)
-P_\star^{I,\mathrm{cl}}(m)
\right|
\le C_r(h_{T}+V).
\]
Finally, Lemma~\ref{lem:global-recovery-bounds} gives
\[
\begin{aligned}
J_I
&\le
C_r\left(
1+\Delta_{\mathrm{sep}}^{-(2s_{\max}-2)}
\right)h_{T}\\
&\le
C_r\Delta_{\mathrm{sep}}^{-(2s_{\max}-2)}h_{T},
\end{aligned}
\]
where the last inequality uses
$0<\Delta_{\mathrm{sep}}\le R$ and $s_{\max}\ge1$.
This proves~\eqref{C.2}.
\end{proof}

\paragraph{Controlling individual displacements under exact specification.}

To obtain a bound for $D_\mu$ that is linear in the observation-density error, we use the bounds for quadratic displacement and cell displacement vectors established in Lemma~\ref{lem:global-recovery-bounds}. A cell displacement vector controls a single nearby fitted center, but several displacement vectors in the same cell may cancel. Under exact specification, a counting argument controls the contribution of cells containing several nearby fitted centers.

\begin{lemma}[Counting near centers]
\label{lem:counting-near}
Suppose $k=K$, and define
\[
N_j=\{i:c(i)=j,\ r_i\le\Delta_{\mathrm{sep}}/4\},
\qquad
F=\{i:r_i>\Delta_{\mathrm{sep}}/4\}.
\]
Then
\begin{equation}
\label{eq:counting-near-mass}
\max_{t\le T}
\sum_{j=1}^K
\left|
\sum_{i\in N_j}q_t(i)-q_t^\star(j)
\right|
\le
C\Delta_{\mathrm{sep}}^{-2s_{\max}}
h_{T}.
\end{equation}
If some $N_j$ is empty, then
\begin{equation}
\label{eq:counting-empty-near}
\pi_{\min}^\star
\le
C\Delta_{\mathrm{sep}}^{-2s_{\max}}
h_{T}.
\end{equation}
Moreover,
\begin{equation}
\label{eq:counting-global-displacement}
D_\mu
\le
\frac{C}{\pi_{\min}^\star}
\Delta_{\mathrm{sep}}^{-(2s_{\max}-1)}
h_{T}.
\end{equation}
\end{lemma}

\begin{proof}
Write $\delta=\Delta_{\mathrm{sep}}$ and $s=s_{\max}$.
For each $j\in[K]$ and $t\le T$, define
\[
\bar q_t(j)=\sum_{i\in V_j}q_t(i),
\qquad
b_{j,t}=\sum_{i\in V_j}q_t(i)(\mu_i-a_j),
\qquad
V_t=\sum_{i=1}^Kq_t(i)r_i^2.
\]
Lemma~\ref{lem:global-recovery-bounds} gives, uniformly in $t$,
\begin{equation}
\begin{aligned}
V_t\le C\delta^{-(2s-2)}h_T, \quad \sum_{j=1}^K\|b_{j,t}\|_2
\le C\delta^{-(2s-1)}h_T,\quad \sum_{j=1}^K|\bar q_t(j)-q_t^\star(j)|
\le C\delta^{-2s}h_T.
\end{aligned}
\label{eq:counting-input-bounds}
\end{equation}

\medskip
\noindent\textit{Step 1:}
We bound the total probability and displacement of the far
labels. Since $r_i>\delta/4$ for $i\in F$, we have $1\le\frac{16r_i^2}{\delta^2}$, and $r_i\le\frac{4r_i^2}{\delta}$. Multiplying by $q_t(i)$ and summing over $F$ yields
\begin{equation}
\begin{aligned}
\sum_{i\in F}q_t(i)
\le16\delta^{-2}V_t
\le C\delta^{-2s}h_T, \quad \sum_{i\in F}q_t(i)r_i
\le4\delta^{-1}V_t
\le C\delta^{-(2s-1)}h_T.
\end{aligned}
\label{eq:counting-far-displacement}
\end{equation}

\medskip
\noindent\textit{Step 2:}
We compare the near-set masses with the true state
probabilities and deduce the bound for an empty near set.
Write
\[
n_{j,t}=\sum_{i\in N_j}q_t(i),
\qquad
f_{j,t}=\sum_{i\in F\cap V_j}q_t(i).
\]
Because $V_j=N_j\mathbin{\dot\cup}(F\cap V_j)$,
we have $\bar q_t(j)=n_{j,t}+f_{j,t}$. Therefore,
\[
\begin{aligned}
\sum_{j=1}^K|n_{j,t}-q_t^\star(j)| &= \sum_{j=1}^K |\bar q_t(j)-q_t^\star(j)-f_{j,t}| \le
\sum_{j=1}^K|\bar q_t(j)-q_t^\star(j)|
+\sum_{j=1}^Kf_{j,t}\\
&=
D_{\{t\}}+\sum_{i\in F}q_t(i) \le C\delta^{-2s}h_T.
\end{aligned}
\]
Taking the maximum over $t$ proves
\eqref{eq:counting-near-mass}.

If $N_{j_0}=\varnothing$, then $n_{j_0,1}=0$ and
$q_1^\star(j_0)=\pi_{j_0}^\star$. Consequently, $\pi_{\min}^\star
\le\pi_{j_0}^\star
=
|n_{j_0,1}-q_1^\star(j_0)|
\le C\delta^{-2s}h_T$, which proves~\eqref{eq:counting-empty-near}.

\medskip
\noindent\textit{Step 3:}
We control the displacement contributed by singleton near
sets. Define
\[
\mathcal S=\{j:|N_j|=1\},
\qquad
\mathcal M=\{j:|N_j|\ge2\}.
\]
For $j\in\mathcal S$, write $N_j=\{i_j\}$. Then
\[
b_{j,t}
=
q_t(i_j)(\mu_{i_j}-a_j)
+
\sum_{i\in F\cap V_j}q_t(i)(\mu_i-a_j).
\]
Since $c(i_j)=j$, the triangle inequality gives $q_t(i_j)r_{i_j}
\le
\|b_{j,t}\|_2
+
\sum_{i\in F\cap V_j}q_t(i)r_i$.
Summing over $\mathcal S$ and using
\eqref{eq:counting-input-bounds} and
\eqref{eq:counting-far-displacement}, we obtain
\[
\begin{aligned}
\sum_{j\in\mathcal S}\sum_{i\in N_j}q_t(i)r_i
&\le
\sum_{j=1}^K\|b_{j,t}\|_2
+\sum_{i\in F}q_t(i)r_i \le C\delta^{-(2s-1)}h_T.
\end{aligned}
\]

\medskip
\noindent\textit{Step 4:}
We control the remaining near displacements using $k=K$
and combine the contributions.
If $\mathcal M\ne\varnothing$, at least one near set
must be empty. Indeed, if every $N_j$ were nonempty, then
their disjointness would imply $K=k
\ge\sum_{j=1}^K|N_j|
\ge (K-1)+2
=K+1$, a contradiction. Step~2 therefore gives
\[
\pi_{\min}^\star\le C\delta^{-2s}h_T.
\]
Since $r_i\le\delta/4$ on every near set and
$\sum_iq_t(i)=1$,
\[
\begin{aligned}
\sum_{j\in\mathcal M}\sum_{i\in N_j}q_t(i)r_i
\le
\frac{\delta}{4}
\sum_{j\in\mathcal M}\sum_{i\in N_j}q_t(i)\le\frac{\delta}{4}
\le
\frac{C}{\pi_{\min}^\star}
\delta^{-(2s-1)}h_T.
\end{aligned}
\]
If $\mathcal M=\varnothing$, this contribution is zero,
so the same bound holds.

As the sets $F,N_1,\ldots,N_K$ partition the fitted labels, we
combine Steps~1, 3, and~4, and using
$0<\pi_{\min}^\star\le1$, and yield
\[
\begin{aligned}
\sum_{i=1}^Kq_t(i)r_i
=
\sum_{i\in F}q_t(i)r_i
+\sum_{j\in\mathcal S}\sum_{i\in N_j}q_t(i)r_i
+\sum_{j\in\mathcal M}\sum_{i\in N_j}q_t(i)r_i\le
\frac{C}{\pi_{\min}^\star}
\delta^{-(2s-1)}h_T.
\end{aligned}
\]
Taking the maximum over $t\le T$ proves
\eqref{eq:counting-global-displacement}.
\end{proof}

\subsubsection{Direct transport control and global recovery}
\label{app:initial-transport-improved}

We next interpolate a Lipschitz Wasserstein test function directly. This gives the global bounds with sharper separation powers, and also the unstructured one-time bound under an assumed matching.

\begin{lemma}[Moment and density comparisons for finite-mixture transport]
\label{lem:direct-transport}
Let $H=\sum_{j=1}^K v_j\delta_{b_j}$ and $H^\star=\sum_{j=1}^K w_j\delta_{a_j}$ be probability measures on $B_d(0,R)$ with nonnegative weights. Assume $K\ge2$, the $a_j$ are distinct, and put $\Delta=1\wedge\min_{i\ne j}\|a_i-a_j\|_2$. Define the polynomial-moment discrepancy
\begin{equation}
\mathfrak M_D(H,H^\star)
:=\sup_{\substack{\deg Q\le D\\\|Q\|_{\infty,B_d(0,R)}\le1}}
       \left|\int Q\,d(H-H^\star)\right|.
\label{eq:finite-moment-discrepancy}
\end{equation}
If $\|b_j-a_j\|_2\le\Delta/4$ for every $j$, then
\begin{equation}
\label{eqn:dung_lemma_D_4_first_comparison}
    W_1(H,H^\star)\le C(K,R)\Delta^{-(2K-2)}
                  \mathfrak M_{2K-1}(H,H^\star).
\end{equation}
Alternatively, if $w_j\ge p_0>0$ for every $j$, the same inequality holds without matching, with constant $C(K,R,p_0)$. The fitted weights may vanish and the fitted centers may coincide in this latter assertion. In particular, for the identity-covariance Gaussian mixture densities,
\begin{align}
W_1(H,H^\star)&\le C(K,R)\Delta^{-(2K-2)}h(p_H,p_{H^\star})
&&\text{under matching},\label{eq:direct-transport-paired}\\
W_1(H,H^\star)&\le C(K,R,p_0)\Delta^{-(2K-2)}h(p_H,p_{H^\star})
&&\text{if }\min_jw_j\ge p_0.\label{eq:direct-transport-global}
\end{align}
For $K=1$, no separation is needed: $W_1(H,H^\star)\le R\mathfrak M_1(H,H^\star)$.
\end{lemma}
\begin{proof}
Write $M=\mathfrak M_{2K-1}(H,H^\star).$ We divide our proof into following steps. 

\medskip
\noindent\textit{Step 1:}
We reduce the moment comparison to a one-dimensional
interpolation problem.
By Lemma~\ref{lem:separation-projection}, there exist a
unit vector $u$ and a constant $c_K\in(0,1]$ such that
\[
|\langle u,z-z'\rangle|
\ge c_K\|z-z'\|_2
\]
for all distinct points in
$\{a_1,\ldots,a_K,b_1,\ldots,b_K\}$.
Define $\bar H=(\langle u,\cdot\rangle)_\#H$, $\bar H^\star=(\langle u,\cdot\rangle)_\#H^\star$, and
$\nu=\bar H-\bar H^\star.$ The projection $z\mapsto\langle u,z\rangle$ onto the direction
$u$ is injective on
$\operatorname{supp}(H)\cup\operatorname{supp}(H^\star)$ as distinct points in this union have distinct projected values.
Lifting any projected coupling therefore gives $W_1(H,H^\star)
\le c_K^{-1}W_1(\bar H,\bar H^\star).$ 

Let $W=W_1(\bar H,\bar H^\star)$. If $W=0$, the assertion \eqref{eqn:dung_lemma_D_4_first_comparison} is immediate. Otherwise, Kantorovich--Rubinstein duality provides a one-Lipschitz function $\psi$ on $[-R,R]$ such that 
\[
\psi(0)=0,
\qquad
\|\psi\|_{\infty,[-R,R]}\le R,
\qquad
\int\psi\,d\nu=W.
\]
As constant extension beyond the endpoints preserves the
one-Lipschitz property, we can extend $\psi$ to $\mathbb R$
by setting $\psi(x)=\psi(-R)$ for $x<-R$ and
$\psi(x)=\psi(R)$ for $x>R$. 

Having bounded the original Wasserstein distance by its
one-dimensional counterpart, we now show that polynomial
moments of the projected measures are controlled by the
original discrepancy $M = \mathfrak M_{2K-1}(H,H^\star)$. For every univariate polynomial $P$ of degree at most $2K-1$,
the polynomial $Q(x)=P(\langle u,x\rangle)$ satisfies
\[
\|Q\|_{\infty,B_d(0,R)}
\le\|P\|_{\infty,[-R,R]},
\qquad
\int Q\,d(H-H^\star)=\int P\,d\nu.
\]
Consequently,
\begin{equation}
\left|\int P\,d\nu\right|
\le \|P\|_{\infty,[-R,R]}M.
\label{eq:transport-projected-moment}
\end{equation}

\medskip
\noindent\textit{Step 2:}
We bound the interpolating polynomial when every three
consecutive projected nodes have a fixed minimum span.
Let $z_1<\cdots<z_N$, $N\le2K$, be the distinct projected
union nodes, and suppose
\[
z_{i+2}-z_i\ge\gamma>0,
\qquad i=1,\ldots,N-2.
\]
Put $\bar\gamma=1\wedge\gamma$. The Lipschitz property gives
$|\psi[z_i,z_{i+1}]|\le1$. For $q\ge2$, the recursion
\[
\psi[z_i,\ldots,z_{i+q}]
=
\frac{
\psi[z_{i+1},\ldots,z_{i+q}]
-\psi[z_i,\ldots,z_{i+q-1}]
}{
z_{i+q}-z_i
}
\]
therefore implies $|\psi[z_i,\ldots,z_{i+q}]|
\le 2^{q-1}\bar\gamma^{-(q-1)}$. The Newton interpolant
\[
P(z):=\psi(z_1)+
\sum_{q=1}^{N-1}
\psi[z_1,\ldots,z_{q+1}]
\prod_{\ell=1}^q(z-z_\ell)
\]
has degree at most $2K-1$ and satisfies
\[
\begin{aligned}
\|P\|_{\infty,[-R,R]}
\le
R+\sum_{q=1}^{N-1}
2^{q-1}\bar\gamma^{-(q-1)}(2R)^q\le C(K,R)\bar\gamma^{-(2K-2)}.
\end{aligned}
\]
The same bound holds for $N\le2$, without a three-node
condition. Since $P=\psi$ on the support of $\nu$,
\eqref{eq:transport-projected-moment} yields
\begin{equation}
W\le C(K,R)(1\wedge\gamma)^{-(2K-2)}M.
\label{eq:transport-triple-bound}
\end{equation}

\medskip
\noindent\textit{Step 3:}
We verify the three-node condition under the matching
hypothesis. For $i\ne j$, $z\in\{a_i,b_i\}$, and
$z'\in\{a_j,b_j\}$,
\[
\begin{aligned}
\|z-z'\|_2
\ge
\|a_i-a_j\|_2-\|z-a_i\|_2-\|z'-a_j\|_2\geq  
\frac{\Delta}{2}.
\end{aligned}
\]
Thus points from different pairs have projected distance
at least $c_K\Delta/2$.
Let $z_1<\cdots<z_N$ be the distinct projected union nodes.
Any three consecutive nodes include two from different
pairs, since each pair contains at most two distinct points.
Their distance is bounded above by the span of the triple,
so
\[
z_{\ell+2}-z_\ell\ge c_K\Delta/2,
\qquad \ell=1,\ldots,N-2.
\]
Applying~\eqref{eq:transport-triple-bound} with
$\gamma=c_K\Delta/2$ and then Step~1 gives
\[
W_1(H,H^\star)
\le c_K^{-1}W
\le C(K,R)\Delta^{-(2K-2)}
\mathfrak M_{2K-1}(H,H^\star),
\]
where factors depending only on $K$ are absorbed into $C(K,R)$.
This proves~\eqref{eqn:dung_lemma_D_4_first_comparison}.

\medskip
\noindent\textit{Step 4:}
We establish the moment bound for arbitrary fitted centers
when $\min_jw_j\ge p_0$.
Define the projected true and fitted centers by $\alpha_j=\langle u,a_j\rangle$ and $\beta_i=\langle u,b_i\rangle$, so that
\[
\bar H=\sum_{i=1}^K v_i\delta_{\beta_i},
\qquad
\bar H^\star=\sum_{j=1}^K w_j\delta_{\alpha_j},
\qquad
\delta_u=\min_{i\ne j}|\alpha_i-\alpha_j|
\ge c_K\Delta.
\]
We distinguish two cases. If $W<p_0\delta_u/8$, every interval $I_j=(\alpha_j-\delta_u/4,\alpha_j+\delta_u/4)$ contains a projected fitted center $\beta_i$ with $v_i>0$.
Otherwise, every coupling must transport the true mass
$w_j$ at $\alpha_j$ by distance at least $\delta_u/4$, giving $W\ge w_j\delta_u/4\ge p_0\delta_u/4$, a contradiction. The $K$ intervals are disjoint, and there are only $K$
fitted labels. Thus each interval contains exactly one
fitted label $\beta_j$. After relabeling, we can suppose
$|\beta_j-\alpha_j|<\delta_u/4$ for every $j$.
For points $\xi\in\{\alpha_i,\beta_i\}$ and
$\xi'\in\{\alpha_j,\beta_j\}$ with $i\ne j$,
\[
|\xi-\xi'|
\ge|\alpha_i-\alpha_j|
-|\xi-\alpha_i|-|\xi'-\alpha_j|
\ge\delta_u/2.
\]
As in Step~3, every three consecutive distinct union nodes
therefore span at least $\delta_u/2$.
Step~2 and the fact that $\delta_u\ge c_K\Delta$ yield
\[
W\le C(K,R)\Delta^{-(2K-2)}M.
\]

If $W\ge p_0\delta_u/8$, we smooth the Lipschitz dual
function before interpolating, thereby avoiding any
spacing requirement on the fitted centers.
Choose a nonnegative $C^\infty$ density $\zeta$ supported
on $[-1,1]$, and set
\[
\eta=\frac{p_0c_K\Delta}{64},
\qquad
\zeta_\eta(x)=\eta^{-1}\zeta(x/\eta),
\qquad
\psi_\eta=\psi*\zeta_\eta.
\]
Since $\psi$ is one-Lipschitz,
\[
|\psi_\eta(x)-\psi(x)|
\le\int|\psi(x-y)-\psi(x)|\zeta_\eta(y)\,dy
\le\eta.
\]
Moreover, $W\ge p_0c_K\Delta/8=8\eta$ and
$\|\nu\|_{\mathrm{TV}}\le2$, so
\[
\int\psi_\eta\,d\nu
\ge\int\psi\,d\nu-2\eta
=W-2\eta
\ge\frac34W.
\]
The weak derivative of $\psi$ satisfies
$\|\psi'\|_\infty\le1$. Hence
\[
\|\psi_\eta\|_{\infty,[-R,R]}\le R+1,
\qquad
\|\psi_\eta^{(q)}\|_\infty
\le
\|\psi'\|_\infty\|\zeta_\eta^{(q-1)}\|_1
\le C_q\eta^{-(q-1)},
\quad q\ge1.
\]
Lemma~\ref{lem:confluent-interpolation} therefore supplies
a polynomial $P_\eta$ interpolating $\psi_\eta$ at the
distinct projected union nodes, with
\[
\deg P_\eta\le2K-1,
\qquad
\|P_\eta\|_{\infty,[-R,R]}
\le C(K,R,p_0)\Delta^{-(2K-2)}.
\]
Because $\nu$ is supported on these nodes,
\[
\frac34W
\le\int\psi_\eta\,d\nu
=\int P_\eta\,d\nu
\le C(K,R,p_0)\Delta^{-(2K-2)}M.
\]
Combining both cases with Step~1 proves
\[
W_1(H,H^\star)
\le C(K,R,p_0)\Delta^{-(2K-2)}
\mathfrak M_{2K-1}(H,H^\star).
\]

\medskip
\noindent\textit{Step 5:}
We derive the density comparisons before treating $K=1$.
By restricting each test polynomial to the span of the true and
fitted centers, whose dimension is at most $2K$, we use Lemma~\ref{lem:moment-extraction} to achieve
\[
\left|\int Q\,d(H-H^\star)\right|
\le C(K,R)h(p_H,p_{H^\star})
\]
whenever $\deg Q\le2K-1$ and
$\|Q\|_{\infty,B_d(0,R)}\le1$.
Taking the supremum over $Q$ yields $M\le C(K,R)h(p_H,p_{H^\star}).$
Substitution into the moment bounds \eqref{eqn:dung_lemma_D_4_first_comparison} proves
\eqref{eq:direct-transport-paired}
and~\eqref{eq:direct-transport-global}.

If $K=1$, both measures are point masses.
For $R>0$ and $a_1\ne b_1$, take
\[
Q(x)=
\frac{\langle b_1-a_1,x\rangle}
{R\|b_1-a_1\|_2}.
\]
Then $\|Q\|_{\infty,B_d(0,R)}\le1$ and
\[
\mathfrak M_1(H,H^\star)
\ge
\left|\int Q\,d(H-H^\star)\right|
=
\frac{\|b_1-a_1\|_2}{R}.
\]
Thus $W_1(H,H^\star)\le R\mathfrak M_1(H,H^\star)$.
If the centers coincide or $R=0$, both sides vanish.
\end{proof}

\begin{lemma}[Cluster-adaptive global mixture transport]
\label{lem:structured-global-transport}
Under the assumptions of the global assertion of Lemma~\ref{lem:direct-transport}, suppose additionally that the true centers have the structured geometry. Then
\begin{equation}
W_1(H,H^\star)\le C(p_0)\Delta_{\mathrm{sep}}^{-(2s_{\max}-2)}h(p_H,p_{H^\star}).
\label{eq:structured-global-mixture}
\end{equation}
The constant may depend on $K,R$ and the cluster parameters, as well as on $p_0$, but not on the fitted weights or the ambient dimension.
\end{lemma}

\begin{proof}
Write $\Delta=\Delta_{\mathrm{sep}}$, $s=s_{\max}$, and
$h=h(p_H,p_{H^\star})$. All constants below may depend
on $K,R$ and the fixed cluster parameters.

\medskip
\noindent\textit{Step 1:}
We construct disjoint groups containing the true clusters
whose boundaries stay uniformly separated from both supports.
If there is only one cluster, then $s=K$ and
Lemma~\ref{lem:direct-transport} proves the assertion.
That lemma also proves the assertion when $\Delta$ is
bounded below by a fixed positive constant. We may therefore
assume $C_0\Delta<D_0/16.$ 

Let $C_1,\ldots,C_M$ be the true clusters, and choose $a_{j_m}\in C_m$ for each $m$. Partition $[D_0/8,D_0/4]$ into $K+1$ equal subintervals.
For fixed $a_{j_m}$, since there are only $K$ fitted labels, one subinterval
contains no distance among $\{\|b_i-a_{j_m}\|_2,\,1\leq i\leq K\}$ in its interior.
Let $\rho_m$ be its midpoint and define $\mathcal A_m=B_d(a_{j_m},\rho_m).$ The subinterval length is $D_0/[8(K+1)]$, so
\[
\left|\|b_i-a_{j_m}\|_2-\rho_m\right|
\ge\frac{D_0}{16(K+1)}
\qquad\text{for every }i.
\]
Recall from assumption in Equation \eqref{eqn:dung_cluster_assumption} that the true clusters have diameter at most
$C_0\Delta<D_0/16$ and are separated by at least $D_0$.
Since $\rho_m\in[D_0/8,D_0/4]$, it follows that
$a_j\in\mathcal A_m$ if and only if $j\in C_m$, and
the balls $\mathcal A_m$ are pairwise disjoint. For the set of distinct union centers $\mathcal Z=\{a_1,\ldots,a_K,b_1,\ldots,b_K\}$, there is consequently a fixed $g>0$ such that $\|z-z'\|_2\ge g$ whenever $z\in\mathcal Z\cap\mathcal A_m$ and
$z'\in\mathcal Z\setminus\mathcal A_m$.

\medskip
\noindent\textit{Step 2:}
We construct bounded polynomial selectors and compare
the fitted and true masses of each group.
Choose a projection direction $u$ from
Lemma~\ref{lem:separation-projection}, and put $E_m=\{\langle u,z\rangle:z\in\mathcal Z\cap\mathcal A_m\}.$ Projected nodes in $E_m$ are separated from all other projected union nodes by at least $g'=c_Kg$.
Let
\[
O_m=\{t\in\mathbb R:\operatorname{dist}(t,E_m)<g'/3\},
\qquad
f_m=\mathbf1_{O_m}*\zeta_{g'/6},
\]
where $\zeta_\varepsilon$ is a fixed smooth probability
density supported on $[-\varepsilon,\varepsilon]$.
Then $f_m$ equals one on $E_m$ and zero on every other
projected union node, and $\|f_m^{(q)}\|_\infty\le C_q$, where $q\ge0$. Interpolating $f_m$ on the distinct projected union nodes
using Lemma~\ref{lem:confluent-interpolation}, and lifting
the interpolant back to $\mathbb R^d$, gives a polynomial
$F_m$ with
\[
\deg F_m\le2K-1,\qquad
\|F_m\|_{\infty,B_d(0,R)}\le C,
\qquad
F_m(z)=\mathbf1\{z\in\mathcal A_m\}
\quad(z\in\mathcal Z).
\]

Define $v_m=\sum_{i:b_i\in\mathcal A_m}v_i$, and $w_m=\sum_{j\in C_m}w_j$. Moment extraction on the span of the union centers,
whose dimension is at most $2K$, yields
\begin{equation}
|v_m-w_m|
=
\left|\int F_m\,d(H-H^\star)\right|
\le Ch.
\label{eq:macro-mass-comparison}
\end{equation}

\medskip
\noindent\textit{Step 3:}
We show that either the desired bound follows immediately
or each group contains exactly as many fitted labels as
true labels.
Write $k_m=\#\{i:b_i\in\mathcal A_m\}$, counting repeated fitted centers and zero-weight labels, and $s_m=|C_m|$. If $k_m<s_m$ for some $m$, we define
\[
Q_m(x)=F_m(x)
\prod_{i:b_i\in\mathcal A_m}\|x-b_i\|_2^2.
\]
This polynomial has uniformly bounded degree and norm,
and $Q_m(b_i)=0$ for every fitted label. The $s_m$ balls $\{B_d(a_j,\Delta/3),\,j\in C_m\}$ are disjoint. Since $k_m<s_m$, some $a_{j_0}$ has
distance at least $\Delta/3$ from every fitted center
in $\mathcal A_m$. Therefore,
\[
\begin{aligned}
\int Q_m\,dH^\star
&=
\sum_{j\in C_m}w_j
\prod_{i:b_i\in\mathcal A_m}\|a_j-b_i\|_2^2 \ge p_0(\Delta/3)^{2k_m}
\ge c p_0\Delta^{2s-2}.
\end{aligned}
\]
Here the last inequality uses $k_m\le s_m-1\le s-1$
and $0<\Delta\le1$.
Since $\int Q_m\,dH=0$, moment extraction implies $h\ge c p_0\Delta^{2s-2}.$ In that case, the trivial transport bound gives
\[
W_1(H,H^\star)\le2R
\le C(p_0)\Delta^{-(2s-2)}h.
\]

Thus, we can assume that $h<c_0p_0\Delta^{2s-2}$ for a sufficiently small fixed $c_0>0$. In this case, $k_m\ge s_m$ for every $m$. Since
\[
K\ge\sum_m k_m\ge\sum_m s_m=K,
\]
we have $k_m=s_m$ for all $m$, and no fitted labels
lie outside the groups.
By~\eqref{eq:macro-mass-comparison}, choosing $c_0$
smaller if necessary also gives
\[
v_m\ge w_m-Ch\ge w_m/2\ge p_0/2.
\]

\medskip
\noindent\textit{Step 4:}
We bound transport between the normalized restrictions
within each group.
Define
\[
H_m=\frac{H|_{\mathcal A_m}}{v_m},
\qquad
H_m^\star=\frac{H^\star|_{\mathcal A_m}}{w_m}.
\]
For every polynomial $Q$ of degree at most $2s_m-1$,
the exact selector values imply
\[
\int Q\,d(H_m-H_m^\star)
=
\frac1{v_m}\int F_mQ\,d(H-H^\star)
+
\frac{w_m-v_m}{v_m}\int Q\,dH_m^\star.
\]
Since $\deg(F_mQ)\le2K+2s_m-2$ and $\|F_mQ\|_\infty\le C\|Q\|_\infty$, moment extraction and~\eqref{eq:macro-mass-comparison}
yield
\[
\left|\int Q\,d(H_m-H_m^\star)\right|
\le C(p_0)\|Q\|_\infty h.
\]
Consequently,
\[
\mathfrak M_{2s_m-1}(H_m,H_m^\star)\le C(p_0)h.
\]
Each normalized true weight is at least $p_0$, since
$w_j/w_m\ge w_j\ge p_0$.
For $s_m\ge2$, the global moment bound in
Lemma~\ref{lem:direct-transport} therefore gives
\[
W_1(H_m,H_m^\star)
\le C(p_0)\Delta^{-(2s_m-2)}h
\le C(p_0)\Delta^{-(2s-2)}h.
\]
For $s_m=1$, the one-component assertion of that lemma
gives the same conclusion.

\medskip
\noindent\textit{Step 5:}
We combine the within-group transport bounds and account
for the group-mass discrepancies.
For each $m$, multiply an optimal coupling of
$H_m,H_m^\star$ by $\min\{v_m,w_m\}$.
The residual marginals have equal total mass
\[
\tau=\frac12\sum_m|v_m-w_m|
\]
and can be coupled at cost at most $2R\tau$.
Thus
\[
\begin{aligned}
W_1(H,H^\star)
&\le
\sum_m\min\{v_m,w_m\}W_1(H_m,H_m^\star)
+R\sum_m|v_m-w_m|\\
&\le C(p_0)\Delta^{-(2s-2)}h+Ch \le C(p_0)\Delta^{-(2s-2)}h.
\end{aligned}
\]
This proves~\eqref{eq:structured-global-mixture}.
All moment comparisons take place on a subspace of
dimension at most $2K$, so the constants are independent
of the ambient dimension and the fitted weights.

\end{proof}

\begin{proof}[Proof of Theorem~\ref{thm:globalexact}]
Put $p_0=\pi_{\min}^\star$, $\Delta=\Delta_{\mathrm{sep}}$, and $h=h(p_{\theta,T},p_{\theta^\star,T})$. Throughout the proof, set $s=s_{\max}$ under the structured geometry and $s=K$ without cluster assumptions.

\medskip
\noindent\textit{Step 1:} We bound the one-time transport errors $W_1(G_\theta^{\{t\}},G_\star^{\{t\}})$ for $t\leq T$ without imposing positivity of the later-time true state probabilities. For each $t\le T$, define
\[
H_t=\frac12(G_\theta^{\{1\}}+G_\theta^{\{t\}}),\qquad H_t^\star=\frac12(G_\star^{\{1\}}+G_\star^{\{t\}}).
\]
The weight of $a_j$ in $H_t^\star$ is $\frac12(\pi_j^\star+q_t^\star(j))\ge p_0/2$. By joint convexity of squared Hellinger distance and contraction under marginalization,
\[
h^2(p_{H_t},p_{H_t^\star})\le\frac12h^2(p_\theta^{\{1\}},p_\star^{\{1\}})+\frac12h^2(p_\theta^{\{t\}},p_\star^{\{t\}})\le h^2.
\]
Under structured geometry, applying Lemma~\ref{lem:structured-global-transport} separately to the averaged and initial mixtures, then adding the bounds, gives
\[
W_1(H_t,H_t^\star)+W_1(G_\theta^{\{1\}},G_\star^{\{1\}})\le C(p_0)\Delta^{-(2s-2)}h.
\]
Since $G_\theta^{\{t\}}-G_\star^{\{t\}}=2(H_t-H_t^\star)-(G_\theta^{\{1\}}-G_\star^{\{1\}})$, Kantorovich--Rubinstein duality yields
\[
W_1(G_\theta^{\{t\}},G_\star^{\{t\}})\le2W_1(H_t,H_t^\star)+W_1(G_\theta^{\{1\}},G_\star^{\{1\}})\le C(p_0)\Delta^{-(2s-2)}h.
\]
Without clustering, the same argument uses Lemma~\ref{lem:direct-transport} with $s=K$. Thus, in both regimes,
\begin{equation}
\max_{t\le T}W_1(G_\theta^{\{t\}},G_\star^{\{t\}})\le C(p_0)\Delta^{-(2s-2)}h.
\label{eq:all-time-transport-proof}
\end{equation}
Lemma~\ref{lem:transport-assignment-matching}, applied with $p=r=1$, gives
\[
D_\mu\le\max_{t\le T}W_1(G_\theta^{\{t\}},G_\star^{\{t\}}),\qquad D_{\{t\}}\le4\Delta^{-1}W_1(G_\theta^{\{t\}},G_\star^{\{t\}}).
\]
Consequently,
\[
D_\mu\le C(p_0)\Delta^{-(2s-2)}h,\qquad \max_{t\le T}D_{\{t\}}\le C(p_0)\Delta^{-(2s-1)}h,
\]
including the bound for $D_\pi=D_{\{1\}}$.

\medskip
\noindent\textit{Step 2:} We prove that, under the structured geometry, $D_I\le C_r(p_0)\Delta^{-(2s-1)}h$ for every $I\subseteq[T]$ with $|I|=r$. Suppose $h\le c(p_0)\Delta^{2s-1}$. By Step~1, choosing $c(p_0)$ sufficiently small ensures
\[
W_1(G_\theta^{\{1\}},G_\star^{\{1\}})\le p_0\Delta/32.
\]
Each ball $B_d(a_j,\Delta/4)$ must contain a fitted center of positive initial weight; otherwise the true mass at $a_j$ alone would incur transport cost at least $p_0\Delta/4$. The $K$ balls are disjoint and there are exactly $K$ fitted labels. Hence, after relabeling, $\|\mu_j-a_j\|_2<\Delta/4$ for every $j\in[K]$.

For an optimal initial coupling, the off-diagonal transport distances are at least $3\Delta/4$. Therefore,
\[
\sum_j|\pi_j-\pi_j^\star|\le\frac{8}{3\Delta}W_1(G_\theta^{\{1\}},G_\star^{\{1\}})\le p_0/12,
\]
so $\pi_j\ge\pi_j^\star/2$. The independently proved paired bound in Theorem~\ref{thm:paired} now gives
\[
D_I\le C_r(p_0)\Delta^{-(2s-1)}h.
\]
If instead $h>c(p_0)\Delta^{2s-1}$, then $D_I\le2<(2/c(p_0))\Delta^{-(2s-1)}h$. Combining the two cases proves the structured joint bound globally.

\medskip
\noindent\textit{Step 3:} We prove that, without cluster assumptions, $D_I\le C_r(p_0)\Delta^{-(2s-1)}h$. If $h\le c(p_0)\Delta^{2s-1}$, Step~1 and the matching argument in Step~2 give, after relabeling,
\[
\|\mu_j-a_j\|_2<\Delta/4,\qquad \pi_j\ge\pi_j^\star/2,\qquad j\in[K].
\]
The paired bound~\eqref{eq:paired-diameter-joint} in Theorem~\ref{thm:paired}, together with $D_\star\ge\Delta$, yields
\[
D_I\le C_r(p_0)\Delta^{-(2s-1)}h.
\]
If instead $h>c(p_0)\Delta^{2s-1}$, the trivial bound $D_I\le2$ gives the same conclusion. Thus
\begin{equation}
D_I\le C_r(p_0)\Delta^{-(2s-1)}h.
\label{eq:global-diameter-probability}
\end{equation}

\medskip
\noindent\textit{Step 4:} We prove that $D_A\le C(p_0)\Delta^{-(2s-1)}h$ in both geometric regimes. If $h\le c(p_0)\Delta^{2s-1}$, Steps~2 and~3 establish matching and $\pi_j\ge\pi_j^\star/2$. The local transition-flow bound in Theorem~\ref{thm:paired} therefore gives
\[
D_A\le C(p_0)\Delta^{-(2s-1)}h.
\]
If instead $h>c(p_0)\Delta^{2s-1}$, each fitted and true transition-flow array is a probability array, so
\[
D_A\le2<\frac{2}{c(p_0)}\Delta^{-(2s-1)}h.
\]
Combining the two cases proves the global transition-flow bound, with $s=s_{\max}$ under the structured geometry and $s=K$ without cluster assumptions.
\end{proof}

\subsection{Proof of Corollary~\ref{cor:globalpath}}
\label{app:proof_cor_global_path}
\begin{proof}
Fix $I\subseteq[T]$ with $|I|=r\ge1$. Put $p_0=\pi_{\min}^\star$, $\Delta=\Delta_{\mathrm{sep}}$, $D_\star=\max_{i,j}\|a_i-a_j\|_2$, and $h=h(p_{\theta,T},p_{\theta^\star,T})$.

\medskip
\noindent\textit{Step 1:} We prove that, under the structured geometry, $W_{1,I}(G_\theta^I,G_\star^I)\le C_r(p_0)\Delta^{-(2s_{\max}-2)}h.$ Lemma~\ref{lem:hierarchical-transport} gives
\[
W_{1,I}(G_\theta^I,G_\star^I)\le V_{1,I}+\frac12C_0\Delta D_I+RJ_I,\qquad V_{1,I}\le D_\mu.
\]
Theorem~\ref{thm:globalexact} and Lemma~\ref{lem:cluster-mass} imply
\[
D_\mu\le C(p_0)\Delta^{-(2s_{\max}-2)}h,\qquad D_I\le C_r(p_0)\Delta^{-(2s_{\max}-1)}h,\qquad J_I\le C_r\Delta^{-(2s_{\max}-2)}h.
\]
Consequently, we achieve the proof of Equation~\eqref{eq:globalexact-path-wasserstein} as
\[
\begin{aligned}
W_{1,I}(G_\theta^I,G_\star^I)
&\le C(p_0)\Delta^{-(2s_{\max}-2)}h+C_r(p_0)\Delta\Delta^{-(2s_{\max}-1)}h+C_r\Delta^{-(2s_{\max}-2)}h\\
&\le C_r(p_0)\Delta^{-(2s_{\max}-2)}h.
\end{aligned}
\]

\medskip
\noindent\textit{Step 2:} We prove that, without cluster assumptions, $W_{1,I}(G_\theta^I,G_\star^I)\le C_r(p_0)\Delta^{-(2K-2)}h.$ For $r=1$, this follows directly from~\eqref{eq:globalexact-state-wasserstein-unstructured}. For $r\ge2$ and $K\ge2$, we use the same let $\bar G_\theta^I$ be the measure obtained by replacing each fitted center by its assigned true center. The assignment coupling gives
\[
W_{1,I}(G_\theta^I,\bar G_\theta^I)\le V_{1,I}.
\]
Both $\bar G_\theta^I$ and $G_\star^I$ are supported on the true grid, whose diameter under $d_{1,r}$ is at most $D_\star$. Matching their common masses and coupling the residual mass, which equals $D_I/2$, yields
\[
W_{1,I}(\bar G_\theta^I,G_\star^I)\le\frac12D_\star D_I.
\]
Thus, by the triangle inequality,
\[
W_{1,I}(G_\theta^I,G_\star^I)\le V_{1,I}+\frac12D_\star D_I.
\]
Using $V_{1,I}\le D_\mu$, Theorem~\ref{thm:globalexact}, and the diameter-dependent bound $D_I\le C_r(p_0)D_\star^{-2}\Delta^{-(2K-3)}h$, we obtain
\[
\begin{aligned}
W_{1,I}(G_\theta^I,G_\star^I)
&\le C(p_0)\Delta^{-(2K-2)}h+C_r(p_0)D_\star^{-1}\Delta^{-(2K-3)}h\\
&\le C_r(p_0)\Delta^{-(2K-2)}h,
\end{aligned}
\]
where the last inequality uses $D_\star\ge\Delta$. For $K=1$, both location-sequence measures are point masses and their normalized path distance equals their one-time distance. This proves~\eqref{eq:globalexact-path-wasserstein-unstructured}. 
\end{proof}

\subsection{Proof of Corollary~\ref{cor:localization}}
\label{app:proof_corollary3.3}
\begin{proof}
The matching argument in the proof of Theorem~\ref{thm:globalexact} proves the corollary under~\eqref{eq:exact-localization-threshold}. Substituting the radius of Proposition~\ref{prop:density} gives~\eqref{eq:exact-localization-sample-size}. For the explicit alternative~\eqref{eq:alternative-localization}, Lemma~\ref{lem:counting-near} bounds the total discrepancy between fitted near-set masses and true initial masses by $C\Delta_{\mathrm{sep}}^{-2s_{\max}}h$. Under the stated threshold this is at most $\pi_{\min}^\star/4$. Thus all $K$ near sets are nonempty, each contains exactly one fitted label, and every fitted initial weight is at least $3\pi_j^\star/4$. The same reasoning applies without clusters with $s_{\max}=K$.
\end{proof}

\subsection{Proof of Theorem~\ref{thm:paired}}
\label{app:proof_paired_theorem}

We first construct lower-degree state selectors for the structured geometry. Their values at fitted centers need not be exact: their errors will be controlled by the weighted displacement bound. We then compare Markov path distributions through their two-time marginals. Without cluster assumptions, we retain full products of center distances in this comparison.

\begin{lemma}[Cluster-localized cardinal polynomials]
\label{lem:lower-degree-cardinal}
Suppose that the structured geometry holds and $\|b_i-a_i\|_2\le\Delta_{\mathrm{sep}}/4$. For all sufficiently small $\Delta_{\mathrm{sep}}$, with the threshold depending only on $K,C_0,D_0$, there are polynomials $L_j$ of degree at most $3K-2$ such that
\begin{align}
L_j(a_i)&=\mathbf1\{i=j\},\qquad
\|L_j\|_{\infty,B_d(0,R)}\le C\Delta_{\mathrm{sep}}^{-(s_{\max}-1)},
\label{eq:lower-cardinal-norm}\\
\sum_{j=1}^K|L_j(b_i)|&\le C,\qquad
\sum_{j=1}^K|L_j(b_i)-\mathbf1\{i=j\}|
\le C\Delta_{\mathrm{sep}}^{-1}\|b_i-a_i\|_2.
\label{eq:lower-cardinal-values}
\end{align}
The constants depend only on $K,R,C_0,D_0$. In the single-cluster case, the assertion holds for every $0<\Delta_{\mathrm{sep}}\le1$.
\end{lemma}

\begin{proof}
We construct cardinal polynomials within each cluster and then separate the clusters.

\paragraph{Interpolation inside a cluster.}
For $j\in C_m$, put $s_m=|C_m|$ and define
\[
\ell_j(x)=\prod_{\ell\in C_m\setminus\{j\}}
\frac{\langle x-a_\ell,a_j-a_\ell\rangle}{\|a_j-a_\ell\|_2^2}.
\]
The empty product is one. On the true centers of $C_m$, this polynomial equals the cardinal indicator of $j$. Its degree is $s_m-1$, and its norm on $B_d(0,R)$ is at most $C\Delta_{\mathrm{sep}}^{-(s_m-1)}$. If $i\in C_m$ and $x\in[a_i,b_i]$, then
\[
\|x-a_\ell\|_2\le(C_0+1/4)\Delta_{\mathrm{sep}},
\qquad \ell\in C_m.
\]
Each factor in the product is therefore uniformly bounded on this segment, and the norm of its gradient is at most $\Delta_{\mathrm{sep}}^{-1}$. The product rule and the mean-value formula give
\[
|\ell_j(b_i)|\le C,
\qquad
|\ell_j(b_i)-\ell_j(a_i)|
\le C\Delta_{\mathrm{sep}}^{-1}\|b_i-a_i\|_2.
\]
For one cluster, take $L_j=\ell_j$, which proves both assertions.

\paragraph{Separating distinct clusters.}
For several clusters, apply Lemma~\ref{lem:separation-projection} to the distinct true and fitted centers, obtaining a projection $u$ that preserves all their pairwise distances up to a constant depending only on $K$. The union of the true and fitted centers assigned to a cluster has diameter at most $(C_0+1/2)\Delta_{\mathrm{sep}}$. Centers assigned to different clusters have distance at least $D_0-\Delta_{\mathrm{sep}}/2$. Thus, for $\Delta_{\mathrm{sep}}$ below a fixed positive threshold, the projected cluster intervals are separated by a fixed positive gap. Indeed, their diameters are $O(\Delta_{\mathrm{sep}})$, whereas the projected distance between any points from different clusters is bounded below by a fixed multiple of $D_0$; overlapping intervals would contradict these two facts.

For each cluster $m$, choose a smooth scalar function $f_m$ that equals one on its projected interval and zero on all other cluster intervals. The functions can be joined across the fixed gaps with derivatives through order $2K$ bounded independently of $\Delta_{\mathrm{sep}}$. Interpolate $f_m$ on the distinct projected union nodes. Lemma~\ref{lem:confluent-interpolation} gives a polynomial $Q_m$ of degree at most $2K-1$ whose norm on $[-R,R]$ is uniformly bounded. Define $F_m(x)=Q_m(\langle u,x\rangle)$. Then
\[
F_m(a_i)=F_m(b_i)=\mathbf1\{i\in C_m\},
\qquad \|F_m\|_{\infty,B_d(0,R)}\le C.
\]

\paragraph{Combining the two constructions.}
For $j\in C_m$, set $L_j=F_m\ell_j$. Its degree is at most $2K+s_m-2\le3K-2$, and its global norm satisfies~\eqref{eq:lower-cardinal-norm}. At every fitted or true center outside $C_m$, the factor $F_m$ vanishes. At every such center inside $C_m$, that factor equals one. Thus the within-cluster value estimates above apply exactly, with no approximation from $F_m$. Summing them over $j$ proves~\eqref{eq:lower-cardinal-values}.
\end{proof}

\begin{lemma}[Markov path comparison]
\label{lem:markov-path-comparison}
Let $P$ and $P^\star$ be the joint distributions of two Markov chains on a common finite state space at times $1,\ldots,r$, allowing time-inhomogeneous transition matrices. Write $q_\ell,q_\ell^\star$ for their marginal state probabilities and $F_\ell,F_\ell^\star$ for the joint probabilities at times $\ell,\ell+1$. Then
\begin{equation}
\|P-P^\star\|_1
\le \|q_1-q_1^\star\|_1
+2\sum_{\ell=1}^{r-1}\|F_\ell-F_\ell^\star\|_{\mathrm{ent},1},
\label{eq:markov-path-comparison}
\end{equation}
where $\|M\|_{\mathrm{ent},1}=\sum_{i,j}|M_{ij}|$. No lower bound on state probabilities is required.
\end{lemma}

\begin{proof}
Let $B_\ell,B_\ell^\star$ be the transition matrices. For a fixed path $x=(x_1,\ldots,x_r)$ the product identity gives
\[
\begin{aligned}
P(x)-P^\star(x)
&=\{q_1(x_1)-q_1^\star(x_1)\}
     \prod_{a=1}^{r-1}B_a(x_a,x_{a+1})\\
&\quad+q_1^\star(x_1)\sum_{\ell=1}^{r-1}
 \left\{\prod_{a<\ell}B_a^\star(x_a,x_{a+1})\right\}
 (B_\ell-B_\ell^\star)(x_\ell,x_{\ell+1})
 \left\{\prod_{a>\ell}B_a(x_a,x_{a+1})\right\}.
\end{aligned}
\]
Empty products equal one. Taking absolute values and summing over paths leaves total future mass one and true past marginal $q_\ell^\star$, giving
\[
\|P-P^\star\|_1\le\|q_1-q_1^\star\|_1
+\sum_{\ell=1}^{r-1}\sum_iq_\ell^\star(i)
       \|B_\ell(i,\cdot)-B_\ell^\star(i,\cdot)\|_1.
\]
For each $i,j$,
\[
q_\ell^\star(i)\{B_\ell(i,j)-B_\ell^\star(i,j)\}
=F_\ell(i,j)-F_\ell^\star(i,j)
 +\{q_\ell^\star(i)-q_\ell(i)\}B_\ell(i,j).
\]
Summing absolute values and using the row sums of $B_\ell$ yields
\[
\begin{aligned}
\sum_iq_\ell^\star(i)\|B_\ell(i,\cdot)-B_\ell^\star(i,\cdot)\|_1
&\le\|F_\ell-F_\ell^\star\|_{\mathrm{ent},1}
   +\|q_\ell-q_\ell^\star\|_1\\
&\le2\|F_\ell-F_\ell^\star\|_{\mathrm{ent},1}.
\end{aligned}
\]
The final inequality is contraction of $\ell_1$ distance on summing the second coordinate of the two-time distributions. Substitution proves equation~\eqref{eq:markov-path-comparison}. For $r=1$ the sum is empty and the assertion is equality. No step divides by a state or path probability. For nonconsecutive selected HMM times, the required matrices are $A^{t_{\ell+1}-t_\ell}$ and $(A^\star)^{t_{\ell+1}-t_\ell}$; they are stochastic even when some marginal probabilities vanish. The argument applies to aligned exact-fit chains, not to an arbitrary grouped overfitted process.
\end{proof}

\begin{proof}[Proof of Theorem~\ref{thm:paired}]
We use the assumed labeling, so that $\mu_j=b_j$ and $\|b_j-a_j\|_2\le\Delta_{\mathrm{sep}}/4$. This labeling also agrees with the nearest-center assignment. Indeed, for $\ell\ne j$,
\[
\begin{aligned}
\|b_j-a_\ell\|_2
&\ge
\|a_j-a_\ell\|_2-\|b_j-a_j\|_2
\qquad\text{(by the reverse triangle inequality)}\\
&\ge
\Delta_{\mathrm{sep}}-\frac{\Delta_{\mathrm{sep}}}{4}\\
&=
\frac{3\Delta_{\mathrm{sep}}}{4}\\
&>
\frac{\Delta_{\mathrm{sep}}}{4}
\ge \|b_j-a_j\|_2.
\end{aligned}
\]
Thus $c(j)=j$, every assignment set is a singleton, and $r_j=h_j=\|b_j-a_j\|_2$. The aggregate losses can consequently be written using the aligned state labels.

\paragraph{Moment comparison for retained observations.}
Before extracting the different losses, we record the observation-error bound used for each polynomial. Let $I\subseteq[T]$ contain $r$ times, and let $U$ be a polynomial on $(\mathbb R^d)^r$ of degree at most $r(3K-2)$. Both location distributions are supported in the span
\[
S=\operatorname{span}\{a_1,\ldots,a_K,b_1,\ldots,b_K\},
\qquad \dim S\le2K.
\]
The HMM application of Lemma~\ref{lem:moment-extraction} therefore has a constant depending only on $K,r,R$ and the degree. Applying it to the retained observations gives the first inequality below; marginalization then replaces their density error by that of the full observed block:
\begin{equation}
\label{eq:paired-moment-control}
\begin{aligned}
\left|\int U\,d(G_\theta^I-G_\star^I)\right|
&\le
C_r\|U\|_\infty
h(p_\theta^I,p_\star^I)
\qquad\text{(by Lemma~\ref{lem:moment-extraction})}\\
&\le
C_r\|U\|_\infty
h(p_{\theta,T},p_{\theta^\star,T})
\qquad\text{(by Lemma~\ref{lem:marginalization})}.
\end{aligned}
\end{equation}
Here the supremum norm is taken on $B_d(0,R)^r$. The remaining task is to choose $U$ whose moment difference equals the recovery quantity being bounded.

\paragraph{Weighted displacement.}
Fix $j\in[K]$ and $t\le T$. At one time, the fitted and true location distributions are
\[
G_\theta^{\{t\}}
=
\sum_{\ell=1}^K q_t(\ell)\delta_{b_\ell},
\qquad
G_\star^{\{t\}}
=
\sum_{\ell=1}^K q_t^\star(\ell)\delta_{a_\ell}.
\]
Expanding the moment of $\widetilde H_j$ over these supports gives
\[
\begin{aligned}
\int\widetilde H_j\,d(G_\theta^{\{t\}}-G_\star^{\{t\}})
&=
\sum_{\ell=1}^K q_t(\ell)\widetilde H_j(b_\ell)
-
\sum_{\ell=1}^K q_t^\star(\ell)\widetilde H_j(a_\ell)\\
&=
\sum_{\ell=1}^K
q_t(\ell)\mathbf1\{j=\ell\}h_j
-
\sum_{\ell=1}^K q_t^\star(\ell)\,0
\qquad\text{(by Equation~\eqref{eq:paired-displacement-values})}\\
&=q_t(j)h_j.
\end{aligned}
\]
The extracted quantity is nonnegative, so taking an absolute value does not change it. Equation~\eqref{eq:paired-moment-control} now gives
\[
\begin{aligned}
q_t(j)h_j
&=
\left|
\int\widetilde H_j\,d(G_\theta^{\{t\}}-G_\star^{\{t\}})
\right|\\
&\le
C\|\widetilde H_j\|_\infty
h(p_{\theta,T},p_{\theta^\star,T})\\
&\le
C\Delta_{\mathrm{sep}}^{-(2s_{\max}-2)}
h(p_{\theta,T},p_{\theta^\star,T})
\qquad\text{(by Equation~\eqref{eq:paired-selector-norms})}.
\end{aligned}
\]
Summing over $j$ yields
\[
\begin{aligned}
\sum_{j=1}^K q_t(j)h_j
&\le
KC\Delta_{\mathrm{sep}}^{-(2s_{\max}-2)}
h(p_{\theta,T},p_{\theta^\star,T})\le
C\Delta_{\mathrm{sep}}^{-(2s_{\max}-2)}
h(p_{\theta,T},p_{\theta^\star,T}),
\end{aligned}
\]
where the fixed factor $K$ is absorbed into $C$. Taking the maximum over $t$ proves
\begin{equation}
\label{eq:paired-proof-displacement}
D_\mu
\le
C\Delta_{\mathrm{sep}}^{-(2s_{\max}-2)}
h(p_{\theta,T},p_{\theta^\star,T}).
\end{equation}
No lower bound on $q_t(j)$ is used; a state with zero probability contributes zero to the extracted displacement.

\paragraph{One-time probabilities.}
For the cardinal polynomial $\widetilde E_j$, the same finite-support expansion gives
\[
\begin{aligned}
\int\widetilde E_j\,d(G_\theta^{\{t\}}-G_\star^{\{t\}})
&=
\sum_{\ell=1}^K q_t(\ell)\widetilde E_j(b_\ell)
-
\sum_{\ell=1}^K q_t^\star(\ell)\widetilde E_j(a_\ell)\\
&=
\sum_{\ell=1}^K
\{q_t(\ell)-q_t^\star(\ell)\}\mathbf1\{j=\ell\}
\qquad\text{(by Equation~\eqref{eq:paired-cardinal-values})}\\
&=q_t(j)-q_t^\star(j).
\end{aligned}
\]
Applying Equation~\eqref{eq:paired-moment-control} to this identity, we obtain
\[
\begin{aligned}
|q_t(j)-q_t^\star(j)|
&\le
C\|\widetilde E_j\|_\infty
h(p_{\theta,T},p_{\theta^\star,T})\\
&\le
C\Delta_{\mathrm{sep}}^{-(2s_{\max}-1)}
h(p_{\theta,T},p_{\theta^\star,T})
\qquad\text{(by Equation~\eqref{eq:paired-selector-norms})}.
\end{aligned}
\]
We take absolute values before summing over states, as required by the definition of $D_{\{t\}}$:
\[
\begin{aligned}
D_{\{t\}}
&=\sum_{j=1}^K|q_t(j)-q_t^\star(j)|\\
&\le
KC\Delta_{\mathrm{sep}}^{-(2s_{\max}-1)}
h(p_{\theta,T},p_{\theta^\star,T})\\
&\le
C\Delta_{\mathrm{sep}}^{-(2s_{\max}-1)}
h(p_{\theta,T},p_{\theta^\star,T}).
\end{aligned}
\]
This bound is uniform in $t$, proving the one-time probability assertion. In particular, $D_\pi=D_{\{1\}}$ satisfies the same bound.

\paragraph{Two-time probabilities in structured configurations.}
We first treat $\Delta_{\mathrm{sep}}$ below the threshold in Lemma~\ref{lem:lower-degree-cardinal}. Fix arbitrary times $t<u$, and let
\[
F^{t,u}_{i\ell}=\mathbb P_\theta(Z_t=i,Z_u=\ell),
\qquad
F^{t,u}_{\star,i\ell}=\mathbb P_{\theta^\star}(Z_t=i,Z_u=\ell),
\qquad E_{ji}=L_j(b_i).
\]
The row and column sums of $F^{t,u}$ are $q_t$ and $q_u$. Since $L_j(a_i)=\mathbf1\{j=i\}$,
\[
\int L_j(x)L_\ell(y)\,d(G_\theta^{\{t,u\}}-G_\star^{\{t,u\}})
=\{EF^{t,u}E^\top-F_\star^{t,u}\}_{j\ell}.
\]
The polynomial has degree at most $2(3K-2)$ and norm at most $C\Delta_{\mathrm{sep}}^{-(2s_{\max}-2)}$. Applying~\eqref{eq:paired-moment-control} entry by entry and summing over the $K^2$ entries gives
\begin{equation}
\|EF^{t,u}E^\top-F_\star^{t,u}\|_{\mathrm{ent},1}
\le C\Delta_{\mathrm{sep}}^{-(2s_{\max}-2)}
h(p_{\theta,T},p_{\theta^\star,T}).
\label{eq:pair-moment-improved}
\end{equation}

It remains to account for the values of the selectors at fitted centers. Put
\[
e_i=\sum_j|E_{ji}-\mathbf1\{j=i\}|,
\qquad m_i=\sum_j|E_{ji}|.
\]
Lemma~\ref{lem:lower-degree-cardinal} gives
$e_i\le C\Delta_{\mathrm{sep}}^{-1}h_i$ and $m_i\le C$.
The identity
\[
EF^{t,u}E^\top-F^{t,u}
=(E-\Id)F^{t,u}E^\top+F^{t,u}(E^\top-\Id)
\]
and the nonnegativity of $F^{t,u}$ imply
\begin{align*}
\|EF^{t,u}E^\top-F^{t,u}\|_{\mathrm{ent},1}
&\le\sum_{i,\ell}F^{t,u}_{i\ell}(e_i m_\ell+e_\ell)\\
&\le C\Delta_{\mathrm{sep}}^{-1}
\left\{\sum_i q_t(i)h_i+\sum_\ell q_u(\ell)h_\ell\right\}\\
&\le C\Delta_{\mathrm{sep}}^{-1}D_\mu.
\end{align*}
Combining this inequality with~\eqref{eq:pair-moment-improved} and~\eqref{eq:paired-proof-displacement} yields
\begin{equation}
D_{\{t,u\}}=\|F^{t,u}-F_\star^{t,u}\|_{\mathrm{ent},1}
\le C\Delta_{\mathrm{sep}}^{-(2s_{\max}-1)}
h(p_{\theta,T},p_{\theta^\star,T}).
\label{eq:paired-two-time-improved}
\end{equation}

For completeness, when $\Delta_{\mathrm{sep}}$ is above the fixed threshold, use the products $\widetilde E_j(x)\widetilde E_\ell(y)$ from Appendix~\ref{app:paired-interpolation}. Their values on both supports are the exact pair indicators, and their norms are bounded by $C\Delta_{\mathrm{sep}}^{-(4s_{\max}-2)}$. Since separation is now bounded below by a fixed constant, this is at most a constant times $\Delta_{\mathrm{sep}}^{-(2s_{\max}-1)}$. Thus~\eqref{eq:paired-two-time-improved} holds throughout the structured regime. Taking $u=t+1$ and maximizing proves the asserted bound for $D_A$.

\paragraph{State-sequence probabilities.}
Let $I=(t_1,\ldots,t_r)$. Under either fitted or true parameters, the selected states form a Markov chain with transition matrices $A^{t_{\ell+1}-t_\ell}$ or $(A^\star)^{t_{\ell+1}-t_\ell}$. Because the labels are aligned, Lemma~\ref{lem:markov-path-comparison} applies directly to these two chains and gives
\[
D_I\le D_{\{t_1\}}
+2\sum_{\ell=1}^{r-1}D_{\{t_\ell,t_{\ell+1}\}}.
\]
The one-time bound already proved and~\eqref{eq:paired-two-time-improved} establish~\eqref{eq:paired-path-improved}. The number of selected times enters only the constant in this step.

\paragraph{Unstructured probability bounds with distance products.}
Write $\Delta=\Delta_{\mathrm{sep}}$, $D_\star=\max_{i,j}\|a_i-a_j\|_2$, and
$\mathcal D_i=\prod_{j\ne i}\|a_i-a_j\|_2$.
The paired displacement polynomial retains the state-specific estimate
\begin{equation}
q_t(i)h_i\le C\mathcal D_i^{-2}h(p_{\theta,T},p_{\theta^\star,T}).
\label{eq:paired-state-product}
\end{equation}
Indeed, the proof of Lemma~\ref{lem:paired-product} gives $\|R_i\|_\infty\le C\mathcal D_i^{-2}$ before replacing its denominator by a separation power. The selected-segment lower bound for $R_i(b_i)$ is uniform. The definition of $\widetilde H_i$ and moment extraction therefore prove~\eqref{eq:paired-state-product}, including when $h_i=0$.

Use the ordinary cardinal polynomials $L_j$ from Lemma~\ref{lem:diameter-joint-selectors}, and put $E_{ji}=L_j(b_i)$ as before. Since $\|L_j\|_\infty\le C/\mathcal D_j$ and $\mathcal D_j\ge(D_\star/2)\Delta^{K-2}$, the pair-moment comparison gives
\[
\|EF^{t,u}E^\top-F_\star^{t,u}\|_{\mathrm{ent},1}
\le CD_\star^{-2}\Delta^{-(2K-4)}h(p_{\theta,T},p_{\theta^\star,T}).
\]
Also $\sqrt{P(b_i)}\le Ch_i\mathcal D_i$ under matching. Equation~\eqref{eq:diameter-linear-error} gives
\[
e_i=\sum_j|E_{ji}-\mathbf1\{j=i\}|
\le Ca_\star h_i\mathcal D_i,
\qquad a_\star=(D_\star\Delta^{K-1})^{-1}.
\]
Expand $EF^{t,u}E^\top-F^{t,u}$ into its two linear error terms and the product of the two errors. Nonnegativity of the pair probabilities gives
\[
\|EF^{t,u}E^\top-F^{t,u}\|_{\mathrm{ent},1}
\le\sum_i q_t(i)e_i+\sum_\ell q_u(\ell)e_\ell
+\sum_{i,\ell}F^{t,u}_{i\ell}e_i e_\ell.
\]
The first two sums are at most
$CD_\star^{-2}\Delta^{-(2K-3)}h(p_{\theta,T},p_{\theta^\star,T})$
by~\eqref{eq:paired-state-product} and the lower bound on $\mathcal D_i$.
For the last sum, use $h_i,h_\ell\le\Delta/4$ and both marginals of $F^{t,u}$:
\[
F^{t,u}_{i\ell}h_i h_\ell
\le C\Delta\min\{\mathcal D_i^{-2},\mathcal D_\ell^{-2}\}
h(p_{\theta,T},p_{\theta^\star,T}).
\]
Multiplication by $a_\star^2\mathcal D_i\mathcal D_\ell$ leaves the factor
\[
a_\star^2\Delta
\min\{\mathcal D_\ell/\mathcal D_i,\mathcal D_i/\mathcal D_\ell\}
\le D_\star^{-2}\Delta^{-(2K-3)}.
\]
Thus the cross term has the same bound. Combining the moment and evaluation comparisons proves
\begin{equation}
D_{\{t,u\}}\le CD_\star^{-2}\Delta^{-(2K-3)}h(p_{\theta,T},p_{\theta^\star,T}).
\label{eq:paired-diameter-pair}
\end{equation}
Every one-time discrepancy is bounded by a pair discrepancy, since $T\ge2$. Lemma~\ref{lem:markov-path-comparison} consequently gives
\begin{equation}
D_I\le C_rD_\star^{-2}\Delta^{-(2K-3)}h(p_{\theta,T},p_{\theta^\star,T})
\le C_r\Delta^{-(2K-1)}h(p_{\theta,T},p_{\theta^\star,T}).
\label{eq:paired-diameter-joint}
\end{equation}
This proves~\eqref{eq:paired-unstructured-probability} for every nonempty $I$, with constants independent of state probabilities.

\paragraph{Structured transport.}
The transport decomposition in Lemma~\ref{lem:hierarchical-transport}, with $p=1$, gives
\[
W_{1,I}(G_\theta^I,G_\star^I)
\le V_{1,I}+\frac{C_0\Delta_{\mathrm{sep}}}{2}D_I+RJ_I.
\]
The displacement term is bounded by $D_\mu$, and Lemma~\ref{lem:cluster-mass} bounds $J_I$ by
$C_r\Delta_{\mathrm{sep}}^{-(2s_{\max}-2)}h(p_{\theta,T},p_{\theta^\star,T})$.
The state-sequence bound just proved gives
\[
\frac{C_0\Delta_{\mathrm{sep}}}{2}D_I
\le C_r\Delta_{\mathrm{sep}}^{-(2s_{\max}-2)}
h(p_{\theta,T},p_{\theta^\star,T}).
\]
All three contributions therefore have separation exponent $2s_{\max}-2$, proving~\eqref{eq:paired-transport-improved}.

\paragraph{Unstructured transport.}
For a single selected time, apply~\eqref{eq:direct-transport-paired} of Lemma~\ref{lem:direct-transport} to the corresponding location measures, and then use marginalization. This gives~\eqref{eq:paired-unstructured-marginal}, even when some true or fitted state probabilities at that time vanish.
For every nonempty $I$, move fitted locations to their assigned true centers and then couple the two probability measures on the true grid. This gives
\[
W_{1,I}(G_\theta^I,G_\star^I)\le V_{1,I}+\tfrac12D_\star D_I.
\]
The first term is bounded by $D_\mu\le C\Delta^{-(2K-2)}h(p_{\theta,T},p_{\theta^\star,T})$. Equation~\eqref{eq:paired-diameter-joint} bounds the second by $C_rD_\star^{-1}\Delta^{-(2K-3)}h(p_{\theta,T},p_{\theta^\star,T})$, which has separation exponent $2K-2$ because $D_\star\ge\Delta$. Their sum proves~\eqref{eq:paired-unstructured-path}. No state probability has been divided out in any of these comparisons.

\end{proof}

\subsubsection{Sharpness of the exact-fit exponents}

The preceding bounds hold without restrictions on the rank of the transition matrix. This includes HMMs whose hidden states are independent. Moment matching in that submodel gives lower examples for the separation powers, even after the centers have been paired.

\begin{proposition}[Sharpness under exact specification]
\label{prop:matched-sharpness}
Fix $s\ge2$, a dimension $d\ge1$, and an observation horizon $T\ge2$. There are exactly specified Gaussian HMMs with $K=s$, positive initial probabilities bounded away from zero, and a single cluster with $\Delta_{\mathrm{sep}}\asymp\delta\downarrow0$, such that their centers satisfy the matching condition of Theorem~\ref{thm:paired} and
\begin{align}
h(p_{\theta_\delta,T},p_{\theta_\delta^\star,T})&\asymp\delta^{2s-1},
\label{eq:sharp-density}\\
D_\pi\asymp D_A\asymp D_I&\asymp1,\qquad \varnothing\ne I\subseteq[T],
\label{eq:sharp-probabilities}\\
D_\mu\asymp W_{1,I}(G_{\theta_\delta}^I,G_{\theta_\delta^\star}^I)&\asymp\delta.
\label{eq:sharp-transport}
\end{align}
The same construction can be embedded in a largest cluster of size $s$ with additional, fixed, separated clusters. Consequently, the structured probability exponent $2s-1$ and transport exponent $2s-2$ cannot be decreased uniformly over classes containing these collinear configurations. Taking $s=K$ also proves sharpness of the unstructured probability and transport exponents $2K-1$ and $2K-2$. This statement concerns the separation exponents of the deterministic inverse inequalities, rather than exact minimax constants or every individual spatial configuration.
\end{proposition}

\begin{proof}
We construct two fixed moment-matched mixing measures and then scale their supports.

\paragraph{Positive atomic measures with matching lower moments.}
Choose distinct real numbers $t_1<\cdots<t_s$ with minimum gap one and strictly positive weights $w_j$ summing to one. Consider the moments
\[
M_m(w,t)=\sum_{j=1}^s w_jt_j^m,\qquad m=0,\ldots,2s-1,
\]
initially allowing the $s$ weights to vary freely. The Jacobian of $(M_0,\ldots,M_{2s-1})$ in the $2s$ weight and location coordinates is invertible. Indeed, a linear dependence between its rows would give a polynomial $P$ of degree at most $2s-1$ satisfying
\[
P(t_j)=0,\qquad w_jP'(t_j)=0,\qquad j=1,\ldots,s.
\]
Since all $w_j>0$, this polynomial has $s$ double roots and must vanish identically. The inverse function theorem therefore supplies a smooth local inverse for this complete moment map. Apply that inverse to the moment vector whose first $2s-1$ coordinates equal their original values and whose last coordinate is $M_{2s-1}(w,t)+u$. This defines a smooth curve $(w(u),t(u))$ such that
\[
M_j(w(u),t(u))=M_j(w,t)\quad(0\le j\le2s-2),\qquad
M_{2s-1}(w(u),t(u))-M_{2s-1}(w,t)=u.
\]
In particular, its tangent is nonzero and its first unmatched moment is exactly $u$, not merely a quantity with a potentially vanishing higher-order remainder.

The tangent can be written explicitly. Put $g(z)=\prod_{j=1}^s(z-t_j)$. The value and derivative Hermite cardinal polynomials have respective leading coefficients
\[
 b_j=-\frac{g''(t_j)}{g'(t_j)^3},\qquad
 c_j=\frac1{g'(t_j)^2}.
\]
Indeed, with $\ell_j(z)=g(z)/\{g'(t_j)(z-t_j)\}$, those polynomials are
$(1-2\ell_j'(t_j)(z-t_j))\ell_j(z)^2$ and $(z-t_j)\ell_j(z)^2$; here $2\ell_j'(t_j)=g''(t_j)/g'(t_j)$. Comparing the leading coefficients in the Hermite interpolation identity for polynomials of degree at most $2s-1$ shows that the unique derivative of the moment-inverse curve is
\begin{equation}
 w_j'(0)=b_j,\qquad
 t_j'(0)=\frac{c_j}{w_j}.
 \label{eq:exact-sharpness-tangent}
\end{equation}
Every location derivative is positive. The weight derivative vector is nonzero because
$g''(t_s)/g'(t_s)=2\sum_{j<s}(t_s-t_j)^{-1}>0$.
Thus both kinds of parameters genuinely change along this curve. Choose a sufficiently small fixed $u\ne0$ so that all weights remain positive, the weight vectors differ, the location vectors differ, and $|t_j(u)-t_j|<1/4$. Set
\[
\nu^\star=\sum_{j=1}^s w_j\delta_{t_j},\qquad
\nu=\sum_{j=1}^s w_j(u)\delta_{t_j(u)}.
\]
These measures match moments through order $2s-2$ and have different moments of order $m=2s-1$.

\paragraph{Scaling and Gaussian density comparison.}
Scale all locations by $\delta$ and embed them along the first coordinate of $\mathbb R^d$. For sufficiently small $\delta$, both mixing measures lie in $B_d(0,R)$, the true minimum separation is $\delta$, and the paired displacements are below $\delta/4$. Their cluster diameter is a fixed multiple of $\delta$. Write $p_\delta,p_\delta^\star$ for the scalar Gaussian mixture densities on the first coordinate; the other coordinates have a common independent standard Gaussian density.

Taylor-expand $\varphi(x-\delta t)$ in its shift through degree $m-1$. Every term cancels in $p_\delta-p_\delta^\star$ by the matched moments. The locations $t_j,t_j(u)$ lie in a fixed bounded interval, so the integral remainder and the formula for derivatives of the Gaussian imply
\[
|p_\delta(x)-p_\delta^\star(x)|
\le C\delta^m(1+|x|^m)\varphi(x)e^{C\delta|x|},
\qquad 0<\delta\le\delta_0.
\]
Also, since the mixture weights are positive and the shifts are bounded by $C\delta$,
\[
p_\delta^\star(x)\ge c\varphi(x)e^{-C\delta|x|-C\delta^2}.
\]
Consequently,
\begin{align*}
h^2(p_\delta,p_\delta^\star)
&\le\frac12\int\frac{(p_\delta-p_\delta^\star)^2}{p_\delta^\star}\,dx\\
&\le C\delta^{2m}\int(1+|x|^{2m})\varphi(x)e^{C\delta|x|}\,dx
\le C\delta^{2m}.
\end{align*}
The last integral is bounded uniformly for $\delta\le\delta_0$. Conversely, the difference of the latent $m$th moments is a nonzero fixed constant times $\delta^m$. Applying Lemma~\ref{lem:moment-extraction} to $x^m$ on the fixed ball $B_d(0,R)$ gives $h(p_\delta,p_\delta^\star)\ge c\delta^m$. Its norm is bounded by $R^m$, independently of $\delta$; no norm on a shrinking support is used. Hence their Hellinger distance is of order $\delta^{2s-1}$.

\paragraph{Embedding into HMMs and comparing the latent losses.}
Take the true initial distribution and each row of its transition matrix to be $w$, and take the fitted initial distribution and each fitted transition row to be $w(u)$. Each observation sequence is then independent over time. The product affinity identity gives
\[
h^2(p_{\theta_\delta,T},p_{\theta_\delta^\star,T})
=1-\{1-h^2(p_\delta,p_\delta^\star)\}^{T},
\]
which lies between the one-time squared Hellinger distance and $T$ times that distance. This proves~\eqref{eq:sharp-density}.

The weight vectors differ by a fixed nonzero amount, so $D_\pi$ is bounded above and below by positive constants. Marginalization of any two-time or path distribution gives the same lower bound for $D_A,D_I$, and their universal upper bound is two. The weighted center displacement is a positive fixed constant times $\delta$. Finally,
\[
W_1(G_{\theta_\delta}^{\{1\}},G_{\theta_\delta^\star}^{\{1\}})
=\delta W_1(\nu,\nu^\star)>0.
\]
For the normalized additive path metric $d_{1,r}$, the Wasserstein distance between product measures equals this marginal distance. A product of optimal marginal couplings proves the upper bound; averaging an optimal marginal Lipschitz test function over the $r$ coordinates proves the lower bound. This establishes~\eqref{eq:sharp-probabilities}--\eqref{eq:sharp-transport}.

To include further clusters, assign the moving cluster a fixed positive total mass and add identical fixed atoms outside it in both models, again using identical rows within each model. The density difference is multiplied by this fixed mass, while the preceding upper Hellinger argument remains valid using the moving cluster as a lower bound for the reference density. The $m$th moment difference still gives the matching lower bound. Probability and displacement discrepancies retain their orders. A Lipschitz transport potential on the moving cluster extends to the additional common atoms, whose signed contributions cancel, so the transport order also remains unchanged. This proves the claimed embedding.
\end{proof}

\begin{example}[Aggregation can retain higher-order dependence]
\label{ex:aggregate-memory}
Consider a two-state true HMM with initial probabilities $(1/2,1/2)$, both transition rows $(1/2,1/2)$, and scalar emission centers $-\delta,\delta$. For a fitted eight-state HMM, label the states by $(u,v,w)\in\{0,1\}^3$, use the uniform initial distribution, and take the deterministic transition
\[
(u,v,w)\longmapsto(v,w,u\oplus v\oplus w),
\]
where $\oplus$ denotes addition modulo two. The fitted emission center is $-\delta$ or $\delta$ according to its first bit. Its first four aggregate labels are
\[
U,\quad V,\quad W,\quad U\oplus V\oplus W,
\]
where $U,V,W$ are independent uniform bits. Every marginal on at most three of these four coordinates is uniform and agrees with the corresponding true marginal. However, the fitted joint distribution of the four coordinates is uniform on the eight even-parity sequences, whereas the true joint distribution is uniform on all sixteen sequences. Thus $D_{\{1,2,3,4\}}=1$ although all corresponding two-time discrepancies vanish.

Writing $f_0(x)=\varphi(x+\delta)$ and $f_1(x)=\varphi(x-\delta)$, the four-observation density difference is exactly
\[
p_{\theta,4}-p_{\theta^\star,4}
=\frac1{16}\prod_{\ell=1}^4\{f_0(x_\ell)-f_1(x_\ell)\}.
\]
Its Hellinger distance is of order $\delta^4$. For the upper bound, divide the squared density difference by the true density $\prod_\ell(f_0+f_1)/2$ and integrate; each scalar factor is $O(\delta^2)$ by the mean-value formula for the Gaussian. For the lower bound, apply Lemma~\ref{lem:moment-extraction} to $x_1x_2x_3x_4$: its latent expectation differs by $\delta^4$. This example satisfies the over-specified parameter assumptions, including repeated fitted centers and zero transition entries. It shows why Lemma~\ref{lem:markov-path-comparison} cannot be applied to an arbitrary aggregate fitted process, and why the structured exact-fit probability exponent $3$ for two true states does not extend to all four-time joint distributions under over-specification.
\end{example}

\subsection{Proofs of Proposition~\ref{prop:mle-matched-rates} and Corollary~\ref{cor:rows}}
\label{app:corollary_3_5_3_6_proof}
\begin{proof}[Proof of Proposition~\ref{prop:mle-matched-rates}]
Let $\widehat\theta_n$ be the exact-fit MLE over $\Theta_K(R)$. Use the Borel relabeling convention in Appendix~\ref{app:loss-measurability}, including its identity-permutation extension outside the matching event. Define
\[
\mathcal E_n
=
\left\{
h(p_{\widehat\theta_n,T},p_{\theta^\star,T})
\le\varepsilon_{n,T}
\right\}.
\]
By Proposition~\ref{prop:density}, $\mathbb P_\star(\mathcal E_n)\ge1-\frac{1}{n}$. The sample-size condition in~\eqref{eq:exact-localization-sample-size} ensures that $\varepsilon_{n,T}\le c(\pi_{\min}^\star)\Delta_{\mathrm{sep}}^{2s_{\max}-1}$. On $\mathcal E_n$, we have:
\[
\begin{aligned}
h(p_{\widehat\theta_n,T},p_{\theta^\star,T})
&\le\varepsilon_{n,T}
\le c(\pi_{\min}^\star)\Delta_{\mathrm{sep}}^{2s_{\max}-1}.
\end{aligned}
\]
Corollary~\ref{cor:localization} therefore supplies the matching required by Theorem~\ref{thm:paired}.  Under this matching, its displacement bound gives
\[
\begin{aligned}
D_\mu(\widehat\theta_n,\theta^\star)
&\le
C\Delta_{\mathrm{sep}}^{-(2s_{\max}-2)}
h(p_{\widehat\theta_n,T},p_{\theta^\star,T})\\
&\le
C\Delta_{\mathrm{sep}}^{-(2s_{\max}-2)}
\varepsilon_{n,T}
\qquad\text{(because $\mathcal E_n$ holds)}.
\end{aligned}
\]
Its one-time probability bound gives
\[
\begin{aligned}
D_\pi(\widehat\theta_n,\theta^\star)
&\le
C\Delta_{\mathrm{sep}}^{-(2s_{\max}-1)}
h(p_{\widehat\theta_n,T},p_{\theta^\star,T})\\
&\le
C\Delta_{\mathrm{sep}}^{-(2s_{\max}-1)}
\varepsilon_{n,T},
\end{aligned}
\]
and, in the structured regimes, its consecutive-pair bound gives
\[
\begin{aligned}
D_A(\widehat\theta_n,\theta^\star)
&\le
C\Delta_{\mathrm{sep}}^{-(2s_{\max}-1)}
h(p_{\widehat\theta_n,T},p_{\theta^\star,T})\\
&\le
C\Delta_{\mathrm{sep}}^{-(2s_{\max}-1)}
\varepsilon_{n,T}.
\end{aligned}
\]
The structured joint probability and transport bounds follow in the same way from Theorem~\ref{thm:paired}. In the unstructured regime, take $s_{\max}=K$: the probability exponent is $2K-1$ and the Wasserstein exponent is $2K-2$ for every nonempty set of selected times. All these inequalities hold on $\mathcal E_n$, with probability at least $1-1/n$. Substituting $\varepsilon_{n,T}=C\sqrt{d\log(nTd)/n}$ gives the stated $\mathcal O_{\mathbb P}$ rates.
\end{proof}

\begin{proof}[Proof of Corollary~\ref{cor:rows}]
Work under the localization condition and use the aligned labeling supplied by Corollary~\ref{cor:localization}. In particular,
\[
\pi_j\ge\frac{\pi_j^\star}{2}
\ge\frac{\pi_{\min}^\star}{2},
\qquad j\in[K].
\]

\paragraph{Maximum center displacement.}
Write $h_j=\|b_j-a_j\|_2$, and choose $j_0$ such that $h_{j_0}=\max_j h_j$. The fitted initial probability lower bound allows us to remove the probability weight from this largest displacement:
\[
\begin{aligned}
\frac{\pi_{\min}^\star}{2}\max_j h_j
&=
\frac{\pi_{\min}^\star}{2}h_{j_0}\\
&\le \pi_{j_0}h_{j_0}
\qquad\text{(because $\pi_{j_0}\ge\pi_{\min}^\star/2$)}\\
&\le \sum_{j=1}^K\pi_jh_j
\qquad\text{(because every summand is nonnegative)}\\
&\le D_\mu
\qquad\text{(because $q_1=\pi$)}\\
&\le
C\Delta_{\mathrm{sep}}^{-(2s_{\max}-2)}
h(p_{\theta,T},p_{\theta^\star,T})
\qquad\text{(by Theorem~\ref{thm:paired})}.
\end{aligned}
\]
Dividing by $\pi_{\min}^\star/2>0$ proves
\[
\max_j\|b_j-a_j\|_2
\le
\frac{C}{\pi_{\min}^\star}
\Delta_{\mathrm{sep}}^{-(2s_{\max}-2)}
h(p_{\theta,T},p_{\theta^\star,T}).
\]

\paragraph{Transition-row discrepancies.}
For the first consecutive pair, define
\[
F_{j\ell}=\pi_jA_{j\ell},
\qquad
F_{j\ell}^\star=\pi_j^\star A_{j\ell}^\star.
\]
These are the fitted and true joint probabilities of the states at times $1$ and $2$. Their row sums equal the initial probabilities:
\[
\begin{aligned}
\sum_{\ell=1}^K F_{j\ell}
&=\pi_j\sum_{\ell=1}^K A_{j\ell}
=\pi_j
\qquad\text{(because transition rows sum to one)},\\
\sum_{\ell=1}^K F_{j\ell}^\star
&=\pi_j^\star\sum_{\ell=1}^K A_{j\ell}^\star
=\pi_j^\star.
\end{aligned}
\]
It follows that the initial-probability discrepancy is bounded by the discrepancy in the corresponding flow row:
\begin{equation}
\label{eq:paired-flow-row-mass}
\begin{aligned}
|\pi_j-\pi_j^\star|
&=
\left|
\sum_{\ell=1}^K(F_{j\ell}-F_{j\ell}^\star)
\right|\\
&\le
\sum_{\ell=1}^K|F_{j\ell}-F_{j\ell}^\star|
\qquad\text{(by the triangle inequality)}.
\end{aligned}
\end{equation}

To pass from joint flows to conditional transition probabilities, isolate the transition-row difference entry by entry:
\[
\begin{aligned}
\pi_j^\star(A_{j\ell}-A_{j\ell}^\star)
&=
\pi_j^\star A_{j\ell}-F_{j\ell}^\star\\
&=
\underbrace{\pi_jA_{j\ell}}_{=F_{j\ell}}
-F_{j\ell}^\star
+(\pi_j^\star-\pi_j)A_{j\ell}\\
&=
F_{j\ell}-F_{j\ell}^\star
+(\pi_j^\star-\pi_j)A_{j\ell}.
\end{aligned}
\]
Taking absolute values gives
\[
\pi_j^\star|A_{j\ell}-A_{j\ell}^\star|
\le
|F_{j\ell}-F_{j\ell}^\star|
+
|\pi_j^\star-\pi_j|A_{j\ell},
\]
where we used $\pi_j^\star>0$ and $A_{j\ell}\ge0$. Summing over $\ell$ yields
\[
\begin{aligned}
\pi_j^\star\|A_{j\cdot}-A_{j\cdot}^\star\|_1
&\le
\sum_{\ell=1}^K|F_{j\ell}-F_{j\ell}^\star|
+
|\pi_j^\star-\pi_j|
\sum_{\ell=1}^K A_{j\ell}\\
&=
\sum_{\ell=1}^K|F_{j\ell}-F_{j\ell}^\star|
+
|\pi_j^\star-\pi_j|
\qquad\text{(because $\sum_\ell A_{j\ell}=1$)}\\
&\le
2\sum_{\ell=1}^K|F_{j\ell}-F_{j\ell}^\star|
\qquad\text{(by Equation~\eqref{eq:paired-flow-row-mass})}.
\end{aligned}
\]
We divide by the positive true initial probability to obtain
\begin{equation}
\label{eq:paired-row-division}
\begin{aligned}
\|A_{j\cdot}-A_{j\cdot}^\star\|_1
&\le
\frac{2}{\pi_j^\star}
\sum_{\ell=1}^K|F_{j\ell}-F_{j\ell}^\star|\\
&\le
\frac{2}{\pi_{\min}^\star}
\sum_{\ell=1}^K|F_{j\ell}-F_{j\ell}^\star|\\
&\le
\frac{2}{\pi_{\min}^\star}
\sum_{i=1}^K\sum_{\ell=1}^K
|F_{i\ell}-F_{i\ell}^\star|\\
&=
\frac{2}{\pi_{\min}^\star}D_{\{1,2\}}
\qquad\text{(because the labels are aligned)}.
\end{aligned}
\end{equation}
In the structured regimes, the two-time probability bound of Theorem~\ref{thm:paired} gives
\[
D_{\{1,2\}}
\le
C\Delta_{\mathrm{sep}}^{-(2s_{\max}-1)}
h(p_{\theta,T},p_{\theta^\star,T}).
\]
Substitution into Equation~\eqref{eq:paired-row-division}, followed by taking the maximum over $j$, proves
\[
\max_j\|A_{j\cdot}-A_{j\cdot}^\star\|_1
\le
\frac{C}{\pi_{\min}^\star}
\Delta_{\mathrm{sep}}^{-(2s_{\max}-1)}
h(p_{\theta,T},p_{\theta^\star,T}).
\]
Without cluster assumptions, the same calculation uses the pair exponent $2K-1$, giving the stated unstructured row bound. The row calculation divides only by a true initial probability. Fitted initial-probability positivity was needed for the maximum center bound, but is not needed for this row-division inequality.

For the MLE, assume the same sufficient condition $\varepsilon_{n,T}\le c(\pi_{\min}^\star)\Delta_{\mathrm{sep}}^{2s_{\max}-1}$ used in Proposition~\ref{prop:mle-matched-rates}. On $\mathcal E_n$, localization supplies the aligned labeling and the fitted initial-probability lower bounds. Applying the two deterministic inequalities to $\widehat\theta_n$ and using $h(p_{\widehat\theta_n,T},p_{\theta^\star,T})\le\varepsilon_{n,T}$ proves their statistical versions on this same event.
\end{proof}

\begin{corollary}[Identifiability from two observations]
\label{cor:identifiability-two}
Under the standing assumptions, two exact-fit Gaussian HMMs with $p_{\theta,2}=p_{\theta^\star,2}$ have the same emission centers, initial distribution, and transition matrix, up to one common permutation.
\end{corollary}

\begin{proof}

Apply Corollary~\ref{cor:localization} with $T=2$ and
$h(p_{\theta,2},p_{\theta^\star,2})=0$ to obtain one common
matching permutation and positive fitted initial probabilities.
After relabeling, Theorem~\ref{thm:paired} yields
$D_\mu=D_{\{1\}}=D_{\{1,2\}}=0$. Since
\[
0\le\sum_{j=1}^K\pi_j\|b_j-a_j\|_2\le D_\mu=0
\quad\text{and}\quad \pi_j>0,
\]
every fitted center equals its matched true center.
Also $D_{\{1\}}=\sum_j|\pi_j-\pi_j^\star|=0$, so
$\pi=\pi^\star$. Finally,
$D_{\{1,2\}}=\sum_{j,\ell}|\pi_jA_{j\ell}
-\pi_j^\star A_{j\ell}^\star|=0$; division by
$\pi_j=\pi_j^\star>0$ gives $A=A^\star$.
All equalities hold under the single permutation fixed at the start.

\end{proof}
\section{Proofs under over-specification}
\label{app:over-specification-proofs}

The displacement and aggregate probability bounds follow from Lemma~\ref{lem:global-recovery-bounds}, established within the proof of Theorem~\ref{thm:globalexact} in Appendix~\ref{app:proof_3_1}. The lemma applies without assuming exact specification. We then establish the one-time transport comparison, prove the path transport bounds, give lower examples for the separation exponents in both regimes, and treat the single true state case.

\subsection{Proof of Theorem~\ref{thm:overglobal}}
\label{proof:global_overift}

\begin{proof}[Proof of Theorem~\ref{thm:overglobal}]
Lemma~\ref{lem:global-recovery-bounds} gives the asserted bounds for $V$, $B$, and $D_I$ for every $k\ge K$. Taking $I=\{1\}$ gives the bound for $D_\pi$, and taking the maximum over consecutive pairs $I=\{t,t+1\}$ gives the bound for $D_A$. The structured separation exponents are $2s_{\max}$ and $\beta_2=\max\{2s_{\max},3s_{\max}-2\}$, respectively. The same lemma gives the unstructured probability exponents $2K$ and $\gamma_{K,2}$, without a lower bound on fitted state probabilities or a factor $(\pi_{\min}^\star)^{-1}$.
\end{proof}

\subsection{Proof of Theorem~\ref{thm:overtransport}}
\label{app:proof_43}

\begin{proof}[Proof of Theorem~\ref{thm:overtransport}]
Fix a nonempty set $I\subseteq[T]$ with $|I|=r$. The bounds of Theorem~\ref{thm:overglobal} control the displacement of fitted centers and the discrepancy between aggregate label probabilities. To convert these bounds into a transport guarantee, we also account for the distances over which the probability discrepancies are redistributed. Lemma~\ref{lem:hierarchical-transport} provides this comparison, while Lemma~\ref{lem:cluster-mass} controls the mass that must be transported between cluster sequences.

\paragraph{Transport under the cluster assumptions.}
Recall that $V=\max_{t\le T}\sum_{i=1}^k q_t(i)r_i^2$ and $V_{2,I}=r^{-1}\sum_{t\in I}\sum_{i=1}^k q_t(i)r_i^2$, where $r=|I|$. Since $V$ is the maximum over all observed times and $V_{2,I}$ is the average over the times in $I$, we have $V_{2,I}\le V$. Therefore:
\[
\begin{aligned}
V_{2,I} =\frac1r\sum_{t\in I}\sum_{i=1}^k q_t(i)r_i^2 
\le\frac1r\sum_{t\in I}V =V.
\end{aligned}
\]
Applying Lemma~\ref{lem:hierarchical-transport} with $p=2$, we obtain
\[
\begin{aligned}
W_{2,I}^2(G_\theta^I,G_\star^I)
&\le
\underbrace{2V_{2,I}}_{\text{displacement contribution}}
+
\underbrace{C_0^2\Delta_{\mathrm{sep}}^2D_I}_{\text{within-cluster cost bound}}
+
\underbrace{4R^2J_I}_{\text{between-cluster cost bound}}.
\end{aligned}
\]
The factor $\Delta_{\mathrm{sep}}^2$ in the second term comes from the squared distance between locations in the same cluster at each retained time. The last term allows transport between different cluster sequences, where distances are bounded by $2R$.

We bound these contributions separately before combining them. Theorem~\ref{thm:overglobal} gives
\[
V_{2,I}\le V\le C\Delta_{\mathrm{sep}}^{-(2s_{\max}-2)}h(p_{\theta,T},p_{\theta^\star,T})
\]

and
\[
\Delta_{\mathrm{sep}}^2D_I
\le C_r\Delta_{\mathrm{sep}}^{-(\beta_r-2)}h(p_{\theta,T},p_{\theta^\star,T}).
\]
The squared within-cluster distance reduces the probability exponent by two. Lemma~\ref{lem:cluster-mass} also gives
\[
J_I\le C_r\Delta_{\mathrm{sep}}^{-(2s_{\max}-2)}h(p_{\theta,T},p_{\theta^\star,T}).
\]
Since $\beta_r\ge2s_{\max}$ and $0<\Delta_{\mathrm{sep}}\le1$, both the displacement and between-cluster contributions are bounded by the multiplier $\Delta_{\mathrm{sep}}^{-(\beta_r-2)}$. Adding the three contributions proves
\begin{equation}
W_{2,I}^2(G_\theta^I,G_\star^I)
\le C_r\Delta_{\mathrm{sep}}^{-(\beta_r-2)}h(p_{\theta,T},p_{\theta^\star,T}).
\label{eq:proof-overtransport-structured}
\end{equation}
For $r=1$, $\beta_1=2s_{\max}$, so the one-time squared Wasserstein exponent remains $2s_{\max}-2$. For $s_{\max}\ge2$ and $r\ge2$, the exponent is $r(s_{\max}-1)+s_{\max}-2$.

\paragraph{Transport without cluster assumptions.}
For any nonempty $I$, first move fitted locations to their assigned true centers, then couple the resulting distribution with $G_\star^I$ on the true grid. The first step has squared cost $V_{2,I}$. For the second step, the two grid distributions have $\ell_1$ discrepancy $D_I$, so a coupling can leave their common mass fixed and move the remaining mass, which is $D_I/2$. The grid diameter under $d_{2,r}$ is $D_\star=\max_{i,j}\|a_i-a_j\|_2$. The squared triangle inequality therefore gives
\[
W_{2,I}^2(G_\theta^I,G_\star^I)
\le 2V_{2,I}+D_\star^2D_I.
\]
Since $V_{2,I}\le V$, Lemma~\ref{lem:global-recovery-bounds} and~\eqref{eq:overfit-diameter-probability} yield
\[
\begin{aligned}
W_{2,I}^2(G_\theta^I,G_\star^I)
&\le C_r\left\{\Delta_{\mathrm{sep}}^{-(2K-2)}
+\Delta_{\mathrm{sep}}^{-(\gamma_{K,r}-2)}\right\}
h(p_{\theta,T},p_{\theta^\star,T}).
\end{aligned}
\]
Here the factor $D_\star^2$ cancels the inverse diameter factor in the probability bound. Using $\gamma_{K,r}\ge2K$ and $0<\Delta_{\mathrm{sep}}\le1$, we obtain
\begin{equation}
W_{2,I}^2(G_\theta^I,G_\star^I)
\le C_r\Delta_{\mathrm{sep}}^{-(\gamma_{K,r}-2)}
h(p_{\theta,T},p_{\theta^\star,T}).
\label{eq:proof-overtransport-unstructured}
\end{equation}
For $I=\{t\}$, we have $\gamma_{K,1}=2K$. The constant is uniform in $t$, so taking the maximum gives
\begin{equation}
\max_{t\le T}W_2^2(G_\theta^{\{t\}},G_\star^{\{t\}})
\le C\Delta_{\mathrm{sep}}^{-(2K-2)}
h(p_{\theta,T},p_{\theta^\star,T}).
\label{eq:proof-overtransport-unstructured-marginal}
\end{equation}

\paragraph{Transport bounds for the MLE.}
The preceding inequalities hold for every fitted parameter satisfying the theorem assumptions. We may therefore apply them to $\widehat\theta_n$ on the density event
\[
\mathcal E_n
=
\left\{
h(p_{\widehat\theta_n,T},p_{\theta^\star,T})
\le\varepsilon_{n,T}
\right\},
\qquad
\mathbb P_{\theta^\star}(\mathcal E_n)\ge 1-\frac{1}{n},
\]
where the probability bound follows from Proposition~\ref{prop:density}.

For structured configurations, taking square roots in \eqref{eq:proof-overtransport-structured} gives, on $\mathcal E_n$,
\[
\begin{aligned}
W_{2,I}(G_{\widehat\theta_n}^I,G_\star^I)
&\le
C_r\Delta_{\mathrm{sep}}^{-(\beta_r-2)/2}
h(p_{\widehat\theta_n,T},p_{\theta^\star,T})^{1/2}\\
&\le
C_r\Delta_{\mathrm{sep}}^{-(\beta_r-2)/2}
\varepsilon_{n,T}^{1/2}
\qquad \text{(by the definition of $\mathcal E_n$)}.
\end{aligned}
\]
In particular, for $I=\{t\}$,
\[
W_2(G_{\widehat\theta_n}^{\{t\}},G_\star^{\{t\}})
\le C\Delta_{\mathrm{sep}}^{-(s_{\max}-1)}
\varepsilon_{n,T}^{1/2},
\qquad t\le T.
\]
The constant is uniform in $t$, and each inequality holds on the same event $\mathcal E_n$. Taking the maximum therefore yields
\[
\max_{t\le T}
W_2(G_{\widehat\theta_n}^{\{t\}},G_\star^{\{t\}})
\le C\Delta_{\mathrm{sep}}^{-(s_{\max}-1)}
\varepsilon_{n,T}^{1/2}.
\]

Without cluster assumptions, Equation~\eqref{eq:proof-overtransport-unstructured-marginal} gives
\[
\max_{t\le T}W_2(G_{\widehat\theta_n}^{\{t\}},G_\star^{\{t\}})
\le C\Delta_{\mathrm{sep}}^{-(K-1)}\varepsilon_{n,T}^{1/2}.
\]
For $r\ge2$, taking square roots in~\eqref{eq:proof-overtransport-unstructured} gives

\[
W_{2,I}(G_{\widehat\theta_n}^I,G_\star^I)
\le C_r\Delta_{\mathrm{sep}}^{-(\gamma_{K,r}-2)/2}
\varepsilon_{n,T}^{1/2}
\]
for the joint latent location distributions. Since the deterministic inequalities apply to every nonempty $I\subseteq[T]$, all these MLE bounds hold simultaneously on $\mathcal E_n$.

By Proposition~\ref{prop:density},
\[
\varepsilon_{n,T}^{1/2}
=
\sqrt{C}\left(\frac{d\log(nTd)}{n}\right)^{1/4}.
\]
Substituting this expression into the preceding inequalities gives the stated Wasserstein rates. For fixed model dimensions, emission separation, and cluster parameters, these rates are of order $(\log n/n)^{1/4}$.
\end{proof}

\subsection{Sharpness of the structured overfitted exponents}
\label{app:overfit-sharpness}

The upper bounds under over-specification allow repeated fitted centers and zero fitted initial probabilities. The following lower examples use these possibilities. They establish sharpness over this full class; the required fitted state count for joint losses can exceed $K+1$.

\begin{proposition}[Sharpness under over-specification]
\label{prop:overfit-sharpness}
Fix $s\ge2$ and $d\ge1$. There are true $K=s$ state models with a single collinear cluster, $\Delta_{\mathrm{sep}}=\delta\downarrow0$, and true initial probabilities bounded away from zero, with the following properties.

For one-time recovery, a fitted model with $k=s+1$ states can satisfy, for a positive sequence $\eta_\delta\downarrow0$ and every fixed $T\ge2$,
\begin{align}
h(p_{\theta_\delta,T},p_{\theta_\delta^\star,T})&\asymp\eta_\delta\delta^{2s},\label{eq:over-sharp-one-density}\\
D_\pi\asymp\max_{t\le T}D_{\{t\}}&\asymp\eta_\delta,\qquad
V\asymp\max_{t\le T}W_2^2(G_{\theta_\delta}^{\{t\}},G_\star^{\{t\}})\asymp\eta_\delta\delta^2,\label{eq:over-sharp-one-loss}\\
B&\asymp\eta_\delta\delta,\qquad D_\mu\asymp\sqrt{\eta_\delta}\,\delta.
\label{eq:over-sharp-displacements}
\end{align}
For every fixed $r\ge2$ and $T\ge r$, a fitted model with a fixed number of states depending only on $s,r$ can satisfy, for $I=\{1,\ldots,r\}$,
\begin{align}
h(p_{\theta_\delta,T},p_{\theta_\delta^\star,T})&\asymp\eta_\delta\delta^{r(s-1)+s},\label{eq:over-sharp-joint-density}\\
D_I&\asymp\eta_\delta,\qquad
W_{2,I}^2(G_{\theta_\delta}^I,G_\star^I)\asymp\eta_\delta\delta^2.
\label{eq:over-sharp-joint-loss}
\end{align}
Consequently, the structured exponents $\beta_r=\max\{2s,r(s-1)+s\}$ for $D_I$ and $\beta_r-2$ for $W_{2,I}^2$ cannot be reduced uniformly over the allowed over-specified models. The one-time displacement exponents are also sharp. Taking $s=K$ proves the unstructured exponents $\gamma_{K,r}$ and $\gamma_{K,r}-2$ as well. These conclusions concern deterministic separation powers with the displayed Hellinger powers fixed, not minimax rates in $n$ or sharpness for each prescribed fitted state count.
\end{proposition}

\begin{proof}
Work first in dimension one with fixed nodes $z_1<\cdots<z_s$ of minimum gap one, and put $g(z)=\prod_{j=1}^s(z-z_j)$. All locations below are multiplied by $\delta$. We use perturbations of size $\eta=\eta_\delta$, chosen sufficiently small relative to $\delta$.

\paragraph{A uniform Gaussian comparison for the perturbations.}
We first give the density comparison with both small parameters explicit. In each construction below, the unscaled locations remain in a fixed compact set for $0\le\eta\le\eta_0$, and all coefficients are fixed independently of $\eta,\delta$. Let $q_{\eta,\delta}$ be the resulting $r$-coordinate observation density, and write $\Lambda$ for the first variation of its unscaled mixing distribution at $\eta=0$. In the one-time construction take $r=1$ in this paragraph. Suppose $\Lambda$ annihilates every polynomial of total degree less than $m$ and is nonzero on a monomial $z^\alpha$ of degree $m$.

For fixed $x$, apply $\Lambda$ to $z\mapsto\varphi_r(x-\delta z)$. Taylor expansion in $\delta$ through degree $m-1$ has polynomial coefficients in $z$ of the same total degrees. All of them vanish under $\Lambda$. Gaussian derivative formulas and the integral remainder consequently give, for a fixed integer $L$,
\[
|\dot q_{0,\delta}(x)|
\le C\delta^m(1+\|x\|_2^L)\varphi_r(x)e^{C\delta\|x\|_2},
\qquad 0<\delta\le\delta_0.
\]
This bound also holds when $\Lambda$ includes derivatives at grid points: the finite point-evaluation and derivative functionals commute with Taylor expansion, and their orders and coefficients are fixed.

Ordinary shifts and weight perturbations have second-order remainders bounded uniformly in $\eta,\delta$. For the symmetric split, the same fact follows without differentiating $\sqrt\eta$ at zero. For any smooth scalar $f$, fixed $z,\alpha$, and small $\eta\ge0$, the central Taylor formula gives
\[
\frac{f(z+\alpha\eta+\sqrt\eta)+f(z+\alpha\eta-\sqrt\eta)}2
=f(z)+\eta\{\alpha f'(z)+\tfrac12f''(z)\}
 +O\!\left(\eta^2\max_{1\le j\le4}\|f^{(j)}\|_\infty\right),
\]
where the norms are over one fixed compact interval containing all the evaluation points. This follows by first expanding the central difference to order two with a fourth-derivative remainder, and then expanding about $z$. Applying it to the Gaussian emission and multiplying by the perturbed weight shows that, after increasing $L$ if necessary,
\[
|q_{\eta,\delta}(x)-q_{0,\delta}(x)-\eta\dot q_{0,\delta}(x)|
\le C\eta^2(1+\|x\|_2^L)\varphi_r(x)e^{C\delta\|x\|_2}.
\]
All these constants are uniform for $0\le\eta\le\eta_0$ and $0<\delta\le\delta_0$.
The true grid density satisfies
\[
q_{0,\delta}(x)\ge c\varphi_r(x)
                    e^{-C\delta\|x\|_2-C\delta^2}.
\]
Therefore the inequality $(\sqrt q-\sqrt p)^2\le(q-p)^2/p$ and Gaussian integrability imply
\begin{equation}
\begin{aligned}
h^2(q_{\eta,\delta},q_{0,\delta})
&\le\frac12\int\frac{(q_{\eta,\delta}-q_{0,\delta})^2}{q_{0,\delta}}\,dx
 \le C(\eta\delta^m+\eta^2)^2,\\
h(q_{\eta,\delta},q_{0,\delta})&\le C(\eta\delta^m+\eta^2).
\end{aligned}
\label{eq:over-sharp-density-comparison}
\end{equation}
For the chosen monomial, scaling the mixing locations gives the exact leading moment difference
\[
\int x^\alpha\,d(H_{\eta,\delta}-H_{0,\delta})
=\delta^m\{\eta\Lambda(z^\alpha)+O(\eta^2)\}.
\]
Its coefficient is nonzero. More explicitly, write $\lambda=\Lambda(z^\alpha)\ne0$. The remainder in this moment expansion is bounded by $C\eta^2\delta^m$ with $C$ independent of $\eta,\delta$: the unscaled nodes and perturbation coefficients lie in a fixed compact set. Decrease $\eta_0$ so that $C\eta_0\le|\lambda|/2$. Then the moment difference has absolute value at least $|\lambda|\eta\delta^m/2$ throughout $0<\eta\le\eta_0$. Lemma~\ref{lem:moment-extraction}, applied to this monomial on a fixed ball rather than a shrinking support, gives a uniform lower bound $h(q_{\eta,\delta},q_{0,\delta})\ge c\eta\delta^m$. Combined with equation~\eqref{eq:over-sharp-density-comparison}, this proves
\begin{equation}
 c\eta\delta^m\le h(q_{\eta,\delta},q_{0,\delta})
             \le C\eta\delta^m,
 \quad 0<\delta\le\delta_0,\quad
 0<\eta\le\min\{\eta_0,\delta^m\}.
 \label{eq:over-sharp-uniform-wedge}
\end{equation}
The constants are uniform over this two-parameter region. Choose $\delta_0\le1$ and $\eta_\delta=\delta^{m+1}$, decreasing $\delta_0$ further to ensure $\eta_\delta\le\eta_0$. Then
\[
\frac{\eta_\delta^2}{\eta_\delta\delta^m}=\delta\longrightarrow0,
\]
so the upper and lower bounds have the common order $\eta_\delta\delta^m$. This is a joint choice of the two parameters, not an exchange of iterated limits. In dimension $d>1$, putting the centers on the first coordinate axis appends the same independent Gaussian factor to both models and leaves their Hellinger distance unchanged.

\paragraph{One-time probabilities and displacement.}
Choose fixed positive weights $w_j$ summing to one. There are unique coefficients $u_j,c_j$ such that
\[
\Lambda(f)=\sum_{j=1}^s\{u_j f(z_j)+c_j f'(z_j)\}
+\frac{w_1}{2}f''(z_1)
\]
vanishes for every polynomial of degree at most $2s-1$. Existence and uniqueness follow from the invertibility of the confluent Vandermonde matrix, equivalently from the fact that a polynomial of degree at most $2s-1$ with $s$ double roots is zero. In particular, $\sum_j u_j=0$.

Neither vector $u$ nor $c$ is zero. If $u=0$, take an antiderivative $f$ of $g$, whose degree is $s+1\le2s-1$. Then $f'(z_j)=0$ for all $j$ but $f''(z_1)=g'(z_1)\ne0$, a contradiction. If $c=0$, use $f(z)=(z-z_1)g(z)$, which vanishes on the grid but has $f''(z_1)=2g'(z_1)\ne0$. Also
\[
\Lambda(g^2)=w_1g'(z_1)^2>0.
\]
Since $g^2$ is monic of degree $2s$, the first nonzero moment order of $\Lambda$ is $m=2s$.

Give group $j$ total weight $w_j+\eta u_j$ and center $z_j+\eta c_j/w_j$. For $j=1$, split this group into two equal weights at
\[
z_1+\eta c_1/w_1-\sqrt\eta,
\qquad z_1+\eta c_1/w_1+\sqrt\eta.
\]
The other groups each contain one atom. Choose $\eta_0>0$ so that $w_j+\eta u_j\ge w_j/2$ for every $j$ and $\eta\max_j|c_j|/w_j+\sqrt\eta<1/4$ for $0<\eta\le\eta_0$. The weights then sum to one, and every fitted atom is uniquely assigned to its indicated true center. The choice of $\eta_0$ depends only on the fixed nodes, weights, and coefficients, not on $\delta$. Expanding a smooth test function verifies that the first variation is exactly $\Lambda$. Thus~\eqref{eq:over-sharp-density-comparison}, with $m=2s$ and $\eta_\delta=\delta^{2s+1}$, gives the required density order for one observation. Take the initial distributions and every transition row to be the corresponding mixture weights. Both hidden processes are then independent over time, and the product-affinity identity preserves this Hellinger order at every fixed $T$.

The aggregate mass difference is exactly $\eta u$, so its $\ell_1$ norm is of order $\eta$. The split contributes order $\eta\delta^2$ to $V$ and order $\sqrt\eta\,\delta$ to $D_\mu$. The within-group signed displacements equal
$\delta\eta c_j+O(\delta\eta^2)$,
giving $B\asymp\eta\delta$ because $c\ne0$. Every coupling to the true support has squared cost at least $V$. Conversely, move fitted atoms to their assigned centers and redistribute the $O(\eta)$ mass discrepancy over distances $O(\delta)$. This gives $W_2^2\asymp\eta\delta^2$ and proves~\eqref{eq:over-sharp-one-density}--\eqref{eq:over-sharp-displacements}.

\paragraph{Joint probabilities and an HMM realization.}
Define two linear functionals on univariate polynomials by
\[
\mathcal A(f)=\sum_j\frac{f(z_j)}{g'(z_j)},\qquad
\mathcal B(f)=\sum_j\left\{
-\frac{g''(z_j)}{g'(z_j)^3}f(z_j)
+\frac{f'(z_j)}{g'(z_j)^2}\right\}.
\]
Ordinary Lagrange interpolation shows that $\mathcal A$ extracts the leading coefficient on polynomials of degree at most $s-1$. Hermite interpolation shows that $\mathcal B$ extracts the leading coefficient on polynomials of degree at most $2s-1$: the leading coefficients of the value and derivative cardinal polynomials are respectively $-g''(z_j)/g'(z_j)^3$ and $1/g'(z_j)^2$. Consequently,
\[
\Lambda_r=\mathcal B\otimes\mathcal A^{\otimes(r-1)}
\]
vanishes on every monomial of total degree below
$m=(2s-1)+(r-1)(s-1)=r(s-1)+s$,
and equals one on $x_1^{2s-1}\prod_{\ell=2}^r x_\ell^{s-1}$. This functional contains only values and derivatives in the first coordinate; it contains no mixed derivatives.

Write its coefficients at the true grid as
\[
\Lambda_r(f)=\sum_{\mathbf j\in[s]^r}
\{U_{\mathbf j}f(z_{j_1},\ldots,z_{j_r})
+C_{\mathbf j}\partial_1f(z_{j_1},\ldots,z_{j_r})\}.
\]
The value vector $U$ is nonzero. Indeed, $g''(z_s)/g'(z_s)=2\sum_{j<s}(z_s-z_j)^{-1}>0$, and every coefficient of $\mathcal A$ is nonzero. Also $\sum_{\mathbf j}U_{\mathbf j}=0$. Start with the product grid weights $w_{\mathbf j}=s^{-r}$. Change them to $w_{\mathbf j}+\eta U_{\mathbf j}$, move the first coordinate of grid point $\mathbf j$ by $\eta C_{\mathbf j}/w_{\mathbf j}$, and leave the other coordinates fixed. Choose $\eta_0$ so that $\eta_0\max_{\mathbf j}|U_{\mathbf j}|<s^{-r}/2$ and $\eta_0\max_{\mathbf j}|C_{\mathbf j}|/w_{\mathbf j}<1/4$. For $0<\eta\le\eta_0$, this gives positive normalized weights, unique indicated nearest-center assignments, and the prescribed first variation. Both inequalities involve only $s,r$ and the fixed grid. Equation~\eqref{eq:over-sharp-density-comparison}, with $\eta_\delta=\delta^{m+1}$, gives~\eqref{eq:over-sharp-joint-density} for the first $r$ observations.

These joint mixing measures are realized by time-homogeneous HMMs. The true model has initial probabilities and transition rows all equal to $(1/s,\ldots,1/s)$. For the fitted model, introduce states $(\mathbf j,\ell)$, with $\mathbf j\in[s]^r$ and $1\le\ell\le r$. Put initial mass $w_{\mathbf j}+\eta U_{\mathbf j}$ on $(\mathbf j,1)$ and zero on the other states. Move deterministically from $(\mathbf j,\ell)$ to $(\mathbf j,\ell+1)$ for $\ell<r$, and use the perturbed grid coordinate as that state's emission center. After time $r$, transition to an additional set of $s$ states with the true emission centers. From every phase-$r$ state, enter these additional states with probabilities $(1/s,\ldots,1/s)$, and give each additional state this same transition row. Thus $k=rs^r+s$ suffices. The remaining observations are independent true draws, so they do not change the Hellinger distance for any $T\ge r$.

Nearest-center aggregation gives $D_I=\eta\sum_{\mathbf j}|U_{\mathbf j}|\asymp\eta$. In any coupling, at least $D_I/2$ of the mass crosses different assigned grid labels. The assigned-label argument in Lemma~\ref{lem:transport-assignment-matching} applies without any Markov assumption on the grouped process: the mass on different labels is at least their total variation, $D_I/2$. On each such pair, at least one true coordinate differs from the fitted point's assigned coordinate by at least $\delta$, whereas its perturbation is less than $\delta/4$. Thus the physical distance in that coordinate is at least $3\delta/4$, and in particular its normalized squared cost is at least $\delta^2/(4r)$. This proves $W_{2,I}^2\ge c_r\eta\delta^2$. Moving fitted grid points back to their assigned points costs $O(\eta^2\delta^2)$, and the remaining mass discrepancy costs $O(\eta\delta^2)$. This proves the reverse inequality and~\eqref{eq:over-sharp-joint-loss}. Since $m=\beta_r$ for $s\ge2,r\ge2$, the claimed joint powers are necessary. For $r=2$, this also proves sharpness of the transition-flow exponent.
\end{proof}

To add separated clusters, give the moving cluster a fixed total initial mass $\alpha\in(0,1)$, keep both its true and fitted chains within that cluster, and give identical additional absorbing states the remaining positive initial masses in both models. Every true initial mass is then bounded away from zero along the separation sequence. The observation laws have the form $\alpha p_\eta+(1-\alpha)q$ and $\alpha p_0+(1-\alpha)q$ with the same $q$. Convexity of squared Hellinger distance gives the required upper order from the within-cluster comparison. The fixed polynomial moment used above has difference exactly $\alpha$ times its within-cluster difference, giving the matching lower order by moment extraction. Probability discrepancies and displacement quantities acquire the same fixed factor. For transport, couple the additional identical paths to themselves for an upper bound; the preceding nearest-support or assigned-label arguments give the lower bound. Thus all displayed orders survive, with constants depending on $\alpha$, and establish the claimed structured worst-case powers when $s=s_{\max}$.

\subsection{Emission-Center Recovery with One True State}
\label{app:proof_remark_43}
The main text assumes $K\ge2$ because the paper studies separation among distinct true emission centers. When $K=1$, there is no separation parameter. We can instead bound displacement directly from the sole true center, using a quadratic polynomial for an overfitted model and a linear polynomial under exact specification.

\begin{remark}[Recovery when $K=1$]
\label{remark:single_true_state}
Suppose $K=1$, with true emission center $a$. For any $\theta\in\Theta_k(R)$ and nonempty $I\subseteq[T]$ with $|I|=r$,
\[
W_{2,I}^2\bigl(G_\theta^I,\delta_{(a,\ldots,a)}\bigr)\le C h(p_{\theta,T},p_{\theta^\star,T}).
\]
Under exact specification, $k=1$, and
\[
\|\mu_1-a\|_2\le C h(p_{\theta,T},p_{\theta^\star,T}).
\]
\end{remark}

\begin{proof}
The true initial probability and transition probability are both one, so the true latent location equals $a$ at every time. Thus $G_\star^{\{t\}}=\delta_a$ and $G_\star^I=\delta_{(a,\ldots,a)}$. Every fitted state is assigned to this same center, meaning that $c(i)=1$ for all $i\in[k]$. After aggregation, the fitted and true label paths therefore both equal $(1,\ldots,1)$ with probability one. All aggregate probability and flow discrepancies are zero; we only need to bound displacement from $a$.

\medskip
\paragraph{Quadratic displacement.}
Fix $t\le T$. With $q_t(i)=\mathbb P_\theta(Z_t=i)$, the squared displacement at time $t$ is $V_t=\sum_{i=1}^k q_t(i)\|\mu_i-a\|_2^2$. We use $P_1(x)=\|x-a\|_2^2$ because its fitted expectation is exactly $V_t$ and its true expectation is zero:
\[
\begin{aligned}
\int P_1\,d\bigl(G_\theta^{\{t\}}-\delta_a\bigr)
&=\sum_{i=1}^k q_t(i)P_1(\mu_i)-P_1(a)\\
&=\sum_{i=1}^k q_t(i)\|\mu_i-a\|_2^2
=V_t.
\end{aligned}
\]
Each summand is nonnegative, so displacements of different fitted centers cannot cancel in this expectation.

The polynomial has degree two, and $\|P_1\|_{\infty,B_d(0,R)}\le4R^2$ because $\|x-a\|_2\le2R$ on the ball. Let $p_\theta^{\{t\}}$ and $p_\star^{\{t\}}$ be the one-time observation densities. Lemma~\ref{lem:moment-extraction} bounds the polynomial expectation by their Hellinger distance, and Lemma~\ref{lem:marginalization} bounds this distance by the full-block error. Hence
\begin{equation}
\begin{aligned}
V_t
&\le C\|P_1\|_\infty h\bigl(p_\theta^{\{t\}},p_\star^{\{t\}}\bigr)\\
&\le C h(p_{\theta,T},p_{\theta^\star,T}),
\qquad t\le T.
\end{aligned}
\label{eq:single-state-Vt}
\end{equation}
Here $4R^2$ is absorbed into $C$. As in the preceding proofs, we apply the moment bound on the span of the true and fitted centers. This span has dimension at most $k+1$, so $C$ may depend on $k$ and $R$ but not on the ambient dimension. The constant is uniform in $t$, since the polynomial and its norm bound are the same at every time. The argument also allows zero fitted probabilities and repeated fitted centers.

\medskip
\paragraph{Transport of the location block.}
Fix a nonempty $I\subseteq[T]$ with $|I|=r$. The true location block is the fixed vector $(a,\ldots,a)$. Any coupling with its point mass must therefore pair every fitted block with this vector. There is no choice of destination in the transport problem, and its cost is
\[
\begin{aligned}
W_{2,I}^2\bigl(G_\theta^I,\delta_{(a,\ldots,a)}\bigr)
&=\mathbb E_\theta\left[\frac1r\sum_{t\in I}\|\mu_{Z_t}-a\|_2^2\right]\\
&=\frac1r\sum_{t\in I}\sum_{i=1}^k q_t(i)\|\mu_i-a\|_2^2\\
&=\frac1r\sum_{t\in I}V_t\\
&\le C h(p_{\theta,T},p_{\theta^\star,T}),
\end{aligned}
\]
by~\eqref{eq:single-state-Vt}. The second equality follows from linearity of expectation and the marginal distribution of $Z_t$; no independence across times is needed. Since the normalized cost averages the one-time squared displacements, the constant does not grow with $r$. This proves the first bound.
\paragraph{Exact specification.}
If $k=1$, the fitted location is also deterministic and equals $\mu_1$ at every time. When $\mu_1=a$, the second bound is immediate. Otherwise, take $u=(\mu_1-a)/\|\mu_1-a\|_2$ and $U(x)=\langle u,x-a\rangle$. Choosing the direction of the displacement gives
\[
U(a)=0,
\qquad
U(\mu_1)=\left\langle\frac{\mu_1-a}{\|\mu_1-a\|_2},\mu_1-a\right\rangle
=\|\mu_1-a\|_2.
\]
Thus, unlike the quadratic polynomial, this linear polynomial recovers the displacement itself:
\[
\int U\,d(\delta_{\mu_1}-\delta_a)=\|\mu_1-a\|_2.
\]
Its supremum norm is at most $2R$, since $|U(x)|\le\|u\|_2\|x-a\|_2\le2R$. Lemmas~\ref{lem:moment-extraction} and~\ref{lem:marginalization} now give
\[
\begin{aligned}
\|\mu_1-a\|_2
&\le C\|U\|_\infty h\bigl(p_\theta^{\{t\}},p_\star^{\{t\}}\bigr)\\
&\le C h(p_{\theta,T},p_{\theta^\star,T}),
\end{aligned}
\]
which proves the second bound. In this case, both location blocks are deterministic. Each of their $r$ coordinates has squared displacement $\|\mu_1-a\|_2^2$, so the normalized path cost gives
\[
W_{2,I}\bigl(\delta_{(\mu_1,\ldots,\mu_1)},\delta_{(a,\ldots,a)}\bigr)=\left\{\frac{r\|\mu_1-a\|_2^2}{r}\right\}^{1/2}=\|\mu_1-a\|_2.
\]
Thus the same linear bound also holds for the location-block Wasserstein distance.
\end{proof}

The proof of Proposition~\ref{prop:density} uses only compactness, Gaussian emissions, and inclusion of the true density in the fitted class; it does not use $K\ge2$ or distinctness of two true centers. It therefore applies unchanged when $K=1$. For the MLE, it bounds the Hellinger error by $\varepsilon_{n,T}$ on $\mathcal E_n$. Taking square roots in the first bound of Remark~\ref{remark:single_true_state} gives $W_{2,I}\le C\varepsilon_{n,T}^{1/2}$. Under exact specification, the second bound instead gives $\|\widehat\mu_1-a\|_2\le C\varepsilon_{n,T}$. For fixed model dimensions and radius, the resulting upper convergence rates are $n^{-1/4}$ and $n^{-1/2}$, respectively, up to logarithmic factors.

\subsection{State splitting and aggregate recovery under over-specification}
\label{app:state-splitting}
This appendix verifies Example~\ref{ex:state-splitting} in Section~\ref{sec:over_specification}. The construction explains why the over-specification results are formulated for aggregate quantities. Individual fitted states can split a single true state while leaving the joint observation density unchanged at every sequence length, so their individual probabilities and transition rows are not identifiable from observations alone. Theorem~\ref{thm:overglobal} therefore controls aggregate state probabilities and transition flows, while Theorem~\ref{thm:overtransport} compares the induced emission-center distributions and finite-block latent location laws without requiring a one-to-one correspondence between fitted and true states.

\begin{proof}
For each true state $j$, choose $m_j\ge1$ substates $(j,1),\ldots,(j,m_j)$, all with emission center $a_j$. Write $\widetilde Z_t$ for the fitted substate at time $t$, and let $Y_t=j$ whenever $\widetilde Z_t=(j,u)$ for some $u$. Thus $Y_t$ records the true-state group containing the fitted substate. We first verify the fitted initial and transition probabilities, then show that summing over substates recovers the true path probabilities. Since all substates in a group have the same emission density, this will also give equality of the joint observation densities.

\paragraph{Initial and transition probabilities.}
Choose nonnegative fitted initial probabilities satisfying
\[
\sum_{u=1}^{m_j}\widetilde\pi_{j,u}=\pi_j^\star,
\qquad j\in[K].
\]
For each source substate $(j,u)$ and destination group $\ell$, choose a probability vector $r_{j\ell,u}$ on $\{1,\ldots,m_\ell\}$ and set
\[
\widetilde A_{(j,u),(\ell,v)}=A_{j\ell}^\star r_{j\ell,u}(v),
\qquad
r_{j\ell,u}(v)\ge0,
\qquad
\sum_{v=1}^{m_\ell}r_{j\ell,u}(v)=1.
\]
The initial probabilities sum to $\sum_j\pi_j^\star=1$. Every transition entry is nonnegative, and each fitted row sums to one:
\[
\begin{aligned}
\sum_{\ell=1}^K\sum_{v=1}^{m_\ell}\widetilde A_{(j,u),(\ell,v)}
&=\sum_{\ell=1}^K A_{j\ell}^\star\sum_{v=1}^{m_\ell}r_{j\ell,u}(v)=\sum_{\ell=1}^K A_{j\ell}^\star
=1.
\end{aligned}
\]
Thus these parameters define a valid HMM. More specifically, the total transition probability from any substate $(j,u)$ into group $\ell$ is
\begin{equation}
\sum_{v=1}^{m_\ell}\widetilde A_{(j,u),(\ell,v)}
=A_{j\ell}^\star.
\label{eq:app-splitting-group-transition}
\end{equation}
The allocation among destination substates may depend on $u$, but the total probability of entering their group does not.

\paragraph{Aggregate path probabilities.}
Fix a horizon $L\ge1$ and a group-label path $z=(z_1,\ldots,z_L)\in[K]^L$. At the first time,
\[
\mathbb P_\theta(Y_1=z_1)=\sum_{u=1}^{m_{z_1}}\widetilde\pi_{z_1,u}=\pi_{z_1}^\star.
\]
For $1\le t<L$, partition the event $\{Y_{1:t}=z_{1:t}\}$ according to the current substate $(z_t,u)$, and sum the next transition over all substates of group $z_{t+1}$. The Markov property of $\widetilde Z_t$ gives
\[
\begin{aligned}
\mathbb P_\theta(Y_{1:t+1}=z_{1:t+1})
&=\sum_{u=1}^{m_{z_t}}\mathbb P_\theta\{Y_{1:t}=z_{1:t},\widetilde Z_t=(z_t,u)\}
\sum_{v=1}^{m_{z_{t+1}}}\widetilde A_{(z_t,u),(z_{t+1},v)}\\
&=A_{z_tz_{t+1}}^\star
\sum_{u=1}^{m_{z_t}}\mathbb P_\theta\{Y_{1:t}=z_{1:t},\widetilde Z_t=(z_t,u)\}\\
&=A_{z_tz_{t+1}}^\star\mathbb P_\theta(Y_{1:t}=z_{1:t}).
\end{aligned}
\]
The second equality uses~\eqref{eq:app-splitting-group-transition}; the third uses the fact that the current-substate events partition the specified group path. Starting from the initial probability and applying this recursion for $t=1,\ldots,L-1$ yields
\begin{equation}
\mathbb P_\theta(Y_{1:L}=z)
=\pi_{z_1}^\star\prod_{t=1}^{L-1}A_{z_tz_{t+1}}^\star
=\mathbb P_\star(Z_{1:L}=z).
\label{eq:app-splitting-path-law}
\end{equation}
For $L=1$, the product is empty and equals one. The calculation also covers zero-probability paths, since it never divides by a path probability or conditions on a specified group path.

Equation~\eqref{eq:app-splitting-path-law} proves that the fitted and true models assign the same probability to every group-label sequence. Summing over coordinates not retained in a given marginal therefore gives equality of all aggregate state probabilities and joint probabilities, including the transition flows.

\paragraph{Observation densities.}
Consider a fitted substate path $((z_1,u_1),\ldots,(z_L,u_L))$ whose group labels are $z$. Every substate $(z_t,u_t)$ has center $a_{z_t}$, so the observation density associated with this path is $\prod_{t=1}^L\phi_d(x_t-a_{z_t})$. This product does not depend on $u_1,\ldots,u_L$.

In the finite-path expansion of the fitted observation density, we may therefore group together all substate paths with the same group labels and factor out their common emission density:
\[
\begin{aligned}
p_{\theta,L}(x_{1:L})
&=\sum_{z\in[K]^L}
\left\{
\sum_{u_1=1}^{m_{z_1}}\cdots\sum_{u_L=1}^{m_{z_L}}
\widetilde\pi_{z_1,u_1}
\prod_{t=1}^{L-1}\widetilde A_{(z_t,u_t),(z_{t+1},u_{t+1})}
\right\}
\prod_{t=1}^L\phi_d(x_t-a_{z_t})\\
&=\sum_{z\in[K]^L}\mathbb P_\theta(Y_{1:L}=z)
\prod_{t=1}^L\phi_d(x_t-a_{z_t})\\
&=\sum_{z\in[K]^L}\pi_{z_1}^\star
\prod_{t=1}^{L-1}A_{z_tz_{t+1}}^\star
\prod_{t=1}^L\phi_d(x_t-a_{z_t})\\
&=p_{\theta^\star,L}(x_{1:L}).
\end{aligned}
\]
The expression in braces sums the probabilities of all fitted paths having group labels $z$, which explains the second equality. The third follows from~\eqref{eq:app-splitting-path-law}, and the last is the true HMM path expansion. Since $L$ and $x_{1:L}$ were arbitrary, the joint observation densities agree for every sequence length. This holds for every admissible allocation of the initial probabilities and every choice of the vectors $r_{j\ell,u}$.

\paragraph{Different fitted parameters with identical observation densities.}
To see that the construction changes transition parameters, choose $j,\ell$ with $A_{j\ell}^\star>0$ and give state $\ell$ at least two substates. Fix a source substate $(j,u)$ and vary only its allocation into group $\ell$:
\[
r_{j\ell,u}(1)=\lambda,
\qquad
r_{j\ell,u}(2)=1-\lambda,
\qquad
r_{j\ell,u}(v)=0\quad(v\ge3),
\qquad 0<\lambda<1.
\]
The corresponding transition entries are $A_{j\ell}^\star\lambda$ and $A_{j\ell}^\star(1-\lambda)$. They change with $\lambda$, while their sum remains $A_{j\ell}^\star$. Keeping all other parameters fixed gives continuously many distinct fitted transition matrices, all yielding the same joint observation densities, as proved above.

This ambiguity cannot be accounted for by state permutations alone. For a fixed fitted state count, there are only finitely many permutations, so the continuously many parameter choices cannot all be relabelings of one parameter.

The initial probabilities can likewise vary whenever a state $j$ with $\pi_j^\star>0$ has at least two substates. For example, setting $\widetilde\pi_{j,1}=\lambda\pi_j^\star$ and $\widetilde\pi_{j,2}=(1-\lambda)\pi_j^\star$, with zero initial mass on any remaining substates of that group, preserves its total initial probability for every $0<\lambda<1$. The joint observation densities again remain unchanged.

Thus the observations cannot distinguish the individual substate allocations in this construction. Their aggregate probabilities and flows remain equal to the true quantities, and their latent location measures also coincide with those of the true model because each group retains its true emission center.
\end{proof}
\section{Additional Experimental Details}
\label{app:synthetic-details}
This appendix gives the fitting and evaluation details for Section~\ref{sec:simulation}, followed by additional results on emission geometry, fitted-center grouping, and reference-center accuracy.

\subsection{Model fitting and evaluation}
\label{app:synthetic-protocol}
\label{app:synthetic-losses}
All estimation experiments use independent trajectories of length $T=6$, known identity emission covariance, and centers constrained to $B_d(0,6)$. Models with different fitted state counts use the same data within each repetition. We implement Baum--Welch updates using scaled forward--backward recursions, following \citet{rabiner1989tutorial}. Each fit uses four initializations and retains the one with the largest final training likelihood. Runs stop when the absolute likelihood increment divided by $nT$ is at most $10^{-7}$ or after $400$ iterations. All selected fits enter the summaries. These stopping rules do not certify global likelihood maximization.

\paragraph{Initialization.} In the scalar experiments, the $k$ initial centers are empirical quantiles at equally spaced levels from $0.08$ to $0.92$, with independent Gaussian perturbations of standard deviation $0.25+0.08j$ at restart $j\in\{0,1,2,3\}$. Initial probabilities are drawn from $\operatorname{Dirichlet}(2,\ldots,2)$. Transition row $i$ is drawn from $\operatorname{Dirichlet}(\mathbf 1+\lambda e_i)$, with an independent $\lambda\sim\operatorname{Unif}[1,8]$. Center updates are projected onto $B_d(0,6)$; probabilities are updated from normalized expected counts without a positive floor. If an expected count is zero, the corresponding center or transition row is retained. The master statistical seed is $202610240$, with separate random streams for data, initialization, and evaluation.

\paragraph{Latent losses.} For $I=(t_1,\ldots,t_r)$, we compute the state-tuple probability as
\[
q_{t_1}(i_1)\prod_{\ell=2}^r\left(A^{t_\ell-t_{\ell-1}}\right)_{i_{\ell-1}i_\ell}.
\]
We sum these probabilities within the assigned groups before taking differences from the true law. This also applies to nonconsecutive times. The Wasserstein losses use the normalized block costs from Section~\ref{sec:latent_recovery_loss}. We solve the finite transport problem by linear programming with HiGHS; in the scalar independent-state comparisons, we use monotone matching of the one-time measures. No Markov approximation is imposed on a grouped process.

\paragraph{Observation error and uncertainty.} We evaluate HMM densities by the forward algorithm. For true and fitted densities $p,q$, let $m=(p+q)/2$. We use
\begin{equation}
h^2(p,q)=\mathbb E_{X\sim m}\left[\frac{(\sqrt{p(X)}-\sqrt{q(X)})^2}{p(X)+q(X)}\right].
\label{eq:synthetic-hellinger-bridge}
\end{equation}
The expectation is estimated using $4096$ independent trajectories from each model, independently of training, with equal weight given to the two sample means. Reported Hellinger errors are the square roots of these estimates. Figures~\ref{fig:synthetic-statistical} and~\ref{fig:synthetic-overfit} use means and pointwise percentile bootstrap intervals from $2000$ resamples of whole data sets. The restricted benchmark in Figure~\ref{fig:synthetic-rates-matching} uses $4000$ bootstrap resamples over its $256$ data sets. Center-matching frequencies use $95\%$ Wilson intervals. Restarts are not treated as independent data sets.

\begin{figure}[t]
\centering
\includegraphics[width=.95\textwidth]{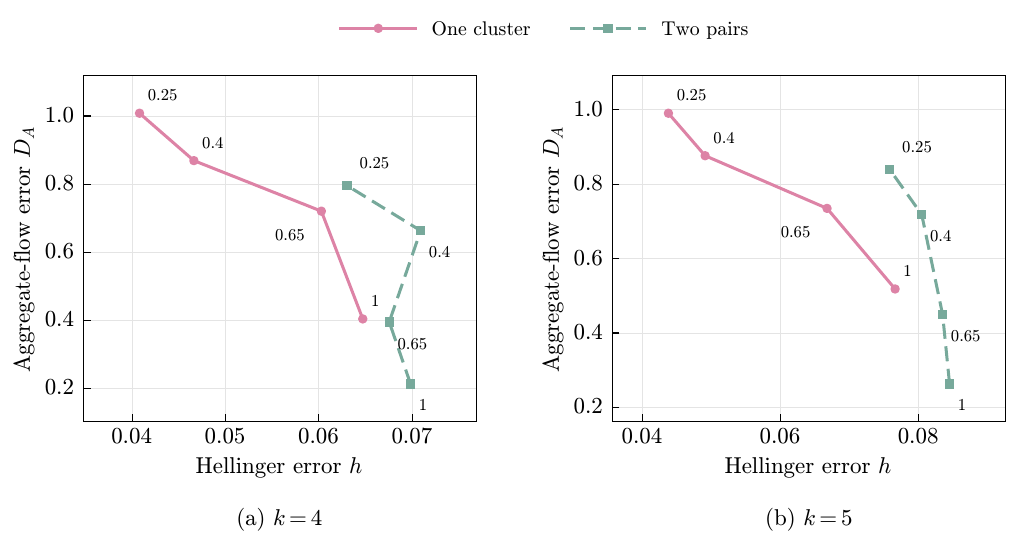}
\caption{\textbf{Observation fit and aggregate recovery under different cluster geometries.} Each point reports mean $h$ and $D_A$ over $12$ independent data sets with $K=4$, $n=512$, and $T=6$. Labels give the minimum separation $\delta$; panels use $k=4$ and $k=5$.}
\label{fig:synthetic-geometry}
\end{figure} 

\subsection{Separation Comparisons and the Global Likelihood Benchmark}
\label{app:synthetic-finite}
\paragraph{Moment-matched pairs.} For $s=2,3,4$, the true unscaled nodes are $z_j=j-(s+1)/2$ with weights $1/s$. In the exact-fit comparison, we fix the last candidate node at $z_s+0.01$ and solve for the remaining nodes and weights by matching moments of orders $0,\ldots,2s-2$. For the over-specified comparison, the first candidate group has two nodes at $y_1\pm0.01$, each with half the group weight; the group centers and weights are chosen to match moments of orders $0,\ldots,2s-1$. These are the constructions used in Sections~\ref{sec:synthetic-exact} and~\ref{sec:synthetic-overconditioning}. The four-state comparison uses the two-pair construction described in Section~\ref{sec:synthetic-exact}.

We scale the within-cluster nodes by each $\delta$ in the main-text grid and evaluate the one-time Hellinger distance $h_1$ by Gauss--Hermite quadrature \citep{golub1969calculation}. The full-block distance follows from $h_T^2=1-(1-h_1^2)^T$. Moment equations and density integrals use $80$-digit arithmetic, with quadrature orders $128$ and $192$ compared for stability. The exponents in Tables~\ref{tab:exact_separation_exponents} and~\ref{tab:overfit_separation_exponents} are the negative least-squares slopes of the log error ratios against $\log\delta$ over $\delta\in\{0.18,0.12,0.08,0.05\}$. The theoretical exponents are used only for comparison after the densities and losses have been evaluated.

\paragraph{Joint over-specified comparison.}
\label{app:synthetic-joint}
We use the perturbations in Proposition~\ref{prop:overfit-sharpness}. With $z_j=j-(s+1)/2$ and $g(z)=\prod_j(z-z_j)$, set $a_j=1/g'(z_j)$, $b_j=-g''(z_j)/g'(z_j)^3$, and $c_j=1/g'(z_j)^2$. Writing $\phi$ for the standard Gaussian density, define
\[
f_\delta(x)=\frac1s\sum_j\phi(x-\delta z_j),\qquad A_\delta(x)=\sum_j a_j\phi(x-\delta z_j),\qquad B_\delta(x)=\sum_j\{b_j+c_j\delta(x-\delta z_j)\}\phi(x-\delta z_j).
\]
The factorized first density variation gives the Hellinger coefficient
\begin{equation}
c_h(\delta)=\frac1{\sqrt8}\left(\int\frac{B_\delta(x)^2}{f_\delta(x)}\,dx\right)^{1/2}\left(\int\frac{A_\delta(x)^2}{f_\delta(x)}\,dx\right)^{(r-1)/2}.
\label{eq:synthetic-tangent-hellinger}
\end{equation}
We evaluate these scalar integrals by the same quadrature procedure. For a grid tuple $\boldsymbol j$, the mass derivative is $U_{\boldsymbol j}=b_{j_1}\prod_{\ell=2}^r a_{j_\ell}$, so $c_D=\sum_{\boldsymbol j}|U_{\boldsymbol j}|$. To compute $c_W$, we solve the grid-flow problem with divergence $U$, allowing mass to pass through occupied intermediate nodes. Directed nearest-neighbor edges have cost $1/r$, and the optimal cost is multiplied by $\delta^2$. The clock-state realization has $k=rs^r+s$ states, as specified in Proposition~\ref{prop:overfit-sharpness}. Figure~\ref{fig:synthetic-joint-conditioning} compares $c_D/c_h$ and $c_W/c_h$ at fixed $s=3$ and $r=2,3,4$.

\paragraph{Restricted likelihood benchmark.}
\label{app:synthetic-single}
Each data set contains $N=6n$ independent standard Gaussian observations. The single-Gaussian MLE is the sample mean projected onto $[-6,6]$. For the equally weighted symmetric mixture with centers $\pm a$, the log likelihood, up to a constant, is
\[
\ell_N(a)=-\frac{Na^2}{2}+\sum_{i=1}^N\log\cosh(aY_i),\qquad a\in[0,6].
\]
For $a>0$, $\ell_N'(a)/(Na)=N^{-1}\sum_iY_i^2\tanh(aY_i)/(aY_i)-1$ is decreasing in $a$, since $\tanh(x)/x$ decreases with $|x|$. The maximizer is zero if $N^{-1}\sum_iY_i^2\le1$. Otherwise it is six if the derivative remains nonnegative there, or the unique positive root of the derivative, found by a bracketed scalar solver. This gives a global maximizer within the stated symmetric submodel. Zero estimates are retained in the means and bootstrap intervals.

\subsection{Estimation under different cluster geometries}
\label{app:synthetic-extra-results}
\label{app:synthetic-geometry}
We supplement the controlled comparisons in Section~\ref{sec:synthetic-exact} by fitting models under two emission geometries with the same $K=4$ and minimum separation $\delta$. The single-cluster centers are $\delta(-1.5,-0.5,0.5,1.5)$, while the two-pair centers are
\[
(-2.5-\delta/2,-2.5+\delta/2,2.5-\delta/2,2.5+\delta/2).
\]
The initial distribution is uniform and both geometries use $A^\star=0.76I_4+0.24\mathbf1\mathbf1^\top/4$. For each $\delta\in\{1,0.65,0.4,0.25\}$, we generate $12$ independent data sets with $n=512$ and fit $k=4,5$ models using the procedure above.

Figure~\ref{fig:synthetic-geometry} shows that smaller observation error need not imply better aggregate recovery. For the single cluster with $k=4$, reducing $\delta$ from $1$ to $0.25$ decreases mean $h$ from $0.065$ to $0.041$, while mean $D_A$ increases from $0.404$ to $1.008$. At $\delta=0.25$, the two-pair model has larger mean $h$ but smaller mean $D_A$. The same pattern occurs for $k=5$. These fitted-model results complement the controlled separation scans, where the loss-to-Hellinger ratios isolate the dependence on cluster size.

\subsection{Grouping based on fitted centers}
\label{app:synthetic-grouping}
We reuse the $k=4,5$ fits from Section~\ref{sec:synthetic-statistical}. We apply complete-linkage clustering to the fitted centers and cut the hierarchy at the known true group count $K=3$. Each group is represented by its mean center, weighted by the fitted states' average occupancies over the six times. To compare the resulting groups with the true states, we match their representatives one-to-one with the true centers by minimizing the total distance between matched centers. The true centers are not used to form the partition. We then compute $D_A$ using this assignment and compare it with the nearest-true-center assignment.

\begin{table}[t]
\centering
\caption{\textbf{Aggregate recovery using fitted-center groups.} Entries are means over the same $24$ independent data sets as Table~\ref{tab:statistical-recovery}. ``Same groups'' counts partitions agreeing with the nearest-true-center assignment after label alignment. The fitted-center partition uses the known group count $K=3$.}
\label{tab:synthetic-grouping}
\setlength{\tabcolsep}{9pt}
\renewcommand{\arraystretch}{1.12}
\arrayrulecolor{black}
\begin{tabular}{ccccc}
\toprule
& & & \multicolumn{2}{c}{Aggregate-flow error} \\
\cmidrule(lr){4-5}
$n$ & $k$ & Same groups & Oracle $D_A$ & Fitted-group $D_A$ \\
\midrule
\multirow{2}{*}{128} & 4 & 23/24 & 0.255 & 0.258 \\
 & 5 & 22/24 & 0.290 & 0.286 \\
\arrayrulecolor{black!12}\midrule[0.3pt]\arrayrulecolor{black}
\multirow{2}{*}{512} & 4 & 23/24 & 0.114 & 0.124 \\
 & 5 & 24/24 & 0.136 & 0.136 \\
\arrayrulecolor{black!12}\midrule[0.3pt]\arrayrulecolor{black}
\multirow{2}{*}{2048} & 4 & 24/24 & 0.061 & 0.061 \\
 & 5 & 24/24 & 0.074 & 0.074 \\
\bottomrule
\end{tabular}
\end{table}

As shown in Table~\ref{tab:synthetic-grouping}, the two assignments agree for all $24$ data sets at $n=2048$. At smaller sample sizes, they differ in up to two data sets per configuration, while the mean flow errors remain close. Fitted-center grouping therefore gives similar aggregate recovery in this experiment. This comparison assumes $K$ is known and does not establish the reference-accuracy condition in Section~\ref{sec:estimated-reference-grouping}.

\subsection{Sensitivity to reference-center error}
\label{app:synthetic-reference-sensitivity}
To examine Proposition~\ref{prop:reference-grouping-stability}, we keep each fitted HMM fixed and perturb its true centers to obtain reference centers $\xi_j=\mu_j^\star+\tau u_j$, where $\tau/\delta_{\mathrm{sep}}\in\{0,0.10,0.20,0.25\}$. Each $u_j$ is an independently sampled unit direction, a random sign in one dimension. At every level, we draw eight sets of references and reassign fitted states to the nearest reference. We compare the resulting grouped distributions at $I=\{1,2\}$ and $I=\{1,3,6\}$ with the true distributions.

This calculation uses $216$ fits from the sample-size experiment, $192$ from the emission-geometry experiment, and $144$ from the matching experiment with $n\le2048$, together with $80$ fits from two additional settings described below. The $48$ matching fits at $n=8192$ are excluded. It contains $632\times4\times8\times2=40{,}448$ comparisons. Reference perturbations use the same master seed $202610240$, task identifiers, and a separate random stream. Table~\ref{tab:synthetic-reference-sensitivity} reports three-time recovery; the assignment-change counts are the same for the two-time calculation.

The first additional setting varies the initial-state probabilities in a two-state Gaussian HMM with emission centers $(-1,1)$ and transition matrix
\[
A^\star=
\begin{pmatrix}
0.8&0.2\\
0.01&0.99
\end{pmatrix}.
\]
For each $\rho\in\{0.3,0.03\}$, we set $\pi^\star=(\rho,1-\rho)$ and generate $16$ independent data sets, each containing $n=1024$ trajectories of length $T=6$. Fitting a two-state model to each data set gives $32$ fits.

The second setting uses three noncollinear emission centers $\delta(-1,0)$, $\delta(0,0.6)$, and $\delta(1,0)$ in two dimensions, with identity emission covariance, initial probabilities $(0.3,0.4,0.3)$, and transition matrix $A^\star=0.75I_3+0.25\mathbf1\mathbf1^\top/3$. For each $\delta\in\{0.4,0.8\}$, we generate $12$ independent data sets of $n=512$ trajectories of length $T=6$ and fit models with $k=3,4$ to each data set, giving another $48$ fits.

\begin{table}[t]
\centering
\setlength{\tabcolsep}{7pt}
\renewcommand{\arraystretch}{1.12}
\arrayrulecolor{black}
\begin{tabular}{ccccc}
\toprule
& & \multicolumn{2}{c}{Three-time probability error} & \\
\cmidrule(lr){3-4}
$\tau/\delta_{\mathrm{sep}}$ & Changed assignments & Oracle & Reference & Max. increase \\
\midrule
0.00 & 0/5056 & 0.425 & 0.425 & 0.000 \\
0.10 & 225/5056 & 0.425 & 0.426 & 0.478 \\
0.20 & 437/5056 & 0.425 & 0.430 & 0.891 \\
0.25 & 546/5056 & 0.425 & 0.435 & 0.891 \\
\bottomrule
\end{tabular}
\caption{\textbf{Grouped recovery under reference-center perturbations.} Oracle and Reference report mean three-time probability errors $D_{\{1,3,6\}}$ under grouping by true and perturbed centers, respectively, over $5056$ reference configurations at each perturbation level. ``Changed assignments'' counts reference configurations in which at least one fitted state changes groups relative to true-center grouping. The last column gives the largest increase in error caused by changing the grouping, with each fitted HMM held fixed.}
\label{tab:synthetic-reference-sensitivity}
\end{table}

\begin{samepage}
All $40{,}448$ comparisons satisfy equation~\eqref{eq:reference-grouping-stability}, including configurations where an assignment changes. At $\tau=0.10\delta_{\mathrm{sep}}$, assignments change in $225$ of $5056$ configurations, increasing to $546$ at $0.25\delta_{\mathrm{sep}}$. Mean probability errors change little, although individual increases can be much larger. The result therefore allows assignments to change within the permitted radius while controlling their effect on grouped recovery. Since the references are constructed from known true centers, this experiment evaluates the conditional stability bound rather than the accuracy of an estimated reference-center procedure.\par
\end{samepage}

\clearpage
\bibliographystyle{plainnat}
\bibliography{references_trang}

\end{document}